\documentclass[12pt,reqno]{amsart}
\def\draft{n}
\usepackage{
  fullpage,graphicx,dbnsymb,amsmath,amssymb,floatflt,picins,
  multicol,mathtools,stmaryrd,xcolor
}
\usepackage[all]{xy}

\usepackage{hyperref}\hypersetup{colorlinks,
   linkcolor={green!50!black},
   citecolor={green!50!black},
   urlcolor=blue
}

\if\draft y
  \usepackage[color]{showkeys}
  \usepackage{everypage,draftwatermark}
  \SetWatermarkScale{4}
  \SetWatermarkLightness{0.9}
\fi

\theoremstyle{plain}
\newtheorem{theorem}{Theorem}[section]

\newtheorem{proposition}[theorem]{Proposition}
\newtheorem{fact}[theorem]{Fact}
\newtheorem{lemma}[theorem]{Lemma}

\newtheorem{claim}[theorem]{Claim}

\newtheorem{conjecture}[theorem]{Conjecture}

\theoremstyle{definition}
\newtheorem{definition}[theorem]{Definition}
\newtheorem{defwarn}[theorem]{Definition and Warning}
\newtheorem{image}[theorem]{Image}

\newtheorem{problem}[theorem]{Problem}

\theoremstyle{remark}
\newtheorem{comment}[theorem]{Comment}
\newtheorem{comments}[theorem]{Comments}
\newtheorem{example}[theorem]{Example}
\newtheorem{exercise}[theorem]{Exercise}

\newtheorem{remark}[theorem]{Remark}
\newtheorem{warning}[theorem]{Warning}

\newlength{\standardunitlength}
\catcode`\@=11
\long\def\@makecaption#1#2{%
    \vskip 10pt
    \setbox\@tempboxa\hbox{
      \small\sf{\bfcaptionfont #1. }\ignorespaces #2}%
    \ifdim \wd\@tempboxa >\captionwidth {%
        \rightskip=\@captionmargin\leftskip=\@captionmargin
        \unhbox\@tempboxa\par}%
      \else
        \hbox to\hsize{\hfil\box\@tempboxa\hfil}%
    \fi}
\font\bfcaptionfont=cmssbx10 scaled \magstephalf
\newdimen\@captionmargin\@captionmargin=2\parindent
\newdimen\captionwidth\captionwidth=\hsize
\catcode`\@=12

\newcommand{\ad}{\operatorname{ad}}

\newcommand{\Aut}{\operatorname{Aut}}

\newcommand{\gr}{\operatorname{gr}}
\newcommand{\im}{\operatorname{im}}

\newcommand{\sign}{\operatorname{sign}}

\newcommand{\tr}{\operatorname{tr}}

\def\qed{{\hfill\text{$\Box$}}}

\newlength{\globalparindent}
\newenvironment{myitemize}{
        \begin{list}{$\bullet$}{\setlength{\leftmargin}{16pt}
        \setlength{\labelwidth}{12pt}
        \setlength{\labelsep}{4pt}}
}{
        \end{list}
}

\def\arXiv#1{{\href{http://front.math.ucdavis.edu/#1}{arXiv:#1}}}

\def\bbN{{\mathbb N}}
\def\bbQ{{\mathbb Q}}
\def\bbR{{\mathbb R}}
\def\bbZ{{\mathbb Z}}
\def\calA{{\mathcal A}}
\def\calB{{\mathcal B}}
\def\calAcf{{\mathcal A}^\text{\it cf}}
\def\calC{{\mathcal C}}
\def\calD{{\mathcal D}}
\def\calF{{\mathcal F}}
\def\calG{{\mathcal G}}
\def\calI{{\mathcal I}}
\def\calK{{\mathcal K}}
\def\calL{{\mathcal L}}
\def\calO{{\mathcal O}}
\def\calP{{\mathcal P}}
\def\calR{{\mathcal R}}
\def\calS{{\mathcal S}}
\def\calT{{\mathcal T}}
\def\calU{{\mathcal U}}
\def\fraka{{\mathfrak a}}
\def\frakg{{\mathfrak g}}
\def\tilE{{\tilde{E}}}

\def\ad{\operatorname{ad}}
\def\Aut{\operatorname{Aut}}
\def\Autop{\operatorname{Aut}^\text{op}}
\def\fil{\operatorname{fil}\,}
\def\gr{\operatorname{gr}\,}
\def\IAM{\mathit{IAM}}
\def\proj{{\operatorname{proj}\,}}
\def\projs{{\operatorname{proj}}}
\newcommand{\Ass}{\operatorname{Ass}}
\def\divop{\operatorname{div}}
\def\CA{\operatorname{CA}}

\def\attr{\operatorname{\mathfrak{tr}}}
\def\lie{\operatorname{\mathfrak{lie}}}
\def\der{\operatorname{\mathfrak{der}}}
\def\sder{\operatorname{\mathfrak{sder}}}
\def\tder{\operatorname{\mathfrak{tder}}}
\def\TAut{\operatorname{TAut}}
\def\SAut{\operatorname{SAut}}

\def\FA{{\mathit F\!A}}
\def\uB{{\mathit u\!B}}
\def\PuB{{\mathit P\!u\!B}}
\def\vB{{\mathit v\!B}}
\def\vT{{\mathit v\!T}}
\def\PvB{{\mathit P\!v\!B}}
\def\wB{{\mathit w\!B}}
\def\swB{{\mathit s\!w\!B}}
\def\uT{{\mathit u\!T}}
\def\wT{{\mathit w\!T}}
\def\wTF{{\mathit w\!T\!F}}
\def\wTFo{{\mathit w\!T\!F^o}}
\def\sKTG{{\mathit s\!K\!T\!G}}
\def\PwB{{\mathit P\!w\!B}}
\def\sl{{\mathit sl}}

\def\aft{{\overrightarrow{4T}}}
\def\aAS{{\overrightarrow{AS}}}
\def\aSTU{{\overrightarrow{STU}}}
\def\aIHX{{\overrightarrow{IHX}}}
\def\Rs{{R$1^{\!s}$}}

\def\up{{\uparrow}}

\def\mathsize#1#2{{\begin{array}{c}\text{#1$\!#2\!$}\end{array}}}
\def\pstex#1{\begin{array}{c}
  \input figs/#1.pstex_t
  \if\draft y
    \smash{\makebox[0in]{\color{labelkey}{\tt figs/#1}}}
  \fi
\end{array}}

\def\draftcut{\if\draft y \newpage \fi}

\def\wClip#1{{\href{http://katlas.math.toronto.edu/drorbn/dbnvp/wClips-#1.php}{wClip:#1}}}

\def\wClipStart#1#2{\marginpar{\bf\footnotesize\begin{center}%
\vskip 12pt
wClip%
\newline\href{http://katlas.math.toronto.edu/drorbn/dbnvp/wClips-#1.php}{#1}%
\vskip 4pt%
\par\includegraphics[width=0.8in]{frames/wClips-#1_#2_200.ps}
\newline starts%
\end{center}}}

\def\wClipAt#1#2#3#4{\marginpar{\bf\footnotesize\begin{center}%
\vskip 12pt
wClip%
\newline\href{http://katlas.math.toronto.edu/drorbn/dbnvp/wClips-#1.php}{#1}%
\vskip 4pt%
\par\includegraphics[width=0.8in]{frames/wClips-#1_#2-#3-#4_200.ps}
\newline at #2:#3:#4%
\end{center}}}

\def\wClipComment#1#2#3{\marginpar{\footnotesize\begin{center}{\bf%
wClip%
\newline\href{http://katlas.math.toronto.edu/drorbn/dbnvp/wClips-#1.php}{#1}}%
\vskip 4pt%
\par\includegraphics[width=0.8in]{frames/wClips-#1_#2_200.ps}
\newline #3%
\end{center}}}

\def\wClipEnd#1{\marginpar{\bf\footnotesize\begin{center}%
\vskip -12mm
wClip%
\newline\href{http://katlas.math.toronto.edu/drorbn/dbnvp/wClips-#1.php}{#1}%
\newline ends%
\end{center}}}

\def\glos#1{\setlength{\fboxsep}{0pt}\colorbox{yellow}{$#1$}}
\def\glost#1{\setlength{\fboxsep}{0pt}\colorbox{yellow}{#1}}

\begin{document}
\newdimen\captionwidth\captionwidth=\hsize
\setcounter{secnumdepth}{4}

\title[Finite Type Invariants of w-Knotted Objects]{Finite Type Invariants
  of w-Knotted Objects: From Alexander to Kashiwara and Vergne}

\author{Dror~Bar-Natan and Zsuzsanna Dancso}
\address{
  Department of Mathematics\\
  University of Toronto\\
  Toronto Ontario M5S 2E4\\
  Canada
}
\email{drorbn@math.toronto.edu, zsuzsi@math.toronto.edu}
\urladdr{
  \url{http://www.math.toronto.edu/~drorbn},
  \url{http://www.math.toronto.edu/zsuzsi}
}

\date{first edition Sep.\ 27, 2013, this edition May.~11,~2014. The
\arXiv{1309.7155} edition may be older}


\thanks{This work was partially supported by NSERC grant RGPIN 262178.
  Electronic version, videos (wClips) and related files at~\cite{WKO},
  \url{http://www.math.toronto.edu/~drorbn/papers/WKO/}.
}

\begin{abstract}
  w-Knots, and more generally, w-knotted objects (w-braids, w-tangles,
etc.)  make a class of knotted objects which is \underline{w}ider but
\underline{w}eaker than their ``\underline{u}sual'' counterparts. To
get (say) w-knots from u-knots, one has to allow non-planar ``virtual''
knot diagrams, hence enlarging the the base set of knots. But then one
imposes a new relation, the ``overcrossings commute'' relation, further
beyond the ordinary collection of Reidemeister moves, making w-knotted
objects a bit weaker once again.

The group of w-braids was studied (under the name
``\underline{w}elded braids'') by Fenn, Rimanyi and
Rourke~\cite{FennRimanyiRourke:BraidPermutation} and was shown to
be isomorphic to the McCool group~\cite{McCool:BasisConjugating}
of ``basis-conjugating'' automorphisms of a free group $F_n$ ---
the smallest subgroup of $\Aut(F_n)$ that contains both braids and
permutations. Brendle and Hatcher~\cite{BrendleHatcher:RingsAndWickets},
in work that traces back to Goldsmith~\cite{Goldsmith:MotionGroups},
have shown this group to be a group of movies of flying rings in
$\bbR^3$. Satoh~\cite{Satoh:RibbonTorusKnots} studied several classes of
w-knotted objects (under the name ``\underline{w}eakly-virtual'') and has
shown them to be closely related to certain classes of knotted surfaces
in $\bbR^4$. So w-knotted objects are algebraically and topologically 
interesting.

In this article we study finite type invariants of several
classes of w-knotted objects.  Following Berceanu and
Papadima~\cite{BerceanuPapadima:BraidPermutation}, we construct
homomorphic universal finite type invariants of w-braids and of
w-tangles. We find that the universal finite type invariant of w-knots
is more or less the Alexander polynomial (details inside).

Much as the spaces $\calA$ of chord diagrams for ordinary knotted
objects are related to metrized Lie algebras, we find that the spaces
$\calA^w$ of ``arrow diagrams'' for w-knotted objects are related to
not-necessarily-metrized Lie algebras. Many questions concerning w-knotted
objects turn out to be equivalent to questions about Lie algebras. Most
notably we find that a homomorphic universal finite type invariant of
w-knotted foams is essentially the same as a solution of
the Kashiwara-Vergne~\cite{KashiwaraVergne:Conjecture} conjecture and
much of the Alekseev-Torossian~\cite{AlekseevTorossian:KashiwaraVergne}
work on Drinfel'd associators and Kashiwara-Vergne can be re-interpreted
as a study of w-knotted trivalent graphs.

The true value of w-knots, though, is likely to emerge later, for we
expect them to serve as a \underline{w}armup example for what we expect
will be even more interesting --- the study of \underline{v}irtual
knots, or v-knots. We expect v-knotted objects to provide the global
context whose projectivization (or ``associated graded structure'')
will be the Etingof-Kazhdan theory of deformation quantization of Lie
bialgebras~\cite{EtingofKazhdan:BialgebrasI}.

  \vskip 3mm
  {\color{red} This paper was split in two and became the first two
  parts of a four-part series (\cite{WKO1}--\cite{WKO4}).  The remaining
  relevance of this paper is due to the series of videotaped lectures
  (wClips) that are linked here.}
\end{abstract}

\maketitle

\clearpage
{\small \tableofcontents}

\draftcut
\section{Introduction} \label{sec:intro}

\subsection{Dreams} \label{subsec:dreams} We \wClipStart{120111-1}{0-00-25}
\marginpar{\bf\footnotesize\raggedright wClips are explained in
Section~\ref{subsec:wClips}.} have a dream\footnote{Understanding
the authors' history and psychology ought never be necessary to
understand their papers, yet it may be helpful. Nothing material
in the rest of this paper relies on Section~\ref{subsec:dreams}.},
at least partially founded on reality, that many of the difficult
algebraic equations in mathematics, especially those that are
written in graded spaces, more especially those that are related in
one way or another to quantum groups~\cite{Drinfeld:QuantumGroups},
and even more especially those related to the work of Etingof and
Kazhdan~\cite{EtingofKazhdan:BialgebrasI}, can be understood, and indeed,
would appear more natural, in terms of finite type invariants of various
topological objects.

We believe this is the case for Drinfel'd's theory
of associators~\cite{Drinfeld:QuasiHopf}, which can be
interpreted as a theory of well-behaved universal finite type
invariants of parenthesized tangles\footnote{``$q$-tangles''
in~\cite{LeMurakami:Universal}, ``non-associative tangles''
in~\cite{Bar-Natan:NAT}.}~\cite{LeMurakami:Universal, Bar-Natan:NAT},
and even more elegantly, as a theory of universal finite type invariants
of knotted trivalent graphs~\cite{Dancso:KIforKTG}.

We believe this is the case for Drinfel'd's ``Grothendieck-Teichmuller
group''~\cite{Drinfeld:GalQQ} which is better understood as a
group of automorphisms of a certain algebraic structure, also
related to universal finite type invariants of parenthesized
tangles~\cite{Bar-Natan:Associators}.

And we're optimistic, indeed we believe, that sooner or later the
work of Etingof and Kazhdan~\cite{EtingofKazhdan:BialgebrasI}
on quantization of Lie bialgebras will be re-interpreted as a
construction of a well-behaved universal finite type invariant of
virtual knots~\cite{Kauffman:VirtualKnotTheory} or of some other class
of virtually knotted objects. Some steps in that direction were taken
by Haviv~\cite{Haviv:DiagrammaticAnalogue}.

We have another dream, to construct a useful ``Algebraic Knot Theory''. As
at least a partial writeup exists~\cite{Bar-Natan:AKT-CFA},
we'll only state that an important ingredient necessary to
fulfil that dream would be a ``closed form''\footnote{The
phrase ``closed form'' in itself requires an explanation. See
Section~\ref{subsec:ClosedForm}. \label{foot:ClosedForm}} formula for an
associator, at least in some reduced sense. Formulae for associators or
reduced associators were in themselves the goal of several studies
undertaken for various other reasons~\cite{LeMurakami:HOMFLY,
Lieberum:gl11, Kurlin:CompressedAssociators, LeeP:ClosedForm}.

\draftcut \subsection{Stories}

Thus the first named author, DBN, was absolutely delighted
when in January 2008 Anton Alekseev described to him his joint
work~\cite{AlekseevTorossian:KashiwaraVergne} with Charles
Torossian --- Anton told DBN that they found a relationship between the
Kashiwara-Vergne conjecture~\cite{KashiwaraVergne:Conjecture},
a cousin of the Duflo isomorphism (which DBN already knew to be
knot-theoretic~\cite{Bar-NatanLeThurston:TwoApplications}), and
associators taking values in a space called $\sder$, which he could
identify as ``tree-level Jacobi diagrams'', also a knot-theoretic
space related to the Milnor invariants~\cite{Bar-Natan:Homotopy,
HabeggerMasbaum:Milnor}. What's more, Anton told DBN that in certain
quotient spaces the Kashiwara-Vergne conjecture can be solved explicitly;
this should lead to some explicit associators!

So DBN spent the following several months trying to
understand~\cite{AlekseevTorossian:KashiwaraVergne}, and this
paper is a summary of these efforts. The main thing we learned is that
the Alekseev-Torossian paper, and with it the Kashiwara-Vergne
conjecture, fit very nicely with our first dream recalled above,
about interpreting algebra in terms of knot theory. Indeed much
of~\cite{AlekseevTorossian:KashiwaraVergne} can be reformulated as a
construction and a discussion of a well-behaved universal finite type
invariant $Z$ of a certain class of knotted objects (which we will call here
w-knotted), a certain natural quotient of the space of virtual knots
(more precisely, virtual trivalent tangles). And our hopes remain high
that later we (or somebody else) will be able to exploit this relationship
in directions compatible with our second dream recalled above, on the
construction of an ``algebraic knot theory''.

The story, in fact, is prettier than we were hoping for, for it has the
following additional qualities:

\begin{myitemize}

\item w-Knotted objects are quite interesting in themselves: as
stated in the abstract, they are related to combinatorial group
theory via ``basis-conjugating'' automorphisms of a free group $F_n$,
to groups of movies of flying rings in $\bbR^3$, and more generally, to
certain classes of knotted surfaces in $\bbR^4$. The references include
\cite{BrendleHatcher:RingsAndWickets, FennRimanyiRourke:BraidPermutation,
Goldsmith:MotionGroups, McCool:BasisConjugating, Satoh:RibbonTorusKnots}.

\item The ``chord diagrams'' for w-knotted objects (really, these are ``arrow
diagrams'') describe formulae for invariant tensors in spaces pertaining to
not-necessarily-metrized Lie algebras in much of the same way as ordinary
chord diagrams for ordinary knotted objects describe formulae for invariant
tensors in spaces pertaining to metrized Lie algebras. This observation is
bound to have further implications.

\item Arrow diagrams also describe the Feynman diagrams of topological BF
theory \cite{CattaneoCotta-RamusinoMartellini:Alexander, CCFM:BF34} and of a
certain class of Chern-Simons theories~\cite{Naot:BF}. Thus it is likely that
our story is directly related to quantum field theory\footnote{Some
non-perturbative relations between BF theory and w-knots was discussed by
Baez, Wise and Crans~\cite{BaezWiseCrans:ExoticStatistics}.}.

\item When composed with the map from knots to w-knots, $Z$ becomes the
Alexander polynomial. For links, it becomes an invariant stronger than the
multi-variable Alexander polynomial which contains the multi-variable
Alexander polynomial as an easily identifiable reduction. On other
w-knotted objects $Z$ has easily identifiable reductions that can be
considered as ``Alexander polynomials'' with good behaviour relative
to various knot-theoretic operations --- cablings, compositions
of tangles, etc. There is also a certain specific reduction of $Z$
that can be considered as the ``ultimate Alexander polynomial'' ---
in the appropriate sense, it is the minimal extension of the Alexander
polynomial to other knotted objects which is well behaved under a whole
slew of knot theoretic operations, including the ones named above.

\end{myitemize}

\begin{figure}
\[
  \def\uT{\parbox[t]{1.875in}{\small
    Ordinary (\underline{u}sual) knotted objects in 3D --- braids,
    knots, links, tangles, knotted graphs, etc.
  }}
  \def\vT{\parbox[t]{1.875in}{\small
    \underline{V}irtual knotted objects --- ``algebraic'' knotted objects,
    or ``not specifically embedded'' knotted objects; knots drawn on a
    surface, modulo stabilization.
  }}
  \def\wT{\parbox[t]{1.875in}{\small
    Ribbon knotted objects in 4D; ``flying rings''. Like v, but also with
    ``overcrossings commute''.
  }}
  \def\uC{\parbox[t]{1.875in}{\small
    Chord diagrams and Jacobi diagrams, modulo $4T$, $STU$, $IHX$, etc.
  }}
  \def\vC{\parbox[t]{1.875in}{\small
   Arrow diagrams and v-Jacobi diagrams, modulo $6T$ and various
   ``directed'' $STU$s and $IHX$s, etc.
  }}
  \def\wC{\parbox[t]{1.875in}{\small
    Like v, but also with ``tails commute''. Only ``two in one out''
    internal vertices.
  }}
  \def\uL{\parbox[t]{1.875in}{\small
    Finite dimensional metrized Lie algebras, representations,  and
    associated spaces.
  }}
  \def\vL{\parbox[t]{1.875in}{\small
    Finite dimensional Lie bi-algebras, representations,  and associated
    spaces.
  }}
  \def\wL{\parbox[t]{1.875in}{\small
    Finite dimensional co-commutative Lie bi-algebras (i.e.,
    $\frakg\ltimes\frakg^\ast$), representations,  and associated
    spaces.
  }}
  \def\uH{\parbox[t]{1.875in}{\small
    The Drinfel'd theory of associators.
  }}
  \def\vH{\parbox[t]{1.875in}{\small
    Likely, quantum groups and the Etingof-Kazhdan theory of quantization
    of Lie bi-algebras.
  }}
  \def\wH{\parbox[t]{1.875in}{\small
    The Kashiwara-Vergne-Alekseev-Torossian theory of convolutions on Lie
    groups and Lie algebras.
  }}
  \input figs/uvw.pstex_t
\]
\caption{The u-v-w Stories} \label{fig:uvw}
\end{figure}

\draftcut \subsection{The Bigger Picture} 
Parallel to the w-story run the possibly more significant u-story
and v-story. The u-story is about u-knots, or more generally,
u-knotted objects (braids, links, tangles, etc.), where ``u'' stands
for \underline{u}sual; hence the u-story is about ordinary knot
theory. The v-story is about v-knots, or more generally, v-knotted
objects, where ``v'' stands for \underline{v}irtual, in the sense of
Kauffman~\cite{Kauffman:VirtualKnotTheory}.

The three stories, u, v, and w, are different from each other. Yet they
can be told along similar lines --- first the knots (topology), then their
finite type invariants and their ``chord diagrams'' (combinatorics),
then those map into certain universal enveloping algebras and similar
spaces associated with various classes of Lie algebras (low algebra),
and finally, in order to construct a ``good'' universal finite type
invariant, in each case one has to confront a certain deeper algebraic
subject (high algebra). These stories are summarized in a table form
in Figure~\ref{fig:uvw}.

u-Knots map into v-knots, and v-knots map into w-knots\footnote{Though the
composition ``$u\to v\to w$'' is not $0$. In fact, the composed map $u\to
w$ is injective.}. The other parts of our stories, the ``combinatorics''
and ``low algebra'' and ``high algebra'' rows of Figure~\ref{fig:uvw},
are likewise related, and this relationship is a crucial part of our
overall theme. Thus we cannot and will not tell the w-story in isolation,
and while it is central to this article, we will necessarily also include
some episodes from the u and v series. \wClipEnd{120111-1}

\subsection{Plans} Our order of proceedings is: w-braids, w-knots,
generalities, w-tangles, w-tangled foams.  For more detailed information
consult the ``Section Summary'' paragraphs below and at the beginning
of each of the sections. An ``odds and ends'' section follows on
page~\pageref{sec:OddsAndEnds}, and a glossary of notation is on
page~\pageref{sec:glossary}.

\def\summarybraids{This section is largely a compilation of existing
literature, though we also introduce the language of arrow diagrams that
we use throughout the rest of the paper. In~\ref{subsec:VirtualBraids}
and~\ref{subsec:wBraids} we define v-braids and then w-braids and
survey their relationship with basis-conjugating automorphisms of free
groups and with ``the group of (horizontal) flying rings in $\bbR^3$''
(really, a group of knotted tubes in $\bbR^4$). In~\ref{subsec:FT4Braids}
we play the usual game of introducing finite type invariants, weight
systems, chord diagrams (arrow diagrams, for this case), and 4T-like
relations. In~\ref{subsec:wBraidExpansion} we define and construct
a universal finite type invariant $Z$ for w-braids --- it turns out
that the only algebraic tool we need to use is the formal exponential
function $\exp(a):=\sum a^n/n!$. In~\ref{subsec:bcomments} we study
some good algebraic properties of $Z$, its injectivity, and its
uniqueness, and we conclude with the slight modifications needed for the
study of non-horizontal flying rings.}

\def\summaryknots{In~\ref{subsec:VirtualKnots} we define v-knots and
w-knots (long v-knots and long w-knots, to be precise) and discuss the map
$v\to w$. In~\ref{subsec:FTforvwKnots} we determine the space of ``chord
diagrams'' for w-knots to be the space $\calA^w(\uparrow)$ of arrow
diagrams modulo $\aft$ and TC relations and in~\ref{subsec:SomeDimensions}
we compute some relevant dimensions. In~\ref{subsec:Jacobi} we show
that $\calA^w(\uparrow)$ can be re-interpreted as a space of trivalent
graphs modulo STU- and IHX-like relations, and is therefore related
to Lie algebras (Sec.~\ref{subsec:LieAlgebras}). This allows us to
completely determine $\calA^w(\uparrow)$.  With no difficulty at
all in~\ref{subsec:Z4Knots} we construct a universal finite type
invariant for w-knots. With a bit of further difficulty we show in
Sec.~\ref{subsec:Alexander} that it is essentially equal to the Alexander
polynomial.}

\def\summaryalg{In this section we define the ``projectivization''
(Sec.~\ref{subsec:Projectivization}) of an arbitrary algebraic structure
(\ref{subsec:AlgebraicStructures}) and introduce the notions of
``expansions'' and ``homomorphic expansions'' (\ref{subsec:Expansions})
for such projectivizations. Everything is so general that practically
anything is an example. The baby-example of quandles is built in into
the section; the braid groups and w-braid groups appeared already in
Section~\ref{sec:w-braids}, yet our main goal is to set the language
for the examples of w-tangles and w-tangled foams, which appear later in
this paper. Both of these examples are types of ``circuit algebras'', and
hence we end this section with a general discussion of circuit algebras
(Sec.~\ref{subsec:CircuitAlgebras}).}

\def\summarytangles{In Sec.~\ref{subsec:vw-tangles} we introduce
v-tangles and w-tangles, the obvious v- and w- counterparts of the
standard knot-theoretic notion of ``tangles'', and briefly discuss their
finite type invariants and their associated spaces of ``arrow diagrams'',
$\calA^v(\uparrow_n)$ and $\calA^w(\uparrow_n)$. We then construct a
homomorphic expansion $Z$, or a ``well-behaved'' universal finite type
invariant for w-tangles. Once again, the only algebraic tool we need to
use is $\exp(a):=\sum a^n/n!$, and indeed, Sec.~\ref{subsec:vw-tangles}
is but a routine extension of parts of Section~\ref{sec:w-knots}.
We break away in Sec.~\ref{subsec:ATSpaces} and show that
$\calA^w(\uparrow_n)\cong\calU(\fraka_n\oplus\tder_n\ltimes\attr_n)$,
where $\fraka_n$ is an Abelian algebra of rank $n$ and where
$\tder_n$ and $\attr_n$, two of the primary spaces used by Alekseev
and Torossian~\cite{AlekseevTorossian:KashiwaraVergne}, have simple
descriptions in terms of words and free Lie algebras. We also show that
some functionals studied in~\cite{AlekseevTorossian:KashiwaraVergne},
$\divop$ and $j$, have a natural interpretation in our language.
In~\ref{subsec:sder} we discuss a subclass of w-tangles called ``special''
w-tangles, and relate them by similar means to Alekseev and Torossian's
$\sder_n$ and to ``tree level'' ordinary Vassiliev theory. Some
conventions are described in Sec.~\ref{subsec:TangleTopology} and the
uniqueness of $Z$ is studied in Sec.\ref{subsec:UniquenessForTangles}.}

\def\summaryfoams{If you have come this far, you must have noticed
the approximate Bolero spirit of this article. In every chapter
a new instrument comes to play; the overall theme remains the
same, but the composition is more and more intricate. In this
chapter we add ``foam vertices'' to w-tangles (and a few lesser
things as well) and ask the same questions we asked before;
primarily, ``is there a homomorphic expansion?''. As we shall
see, in the current context this question is equivalent to the
Alekseev-Torossian~\cite{AlekseevTorossian:KashiwaraVergne} version
of the Kashiwara-Vergne~\cite{KashiwaraVergne:Conjecture} problem and
explains the relationship between these topics and Drinfel'd's theory
of associators.}

\noindent
{\small \begin{multicols}{2}

{\bf Section~\ref{sec:w-braids}, w-Braids.} (page~\pageref{sec:w-braids})
\summarybraids

{\bf Section~\ref{sec:w-knots}, w-Knots.} (page~\pageref{sec:w-knots})
\summaryknots

{\bf Section~\ref{sec:generalities}, Algebraic Structures,
Projectivizations, Expansions, Circuit Algebras.}
(page~\pageref{sec:generalities}) \summaryalg

{\bf Section~\ref{sec:w-tangles}, w-Tangles.}
(page~\pageref{sec:w-tangles}) \summarytangles

{\bf Section~\ref{sec:w-foams}, w-Tangled Foams.}
(page~\pageref{sec:w-foams}) \summaryfoams

\end{multicols}}

\subsection{Acknowledgement} We wish to thank Anton Alekseev, Jana
Archibald, Scott Carter, Karene Chu, Iva Halacheva, Joel Kamnitzer,
Lou Kauffman, Peter Lee, Louis Leung, Dylan Thurston, Lucy Zhang, and Jean-Baptiste Meilhan
for comments and suggestions.

\subsection{wClips} \label{subsec:wClips} \wClipAt{120118-1}{0}{03}{10}
Alongside this paper there is a series of video clips explaining parts
of it. The series as a whole can be found at~\cite{WKO}; references to
specific clips and specific times within clips appear at the margin of
this paper. We thank Peter Lee for contributing \wClip{120201} and Karene
Chu for contributing \wClip{120314}.

\clearpage\draftcut
\section{w-Braids} \label{sec:w-braids}

\begin{quote} \small {\bf Section Summary. }
  \summarybraids
\end{quote}

\subsection{Preliminary: Virtual Braids, or v-Braids.}
\label{subsec:VirtualBraids} \wClipStart{120111-2}{0-00-00}
Our main object of study for this section, w-braids, are best
viewed as ``virtual braids''~\cite{Bardakov:VirtualAndUniversal,
KauffmanLambropoulou:VirtualBraids, BardakovBellingeri:VirtualBraids},
or v-braids, modulo one additional relation. Hence we start with v-braids.

It is simplest to define v-braids in
terms of generators and relations, either algebraically or pictorially.
This can be done in at least two ways --- the easier-at-first but
philosophically-less-satisfactory ``planar'' way, and the harder to
digest but morally more correct ``abstract'' way.\footnote{Compare with a
similar choice that exists in the definition of manifolds, as either
appropriate subsets of some ambient Euclidean spaces (module some
equivalences) or as abstract gluings of coordinate patches (modulo some
other equivalences). Here in the ``planar'' approach of
Section~\ref{subsubsec:Planar} we consider v-braids
as ``planar'' objects, and in the ``abstract approach'' of
Section~\ref{subsubsec:Abstract} they are just ``gluings'' of abstract
``crossings'', not drawn anywhere in particular.}

\subsubsection{The ``Planar'' Way} \label{subsubsec:Planar} For a
natural number $n$ set $\glos{\vB_n}$ to be the group generated by
symbols $\glos{\sigma_i}$ ($1\leq i\leq n-1$), called ``crossings''
and graphically represented by an overcrossing $\overcrossing$ ``between
strand $i$ and strand $i+1$'' (with inverse $\undercrossing$)\footnote{We
sometimes refer to $\overcrossing$ as a ``positive crossing'' and to
$\undercrossing$ as a ``negative crossing''.}, and $\glos{s_i}$, called
``virtual crossings'' and graphically represented by a non-crossing,
$\virtualcrossing$, also ``between strand $i$ and strand $i+1$'',
subject to the following relations:

\begin{myitemize}

\item The subgroup of $\vB_n$ generated by the virtual crossings $s_i$
is the symmetric group $\glos{S_n}$, and the $s_i$'s correspond to the
transpositions $(i,i+1)$. That is, we have
\begin{equation} \label{eq:sRelations}
  s_i^2=1,
  \qquad s_is_{i+1}s_i = s_{i+1}s_is_{i+1},
  \qquad\text{and if $|i-j|>1$ then}
  \qquad s_is_j=s_js_i.
\end{equation}
In pictures, this is
\begin{equation} \label{eq:sRels}
  \def\i{{$i$}}
  \def\ip{{$i\!+\!1$}}
  \def\ipp{{$i\!+\!2$}}
  \def\j{{$j$}}
  \def\jp{{$j\!+\!1$}}
  \input figs/sRels.pstex_t
\end{equation}
Note that we read our braids from bottom to top.

\item The subgroup of $\vB_n$ generated by the crossings $\sigma_i$'s is
the usual braid group $\glos{\uB_n}$, and $\sigma_i$ corresponds to the braiding
of strand $i$ over strand $i+1$. That is, we have
\begin{equation} \label{eq:R3}
  \sigma_i\sigma_{i+1}\sigma_i
    = \sigma_{i+1}\sigma_i\sigma_{i+1},
  \qquad\text{and if $|i-j|>1$ then}
  \qquad \sigma_i\sigma_j=\sigma_j\sigma_i.
\end{equation}
In pictures, dropping the indices, this is
\begin{equation} \label{eq:sigmaRels} \begin{picture}(0,0)%
\includegraphics{figs/sigmaRels.pstex}%
\end{picture}%
%
%
\setlength{\unitlength}{3158sp}%
\begingroup\makeatletter\ifx\SetFigFont\undefined%
\gdef\SetFigFont#1#2#3#4#5{%
  \reset@font\fontsize{#1}{#2pt}%
  \fontfamily{#3}\fontseries{#4}\fontshape{#5}%
  \selectfont}%
\fi\endgroup%
\begin{picture}(5154,924)(1639,-673)
\put(5551,-286){\makebox(0,0)[b]{\smash{{\SetFigFont{10}{12.0}{\rmdefault}{\mddefault}{\updefault}{\color[rgb]{0,0,0}$=$}%
}}}}
\put(2401,-286){\makebox(0,0)[b]{\smash{{\SetFigFont{10}{12.0}{\rmdefault}{\mddefault}{\updefault}{\color[rgb]{0,0,0}$=$}%
}}}}
\end{picture}%
 \end{equation}
The first of these relations is the ``Reidemeister 3 move''\footnote{The
Reidemeister 2 move is the relations $\sigma_i\sigma_i^{-1}=1$ which is
part of the definition of ``a group''. There is no Reidemeister 1 move in
the theory of braids.} of knot theory. The second is sometimes called
``locality in space''~\cite{Bar-Natan:NAT}.

\item Some ``mixed relations'',
\begin{equation} \label{eq:MixedRelations}
  s_i\sigma^{\pm 1}_{i+1}s_i = s_{i+1}\sigma^{\pm 1}_is_{i+1},
  \qquad\text{and if $|i-j|>1$ then}
  \qquad s_i\sigma_j=\sigma_js_i.
\end{equation}
In pictures, this is
\begin{equation} \label{eq:MixedRels}
  \input figs/MixedRels.pstex_t
\end{equation}

\end{myitemize}

\begin{remark} \label{rem:Skeleton} The ``skeleton'' of a v-braid $B$
is the set of strands appearing in it, retaining the association
between their beginning and ends but ignoring all the crossing
information. More precisely, it is the permutation induced by
tracing along $B$, and even more precisely it is the image of $B$
via the ``skeleton morphism'' $\glos{\varsigma}\colon\vB_n\to S_n$ defined
by $\varsigma(\sigma_i)=\varsigma(s_i)=s_i$ (or pictorially, by
$\varsigma(\overcrossing)=\varsigma(\virtualcrossing)=\virtualcrossing$).
Thus the symmetric group $S_n$ is both a subgroup and a quotient group
of $\vB_n$.
\end{remark}

Like there are pure braids to accompany braids, there are pure virtual
braids as well:

\begin{definition} A pure v-braid is a v-braid whose skeleton is the
identity permutation; the group $\glos{\PvB_n}$ of all pure v-braids is
simply the kernel of the skeleton morphism $\varsigma\colon\vB_n\to S_n$.
\end{definition}

We note the sequence of group homomorphisms
\begin{equation} \label{eq:ExcatSeqForPvB}
  1\longrightarrow\PvB_n\xhookrightarrow{\quad}\vB_n
  \overset{\varsigma}{\longrightarrow}S_n
  \longrightarrow 1.
\end{equation}
This sequence is exact and split, with the splitting given by the inclusion
$S_n\hookrightarrow\vB_n$ mentioned above~\eqref{eq:sRelations}. Therefore
we have that
\begin{equation} \label{eq:vBSemiDirect}
  \vB_n=\PvB_n\rtimes S_n.
\end{equation}

\subsubsection{The ``Abstract'' Way} \label{subsubsec:Abstract}
The relations~\eqref{eq:sRels} and~\eqref{eq:MixedRels} that govern
the behaviour of virtual crossings precisely say that virtual crossings
really are ``virtual'' --- if a piece of strand is routed within a braid
so that there are only virtual crossings around it, it can be rerouted
in any other ``virtual only'' way, provided the ends remain fixed
(this is Kauffman's ``detour move''~\cite{Kauffman:VirtualKnotTheory,
KauffmanLambropoulou:VirtualBraids}). Since a v-braid $B$ is independent
of the routing of virtual pieces of strand, we may as well never supply
this routing information.

\parpic[r]{$\input figs/PvBExample.pstex_t$}
Thus for example, a perfectly fair verbal description of the
(pure!) v-braid on the right is ``strand 1 goes over strand 3 by a
positive crossing then likewise positively over strand 2 then negatively
over 3 then 2 goes positively over 1''. We don't need to specify how
strand 1 got to be near strand 3 so it can go over it --- it got there
by means of virtual crossings, and it doesn't matter how. Hence we arrive
at the following ``abstract'' presentation of $\PvB_n$ and $\vB_n$:

\begin{proposition} (E.g.~\cite{Bardakov:VirtualAndUniversal})
\begin{enumerate}
\item The \wClipStart{120118-2}{0-25-21} group $\PvB_n$ of pure v-braids
is isomorphic to the group generated by symbols $\glos{\sigma_{ij}}$ for
$1\leq i\neq j\leq n$ (meaning ``strand $i$ crosses over strand $j$
at a positive crossing''\footnote{The inverse, $\sigma_{ij}^{-1}$,
is ``strand $i$ crosses over strand $j$ at a negative crossing''}),
subject to the third Reidemeister move and to locality in space (compare
with~\eqref{eq:R3} and~\eqref{eq:sigmaRels}):

\begin{align*}
  \sigma_{ij}\sigma_{ik}\sigma_{jk} &= \sigma_{jk}\sigma_{ik}\sigma_{ij}
    & \text{whenever}\qquad & |\{i,j,k\}|=3, \\
  \sigma_{ij}\sigma_{kl} &= \sigma_{kl}\sigma_{ij}
     & \text{whenever}\qquad & |\{i,j,k,l\}|=4.
\end{align*}
\item If $\tau\in S_n$, then with the action
$\sigma_{ij}^\tau:=\sigma_{\tau i,\tau j}$ we recover the semi-direct
product decomposition $\vB_n=\PvB_n\rtimes S_n$. \qed \wClipEnd{120111-2}
\end{enumerate}
\end{proposition}

\draftcut \subsection{On to w-Braids} \label{subsec:wBraids}
To define w-braids, we break the symmetry between over crossings and
under crossings by imposing one of the ``forbidden moves'' virtual knot
theory, but not the other:
\begin{equation} \label{eq:OvercrossingsCommute}
  \sigma_i\sigma_{i+1}s_i = s_{i+1}\sigma_i\sigma_{i+1},
  \qquad\text{yet}\qquad
  s_i\sigma_{i+1}\sigma_i \neq \sigma_{i+1}\sigma_is_{i+1}.
\end{equation}
Alternatively,
\[ \sigma_{ij}\sigma_{ik} = \sigma_{ik}\sigma_{ij},
  \qquad\text{yet}\qquad
  \sigma_{ik}\sigma_{jk} \neq \sigma_{jk}\sigma_{ik}.
\]
In pictures, this is
\begin{equation} \label{eq:OC} \input figs/OCUC.pstex_t \end{equation}

The relation we have just imposed may be called the ``unforbidden
relation'', or, perhaps more appropriately, the ``overcrossings commute''
relation (\glost{OC}). Ignoring the non-crossings\footnote{Why this is
appropriate was explained in the previous section.} $\virtualcrossing$,
the OC relation says that it is the same if strand $i$ first crosses
over strand $j$ and then over strand $k$, or if it first crosses over
strand $k$ and then over strand $j$. The ``undercrossings commute''
relation \glost{UC}, the one we do not impose
in~\eqref{eq:OvercrossingsCommute}, would say the same except with
``under'' replacing ``over''.

\begin{definition} The group of w-braids is $\glos{\wB_n}:=\vB_n/OC$. Note
that $\varsigma$ descends to $\wB_n$ and hence we can define the group of
pure w-braids to be $\glos{\PwB_n}:=\ker\varsigma\colon\wB_n\to S_n$. We
still have a split exact sequence as at~\eqref{eq:ExcatSeqForPvB} and
a semi-direct product decomposition $\wB_n=\PwB_n\rtimes S_n$.
\end{definition}

\begin{exercise} Show that the OC relation is equivalent to the relation
\[
  \sigma_i^{-1}s_{i+1}\sigma_i = \sigma_{i+1}s_i\sigma_{i+1}^{-1}
  \qquad\text{or}\qquad
  \parbox[m]{1.5in}{\begin{picture}(0,0)%
\includegraphics{figs/AltOC.pstex}%
\end{picture}%
%
%
\setlength{\unitlength}{3158sp}%
\begingroup\makeatletter\ifx\SetFigFont\undefined%
\gdef\SetFigFont#1#2#3#4#5{%
  \reset@font\fontsize{#1}{#2pt}%
  \fontfamily{#3}\fontseries{#4}\fontshape{#5}%
  \selectfont}%
\fi\endgroup%
\begin{picture}(1524,924)(1639,-673)
\put(2401,-286){\makebox(0,0)[b]{\smash{{\SetFigFont{10}{12.0}{\rmdefault}{\mddefault}{\updefault}{\color[rgb]{0,0,0}$=$}%
}}}}
\end{picture}%
 }
\]
\end{exercise}

While mostly in this paper the pictorial / algebraic definition of w-braids
(and other w-knotted objects) will suffice, we ought describe at least
briefly 2-3 further interpretations of $\wB_n$:

\subsubsection{The group of flying rings} \label{subsubsec:FlyingRings}
Let \glos{X_n} be the space of all placements of $n$ numbered disjoint
geometric circles in $\bbR^3$, such that all circles are parallel to
the $xy$ plane. Such placements will be called horizontal\footnote{
For the group of non-horizontal flying rings see
Section \ref{subsubsec:NonHorRings}}. A horizontal
placement is determined by the centres in $\bbR^3$ of the $n$ circles
and by $n$ radii, so $\dim X_n=3n+n=4n$. The permutation group $S_n$
acts on $X_n$ by permuting the circles, and one may think of the quotient
$\glos{\tilde{X}_n}:=X_n/S_n$ as the space of all horizontal
placements of $n$ unmarked circles in $\bbR^3$. The fundamental group
$\pi_1(\tilde{X}_n)$ is a group of paths traced by $n$ disjoint horizontal
circles (modulo homotopy), so it is fair to think of it as ``the group
of flying rings''.

\begin{theorem} The group of pure w-braids $\PwB_n$ is isomorphic to the group
of flying rings $\pi_1(X_n)$. The group $\wB_n$ is isomorphic to the group
of unmarked flying rings $\pi_1(\tilde{X}_n)$.
\end{theorem}

For the proof of this theorem, see~\cite{Goldsmith:MotionGroups,
Satoh:RibbonTorusKnots} and
especially~\cite{BrendleHatcher:RingsAndWickets}. Here we will contend
ourselves with pictures describing the images of the generators of $\wB_n$
in $\pi_1(\tilde{X}_n)$ and a few comments:

\[ \input figs/FlyingRings.pstex_t \]

Thus we map the permutation $s_i$ to the movie clip in which ring
number $i$ trades its place with ring number $i+1$ by having the two
flying around each other. This acrobatic feat is performed in $\bbR^3$
and it does not matter if ring number $i$ goes ``above'' or ``below''
or ``left'' or ``right'' of ring number $i+1$ when they trade places,
as all of these possibilities are homotopic. More interestingly, we
map the braiding $\sigma_i$ to the movie clip in which ring $i+1$
shrinks a bit and flies through ring $i$. It is a worthwhile
exercise for the reader to verify that the relations in the definition
of $\wB_n$ become homotopies of movie clips. Of these relations it
is most interesting to see why the ``overcrossings commute'' relation
$\sigma_i\sigma_{i+1}s_i = s_{i+1}\sigma_i\sigma_{i+1}$ holds, yet the
``undercrossings commute'' relation $\sigma^{-1}_i\sigma^{-1}_{i+1}s_i =
s_{i+1}\sigma^{-1}_i\sigma^{-1}_{i+1}$ doesn't. \wClipEnd{120118-2} 

\begin{exercise}\label{ex:swBn}
To be perfectly precise, we have to
specify the fly-through direction. In our notation, $\sigma_i$ means
that the ring corresponding to the under-strand approaches the bigger
ring representing the over-strand from below, flies through it and exists
above.  For $\sigma_i^{-1}$ we are ``playing the movie backwards'', i.e.,
the ring of the under-strand comes from above and exits below the ring
of the over-strand.

Let ``the signed $w$ braid group'', $\swB_n$, be the group
of horizontal flying rings where both fly-through 
directions are allowed. This introduces a ``sign'' for
each crossing $\sigma_i$:
\begin{center}
 \input figs/FlyingRings2.pstex_t
\end{center}
In other words, $\swB_n$ is generated by $s_i$, 
$\sigma_{i+}$ and $\sigma_{i-}$, for $i=1,...,n$. Check that in $\swB_n$
$\sigma_{i-}=s_i\sigma_{i+}^{-1}s_i$, and this, along with
the other obvious relations implies $\swB_n \cong \wB_n$.

For a rigorous discussion of orientations and signs, see Section \ref{subsec:TangleTopology}.
\end{exercise}

\subsubsection{Certain ribbon tubes in $\bbR^4$} \label{subsubsec:ribbon}
With
\wClipStart{120118-3}{0-00-00}
time as the added dimension, a flying ring in $\bbR^3$ traces a tube
(an annulus) in $\bbR^4$, as shown in the picture below:
\[ \input figs/RibbonTubes.pstex_t \]
Note that we adopt here the drawing conventions of Carter and
Saito~\cite{CarterSaito:KnottedSurfaces} --- we draw surfaces as if they
were projected from $\bbR^4$ to $\bbR^3$, and we cut them open whenever
they are ``hidden'' by something with a higher fourth coordinate.

Note
\wClipComment{120118-3}{0-09-40}{continues Sec.~\ref{subsubsec:NonHorRings}}
also that the tubes we get in $\bbR^4$ always bound natural 3D
``solids'' --- their ``insides'', in the pictures above. These solids
are disjoint in the case of $s_i$ and have a very specific kind
of intersection in the case of $\sigma_i$ --- these are transverse
intersections with no triple points, and their inverse images are a
meridional disk on the ``thin'' solid tube and an interior disk on the
``thick'' one. By analogy with the case of ribbon knots and ribbon
singularities in $\bbR^3$ (e.g.~\cite[Chapter V]{Kauffman:OnKnots}) and
following Satoh~\cite{Satoh:RibbonTorusKnots}, we call this kind if
intersections of solids in $\bbR^4$ ``ribbon singularities'' and thus our
tubes in $\bbR^4$ are always ``ribbon tubes''.

\subsubsection{Basis conjugating automorphisms of $F_n$}
\label{subsubsec:McCool}
Let
\wClipStart{120125-1}{0-00-00}
$\glos{F_n}$ be the free (non-Abelian) group with generators
$\glos{\xi_1,\ldots,\xi_n}$. Artin's theorem (Theorems 15 and 16
of~\cite{Artin:TheoryOfBraids}) says that the (usual) braid
group $\uB_n$ (equivalently, the subgroup of $\wB_n$ generated by
the $\sigma_i$'s) has a faithful right action on $F_n$. In other
words, $\uB_n$ is isomorphic to a subgroup $H$ of $\Autop(F_n)$
(the group of automorphisms of $F_n$ with opposite multiplication;
$\psi_1\psi_2:=\psi_2\circ\psi_1$). Precisely, using $(\xi,
B)\mapsto\xi\glos{\sslash}B$ to denote the right action of $\Autop(F_n)$ on
$F_n$, the subgroup $H$ consists of those automorphisms $B\colon F_n\to F_n$
of $F_n$ that satisfy the following two conditions:
\begin{enumerate}
\item $B$ maps any generator $\xi_i$ to a
  conjugate of a generator (possibly different). That is, there is a
  permutation $\beta\in S_n$ and elements $a_i\in F_n$ so that for every $i$,
  \begin{equation} \label{eq:BasisConjugating}
    \xi_i \sslash B = a_i^{-1}\xi_{\beta i}a_i.
  \end{equation}
\item $B$ fixes the ordered product of the generators of $F_n$,
  \[ \xi_1\xi_2\cdots \xi_n \sslash B = \xi_1\xi_2\cdots \xi_n. \]
\end{enumerate}

McCool's theorem~\cite{McCool:BasisConjugating} says that the same
holds true\footnote{Though see Warning~\ref{warn:NoArtin}.} if one
replaces the braid group $\uB_n$ with the bigger group $\wB_n$ and
drops the second condition above. So $\wB_n$ is precisely the group
of ``basis-conjugating'' automorphisms of the free group $F_n$,
the group of those automorphisms which map any ``basis element''
in $\{\xi_1,\ldots,\xi_n\}$ to a conjugate of a (possibly different)
basis element.

The relevant action is explicitly defined on the generators of $\wB_n$
and $F_n$ as follows (with the omitted generators of $F_n$ always fixed):
\begin{equation} \label{eq:ExplicitPsi}
  (\xi_i, \xi_{i+1})\sslash s_i = (\xi_{i+1}, \xi_i)
  \qquad
  (\xi_i, \xi_{i+1})\sslash \sigma_i = (\xi_{i+1},
    \xi_{i+1}\xi_i\xi_{i+1}^{-1})
  \qquad
  \xi_j\sslash \sigma_{ij} = \xi_i\xi_j\xi_i^{-1}
\end{equation} 
It is a worthwhile exercise to verify that $\sslash$ respects the
relations in the definition of $\wB_n$ and that the permutation $\beta$
in~\eqref{eq:BasisConjugating} is the skeleton $\varsigma(B)$.

There is a more conceptual description of $\sslash$, in terms of the
structure of $\wB_{n+1}$. Consider the inclusions
\begin{equation} \label{eq:inclusions}
  \wB_n \xhookrightarrow{\iota} \wB_{n+1} \xhookleftarrow{i_u} F_n.
\end{equation}

\parpic[r]{$\input figs/xi.pstex_t$}
Here $\glos{\iota}$ is the inclusion of $\wB_n$ into $\wB_{n+1}$ by adding an
inert $(n+1)-$st strand (it is injective as it has a well defined
one sided inverse --- the deletion of the $(n+1)$-st strand). The
inclusion $\glos{i_u}$ of the free group $F_n$ into $\wB_{n+1}$ is defined by
$i_u(\xi_i):=\sigma_{i,n+1}$.  The image $i_u(F_n)\subset\wB_{n+1}$ is the
set of all w-braids whose first $n$ strands are straight and vertical,
and whose $(n+1)$-st strand wanders among the first $n$ strands mostly
virtually (i.e., mostly using virtual crossings), occasionally slipping
under one of those $n$ strands, but never going over anything.  In the
``flying rings'' picture of Section~\ref{subsubsec:FlyingRings}, the image
$i_u(F_n)\subset\wB_{n+1}$ can be interpreted as the fundamental group
of the complement in $\bbR^3$ of $n$ stationary rings (which is indeed
$F_n$) --- in $i_u(F_n)$ the only ring in motion is the last, and it only
goes under, or ``through'', other rings, so it can be replaced by a point
object whose path is an element of the fundamental group. The injectivity
of $i_u$ follows from this geometric picture.

\parpic[r]{$\input figs/Bgamma.pstex_t$} \picskip{4}
One may explicitly verify that $i_u(F_n)$ is normalized by $\iota(\wB_n)$
in $\wB_{n+1}$ (that is, the set $i_u(F_n)$ is preserved by conjugation
by elements of $\iota(\wB_n)$). Thus the following definition (also shown
as a picture on the right) makes sense, for $B\in\wB_n\subset\wB_{n+1}$
and for $\gamma\in F_n\subset\wB_{n+1}$:

\begin{equation} \label{eq:ConceptualPsi}
  \gamma\sslash B := i_u^{-1}(B^{-1}\gamma B)
\end{equation}

It is a worthwhile exercise to recover the explicit formulae
in~\eqref{eq:ExplicitPsi} from the above definition.

\begin{warning} \label{warn:NoArtin} People familiar with the Artin story for
ordinary braids should be warned that even though $\wB_n$ acts on $F_n$ and
the action is induced from the inclusions in~\eqref{eq:inclusions} in much
of the same way as the Artin action is induced by inclusions $\uB_n
\xhookrightarrow{\iota} \uB_{n+1} \xhookleftarrow{i} F_n$, there are also
some differences, and some further warnings apply:
\begin{myitemize}
\item In the ordinary Artin story, $i(F_n)$ is the set of braids in
$\uB_{n+1}$ whose first $n$ strands are unbraided (that is, whose image in
$\uB_n$ via ``dropping the last strand'' is the identity). This is not true
for w-braids. For w-braids, in $i_u(F_n)$ the last strand always goes
``under'' all other strands (or just virtually crosses them), but never
over.
\item Thus unlike the isomorphism $\PuB_{n+1}\cong \PuB_n\ltimes F_n$,
it is not true that $\PwB_{n+1}$ is isomorphic to $\PwB_n\ltimes F_n$.
\item The Overcrossings Commute relation imposed in $\wB$ breaks
the symmetry between overcrossings and undercrossings. Thus let
$i_o\colon F_n\to\wB_n$ be the ``opposite'' of $i_u$, mapping into braids in
which the last strand is always ``over'' or virtual. Then $i_o$ is not
injective (its image is in fact Abelian) and its image is not normalized
by $\iota(\wB_n)$. So there is no ``second'' action of $\wB_n$ on $F_n$
defined using $i_o$.
\item For v-braids, both $i_u$ and $i_o$ are injective and there are two
actions of $\vB_n$ on $F_n$ --- one defined by first projecting into
w-braids, and the other defined by first projecting into v-braids modulo
``Undercrossings Commute''. Yet v-braids contain more information than
these two actions can see. The ``Kishino'' v-braid below, for example,
is visibly trivial if either overcrossings or undercrossings are made to
commute, yet by computing its Kauffman bracket we know it is non-trivial
as a v-braid~\cite[``The Kishino Braid'']{WKO}:
\[
  \begin{picture}(0,0)%
\includegraphics{figs/KishinoBraid.pstex}%
\end{picture}%
%
%
\setlength{\unitlength}{2763sp}%
\begingroup\makeatletter\ifx\SetFigFont\undefined%
\gdef\SetFigFont#1#2#3#4#5{%
  \reset@font\fontsize{#1}{#2pt}%
  \fontfamily{#3}\fontseries{#4}\fontshape{#5}%
  \selectfont}%
\fi\endgroup%
\begin{picture}(7074,1027)(-11,-176)
\put(976,-136){\makebox(0,0)[b]{\smash{{\SetFigFont{7}{8.4}{\rmdefault}{\mddefault}{\updefault}{\color[rgb]{0,0,0}$a$}%
}}}}
\put(6076,-136){\makebox(0,0)[b]{\smash{{\SetFigFont{7}{8.4}{\rmdefault}{\mddefault}{\updefault}{\color[rgb]{0,0,0}$b$}%
}}}}
\end{picture}%
 \quad \left(\parbox{1.6in}{\footnotesize
    The commutator $ab^{-1}a^{-1}b$ of v-braids $a,b$ annihilated by OC/UC,
    respectively, with a minor cancellation.
  }\right)
\]
\wClipEnd{120125-1}
\end{myitemize}
\end{warning}

\begin{problem} \label{prob:wCombing}
Is $\PwB_n$ a semi-direct product of free groups? Note that both
$\PuB_n$ and $\PvB_n$ are such semi-direct products: For $\PuB_n$, this
is the well known ``combing of braids''; it follows from $\PuB_n\cong
\PuB_{n-1}\ltimes F_{n-1}$ and induction. For $\PvB_n$, it is a result
stated in~\cite{Bardakov:VirtualAndUniversal} (though our own understanding
of~\cite{Bardakov:VirtualAndUniversal} is incomplete).
\end{problem}

\begin{remark} \label{rem:GutierrezKrstic}
Note that Guti\'errez and Krsti\'c~\cite{GutierrezKrstic:NormalForms}
find ``normal forms'' for the elements of $\PwB_n$, yet they
do not decide whether $\PwB_n$ is ``automatic'' in the sense
of~\cite{Epstein:WordProcessing}.
\end{remark}

\draftcut \subsection{Finite Type Invariants of v-Braids and w-Braids}
\label{subsec:FT4Braids}

Just as we had two definitions for v-braids (and thus w-braids)
in Section~\ref{subsec:VirtualBraids}, we will give two (obviously
equivalent) developments of the theory of finite type invariants of
v-braids and w-braids --- a pictorial/topological version in
Section~\ref{subsubsec:FTPictorial}, and a more abstract algebraic version
in Section~\ref{subsubsec:FTAlgebraic}.

\subsubsection{Finite Type Invariants, the Pictorial Approach}
\label{subsubsec:FTPictorial}

In
\wClipComment{120125-2}{0-18-54}{describes the standard theory, briefly}
the standard theory of finite type invariants of knots (also known as
Vassiliev or Goussarov-Vassiliev invariants)~\cite{Goussarov:nEquivalence,
Vassiliev:CohKnot, Bar-Natan:OnVassiliev, Bar-Natan:EMP}
one progresses from the definition of finite type via iterated
differences to chord diagrams and weight systems, to $4T$ (and other)
relations, to the definition of universal finite type invariants,
and beyond. The exact same progression (with different objects
playing similar roles, and sometimes, when yet insufficiently
studied, with the last step or two missing) is also seen in the theories
of finite type invariants of braids~\cite{Bar-Natan:Braids},
3-manifolds~\cite{Ohtsuki:IntegralHomology, LeMurakamiOhtsuki:Universal,
Le:UniversalIHS}, virtual knots~\cite{GoussarovPolyakViro:VirtualKnots,
Polyak:ArrowDiagrams} and of several other classes of objects. We thus
assume that the reader has familiarity with these basic ideas, and we
only indicate briefly how they are implemented in the case of v-braids
and w-braids.

\begin{figure}
\[ \input figs/Dvh1.pstex_t \]
  \caption{
    On the left, a 3-singular v-braid and its corresponding 3-arrow
    diagram. A self-explanatory algebraic notation for this arrow
    diagram is $(\glos{a_{12}a_{41}a_{23}},\,3421)$.  picture and in algebraic
    notation. Note that we regard arrow diagrams as graph-theoretic
    objects, and hence the two arrow diagrams on the right, whose
    underlying graphs are the same, are regarded as equal. In
    algebraic notation this means that we always impose the relation
    $a_{ij}a_{kl}=a_{kl}a_{ij}$ when the indices $i$, $j$, $k$, and $l$
    are all distinct.
  } \label{fig:Dvh1}
\end{figure}

Much
\wClipStart{120208}{0-04-23}
like the formula $\doublepoint\to\overcrossing-\undercrossing$ of
the Vassiliev-Goussarov fame, given a v-braid invariant $\glos{V}\colon
\vB_n\to A$ valued in some Abelian group $A$, we extend it to
``singular'' v-braids, braids that contain ``semi-virtual crossings''
like $\glos{\semivirtualover}$ and $\glos{\semivirtualunder}$ using the
formulae $V(\semivirtualover):=V(\overcrossing)-V(\virtualcrossing)$
and $V(\semivirtualunder):=V(\undercrossing)-V(\virtualcrossing)$
(see~\cite{GoussarovPolyakViro:VirtualKnots, Polyak:ArrowDiagrams,
Bar-NatanHalachevaLeungRoukema:v-Dims}). We
say that ``$V$ is of type $m$'' if its extension vanishes on singular
v-braids having more than $m$ semi-virtual crossings. Up to invariants
of lower type, an invariant of type $m$ is determined by its ``weight
system'', which is a functional $W=\glos{W_m}(V)$ defined on ``$m$-singular
v-braids modulo $\overcrossing=\virtualcrossing=\undercrossing$''. Let
us denote the vector space of all formal linear combinations of such
equivalence classes by $\glos{\calG_m}\calD^v_n$. Much as $m$-singular knots
modulo $\overcrossing=\undercrossing$ can be identified with chord
diagrams, the basis elements of $\calG_m\calD^v_n$ can be identified with
pairs $(D,\beta)$, where $D$ is a horizontal arrow diagram and $\beta$
is a ``skeleton permutation''. See Figure~\ref{fig:Dvh1}.

We assemble the spaces $\calG_m\calD^v_n$ together to form a single
graded space, $\glos{\calD^v_n}:=\oplus_{m=0}^\infty\calG_m\calD^v_n$. Note that
throughout this paper, whenever we write an infinite direct sum,
we automatically complete it. Thus in $\calD^v_n$ we allow infinite sums
with one term in each homogeneous piece $\calG_m\calD^v_n$.

\begin{figure}
\[ \input figs/6T.pstex_t \]
\[
  a_{ij}a_{ik}+a_{ij}a_{jk}+a_{ik}a_{jk}
  = a_{ik}a_{ij}+a_{jk}a_{ij}+a_{jk}a_{ik}
\]
\[
  \text{or}\qquad
  [a_{ij}, a_{ik}] + [a_{ij}, a_{jk}] + [a_{ik}, a_{jk}] = 0
\]
\caption{The $6T$ relation. Standard knot theoretic conventions apply ---
only the relevant parts of each diagram is shown; in reality each diagram
may have further vertical strands and horizontal arrows, provided the
extras are the same in all 6 diagrams. Also, the vertical strands are in no
particular order --- other valid $6T$ relations are obtained when those
strands are permuted in other ways.} \label{fig:6T}
\end{figure}

\begin{figure}
\[ \begin{array}{ccc}
  \input figs/TC.pstex_t & \qquad & \input figs/4TArrow.pstex_t \\
  a_{ij}a_{ik} = a_{ik}a_{ij} &&
  a_{ij}a_{jk} + a_{ik}a_{jk} = a_{jk}a_{ij} + a_{jk}a_{ik} \\
  \text{or} \quad [a_{ij}, a_{ik}] = 0 &&
  \text{or} \quad [a_{ij} + a_{ik}, a_{jk}] = 0
\end{array} \]
\caption{The TC and the $\protect\aft$ relations.} \label{fig:TCand4T}
\end{figure}

In the standard finite-type theory for knots, weight
systems always satisfy the $4T$ relation, and are therefore
functionals on $\calA:=\calD/4T$. Likewise, in the case of
v-braids, weight systems satisfy the ``$\glos{6T}$ relation''
of~\cite{GoussarovPolyakViro:VirtualKnots, Polyak:ArrowDiagrams,
Bar-NatanHalachevaLeungRoukema:v-Dims},
shown in Figure~\ref{fig:6T}, and are therefore functionals on
$\glos{\calA^v_n}:=\calD^v_n/6T$. In the case of w-braids, the ``overcrossings
commute'' relation~\eqref{eq:OvercrossingsCommute} implies the
``Tails Commute'' (\glost{TC}) relation on the level of arrow diagrams,
and in the presence of the TC relation, two of the terms in the $6T$
relation drop out, and what remains is the ``$\glos{\aft}$'' relation. These
relations are shown in Figure~\ref{fig:TCand4T}. Thus weight systems
of finite type invariants of w-braids are linear functionals on
$\glos{\calA^w_n}:=\calD^v_n/TC,\aft$.

The next question that arises is whether we have already found {\em
all} the relations that weight systems always satisfy. More precisely,
given a degree $m$ linear functional on $\calA^v_n=\calD^v_n/6T$ (or
on $\calA^w_n=\calD^v_n/TC,\aft$), is it always the weight system
of some type $m$ invariant $V$ of v-braids (or w-braids)? As in every
other theory of finite type invariants, the answer to this question
is affirmative if and only if there exists a ``universal finite type
invariant'' (or simply, an ``expansion'') of v-braids (w-braids):

\begin{definition} \label{def:vwbraidexpansion} An expansion
for v-braids (w-braids) is an invariant $Z\colon \vB_n\to\calA^v_n$
(or $Z\colon \wB_n\to\calA^w_n$) satisfying the following ``universality
condition'':
\begin{itemize}
\item If $B$ is an $m$-singular v-braid (w-braid) and
$D\in\calG_m\calD^v_n$ is its underlying arrow diagram as in
Figure~\ref{fig:Dvh1}, then
\[ Z(B)=D+(\text{terms of degree\,}>m). \]
\end{itemize}
\end{definition}

Indeed if $Z$ is an expansion and $W\in\calG_m\calA^\star$,\footnote{$\calA$
here denotes either $\calA^v_n$ or $\calA^w_n$, or in fact, any of many
similar spaces that we will discuss later on.} the universality condition
implies that $W\circ Z$ is a finite type invariant whose weight system
is $W$. To go the other way, if $(D_i)$ is a basis of $\calA$ consisting
of homogeneous elements, if $(W_i)$ is the dual basis of $\calA^\star$ and
$(V_i)$ are finite type invariants whose weight systems are the $W_i$'s,
then $Z(B):=\sum_iD_iV_i(B)$ defines an expansion.

In general, constructing a universal finite type invariant is a
hard task. For knots, one uses either the Kontsevich integral or
perturbative Chern-Simons theory (also known as ``configuration
space integrals''~\cite{BottTaubes:SelfLinking} or ``tinker-toy
towers''~\cite{Thurston:IntegralExpressions}) or the rather fancy
algebraic theory of ``Drinfel'd associators'' (a summary of all those
approaches is at~\cite{Bar-NatanStoimenow:Fundamental}). For homology
spheres, this is the ``LMO invariant''~\cite{LeMurakamiOhtsuki:Universal,
Le:UniversalIHS} (also the ``\AA{}rhus
integral''~\cite{Bar-NatanGaroufalidisRozanskyThurston:Aarhus}). For
v-braids, we still don't know if an expansion exists. As we shall see
below, the construction of an expansion for w-braids is quite easy.

\subsubsection{Finite Type Invariants, the Algebraic Approach}
\label{subsubsec:FTAlgebraic}
For
\wClipStart{120201}{0-04-37}
any group $G$, one can form the group algebra ${\mathbb F}G$ for some
field $\mathbb F$ by allowing formal linear combinations of group elements
and extending multiplication linearly. The {\it augmentation ideal} is
the ideal generated by differences, or equivalently, the set of linear
combinations of group elements whose coefficients sum to zero:
\[ \glos{\calI} := \left\{\sum_{i=1}^k a_ig_i\colon
  a_i \in {\mathbb F}, g_i \in G, \sum_{i=1}^k a_i=0\right\}.
\]
Powers of the augmentation ideal provide a filtration of the group
algebra. Let $\glos{\calA(G)}:= \bigoplus_{m\geq 0} \calI^m/\calI^{m+1}$
be the associated graded space corresponding to this filtration.

\begin{definition}\label{def:grpexpansion} An
\wClipComment{120201}{0-52-41}{is much more detailed on these matters}
expansion for the group
$G$ is a map $Z\colon G \to \calA(G)$, such that the linear extension
$Z\colon  {\mathbb F}G \to \calA(G)$ is filtration preserving and
the induced map $$\gr Z\colon  (\gr {\mathbb F}G=\calA(G)) \to (\gr
\calA(G)=\calA(G))$$ is the identity. An equivalent way to phrase this
is that the degree $m$ piece of $Z$ restricted to $\calI^m$ is the
projection onto $\calI^m/\calI^{m+1}$.

\begin{exercise}\label{ex:BraidsAlgApproach}
Verify that for the groups $\PvB_n$ and $\PwB_n$ the m-th power of the
augmentation ideal coincides with the span of all resolutions of
$m$-singular $v$- or $w$-braids (by a resolution we mean the formal
linear combination where each semivirtual crossing is replaced by
the appropriate difference of a virtual and a regular crossing). Then
check that the notion of expansion defined above is the same as that of
Definition \ref{def:vwbraidexpansion}, restricted to pure braids.
\end{exercise}

Finally, note the functorial nature of the construction above. What we have 
described is a functor, called ``projectivization''
$\proj\colon  Groups \to GradedAlgebras$, which assigns to each group
$G$ the graded algebra $\calA(G)$. To each homomorphism $\phi\colon  G \to H$,
$\proj$ assigns
the induced map $\gr \phi\colon  (\calA(G)=\gr {\mathbb F}G) \to (\calA(H)= \gr {\mathbb F}H)$.
 
\end{definition}

\draftcut \subsection{Expansions for w-Braids}\label{subsec:wBraidExpansion}

The space $\calA^w_n$ of arrow diagrams on $n$ strands is an associative
algebra in an obvious manner: If the permutations underlying two arrow
diagrams are the identity permutations, we simply juxtapose the diagrams.
Otherwise we ``slide'' arrows through permutations in the obvious manner
--- if $\tau$ is a permutation, we declare that $\tau a_{(\tau i)(\tau
j)}=a_{ij}\tau$. Instead of seeking an expansion $\wB_n\to\calA^w_n$, we
set the bar a little higher and seek a ``homomorphic expansion'':

\begin{definition} \label{def:Universallity} A homomorphic expansion
$Z\colon \wB_n\to\calA^w_n$ is an expansion that carries products in $\wB_n$
to products in $\calA^w_n$.
\end{definition}

To find a homomorphic expansion, we just need to define it
on the generators of $\wB_n$ and verify that it satisfies the
relations defining $\wB_n$ and the universality condition.
Following~\cite[Section~5.3]{BerceanuPapadima:BraidPermutation}
and~\cite[Section~8.1]{AlekseevTorossian:KashiwaraVergne} we set
$Z(\virtualcrossing)=\virtualcrossing$ (that is, a transposition in
$\wB_n$ gets mapped to the same transposition in $\calA^w_n$, adding no
arrows) and $Z(\overcrossing)=\exp(\rightarrowdiagram)\virtualcrossing$.
This last formula is important so deserves to be magnified, explained
and replaced by some new notation:

\begin{equation} \label{eq:reservoir}
  Z\left(\!\mathsize{\Huge}{\overcrossing}\!\right)\! =
  \exp\left(\!\mathsize{\Huge}{\rightarrowdiagram}\!\right)
    \cdot\mathsize{\Huge}{\virtualcrossing}
  = \input figs/ZIsExp.pstex_t+\ldots =: \input figs/ArrowReservoir.pstex_t.
\end{equation}

Thus the new notation $\overset{e^a}{\longrightarrow}$ stands
for an ``exponential reservoir'' of parallel arrows, much like
$e^a=1+a+aa/2+aaa/3!+\ldots$ is a ``reservoir'' of $a$'s. With
the obvious interpretation for $\overset{e^{-a}}{\longrightarrow}$
(the $-$ sign indicates that the terms should have alternating signs,
as in $e^{-a}=1-a+a^2/2-a^3/3!+\ldots$), the second Reidemeister move
$\overcrossing\undercrossing=1$ forces that we set
\[ Z\left(\mathsize{\Huge}{\undercrossing}\right) =
  \mathsize{\Huge}{\virtualcrossing}
  \cdot\exp\left(-\mathsize{\Huge}{\rightarrowdiagram}\right)
  = \begin{picture}(0,0)%
\includegraphics{figs/NegReservoir1.pstex}%
\end{picture}%
%
%
\setlength{\unitlength}{2960sp}%
\begingroup\makeatletter\ifx\SetFigFont\undefined%
\gdef\SetFigFont#1#2#3#4#5{%
  \reset@font\fontsize{#1}{#2pt}%
  \fontfamily{#3}\fontseries{#4}\fontshape{#5}%
  \selectfont}%
\fi\endgroup%
\begin{picture}(739,1374)(1281,-523)
\put(1651,239){\makebox(0,0)[b]{\smash{{\SetFigFont{9}{10.8}{\rmdefault}{\mddefault}{\updefault}{\color[rgb]{0,0,0}$e^{-a}$}%
}}}}
\end{picture}%
 = \input figs/NegReservoir2.pstex_t.
\]
\wClipEnd{120201}

\begin{theorem} \label{thm:RInvariance} The above formulae define
an invariant $Z\colon \wB_n\to\calA^w_n$ (that is, $Z$ satisfies all the
defining relations of $\wB_n$). The resulting $Z$ is a homomorphic
expansion (that is, it satisfies the universality property of
Definition~\ref{def:Universallity}).
\end{theorem}

\begin{proof} (Following~\cite{BerceanuPapadima:BraidPermutation,
AlekseevTorossian:KashiwaraVergne}) For the invariance of $Z$, the
only interesting relations to check are the Reidemeister 3 relation
of~\eqref{eq:sigmaRels} and the Overcrossings Commute relation
of~\eqref{eq:OC}. For Reidemeister 3, we have
\[ \input figs/R3Left.pstex_t
  = e^{a_{12}}e^{a_{13}}e^{a_{23}}\tau
  \overset{1}{=} e^{a_{12}+a_{13}}e^{a_{23}}\tau
  \overset{2}{=} e^{a_{12}+a_{13}+a_{23}}\tau,
\]
where $\tau$ is the permutation $321$ and equality 1 holds because
$[a_{12},a_{13}]=0$ by a Tails Commute (TC) relation and equality 2 holds
because $[a_{12}+a_{13}, a_{23}]=0$ by a $\aft$ relation.
Likewise, again using TC and $\aft$,
\[ \input figs/R3Right.pstex_t
  = e^{a_{23}}e^{a_{13}}e^{a_{12}}\tau
  = e^{a_{23}}e^{a_{13}+a_{12}}\tau
  = e^{a_{23}+a_{13}+a_{12}}\tau,
\]
and so Reidemeister 3 holds. An even simpler proof using just the Tails
Commute relation shows that the Overcrossings Commute relation also holds.
Finally, since $Z$ is homomorphic, it is enough to check the universality
property at degree $1$, where it is very easy:
\[ Z\left(\mathsize{\Huge}{\semivirtualover}\right) =
  \exp\left(\mathsize{\Huge}{\rightarrowdiagram}\right)
    \cdot\mathsize{\Huge}{\virtualcrossing}
    - \mathsize{\Huge}{\virtualcrossing}
  = \mathsize{\Huge}{\rightarrowdiagram}\cdot\mathsize{\Huge}{\virtualcrossing}
    + (\text{terms of degree\,}>1),
\]
and a similar computation manages the $\semivirtualunder$ case. \qed
\end{proof}

\begin{remark} \label{rem:YangBaxter} Note that the main ingredient
of the above proof was to show that $\glos{R}:=Z(\sigma_{12})=e^{a_{12}}$
satisfies the famed Yang-Baxter equation,
\[ R^{12}R^{13}R^{23} = R^{23}R^{13}R^{12}, \]
where $R^{ij}$ means ``place $R$ on strands $i$ and $j$''.
\wClipEnd{120208}
\end{remark}

\draftcut
\newpage
\subsection{Some Further Comments} \label{subsec:bcomments}
\subsubsection{Compatibility with Braid Operations}
  \label{subsubsec:BraidCompatibility}
As with
\wClipStart{120215}{0-00-00}
any new gadget, we would like to know how compatible the expansion
$Z$ of the previous section is with the gadgets we already have; namely,
with various operations that are available on w-braids and with the action
of w-braids on the free group $F_n$ (Section~\ref{subsubsec:McCool}).

\parpic[r]{$\xymatrix{
  \wB_n \ar[r]^\theta \ar[d]_Z   & \wB_n \ar[d]^Z    \\
  \calA^w_n \ar[r]_\theta        & \calA^w_n
    \ar@{}[ul]|{\text{\huge$\circlearrowleft$}}
}$}
\paragraph{$Z$ is Compatible with Braid Inversion} \label{par:theta}
Let $\theta$ denote both the ``braid inversion'' operation
$\glos{\theta}\colon \wB_n\to\wB_n$ defined by $B\mapsto B^{-1}$ and
the ``antipode'' anti-automorphism $\theta\colon \calA^w_n\to\calA^w_n$
defined by mapping permutations to their inverses and arrows to their
negatives (that is, $a_{ij}\mapsto-a_{ij}$). Then the diagram on the
right commutes.

\pagebreak[2]

\parpic[r]{$\xymatrix{
  \wB_n \ar[r]^<>(0.5)\Delta \ar[d]_Z   & \wB_n\times\wB_n \ar[d]^{Z\times Z} \\
  \calA^w_n \ar[r]_<>(0.5)\Delta        & \calA^w_n\otimes\calA^w_n
    \ar@{}[ul]|{\text{\huge$\circlearrowleft$}}
}$}
\paragraph{Braid Cloning and the Group-Like Property} \label{par:Delta}
Let $\glos{\Delta}$ denote both the ``braid cloning''
operation $\Delta\colon \wB_n\to\wB_n\times\wB_n$ defined by
$B\mapsto (B,B)$ and the ``co-product'' algebra morphism
$\Delta\colon \calA^w_n\to\calA^w_n\otimes\calA^w_n$ defined by cloning
permutations (that is, $\tau\mapsto\tau\otimes\tau$) and by treating
arrows as primitives (that is, $a_{ij}\mapsto a_{ij}\otimes 1+1\otimes
a_{ij}$). Then the diagram on the right commutes. In formulae, this is
$\Delta(Z(B))=Z(B)\otimes Z(B)$, which is the statement ``$Z(B)$ is
group-like''.

\parpic[r]{$\xymatrix{
  \wB_n \ar[r]^<>(0.5)\iota \ar[d]_Z	& \wB_{n+1} \ar[d]^Z	\\
  \calA^w_n \ar[r]_<>(0.5)\iota	& \calA^w_{n+1}
    \ar@{}[ul]|{\text{\huge$\circlearrowleft$}}
}$}
\paragraph{Strand Insertions} \label{par:iota}
Let $\iota\colon \wB_n\to\wB_{n+1}$ be an operation of ``inert strand
insertion''. Given $B\in\wB_n$, the resulting $\iota B\in\wB_{n+1}$
will be $B$ with one strand $S$ added at some location chosen in
advance --- to the left of all existing strands, or to the right, or
starting from between the 3rd and the 4th strand of $B$ and ending
between the 6th and the 7th strand of $B$; when adding $S$, add it
``inert'', so that all crossings on it are virtual (this is well
defined). There is a corresponding inert strand addition operation
$\iota\colon \calA^w_n\to\calA^w_{n+1}$, obtained by adding a strand at the
same location as for the original $\iota$ and adding no arrows. It is
easy to check that $Z$ is compatible with $\iota$; namely, that the
diagram on the right is commutative.

\parpic[r]{$\xymatrix{
  \wB_n \ar[r]^<>(0.5){d_k} \ar[d]_Z   & \wB_{n-1} \ar[d]^Z    \\
  \calA^w_n \ar[r]_<>(0.5){d_k}        & \calA^w_{n-1}
    \ar@{}[ul]|{\text{\huge$\circlearrowleft$}}
}$}
\paragraph{Strand Deletions} \label{par:deletions} Given $k$ between $1$
and $n$, let $\glos{d_k}\colon \wB_n\to\wB_{n-1}$ the operation of
``removing the $k$th strand''.  This operation induces a homonymous
operation $d_k\colon \calA^w_n\to\calA^w_{n-1}$: if $D\in\calA^w_n$ is an
arrow diagram, $d_kD$ is $D$ with its $k$th strand removed if no arrows
in $D$ start or end on the $k$th strand, and it is $0$ otherwise. It
is easy to check that $Z$ is compatible with $d_k$; namely, that the
diagram on the right is commutative.\footnote{
Section~\ref{subsec:Projectivization}, ``$d_k\colon \wB_n\to\wB_{n-1}$''
is an algebraic structure made of two spaces ($\wB_n$ and $\wB_{n-1}$),
two binary operations (braid composition in $\wB_n$ and in $\wB_{n-1}$),
and one unary operation, $d_k$. After projectivization we get the
algebraic structure $d_k\colon \calA^w_n\to\calA^w_{n-1}$ with $d_k$
as described above, and an alternative way of stating our assertion is
to say that $Z$ is a morphism of algebraic structures. A similar remark
applies (sometimes in the negative form) to the other operations discussed
in this section.}

\parpic[r]{$\xymatrix{
  F_n \ar@{}[r]|{\mathsize{\Huge}{\actsonright}} \ar[d]_Z & \wB_n \ar[d]^Z \\
  \FA_n \ar@{}[r]|{\mathsize{\Huge}{\actsonright}} & \calA^w_n 
    \ar@{}[ul]|{\text{\huge$\circlearrowleft$}}
}$}
\paragraph{Compatibility with the action on $F_n$} \label{par:action}
Let $\glos{\FA_n}$ denote the (degree-completed) free associative (but
not commutative) algebra on generators $\glos{x_1,\dots,x_n}$. Then
there is an ``expansion'' $Z\colon F_n\to \FA_n$ defined by $\xi_i\mapsto
e^{x_i}$ (see~\cite{Lin:Expansions} and the related ``Magnus Expansion''
of~\cite{MagnusKarrasSolitar:CGT}). Also, there is a right action of
$\calA^w_n$ on $\FA_n$ defined on generators by $x_i\tau=x_{\tau i}$
for $\tau\in S_n$ and by $x_ja_{ij}=[x_i,x_j]$ and $x_ka_{ij}=0$ for
$k\neq j$ and extended by the Leibniz rule to the rest of $\FA_n$ and
then multiplicatively to the rest of $\calA^w_n$.

\begin{exercise} Using the language of
Section~\ref{subsec:Projectivization}, verify that $\FA_n=\proj F_n$ and
that when the actions involved are regarded as instances of the algebraic
structure ``one monoid acting on another'', we have that
$\left(\FA_n\actsonright\calA^w_n\right)=\proj\left(F_n\actsonright
\wB_n\right)$. Finally, use the definition of the action in
\eqref{eq:ConceptualPsi} and the commutative diagrams of paragraphs
\ref{par:theta}, \ref{par:Delta} and~\ref{par:iota} to show that the
diagram of paragraph~\ref{par:action} is also commutative.
\end{exercise}

\pagebreak[2]

\parpic[r]{$\begin{array}{c}
  \input figs/StrandDoubling.pstex_t \\
  \xymatrix{
    \wB_n \ar[r]^<>(0.5){u_k} \ar[d]_Z   & \wB_{n+1} \ar[d]^Z    \\
    \calA^w_n \ar[r]_<>(0.5){u_k}        & \calA^w_{n+1}
      \ar@{}[ul]|{\text{\huge$\not\circlearrowleft$}}
  }
\end{array}$}
\paragraph{Unzipping a Strand} \label{par:unzip} Given $k$ between $1$ and
$n$, let $\glos{u_k}\colon \wB_n\to\wB_{n+1}$ the operation of ``unzipping
the $k$th strand'', briefly defined on the right\footnote{Unzipping
a knotted zipper turns a single band into two parallel ones. This
operation is also known as ``strand doubling'', but for compatibility with
operations that will be introduced later, we prefer ``unzipping''.}. The
induced operation $u_k\colon \calA^w_n\to\calA^w_{n+1}$ is also shown on
the right --- if an arrow starts (or ends) on the strand being doubled,
it is replaced by a sum of two arrows that start (or end) on either
of the two ``daughter strands'' (and this is performed for each arrow
independently; so if there are $t$ arrows touching the $k$th strands in
a diagram $D$, then $u_kD$ will be a sum of $2^t$ diagrams).

In some sense, this whole paper as well as the work of
Kashiwara and Vergne~\cite{KashiwaraVergne:Conjecture} and Alekseev
and Torossian~\cite{AlekseevTorossian:KashiwaraVergne} is about coming
to grips with the fact that $Z$ is {\bf not} compatible with $u_k$
(that the diagram on the right is {\bf not} commutative). Indeed,
let $x:=a_{13}$ and $y:=a_{23}$ be as on the right, and let
$s$ be the permutation $21$ and $\tau$ the permutation $231$.
Then $d_1Z(\overcrossing)=d_1(e^{a_{12}}s)=e^{x+y}\tau$ while
$Z(d_1\overcrossing)=e^ye^x\tau$. So the failure of $d_1$ and
$Z$ to commute is the ill-behaviour of the exponential function
when its arguments are not commuting, which is measured by the
BCH formula, central to both~\cite{KashiwaraVergne:Conjecture}
and~\cite{AlekseevTorossian:KashiwaraVergne}.

\subsubsection{Power and Injectivity} The following theorem is due to
Berceanu and
Papadima~\cite[Theorem~5.4]{BerceanuPapadima:BraidPermutation}; a variant of
this theorem are also true for ordinary braids~\cite{Bar-Natan:Homotopy,
Kohno:deRham, HabeggerMasbaum:Milnor}, and can be proven by similar means:

\begin{theorem} $Z\colon \wB_n\to\calA^w_n$ is injective. In other words, finite
type invariants separate w-braids.
\end{theorem}

\begin{proof} Follows immediately from the faithfulness of the action
$F_n\actsonright\wB_n$, from the compatibility of $Z$ with this action,
and from the injectivity of $Z\colon F_n\to\FA_n$ (the latter is well known,
see e.g.~\cite{MagnusKarrasSolitar:CGT, Lin:Expansions}). Indeed
if $B_1$ and $B_2$ are w-braids and $Z(B_1)=Z(B_2)$, then
$Z(\xi)Z(B_1)=Z(\xi)Z(B_2)$ for any $\xi\in F_n$, therefore $\forall\xi\,
Z(\xi\sslash B_1)=Z(\xi\sslash B_2)$, therefore $\forall\xi\,\xi\sslash
B_1=\xi\sslash B_2$, therefore $B_1=B_2$.
\end{proof}

\begin{remark} Apart from the obvious, that $\calA^w_n$ can be computed
degree by degree in exponential time, we do not know a simple formula for
the dimension of the degree $m$ piece of $\calA^w_n$ or a natural basis of
that space. This compares unfavourably with the situation for ordinary
braids (see e.g.~\cite{Bar-Natan:Braids}). Also compare with
Problem~\ref{prob:wCombing} and with Remark~\ref{rem:GutierrezKrstic}.
\end{remark}

\subsubsection{Uniqueness} There is certainly not a unique expansion for
w-braids --- if $Z_1$ is an expansion and and $P$ is any degree-increasing
linear map $\calA^w\to\calA^w$ (a ``pollution'' map), then $Z_2:=(I+P)\circ
Z_1$ is also an expansion, where $I\colon \calA^w\to\calA^w$ is the
identity. But that's all, and if we require a bit more, even that
freedom disappears.

\begin{proposition} If $Z_{1,2}\colon \wB_n\to\calA^w_n$ are expansions then
there exists a degree-increasing linear map $P\colon \calA^w\to\calA^w$ so
that $Z_2:=(I+P)\circ Z_1$.
\end{proposition}

\begin{proof} (Sketch).  Let $\widehat{\wB_n}$ be the unipotent completion
of $\wB_n$. That is, let $\bbQ\wB_n$ be the algebra of formal linear
combinations of w-braids, let $\calI$ be the ideal in $\bbQ\wB_n$ be the
ideal generated by $\semivirtualover=\overcrossing-\virtualcrossing$ and by
$\semivirtualunder=\virtualcrossing-\undercrossing$, and set
\[ \widehat{\wB_n}:=
  \underleftarrow{\lim}_{m\to\infty} \bbQ\wB_n \left/\calI^m\right..
\]
$\widehat{\wB_n}$ is filtered with
$\calF_m\widehat{\wB_n}:=\underleftarrow{\lim}_{m'>m} \calI^m
\left/\calI^{m'}\right..$ An ``expansion'' can be re-interpreted as an
``isomorphism of $\widehat{\wB_n}$ and $\calA^w_n$ as filtered vector
spaces''. Always, any two isomorphisms differ by an automorphism of the
target space, and that's the essence of $I+P$. \qed
\end{proof}

\begin{proposition} If $Z_{1,2}\colon \wB_n\to\calA^w_n$ are homomorphic
expansions that commute with braid cloning (paragraph~\ref{par:Delta}) and
with strand insertion (paragraph~\ref{par:iota}), then $Z_1=Z_2$.
\end{proposition}

\begin{proof} (Sketch). A homomorphic expansion that commutes with strand
insertions is determined by its values on the generators $\overcrossing$,
$\undercrossing$ and $\virtualcrossing$ of $\wB_2$. Commutativity
with braid cloning implies that these values must be (up to permuting
the strands) group like, that is, the exponentials of primitives. But
the only primitives are $a_{12}$ and $a_{21}$, and one may easily
verify that there is only one way to arrange these so that $Z$
will respect $\virtualcrossing^2=\overcrossing\cdot\undercrossing=1$ and
$\semivirtualover\mapsto\rightarrowdiagram+(\text{higher degree terms})$. \qed
\wClipEnd{120215}
\end{proof}

\subsubsection{The group of non-horizontal flying rings}
\label{subsubsec:NonHorRings}
Let
\wClipAt{120118-3}{0}{14}{20}
$\glos{Y_n}$ denote the space of all placements of $n$ numbered disjoint
oriented unlinked geometric circles in $\bbR^3$. Such a placement
is determined by the centres in $\bbR^3$ of the circles, the radii,
and a unit normal vector for each circle pointing in the positive
direction, so $\dim Y_n=3n+n+3n=7n$.  $S_n \ltimes \bbZ_2^n$ acts on
$Y_n$ by permuting the circles and mapping each circle to its image in
either an orientation-preserving or an orientation-reversing way. Let
$\glos{\tilde{Y}_n}$ denote the quotient $Y_n/S_n \ltimes \bbZ_2^n$.
The fundamental group $\pi_1(\tilde{Y}_n)$ can be thought of as the
``group of flippable flying rings''.  Without loss of generality, we
can assume that the basepoint is chosen to be a horizontal placement.
We want to study the relationship of this group to $\wB_n$.

\begin{theorem}
$\pi_1(\tilde{Y}_n)$ is a $\bbZ_2^n$-extension of $\wB_n$, generated
by $s_i$, $\sigma_{i}$ ($1\leq i \leq n-1)$, and $\glos{w_i}$ (``flips''),
for $1\leq i \leq n$; with the relations as above, and in addition:
\[
 w_i^2=1; \qquad w_iw_j=w_jw_i; \qquad w_js_i=s_iw_j \quad \text{when } i\neq j, j+1;
\]
\[
 w_is_i=s_iw_{i+1}; \qquad w_{i+1}s_i=s_iw_i;
\]
\[
w_j\sigma_{i}=\sigma_{i}w_j \quad \text{if } j \neq i, i+1; \quad
w_{i+1}\sigma_{i}=\sigma_{i}w_{i}; \quad \text{yet} \quad
w_i\sigma_{i}=s_i\sigma_i^{-1}s_iw_{i+1}.
\]
\end{theorem}
The two most interesting flip relations in pictures:
\begin{equation}\label{eq:FlipRels}
  \raisebox{-10mm}{\input figs/FlipRels.pstex_t}
\end{equation}

\parpic[r]{\input{figs/FlippingRing.pstex_t}}
Instead of a proof, we provide some heuristics.
Since each circle starts out in a horizontal position and returns 
to a horizontal position, there is an integer number of 
``flips'' they do in between, these are the generators $w_i$, as 
shown on the right.

The first relation says that a double flip is homotopic to doing nothing.
Technically, there are two different directions of flips, and they are the 
same via this (non-obvious but true) relation. The rest of the first line is 
obvious: flips of different rings commute, and if
two rings fly around each other while another one flips, the order of these
events can be switched by homotopy. The second line says that if two rings trade
places with no interaction while one flips, the order of these events can be 
switched as well. However, we have to re-number the flip to conform to the
strand labelling convention.

The only subtle point is how flips interact with crossings. First of all,
if one ring flies through another while a third one flips, the order clearly
does not matter. If a ring flies through another and also flips, the 
order can be switched. However, if ring $A$ flips and then ring $B$ flies 
through it, this is homotopic to ring $B$ flying through ring $A$
from the other direction and then ring $A$ flipping. In other words, commuting
$\sigma_i$ with $w_i$ changes the ``sign of the crossing'' in the sense of
Exercise \ref{ex:swBn}. This gives the last, and the only truly non-commutative flip 
relation.

\parpic[r]{\begin{picture}(0,0)%
\includegraphics{figs/Wen.pstex}%
\end{picture}%
%
%
\setlength{\unitlength}{4934sp}%
\begingroup\makeatletter\ifx\SetFigFont\undefined%
\gdef\SetFigFont#1#2#3#4#5{%
  \reset@font\fontsize{#1}{#2pt}%
  \fontfamily{#3}\fontseries{#4}\fontshape{#5}%
  \selectfont}%
\fi\endgroup%
\begin{picture}(545,1337)(-382,-464)
\end{picture}%
}
\wClipAt{120118-3}{0}{19}{30}
To explain why the flip is denoted by $w$, let us consider the alternative
description by ribbon tubes. A flipping ring traces a so called 
wen\footnote{The term wen was coined by Kanenobu and Shima in 
\cite{KanenobuShima:TwoFiltrationsR2K}}
in $\bbR^4$. A wen is a Klein bottle cut along a meridian circle, 
as shown. The wen is embedded in $\bbR^4$.

Finally, note that $\pi_1Y_n$ is exactly the pure $w$-braid group
$\PwB_n$: since each ring has to return to its original position and
orientation, each does an even number of flips.  The flips (or wens)
can all be moved to the bottoms of the braid diagram strands (to the
bottoms of the tubes, to the beginning of words), at a possible cost, as
specified by Equation~\eqref{eq:FlipRels}. Once together, they pairwise
cancel each other.  As a result, this group can be thought of as not
containing wens at all.
\wClipEnd{120118-3}

\subsubsection{The Relationship with u-Braids} \label{subsubsec:RelWithu} 
For
\wClipStart{120222}{0-00-21}
the sake of ignoring strand permutations, we restrict our
attention to pure braids.

\parpic[r]{$\xymatrix{
  \PuB \ar@{.>}[r]^{Z^u} \ar[d]^a & \calA^u \ar[d]^\alpha \\
  \PwB \ar[r]^{Z^w} & \calA^w
}$}
By Section \ref{subsubsec:FTAlgebraic}, for any expansion $Z^u\colon
\PuB_n \to \calA^u_n$ (where $\PuB_n$ is the ``usual'' braid group
and $\calA^u_n$ is the algebra of horizontal chord diagrams on $n$
strands), there is a square of maps as shown on the right. Here $Z^w$
is the expansion constructed in Section~\ref{subsec:wBraidExpansion},
the left vertical map $\glos{a}$ is the composition of the inclusion
and projection maps $\PuB_n \to \PvB_n \to \PwB_n$. The map $\glos{\alpha}$
is the induced map by the functoriality of projectivization, as noted
after Exercise \ref{ex:BraidsAlgApproach}.  The reader can verify that
$\alpha$ maps each chord to the sum of its two possible directed versions.

Note that this square is {\it not} commutative for any choice of $Z^u$ even
in degree 2: the image of a crossing under $Z^w$ is outside the image 
of $\alpha$.

\parpic[r]{\input{figs/uwsquare2.pstex_t}}
More specifically, for any choice $c$ of a ``parenthesization'' of $n$
points, the KZ-construction / Kontsevich integral (see for example
\cite{Bar-Natan:NAT}) returns an expansion $Z_c^u$ of $u$-braids. As
we shall see in Proposition~\ref{prop:uwBT}, for any choice of $c$, the two
compositions $\alpha \circ Z_c^u$ and $Z^w \circ a$ are ``conjugate in a
bigger space'': there is a map $i$ from $\calA^w$ to a larger space of
``non-horizontal arrow diagrams'', and in this space the images of the
above composites are conjugate.  However, we are not certain that $i$
is an injection, and whether the conjugation leaves the $i$-image of
$\calA^w$ invariant, and so we do not know if the two compositions differ
merely by an outer automorphism of $\calA^w$.

\clearpage
\draftcut
\section{w-Knots} \label{sec:w-knots}

\begin{quote} \small {\bf Section Summary. }
  \summaryknots
\end{quote}

{\bf Knots are the wrong objects for study in knot theory,}
\wClipAt{120222}{0}{11}{55}
v-knots are the wrong objects for study in the theory of v-knotted
objects and w-knots are the wrong objects for study in the theory of
w-knotted objects. Studying uvw-knots on their own is the parallel of
studying cakes and pastries as they come out of the bakery --- we sure
want to make them our own, but the theory of desserts is more about the
ingredients and how they are put together than about the end products. In
algebraic knot theory this reflects through the fact that knots are
not finitely generated in any sense (hence they must be made of some
more basic ingredients), and through the fact that there are very few
operations defined on knots (connected sums and satellite operations
being the main exceptions), and thus most interesting properties of
knots are transcendental, or non-algebraic, when viewed from within the
algebra of knots and operations on knots~\cite{Bar-Natan:AKT-CFA}.

The right objects for study in knot theory, or v-knot theory or w-knot
theory, are thus the ingredients that make up knots and that permit a
richer algebraic structure. These are braids, studied in the previous
section, and even more so tangles and tangled graphs, studied in the
following sections.  Yet tradition has its place and the sweets are
tempting, and we feel compelled to introduce some of the tools we will
use in the deeper and healthier study of w-tangles and w-tangled foams
in the limited but tasty arena of the baked goods of knot theory,
the knots themselves.

\draftcut \subsection{v-Knots and w-Knots} \label{subsec:VirtualKnots}
v-Knots may be understood either as knots drawn on surfaces modulo the
addition or removal of empty handles~\cite{Kauffman:VirtualKnotTheory,
Kuperberg:VirtualLink} or as ``Gauss diagrams'' (Remark~\ref{rem:GD}),
or simply ``unembedded but wired together'' crossings modulo the
Reidemeister moves (\cite{Kauffman:VirtualKnotTheory, Roukema:GPV}
and Section~\ref{subsec:CircuitAlgebras}). But right now we forgo the
topological and the abstract and give only the ``planar'' (and somewhat
less philosophically satisfying) definition of v-knots.

\begin{figure}[h]
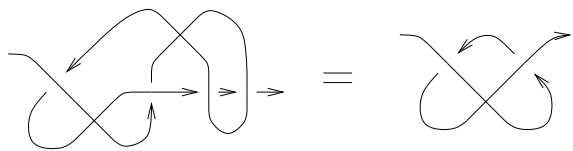

\[ \input figs/VKnot.pstex_t \]
\caption{
  A long v-knot diagram with 2 virtual crossings, 2 positive crossings and
  2 negative crossings. A positive-negative pair can easily be cancelled
  using R2, and then a virtual crossing can be cancelled using VR1, and it
  seems that the rest cannot be simplified any further.
} \label{fig:VKnot}
\end{figure}

\begin{figure}
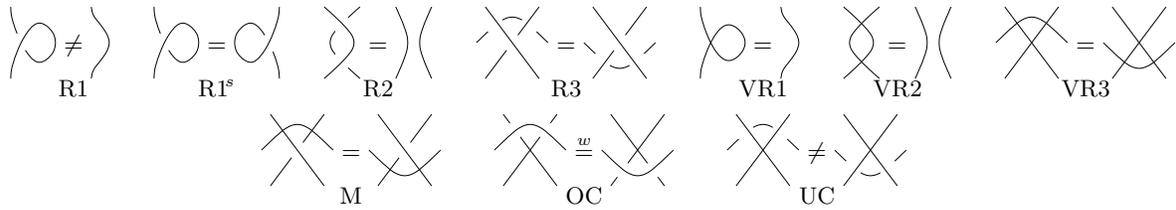

\[ \input figs/VKnotRels.pstex_t \]
\caption{
  The relations defining v-knots and w-knots, along with two relations that
  are {\em not} imposed.
} \label{fig:VKnotRels}
\end{figure}

\begin{definition} A ``long v-knot diagram'' is an arc smoothly
drawn in the plane from $-\infty$ to $+\infty$, with finitely many
self-intersections, divided into ``virtual crossings'' $\virtualcrossing$
and over- and under-crossings, $\overcrossing$ and $\undercrossing$,
and regarded up to planar isotopy. A picture is worth more than a more
formal definition, and one appears in Figure~\ref{fig:VKnot}. A ``long
v-knot'' is an equivalence class of long v-knot diagrams, modulo the
equivalence generated by the Reidemeister $1^{\!s}$, 2 and 3 moves
(\glost{\Rs}, \glost{R2} and \glost{R3})\footnote{
\Rs\ is the ``spun'' version of R1 --- kinks can
be spun around, but not removed outright. See
Figure~\ref{fig:VKnotRels}.}, the virtual Reidemeister 1 through 3 moves
(\glost{VR1}, \glost{VR2}, \glost{VR3}), and by the mixed relations
(\glost{M}); all these are shown in Figure~\ref{fig:VKnotRels}. Finally,
``long w-knots'' are obtained from long v-knots by also dividing
by the Overcrossings Commute (OC) relations, also shown in
Figure~\ref{fig:VKnotRels}.  Note that we never mod out by the
Reidemeister 1 (\glost{R1}) move nor by the Undercrossings Commute
relation (UC).
\end{definition}

\begin{defwarn} A ``circular v-knot'' is like a long v-knot, except
parametrized by a circle rather than by a long line. Unlike the case of
ordinary knots, circular v-knots are {\bf not} equivalent to long v-knots.
The same applies to w-knots.
\end{defwarn}

\begin{defwarn} Long v-knots form a monoid using the concatenation operation
$\#$. Unlike the case of ordinary knots, the resulting monoid is {\bf not}
Abelian. The same applies to w-knots.
\wClipEnd{120222}
\end{defwarn}

\begin{remark} \label{rem:GD} A
\wClipStart{120229}{0-02-53}
``Gauss diagram'' is a straight ``skeleton
line'' along with signed directed chords (signed ``arrows'') marked along
it (more at~\cite{Kauffman:VirtualKnotTheory,
GoussarovPolyakViro:VirtualKnots}). Gauss diagrams are in an obvious
bijection with long v-knot diagrams; the skeleton line of a Gauss diagram
corresponds to the parameter space of the v-knot, and the arrows
correspond to the crossings, with each arrow heading from the upper strand
to the lower strand, marked by the sign of the relevant crossing:
\[ \input figs/GDExample.pstex_t \]
One may also describe the relations in Figure~\ref{fig:VKnotRels} as well
as circular v-knots and other types of v-knots (as we will encounter later)
in terms of Gauss diagrams with varying skeletons.
\end{remark}

\begin{figure}
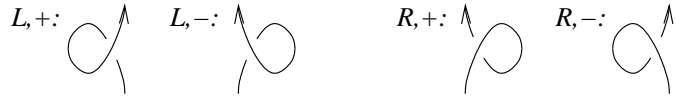

\[ \input figs/Kinks.pstex_t \]
\caption{
  The positive and negative under-then-over kinks (left), and the positive
  and negative over-then-under kinks (right). In each pair the
  negative kink is the $\#$-inverse of the positive kink.
\label{fig:Kinks}}
\end{figure}

\begin{remark}\label{rem:Framing} 
Since we do not mod out by R1, it is perhaps more
appropriate to call our class of v/w-knots ``framed long v/w-knots'',
but since we care more about framed v/w-knots than about unframed ones,
we reserve the unqualified name for the framed case, and when we do wish to
mod out by R1 we will explicitly write ``unframed long v/w-knots''.

Recall that in the case of ``usual knots'', or u-knots, dropping the R1
relation altogether also results in a $\bbZ^2$-extension of unframed
knot theory, where the two factors of $\bbZ$ are framing and rotation
number. If one wants to talk about ``true'' framed knots, one mods out
by the spun Reidemeister 1 relation (\Rs\ of Figure~\ref{fig:VKnotRels}),
which preserves the blackboard framing but does not preserve the rotation
number. We take the analogous approach here, including the \Rs\ relation
but not R1 also in the v and w cases.

This said, note that the monoid of long v-knots is just a central extension
by $\bbZ$ of the monoid of unframed long v-knots, and so studying the
framed case is not very different from studying the unframed case. Indeed
the four ``kinks'' of Figure~\ref{fig:Kinks} generate a central $\bbZ$ within
long v-knots, and it is not hard to show that the sequence
\begin{equation} \label{eq:FramedAndUnframed}
   1\longrightarrow
   \bbZ \longrightarrow
   \{\text{long v-knots}\} \longrightarrow
   \{\text{unframed long v-knots}\} \longrightarrow 1
\end{equation}
is split and exact. The same can be said for w-knots.
\end{remark}

\begin{exercise} \label{ex:sl} Show that a splitting of the
sequence~\eqref{eq:FramedAndUnframed} is given by the ``self-linking''
invariant $\glos{\sl}\colon \{\text{long v-knots}\}\to\bbZ$ defined by
\[
  \sl(K):=\sum_{\text{crossings}\atop x\text{ in }K}\sign x ,
\]
where $K$ is a v-knot diagram, and the sign of a crossing $x$ is defined
so as to agree with the signs in Figure~\ref{fig:Kinks}.
\end{exercise}

\begin{remark} w-Knots are strictly weaker than v-knots --- a notorious
example is the Kishino knot (e.g.~\cite{Dye:Kishinos}) which is non-trivial
as a v-knot yet both it and its mirror are trivial as w-knots. Yet ordinary
knots inject even into w-knots, as the Wirtinger presentation makes sense
for w-knots and therefore w-knots have a ``fundamental quandle'' which
generalizes the fundamental quandle of ordinary
knots~\cite{Kauffman:VirtualKnotTheory}, and as the fundamental
quandle of ordinary knots separates ordinary
knots~\cite{Joyce:TheKnotQuandle}.
\end{remark}

\subsubsection{A topological construction of Satoh's tubing map}
\label{subsubsec:TopTube}
Following Satoh~\cite{Satoh:RibbonTorusKnots}
\wClipComment{120229}{0-36-11}{has more pictures, less formalism}
and using the same constructions as in Section~\ref{subsubsec:ribbon}, we
can map w-knots to (``long'') ribbon tubes in $\bbR^4$ (and the relations
in Figure~\ref{fig:VKnotRels} still hold). It is natural to expect that
this ``tubing'' map is an isomorphism; in other words, that the theory
of w-knots provides a ``Reidemeister framework'' for long ribbon tubes
in $\bbR^4$ --- that every long ribbon tube is in the image of this map
and that two ``w-knot diagrams'' represent the same long ribbon tube iff
they differ by a sequence of moves as in Figure~\ref{fig:VKnotRels}. This
remains unproven.

Let $\glos{\delta}\colon\{\text{v-knots}\} \to \{\text{Ribbon
tori in } \bbR^4\}$ denote the tubing map described in
Section~\ref{subsubsec:ribbon}.  In Satoh's~\cite{Satoh:RibbonTorusKnots}
$\delta$ is called ``Tube''.  It is worthwhile to give a completely
``topological'' definition of $\delta$. To do this we must start with
a topological interpretation of v-knots.

The standard topological interpretation of v-knots
(e.g.~\cite{Kuperberg:VirtualLink}) is that they are oriented framed knots
drawn\footnote{Here and below, ``drawn on $\Sigma$'' means ``embedded in
$\Sigma\times[-\epsilon,\epsilon]$''.} on an oriented surface $\Sigma$,
modulo ``stabilization'', which is the addition and/or removal of empty
handles (handles that do not intersect with the knot). We prefer an
equivalent, yet even more bare-bones approach. For us, a virtual knot is an
oriented framed knot $\gamma$ drawn on a ``virtual surface $\glos{\Sigma}$ for
$\gamma$''. More precisely, $\Sigma$ is an oriented surface that may have
a boundary, $\gamma$ is drawn on $\Sigma$, and the pair $(\Sigma,\gamma)$
is taken modulo the following relations:
\begin{itemize}
\item Isotopies of $\gamma$ on $\Sigma$ (meaning, in
  $\Sigma\times[-\epsilon,\epsilon]$).
\item Tearing and puncturing parts of $\Sigma$ away from $\gamma$:
\end{itemize}
\[ \input{figs/TearingAndPuncturing.pstex_t} \]
(We call $\Sigma$ a ``virtual surface'' because tearing and puncturing
imply that we only care about it in the immediate vicinity of $\gamma$).

We can now define\footnote{Following a private discussion with Dylan
Thurston.} a map $\delta$, defined on v-knots and taking values
in ribbon tori in $\bbR^4$: given $(\Sigma,\gamma)$, embed $\Sigma$
arbitrarily in $\bbR^3_{xzt}\subset\bbR^4$. Note that the unit normal
bundle of $\Sigma$ in $\bbR^4$ is a trivial circle bundle and it has a
distinguished trivialization, constructed using its positive-$y$-direction
section and the orientation that gives each fibre a linking number
$+1$ with the base $\Sigma$.  We say that a normal vector to $\Sigma$
in $\bbR^4$ is ``near unit'' if its norm is between $1-\epsilon$ and
$1+\epsilon$. The near-unit normal bundle of $\Sigma$ has as fibre
an annulus that can be identified with $[-\epsilon,\epsilon]\times
S^1$ (identifying the radial direction $[1-\epsilon,1+\epsilon]$
with $[-\epsilon,\epsilon]$ in an orientation-preserving manner), and
hence the near-unit normal bundle of $\Sigma$ defines an embedding
of $\Sigma\times[-\epsilon,\epsilon]\times S^1$ into $\bbR^4$. On the
other hand, $\gamma$ is embedded in $\Sigma\times[-\epsilon,\epsilon]$ so
$\gamma\times S^1$ is embedded in $\Sigma\times[-\epsilon,\epsilon]\times
S^1$, and we can let $\delta(\gamma)$ be the composition
\[ \gamma\times S^1
  \hookrightarrow\Sigma\times[-\epsilon,\epsilon]\times S^1
  \hookrightarrow\bbR^4,
\]
which is a torus in $\bbR^4$, oriented using the given orientation of
$\gamma$ and the standard orientation of $S^1$.

A framing of a knot (or a v-knot) $\gamma$ can be thought of as a
``nearby companion'' to $\gamma$. Applying the above procedure to a knot
and a nearby companion simultaneously, we find that $\delta$ takes framed
v-knots to framed ribbon tori in $\bbR^4$, where a framing of a tube in
$\bbR^4$ is a continuous up-to-homotopy choice of unit normal vector at
every point of the tube. Note that from the perspective of flying rings as
in Section~\ref{subsubsec:FlyingRings} a framing is a ``companion ring''
to a flying ring. In the framing of $\delta(\gamma)$ the companion ring
is never linked with the main ring, but can fly parallel inside, outside,
above or below it and change these positions, as shown below.

\vspace{2mm}

\begin{center}
 \begin{picture}(0,0)%
\includegraphics{figs/CompanionRing.pstex}%
\end{picture}%
%
%
\setlength{\unitlength}{3947sp}%
\begingroup\makeatletter\ifx\SetFigFont\undefined%
\gdef\SetFigFont#1#2#3#4#5{%
  \reset@font\fontsize{#1}{#2pt}%
  \fontfamily{#3}\fontseries{#4}\fontshape{#5}%
  \selectfont}%
\fi\endgroup%
\begin{picture}(3766,466)(893,-1044)
\end{picture}%

\end{center}

We leave it to the reader to verify that $\delta(\gamma)$ is ribbon, that
it is independent of the choices made within its construction, that it is
invariant under isotopies of $\gamma$ and under tearing and puncturing of
$\Sigma$, that it is also invariant under the ``overcrossings commute''
relation of Figure~\ref{fig:VKnotRels} and hence the true domain of
$\delta$ is w-knots, and that it is equivalent to Satoh's tubing map.

\draftcut
\subsection{Finite Type Invariants of v-Knots and w-Knots}
\label{subsec:FTforvwKnots}

Much as for v-braids and w-braids (Section~\ref{subsec:FT4Braids}) and
much as for ordinary knots (e.g.~\cite{Bar-Natan:OnVassiliev}) we define
finite type invariants for v-knots and for w-knots using an alternation
scheme with $\semivirtualover\to\overcrossing-\virtualcrossing$
and $\semivirtualunder\to\virtualcrossing-\undercrossing$. That is,
we extend any Abelian-group-valued invariant of v- or w-knots to v- or
w-knots also containing ``semi-virtual crossings'' like $\semivirtualover$
and $\semivirtualunder$ using the above assignments, and we declare an
invariant to be ``of type $m$'' if it vanishes on v- or w-knots with more
than $m$ semi-virtuals. As for v- and w-braids and as for ordinary knots,
such invariants have an ``$m$th derivative'', their ``weight system'',
which is a linear functional on the space $\calA^{sv}(\uparrow)$ (for
v-knots) or $\calA^{sw}(\uparrow)$ (for w-knots). We turn to the definitions
of these spaces, following~\cite{GoussarovPolyakViro:VirtualKnots,
Bar-NatanHalachevaLeungRoukema:v-Dims}:

\begin{definition} \label{def:ArrowDiagrams} An ``arrow diagram''
is a chord diagram along a long line (called ``the skeleton''),
in which the chords are oriented (hence ``arrows''). An example is
in Figure~\ref{fig:ADand6T}. Let $\glos{\calD^v(\uparrow)}$ be the space of
formal linear combinations of arrow diagrams.  Let $\glos{\calA^v(\uparrow)}$
be $\calD^v(\uparrow)$ modulo all ``6T relations''. Here a 6T relation is
any (signed) combination of arrow diagrams obtained from the diagrams in
Figure~\ref{fig:6T} by placing the 3 vertical strands there along a long
line in any order, and possibly adding some further arrows in between. An
example is in Figure~\ref{fig:ADand6T}. Let $\glos{\calA^{sv}(\uparrow)}$
be the further quotient of $\calA^v(\uparrow)$ by the \glost{RI} relation,
where the RI (for Rotation number Independence) relation asserts that an
isolated arrow pointing to the right equals an isolated arrow pointing
to the left, as shown in Figure~\ref{fig:ADand6T}\footnote{
  The XII relation of~\cite{Bar-NatanHalachevaLeungRoukema:v-Dims} follows
  from RI and need not be imposed.
}.

\begin{figure}
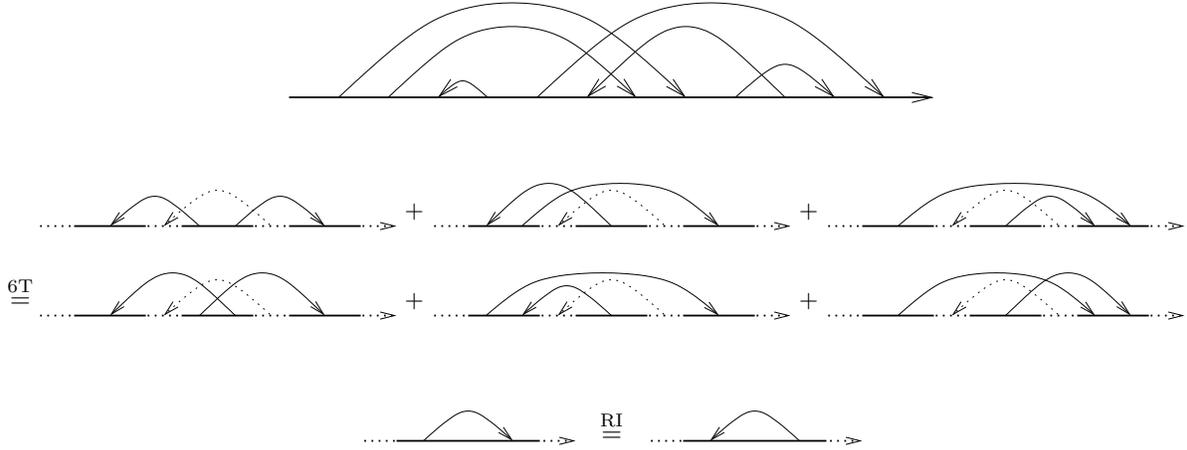

\[ \input figs/ADand6T.pstex_t \]
\caption{
  An arrow diagram of degree 6, a 6T relation, and an RI relation.
} \label{fig:ADand6T}
\end{figure}

Let $\glos{\calA^w(\uparrow)}$ be the further quotient of $\calA^v(\uparrow)$
by the ``Tails Commute'' (TC) relation, first displayed
in Figure~\ref{fig:TCand4T} and reproduced for the case of a
long-line skeleton in Figure~\ref{fig:TCand4TForKnots}. Likewise, let
$\glos{\calA^{sw}(\uparrow)}:=\calA^{sv}(\uparrow)/TC=\calA^w(\uparrow)/RI$.
Alternatively, noting that given TC two of the terms in 6T drop out,
$\calA^w(\uparrow)$ is the space of formal linear combinations
of arrow diagrams modulo TC and $\aft$ relations, displayed in
Figures~\ref{fig:TCand4T} and~\ref{fig:TCand4TForKnots}. Likewise,
$\calA^{sw}=\calD^v/TC,\aft,RI$. Finally, grade $\calD^v(\uparrow)$ and
all of its quotients by declaring that the degree of an arrow diagram
is the number of arrows in it.

\begin{figure}
\[ \input figs/TCand4TForKnots.pstex_t \]
\caption{The TC and the $\protect\aft$ relations for knots.}
\label{fig:TCand4TForKnots}
\end{figure}

\end{definition}

As an example, the spaces $\calA^{v,sv,w,sw}(\uparrow)$
restricted to degrees up to 2 are studied in detail in
Section~\ref{subsec:ToTwo}.

In the same manner as in the theory of finite type invariants of ordinary
knots (see especially~\cite[Section~3]{Bar-Natan:OnVassiliev}, the spaces
$\calA^{-}(\uparrow)$ carry much algebraic structure.  The obvious
juxtaposition product makes them into graded algebras. The product of two
finite type invariants is a finite type invariant (whose type is the sum
of the types of the factors); this induces a product on weight systems,
and therefore a co-product $\Delta$ on arrow diagrams. In brief (and much
the same as in the usual finite type story), the co-product $\Delta D$
of an arrow diagram $D$ is the sum of all ways of dividing the arrows
in $D$ between a ``left co-factor'' and a ``right co-factor''. In summary,

\begin{proposition} \label{prop:CoarseStructure} $\calA^v(\uparrow)$,
$\calA^{sv}(\uparrow)$, $\calA^w(\uparrow)$, and
$\calA^{sw}(\uparrow)$ are co-commutative graded bi-algebras.
\end{proposition}

By the Milnor-Moore theorem~\cite{MilnorMoore:Hopf} we find that
$\calA^{v,sv,w,sw}(\uparrow)$ are the
universal enveloping algebras of their Lie algebras of primitive
elements. Denote these (graded) Lie algebras by
$\glos{\calP^{v,sv,w,sw}(\uparrow)}$, respectively.

When we grow up we'd like to understand $\calA^v(\uparrow)$ and
$\calA^{sv}(\uparrow)$. At the moment we know only very little about these
spaces beyond the generalities of Proposition~\ref{prop:CoarseStructure}.
In the next section some dimensions of low degree parts of
$\calA^{v,sv}(\uparrow)$ are displayed.  Also, given a finite dimensional
Lie bialgebra and a finite dimensional representation thereof, we know
how to construct linear functionals on $\calA^v(\uparrow)$ (one in each
degree)~\cite{Haviv:DiagrammaticAnalogue, Leung:CombinatorialFormulas}
(but not on $\calA^{sv}(\uparrow)$). But we don't even know which degree
$m$ linear functionals on $\calA^{sv}(\uparrow)$ are the weight systems
of degree $m$ invariants of v-knots (that is, we have not solved the
``Fundamental Problem''~\cite{Bar-NatanStoimenow:Fundamental} for
v-knots).

As we shall see below, the situation is much brighter for
$\calA^{w,sw}(\uparrow)$.

\draftcut
\subsection{Some Dimensions} \label{subsec:SomeDimensions}

The table below lists what we could find about $\calA^v$ and $\calA^w$ by
crude brute force computations in low degrees. We list degrees 0 through
7. The spaces we study are $\calA^-(\uparrow)$, $\calA^{s-}(\uparrow)$,
$\calA^{r-}(\uparrow)$ which is $\calA^-(\uparrow)$ moded out by
``isolated'' arrows~\footnote{That is, $\calA^{r-}(\uparrow)$
is $\calA^-(\uparrow)$ modulo ``Framing Independence'' (\glost{FI})
relations~\cite{Bar-Natan:OnVassiliev}, with the isolated arrow
taken with either orientation. It is the space related to finite type
invariants of unframed knots, on which the first Reidemeister move
is also imposed, in the same way as $\calA^-(\uparrow)$ is related
to framed knots.}, $\calP^-(\uparrow)$ which is the space of
primitives in $\calA^-(\uparrow)$, and $\glos{\calA^-(\bigcirc)}$,
$\glos{\calA^{s-}(\bigcirc)}$, and $\glos{\calA^{r-}(\bigcirc)}$,
which are the same as $\calA^-(\uparrow)$, $\calA^{s-}(\uparrow)$,
and $\calA^{r-}(\uparrow)$ except with closed knots (knots with
a circle skeleton) replacing long knots. Each of these spaces we
study in three variants: the ``v'' and the ``w'' variants, as well
as the \underline{u}sual knots ``u'' variant which is here just for
comparison. We also include a row ``$\dim\calG_m\calL ie^-(\uparrow)$''
for the dimensions of ``Lie-algebraic weight systems''. Those
are explained in the u and v cases in~\cite{Bar-Natan:OnVassiliev,
Haviv:DiagrammaticAnalogue, Leung:CombinatorialFormulas}, and in the w
case in Section~\ref{subsec:LieAlgebras}.

{
\def\uvw#1#2#3{{\hspace{-2.5mm}\text{\small
  $\begin{array}{c}#1\mid#2\\#3\end{array}$}\hspace{-3mm}
}}
\begin{center}\begin{tabular}{||c|c||c|c|c|c|c|c|c|c|c||}
\hline \hline
&& \multicolumn{3}{c|}{\footnotesize See Section~\ref{subsec:ToTwo}} &&&&&& \\
$m$ && 0 & 1 & 2 & 3 & 4 & 5 & 6 & 7 & \footnotesize Comments \\
\hline
$\dim\calG_m\calA^-(\uparrow)$ & \uvw{u}{v}{w} &
  \uvw{1}{1}{1} & \uvw{1}{2}{2} & \uvw{2}{7}{4} & \uvw{3}{27}{7} &
  \uvw{6}{139}{12} & \uvw{10}{813}{19} & \uvw{19}{?}{30} & \uvw{33}{?}{45}
  & \uvw{\ref{com:uknots}}{\ref{com:longv}}{\ref{com:wknots},
    \ref{com:longw}, \ref{com:nextfew}} \\
\hline
$\dim\calG_m\calL ie^-(\uparrow)$ & \uvw{u}{v}{w} &
  \uvw{1}{1}{1} & \uvw{1}{2}{2} & \uvw{2}{7}{4} & \uvw{3}{27}{7} &
  \uvw{6}{\,\geq\!128}{12} & \uvw{10}{?}{19} & \uvw{19}{?}{30} &
  \uvw{33}{?}{45}
  & \uvw{\ref{com:uknots}}{\ref{com:Lie}}{\ref{com:nextfew}} \\
\hline
$\dim\calG_m\calA^{s-}(\uparrow)$ & \uvw{u}{v}{w} &
  \uvw{-}{1}{1} & \uvw{-}{1}{1} & \uvw{-}{3}{2} & \uvw{-}{10}{3} &
  \uvw{-}{52}{5} & \uvw{-}{298}{7} & \uvw{-}{?}{11} & \uvw{-}{?}{15}
  & \uvw{\ref{com:su}}{\ref{com:longv}}{\ref{com:wknots}, \ref{com:nextfews}} \\
\hline
$\dim\calG_m\calA^{r-}(\uparrow)$ & \uvw{u}{v}{w} &
  \uvw{1}{1}{1} & \uvw{0}{0}{0} & \uvw{1}{2}{1} & \uvw{1}{7}{1} &
  \uvw{3}{42}{2} & \uvw{4}{246}{2} & \uvw{9}{?}{4} & \uvw{14}{?}{4}
  & \uvw{\ref{com:uknots}}{\ref{com:fiwarning}}{\ref{com:wknots},
    \ref{com:nextfewr}} \\
\hline
$\dim\calG_m\calP^-(\uparrow)$ & \uvw{u}{v}{w} &
  \uvw{0}{0}{0} & \uvw{1}{2}{2} & \uvw{1}{4}{1} & \uvw{1}{15}{1} &
  \uvw{2}{82}{1} & \uvw{3}{502}{1} & \uvw{5}{?}{1} & \uvw{8}{?}{1}
  & \uvw{\ref{com:uknots}}{\ref{com:Pv}}{\ref{com:wknots}} \\
\hline
$\dim\calG_m\calA^-(\bigcirc)$ & \uvw{u}{v}{w} &
  \uvw{1}{1}{1} & \uvw{1}{1}{1} & \uvw{2}{2}{1} & \uvw{3}{5}{1} &
  \uvw{6}{19}{1} & \uvw{10}{77}{1} & \uvw{19}{?}{1} & \uvw{33}{?}{1}
  & \uvw{\ref{com:uknots}}{\ref{com:closedv}}{\ref{com:wknots}} \\
\hline
$\dim\calG_m\calA^{s-}(\bigcirc)$ & \uvw{u}{v}{w} &
  \uvw{-}{1}{1} & \uvw{-}{1}{1} & \uvw{-}{1}{1} & \uvw{-}{2}{1} &
  \uvw{-}{6}{1} & \uvw{-}{23}{1} & \uvw{-}{?}{1} & \uvw{-}{?}{1}
  & \uvw{\ref{com:su}}{\ref{com:longv}}{\ref{com:wknots}} \\
\hline
$\dim\calG_m\calA^{r-}(\bigcirc)$ & \uvw{u}{v}{w} &
  \uvw{1}{1}{1} & \uvw{0}{0}{0} & \uvw{1}{0}{0} & \uvw{1}{1}{0} &
  \uvw{3}{4}{0} & \uvw{4}{17}{0} & \uvw{9}{?}{0} & \uvw{14}{?}{0}
  & \uvw{\ref{com:uknots}}{\ref{com:closedv}}{\ref{com:wknots}} \\
\hline \hline
\end{tabular}\end{center}
}

\begin{comments} \begin{enumerate}
\item \label{com:uknots} Much more is known computationally on the u-knots
  case.  See especially~\cite{Bar-Natan:OnVassiliev, Bar-Natan:Computations,
  Kneissler:Twelve, Amir-KhosraviSankaran:VasCalc}.
\item \label{com:longv} These dimensions were computed by Louis Leung and
  DBN using a program available at~\cite[``Dimensions'']{WKO}.
\item \label{com:wknots} As we shall see in Section~\ref{subsec:Jacobi},
  the spaces associated with w-knots are understood to all degrees.
\item \label{com:longw} To degree 4, these numbers were also verified
  by~\cite[``Dimensions'']{WKO}.
\item \label{com:nextfew} The next few numbers in these sequences are 67, 97,
  139, 195, 272.
\item \label{com:Lie} These dimensions were computed by Louis Leung and
  DBN using a program available at~\cite[``Arrow Diagrams and
  $gl(N)$'']{WKO}. Note the match with the row above.
\item \label{com:su} There is no ``s'' quotient in the ``u'' case.
\item \label{com:nextfews} The next few numbers in this sequence are 22, 30,
  42, 56, 77.
\item \label{com:fiwarning} These numbers were computed
  by~\cite[``Dimensions'']{WKO}. Contrary to the $\calA^u$
  case, $\calA^{rv}$ is {\em not} the quotient of $\calA^{v}$ by the
  ideal generated by degree 1 elements, and therefore the dimensions
  of the graded pieces of these two spaces cannot be deduced from each
  other using the Milnor-Moore theorem.
\item \label{com:nextfewr} The next few numbers in this sequence are
  7,8,12,14,21.
\item \label{com:Pv} These dimensions were deduced from the dimensions of
  $\calG_m\calA^v(\uparrow)$ using the Milnor-Moore theorem.
\item \label{com:closedv} Computed
  by~\cite[``Dimensions'']{WKO}. Contrary to the $\calA^u$
  case, $\calA^v(\bigcirc)$, $\calA^{sv}(\bigcirc)$, and
  $\calA^{rv}(\bigcirc)$ are {\em not} isomorphic to $\calA^v(\uparrow)$,
  $\calA^{sv}(\uparrow)$, and $\calA^{rv}(\uparrow)$ and separate
  computations are required.
\end{enumerate}
\end{comments}

\draftcut
\subsection{Expansions for w-Knots} \label{subsec:Z4Knots}
The notion of ``an expansion'' (or ``a universal finite type invariant'')
for w-knots (or v-knots) is defined in complete analogy with the
parallel notion for ordinary knots (e.g.~\cite{Bar-Natan:OnVassiliev}),
except replacing double points ($\doublepoint$) with semi-virtual
crossings ($\semivirtualover$ and $\semivirtualunder$) and replacing
chord diagrams by arrow diagrams. Alternatively, it is the same as
an expansion for w-braids (Definition~\ref{def:vwbraidexpansion}),
with the obvious replacement of w-braids by w-knots. Just as in the
cases of ordinary knots and/or w-braids, the existence of an expansion
$Z\colon \{\text{w-knots}\}\to\calA^{sw}(\uparrow)$ is equivalent to the
statement ``every weight system integrates'', i.e., ``every degree $m$
linear functional on $\calA^{sw}(\uparrow)$ is the $m$th derivative of
a type $m$ invariant of long w-knots''.

\begin{theorem} \label{thm:ExpansionForKnots}
There exists an expansion $Z\colon \{\text{w-knots}\}\to\calA^{sw}(\uparrow)$.
\end{theorem}
\wClipEnd{120229}

\begin{proof} It is best to define $Z$ by an example, and it is best to
display the example only as a picture:
\wClipStart{120307}{0-00-00}
\[ \input figs/ZwKnotsExample.pstex_t \]
It is clear how to define $Z(K)$ in the general case --- for every crossing
in $K$ place an exponential reservoir of arrows (compare
with~\eqref{eq:reservoir}) next to that crossing, with
the arrows heading from the upper strand to the lower strand, taking
positive reservoirs ($e^a$, with $a$ symbolizing the arrow) for positive
crossings and negative reservoirs ($e^{-a}$) for negative crossings, and
then tug the skeleton until it looks like a straight line. Note that the
Tails Commute relation in $\calA^{sw}$ is used to show that all reasonable
ways of placing an arrow reservoir at a crossing (with its heading and sign
fixed) are equivalent:
\[ \input figs/FourWays.pstex_t \]

The same proof that shows the invariance of $Z$ in the braids case
(Theorem~\ref{thm:RInvariance}) works here as well\footnote{A tiny bit of
extra care is required for invariance under \Rs: it easily follows from
RI.}, and the same argument as in the braids case shows the universality
of $Z$. \qed
\end{proof}

\begin{remark} \label{rem:ZwForGD} Using the language of Gauss diagrams
(Remark~\ref{rem:GD}) the definition of $Z$ is even simpler. Simply map
every positive arrow in a Gauss diagram to a positive ($e^a$) reservoir,
and every negative one to a negative ($e^{-a}$) reservoir:
\[ \input figs/ZwForGD.pstex_t \]
\end{remark}

An expansion (a universal finite type invariant) is as interesting as its
target space, for it is just a tool that takes linear functionals on the
target space to finite type invariants on its domain space. The purpose
of the next section is to find out how interesting are our present target
space, $\calA^{sw}(\uparrow)$, and its ``parent'', $\calA^w(\uparrow)$.

\draftcut
\subsection{Jacobi Diagrams, Trees and Wheels} \label{subsec:Jacobi}

In studying $\calA^w(\uparrow)$ we again follow the model set
by ordinary knots. Compare the following definitions and theorem
with~\cite[Section~3]{Bar-Natan:OnVassiliev}.

\begin{figure}
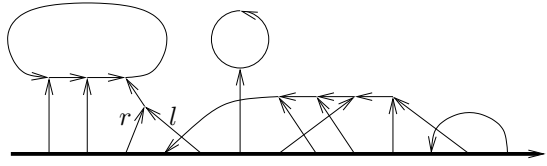

\[ \input figs/wJacDiag.pstex_t \]
\caption{A degree 11 w-Jacobi diagram on a long line skeleton. It has a
skeleton line at the bottom, 13 vertices along the skeleton (of which 2 are
incoming and 11 are outgoing), 9 internal vertices (with only one
explicitly marked with ``left'' ($l$) and ``right'' ($r$)) and one
bubble. The five quadrivalent vertices that seem to appear in the diagram
are just projection artifacts and graph-theoretically, they don't exist.}
\label{fig:wJacDiag}
\end{figure}

\begin{definition} \label{def:wJac} A ``w-Jacobi diagram on a long line
skeleton''\footnote{What a mouthful! We usually short this to
``w-Jacobi diagram'', or sometimes ``arrow diagram'' or just ``diagram''.}
is a connected graph made of the following ingredients:
\begin{itemize}
\item A ``long'' oriented ``skeleton'' line. We usually draw the skeleton
  line a bit thicker for emphasis.
\item Other directed edges, usually called ``arrows''.
\item Trivalent ``skeleton vertices'' in which an arrow starts or ends on
  the skeleton line.
\item Trivalent ``internal vertices'' in which two arrows end and one arrow
  begins. The internal vertices are ``oriented'' --- of the two arrows that
  end in an internal vertices, one is marked as ``left'' and the other is
  marked as ``right''. In reality when a diagram is drawn in the plane, we
  almost never mark ``left'' and ``right'', but instead assume the
  ``left'' and ``right'' inherited from the plane, as seen from the
  outgoing arrow from the given vertex.
\end{itemize}
Note that we allow multiple arrows connecting the same two vertices
(though at most two are possible, given connectedness and trivalence)
and we allow ``bubbles'' --- arrows that begin and end in the same
vertex. Note that for the purpose of determining equality of diagrams
the skeleton line is distinguished.
The ``degree'' of a w-Jacobi diagram is half the number of
trivalent vertices in it, including both internal and skeleton vertices.
An example of a w-Jacobi diagram is in Figure~\ref{fig:wJacDiag}.
\end{definition}

\begin{figure}
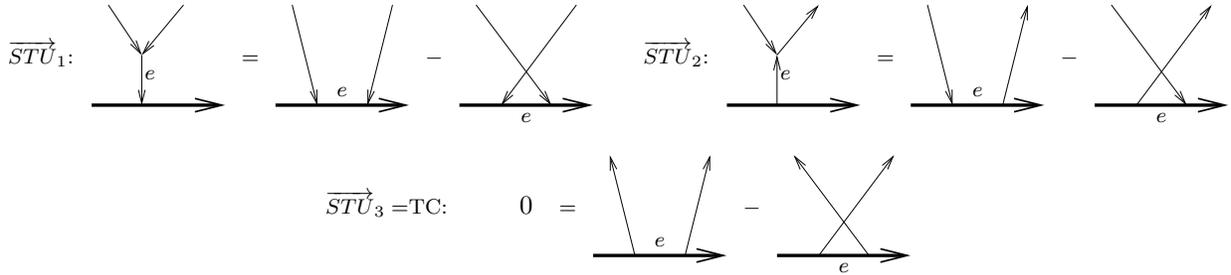

\[ \input figs/aSTU.pstex_t \]
\caption{The $\protect\aSTU_{1,2}$ and TC relations with
their ``central edges'' marked $e$.}
\label{fig:aSTU}
\end{figure}

\begin{definition}
Let $\glos{\calD^{wt}}(\uparrow)$ be the graded vector space of formal
linear combinations of w-Jacobi diagrams on a long line skeleton,
and let $\glos{\calA^{wt}}(\uparrow)$ be $\calD^{wt}(\uparrow)$
modulo the $\glos{\aSTU_1}$, $\glos{\aSTU_2}$, and TC relations of
Figure~\ref{fig:aSTU}. Note that that each diagram appearing in each
$\aSTU$ relation has a ``central edge'' $e$ which can serve as an
``identifying name'' for that $\aSTU$. Thus given a diagram $D$ with
a marked edge $e$ which is either on the skeleton or which contacts
the skeleton, there is an unambiguous $\aSTU$ relation ``around'' or
``along'' the edge $e$.
\end{definition}

\begin{figure}
\[ \input figs/aIHX.pstex_t \]
\caption{The $\protect\aAS$ and $\protect\aIHX$ relations.}
\label{fig:aIHX}
\end{figure}

We like to call the following theorem ``the bracket-rise theorem'',
for it justifies the introduction of internal vertices, and as
should be clear from the $\aSTU$ relations and as will become even
clearer in Section~\ref{subsec:LieAlgebras}, internal vertices can be
viewed as ``brackets''. Two other bracket-rise theorems are Theorem~6
of~\cite{Bar-Natan:OnVassiliev} and Ohtsuki's theorem, Theorem~4.9
of~\cite{Polyak:ArrowDiagrams}.

\begin{theorem}[bracket-rise] \label{thm:BracketRise} The obvious inclusion
$\iota\colon \calD^v(\uparrow)\to\calD^{wt}(\uparrow)$ of arrow diagrams
(Definition~\ref{def:ArrowDiagrams}) into w-Jacobi diagrams descends
to the quotient $\calA^w(\uparrow)$ and induces an isomorphism
$\bar\iota\colon \calA^w(\uparrow)\stackrel{\sim}{\longrightarrow}
\calA^{wt}(\uparrow)$.  Furthermore, the $\glos{\aAS}$ and $\glos{\aIHX}$
relations of Figure~\ref{fig:aIHX} hold in $\calA^{wt}(\uparrow)$.
\end{theorem}

\begin{proof} The proof, joint with D.~Thurston, is modelled after
the proof of Theorem~6 of~\cite{Bar-Natan:OnVassiliev}. To show that
$\iota$ descends to $\calA^w(\uparrow)$ we just need to show that in
$\calA^{wt}(\uparrow)$, $\aft$ follows from $\aSTU_{1,2}$. Indeed,
applying $\aSTU_1$ along the edge $e_1$ and $\aSTU_2$ along $e_2$ in
the picture below, we get the two sides of $\aft$:
\begin{equation} \label{eq:STUto4T}
  \input figs/STUto4T.pstex_t
\end{equation}

The fact that $\bar\iota$ is surjective is obvious; indeed, for diagrams
in $\calA^{wt}(\uparrow)$ that have no internal vertices there is
nothing to show, for they are really in $\calA^w(\uparrow)$. Further,
by repeated use of $\aSTU_{1,2}$ relations, all internal vertices in
any diagram in $\calA^{wt}(\uparrow)$ can be removed (remember that
the diagrams in $\calA^{wt}(\uparrow)$ are always connected, and in
particular, if they have an internal vertex they must have an internal
vertex connected by an edge to the skeleton, and the latter 
vertex can be removed first).

To complete the proof that $\bar\iota$ is an isomorphism it is enough
to show that the ``elimination of internal vertices'' procedure of
the last paragraph is well defined --- that its output is independent
of the order in which $\aSTU_{1,2}$ relations are applied in order to
eliminate internal vertices. Indeed, this done, the elimination map would
by definition satisfy the $\aSTU_{1,2}$ relations and thus descend to
a well defined inverse for $\bar\iota$.

On diagrams with just one internal vertex, Equation~\eqref{eq:STUto4T}
shows that all ways of eliminating that vertex are equivalent modulo $\aft$
relations, and hence the elimination map is well defined on such diagrams.

Now assume that we have shown that the elimination map is well defined on
all diagrams with at most 7 internal vertices, and let $D$ be a diagram
with 8 internal vertices\footnote{``7'' here is a symbol that denotes an
arbitrary natural number greater than 1 and ``8'' denotes $7+1$.}. Let
$e$ and $e'$ be edges in $D$ that connect the skeleton of $D$ to an
internal vertex. We need to show that any elimination process that
begins with eliminating $e$ yields the same answer, modulo $\aft$, as
any elimination process that begins with eliminating $e'$. There are
several cases to consider.

\parpic[r]{$\input figs/CaseI.pstex_t$}
{\bf Case I.} $e$ and $e'$ connect the skeleton to {\em different} internal
vertices of $D$. In this case, after eliminating $e$ we get a signed sum
of two diagrams with exactly 7 internal vertices, and since the elimination
process is well defined on such diagrams, we may as well continue by
eliminating $e'$ in each of those, getting a signed sum of 4 diagrams with
6 internal vertices each. On the other hand, if we start by eliminating
$e'$ we can continue by eliminating $e$, and we get the {\em same} signed
sum of 4 diagrams with 6 internal vertices.

\parpic[r]{$\input figs/CaseII.pstex_t$}
{\bf Case II.} $e$ and $e'$ are connected to the same internal vertex $v$
of $D$, yet some other edge $e''$ exists in $D$ that connects the skeleton
of $D$ to some other internal vertex $v'$ in $D$. In that case, use the
previous case and the transitivity of equality: (elimination starting with
$e$)=(elimination starting with $e''$)=(elimination starting with $e'$).

\parpic[r]{$\input figs/CaseIII.pstex_t$}
{\bf Case III.} Case III is what remains if neither Case I nor Case II
hold. In that case, $D$ must have a schematic form as on the right,
with the ``blob'' not connected to the skeleton other than via $e$
or $e'$, yet further arrows may exist outside of the blob. Let $f$
denote the edge connecting the blob to $e$ and $e'$. The ``two in one
out'' rule for vertices implies that any part of a diagram must have
an excess of incoming edges over outgoing edges, equal to the total
number of vertices in that diagram part. Applying this principle to
the blob, we find that it must contain exactly one vertex, and that $f$
and therefore $e$ and $e'$ must all be oriented upwards.

\parpic[r]{$\input figs/CaseIIIa.pstex_t$}
We leave it to the reader to verify that in this case the two ways of
applying the elimination procedure, $e$ and then $f$ or $e'$ and then $f$,
yield the same answer modulo $\aft$ (in fact, that answer is $0$).

We also leave it to the reader to verify that $\aSTU_1$ implies $\aAS$
and $\aIHX$.  Algebraically, these are restatements of the anti-symmetry
of the bracket and of Jacobi's identity: if $[x,y]:=xy-yx$, then
$0=[x,y]+[y,x]$ and $[x,[y,z]]=[[x,y],z]-[[x,z],y]$. \qed
\wClipEnd{120307}
\end{proof}

Note that $\calA^{wt}(\uparrow)$ inherits algebraic structure from
$\calA^w(\uparrow)$: it is an algebra by concatenation of diagrams,
and a co-algebra with $\Delta(D)$, for $D\in\calD^{wt}(\uparrow)$,
being the sum of all ways of dividing $D$ between a ``left co-factor''
and a ``right co-factor'' so that connected components of $D-S$
are kept intact, where $S$ is the skeleton line of $D$ (compare
with~\cite[Definition~3.7]{Bar-Natan:OnVassiliev}).

As $\calA^w(\uparrow)$ and $\calA^{wt}(\uparrow)$ are canonically
isomorphic, from this point on we will not keep the distinction between
the two spaces.

One may add the RI relation to the definition of $\calA^{wt}(\uparrow)$
to get a space $\calA^{swt}(\uparrow)$, or the FI relation to
get $\calA^{rwt}(\uparrow)$. The statement and proof of the
bracket rise theorem adapt with no difficulty, and we find that
$\calA^{sw}(\uparrow)\cong\calA^{swt}(\uparrow)$ and
$\calA^{rw}(\uparrow)\cong\calA^{rwt}(\uparrow)$.

\begin{figure}
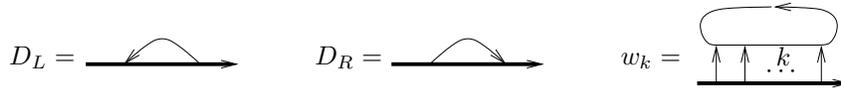

\[ \input figs/AwGenerators.pstex_t \]
\caption{The left-arrow diagram $D_L$, the right-arrow diagram $D_R$ and
  the $k$-wheel $w_k$.}
\label{fig:AwGenerators}
\end{figure}

\begin{theorem} \label{thm:Aw}
The bi-algebra $\calA^w(\uparrow)$ is the bi-algebra of polynomials
in the diagrams $\glos{D_L}$, $\glos{D_R}$ and $\glos{w_k}$
(for $k\geq 1$) shown in Figure~\ref{fig:AwGenerators}, where
$\deg D_L=\deg D_R=1$ and $\deg w_k=k$, subject to the one relation
$w_1=D_L-D_R$. Thus $\calA^w(\uparrow)$ has two generators in degree
1 and one generator in every degree greater than 1, as stated in
Section~\ref{subsec:SomeDimensions}.
\end{theorem}

\begin{proof} (sketch). Readers familiar with the diagrammatic PBW
theorem~\cite[Theorem~8]{Bar-Natan:OnVassiliev} will note that it has
an obvious analogue for the $\calA^w(\uparrow)$ case, and that the proof
in~\cite{Bar-Natan:OnVassiliev} carries through almost verbatim. Namely,
the space $\calA^w(\uparrow)$ is isomorphic to a space $\glos{\calB^w}$
of ``unitrivalent diagrams'' with symmetrized univalent ends modulo
$\aAS$ and $\aIHX$. Given the ``two in one out'' rule for arrow
diagrams in $\calA^w(\uparrow)$ (and hence in $\calB^w$)
the connected components of diagrams in $\calB^w$ can only be
trees or wheels. Trees vanish if they have more than one leaf, as their
leafs are symmetric while their internal vertices are anti-symmetric,
so $\calB^w$ is generated by wheels (which become the $w_k$'s in
$\calA^w(\uparrow)$) and by the one-leaf-one-root tree, which is simply
a single arrow, and which becomes the average of $D_L$ and $D_R$. The
relation $w_1=D_L-D_R$ is then easily verified using $\aSTU_2$.

One may also argue directly, without using sophisticated tools. In
short, let $D$ be a diagram in $\calA^w(\uparrow)$ and $S$ is its
skeleton. Then $D-S$ may have several connected components, whose ``legs''
are intermingled along $S$. Using the $\aSTU$ relations these legs can
be sorted (at a cost of diagrams with fewer connected components, which
could have been treated earlier in an inductive proof). At the end of the
sorting procedure one can see that the only diagrams that remain are our
declared generators. It remains to show that our generators are linearly
independent (apart for the relation $w_1=D_L-D_R$). For the generators
in degree 1, simply write everything out explicitly in the spirit of
Section~\ref{subsubsec:DegreeOne}. In higher degrees there is only one
primitive diagram in each degree, so it is enough to show that $w_k\neq
0$ for every $k$. This can be done ``by hand'', but it is more easily
done using Lie algebraic tools in Section~\ref{subsec:LieAlgebras}. \qed
\end{proof}

\begin{exercise} \label{exe:Asw} Show that the bi-algebra
$\calA^{rw}(\uparrow)$ (see Section~\ref{subsec:SomeDimensions}) is
the bi-algebra of polynomials in the wheel diagrams $w_k$ ($k\geq 2$),
and that $\calA^{sw}(\uparrow)$ is the bi-algebra of polynomials in the
same wheel diagrams and an additional generator $\glos{D_A}:=D_L=D_R$.
\end{exercise}

\begin{theorem} \label{thm:AwCirc} In $\calA^w(\bigcirc)$ all wheels
vanish and hence the bi-algebra $\calA^w(\bigcirc)$ is the bi-algebra
of polynomials in a single variable $D_L=D_R$.
\end{theorem}

\begin{proof} This is Lemma~2.7 of~\cite{Naot:BF}. In short, a wheel in
$\calA^w(\bigcirc)$ can be reduced using $\aSTU_2$ to a difference of
trees. One of these trees has two adjoining leafs and hence is 0 by TC and
$\aAS$. In the other two of the leafs can be commuted ``around the circle''
using TC until they are adjoining and hence vanish by TC and $\aAS$. A
picture is worth a thousand words, but sometimes it takes up more space.
\qed
\end{proof}

\begin{exercise} Show that $\calA^{sw}(\bigcirc)\cong\calA^w(\bigcirc)$
yet $\calA^{rw}(\bigcirc)$ vanishes except in degree $0$.
\end{exercise}

The following two exercises may help the reader to develop a better
``feel'' for $\calA^w(\uparrow)$ and will be needed, within the discussion
of the Alexander polynomial (especially within
Definition~\ref{def:InterpretationMap}).

\parpic[r]{\raisebox{-12mm}{$\input figs/CC.pstex_t$}}
\begin{exercise} Show
\wClipAt{120404}{0}{58}{42}
that the ``Commutators Commute'' (\glost{CC}) relation,
shown on the right, holds in $\calA^w(\uparrow)$. (Interpreted in
Lie algebras as in the next section, this relation becomes $[[x,y],
[z,w]]=0$, and hence the name ``Commutators Commute''). Note that the
proof of CC depends on the skeleton having a single component; later,
when we will work with $\calA^w$-spaces with more complicated skeleta,
the CC relation will not hold.
\end{exercise}

\parpic[r]{\raisebox{-2mm}{$\input figs/Hair.pstex_t$}}
\begin{exercise} \label{ex:Hair} Show that ``detached wheels'' and
``hairy $Y$'s'' make sense in $\calA^w(\uparrow)$. As on the right, a
detached wheel is a wheel with a number of spokes, and a hairy $Y$ is a
combinatorial $Y$ shape with further ``hair'' on its trunk (its outgoing
arrow). It is specified where the trunk and the leafs of the $Y$ connect
to the skeleton, but it is not specified where the spokes of the wheel
and where the hair on the $Y$ connect to the skeleton. The content of the
exercise is to show that modulo the relations of $\calA^w(\uparrow)$,
it is not necessary to specify this further information: all ways of
connecting the spokes and the hair to the skeleton are equivalent. Like
the previous exercise, this result depends on the skeleton having a
single component.
\end{exercise}

\begin{remark} In the case of classical knots and classical chord diagrams,
Jacobi diagrams have a topological interpretation using the
Goussarov-Habiro calculus of claspers~\cite{Goussarov:3Manifolds,
Habiro:Claspers}. In the w case a similar such calculus was developed by 
Watanabe in~\cite{Watanabe:ClasperMoves}. Various related results are 
at~\cite{HabiroKanenobuShima:R2K, HabiroShima:R2KII}.
\end{remark}

\draftcut
\subsection{The Relation with Lie Algebras} \label{subsec:LieAlgebras}
The
\wClipStart{120314}{0-04-47}
theory of finite type invariants of knots is related to the theory
of metrized Lie algebras via the space $\calA$ of chord diagrams, as
explained in~\cite[Theorem~4, Exercise~5.1]{Bar-Natan:OnVassiliev}. In
a similar manner the theory of finite type invariants of w-knots is
related to arbitrary finite-dimensional Lie algebras (or equivalently, to
doubles of co-commutative Lie bialgebra) via the space $\calA^w(\uparrow)$
of arrow diagrams.

\subsubsection{Preliminaries} Given a finite dimensional Lie algebra
$\glos{\frakg}$ let $\glos{I\frakg}:=\frakg^\ast\rtimes\frakg$ be the
semi-direct product of the dual $\frakg^\ast$ of $\frakg$ with $\frakg$,
with $\frakg^\ast$ taken as an Abelian algebra and with $\frakg$ acting
on $\frakg^\ast$ by the usual coadjoint action. In formulae,
\[ I\frakg=\{(\varphi, x)\colon \,\varphi\in\frakg^\ast,\,x\in\frakg\}, \]
\[ [(\varphi_1,x_1), (\varphi_2,x_2)]
  = (x_1\varphi_2-x_2\varphi_1, [x_1,x_2]).
\]

In the case where $\frakg$ is the algebra $so(3)$ of infinitesimal
symmetries of $\bbR^3$, its dual $\frakg^\ast$ is itself $\bbR^3$ with the
usual action of $so(3)$ on it, and $I\frakg$ is the algebra $\bbR^3\rtimes
so(3)$ of infinitesimal affine isometries of $\bbR^3$. This is the
Lie algebra of the Euclidean group of isometries of $\bbR^3$, which is
often denoted $ISO(3)$. This explains our choice of the name $I\frakg$.

Note that if $\frakg$ is a co-commutative Lie bialgebra then $I\frakg$
is the ``double'' of $\frakg$~\cite{Drinfeld:QuantumGroups}. This
is a significant observation, for it is a part of the relationship
between this paper and the Etingof-Kazhdan theory of quantization of
Lie bialgebras~\cite{EtingofKazhdan:BialgebrasI}. Yet we will make no
explicit use of this observation below.

In the construction that follows we are going to construct a map
from $\calA^w$ to $\glos{\calU}(I\frakg)$, the universal enveloping algebra of
$I\frakg$.  Note that a map ${\calA}^w
\to \calU(I\frakg)$ is ``almost the same'' as a map $\calA^{sw} \to
\calU(I\frakg)$, in the following sense.  There is an obvious quotient
map $p\colon {\calA}^w \to \calA^{sw}$, and $p$ has a one-sided inverse
$F\colon  \calA^{sw} \to {\calA}^w$ defined by
\[ F(D)= \sum_{k=0}^\infty \frac{(-1)^k}{k!} S_L^k(D)\cdot w_1^k. \]
Here $S_L$ denotes the map that sends an arrow diagram to the sum of
all ways of deleting a left-going arrow, and $w_1$ denotes the 1-wheel,
as shown in Figure~\ref{fig:AwGenerators}.  The reader can verify that
$F$ is well-defined, an algebra- and co-algebra homomorphism, and that
$p\circ F= id_{\calA^{sw}}$.

\subsubsection{The Construction} Fixing a finite dimensional Lie algebra
$\frakg$ we construct a map $\glos{\calT^w_\frakg}\colon
{\calA}^w\to\calU(I\frakg)$ which assigns to every arrow diagram $D$
an element of the universal enveloping algebra $\calU(I\frakg)$. As is
often the case in our subject, a picture of a typical example is worth
more than a formal definition:
\[
  \def\I{{$I$}}
  \def\B{{$B$}}
  \def\g{{$\frakg$}}
  \def\d{{$\frakg^\ast$}}
  \def\p{{$\frakg^\ast\otimes\frakg^\ast\otimes\frakg\otimes\frakg
    \otimes\frakg^\ast\otimes\frakg^\ast$}}
  \def\u{{$\calU(I\frakg)$}}
  \input figs/Twg.pstex_t
\]

In short, we break up the diagram $D$ into its constituent
pieces and assign a copy of the structure constants tensor
$B\in\frakg^\ast\otimes\frakg^\ast\otimes\frakg$ to each internal
vertex $v$ of $D$ (keeping an association between the tensor factors
in $\frakg^\ast\otimes\frakg^\ast\otimes\frakg$ and the edges emanating
from $v$, as dictated by the orientations of the edges and of the vertex
$v$ itself). We assign the identity tensor in
$\frakg^\ast\otimes\frakg$ to every arrow in $D$ that is not connected to an
internal vertex, and contract any pair of factors connected by a fully
internal arrow. The remaining tensor factors
($\frakg^\ast\otimes\frakg^\ast\otimes\frakg\otimes\frakg
\otimes\frakg^\ast\otimes\frakg^\ast$ in our examples) are all along the
skeleton and can thus be ordered by the skeleton. We then multiply these
factors to get an output $\calT^w_\frakg(D)$ in $\calU(I\frakg)$.

It is also useful to restate this construction given a choice of a basis.
Let $\glos{(x_j)}$ be a basis of $\frakg$ and let $\glos{(\varphi^i)}$
be the dual basis of $\frakg^\ast$, so that $\varphi^i(x_j)=\delta^i_j$,
and let $\glos{b_{ij}^k}$ denote the structure constants of $\frakg$ in
the chosen basis: $[x_i,x_j]=\sum b_{ij}^kx_k$. Mark every arrow in $D$
with lower case Latin letter from within $\{i,j,k,\dots\}$\footnote{The
supply of these can be made inexhaustible by the addition of numerical
subscripts.}. Form a product $P_D$ by taking one $b_{\alpha\beta}^\gamma$
factor for each internal vertex $v$ of $D$ using the letters marking the
edges around $v$ for $\alpha$, $\beta$ and $\gamma$ and by taking one
$x_\alpha$ or $\varphi^\beta$ factor for each skeleton vertex of $D$,
taken in the order that they appear along the skeleton, with the indices
$\alpha$ and $\beta$ dictated by the edge markings and with the choice
between factors in $\frakg$ and factors in $\frakg^\ast$ dictated by the
orientations of the edges. Finally let $\calT^w_\frakg(D)$ be the sum
of $P_D$ over the indices $i,j,k,\dots$ running from $1$ to $\dim\frakg$:

\begin{equation} \label{eq:Twb}
  \def\P{{$\displaystyle
    \sum_{i,j,k,l,m,n=1}^{\dim\frakg}
    \hspace{-4mm} b_{ij}^kb_{kl}^m
    \varphi^i\varphi^jx_nx_m\varphi^l\in\calU(I\frakg)
  $}}
  \input figs/Twb.pstex_t
\end{equation}

The following is easy to verify (compare with~\cite[Theorem~4,
Exercise~5.1]{Bar-Natan:OnVassiliev}):

\begin{proposition} The above two definitions of $T^w_\frakg$ agree, are
independent of the choices made within them, and respect all the relations
defining ${\calA}^w$. \qed
\end{proposition}

While we do not provide a proof of this proposition here, it is worthwhile
to state the correspondence between the relations defining ${\calA}^w$
and the Lie algebraic information in $\calU(I\frakg)$: $\aAS$ is
the antisymmetry of the bracket of $\frakg$, $\aIHX$ is the Jacobi
identity of $\frakg$, $\aSTU_1$ and  $\aSTU_2$ are the relations
$[x_i,x_j]=x_ix_j-x_jx_i$ and $[\varphi^i,x_j]=\varphi^ix_j-x_j\varphi^i$
in $\calU(I\frakg)$, $TC$ is the fact that $\frakg^\star$ is taken as an
Abelian algebra, and $\aft$ is the fact that the identity tensor in
$\frakg^\ast\otimes\frakg$ is $\frakg$-invariant.

\subsubsection{Example: The 2 Dimensional Non-Abelian Lie Algebra}
Let $\frakg$ be the Lie algebra with two generators $x_{1,2}$
satisfying $[x_1,x_2]=x_2$, so that the only non-vanishing structure
constants $b_{ij}^k$ of $\frakg$ are $b_{12}^2=-b_{21}^2=1$. Let
$\varphi^i\in\frakg^\ast$ be the dual basis of $x_i$; by an easy
calculation, we find that in $I\frakg$ the element $\varphi^1$ is
central, while $[x_1,\varphi^2]=-\varphi^2$ and
$[x_2,\varphi^2]=\varphi^1$. We calculate $\calT^w_\frakg(D_L)$,
$\calT^w_\frakg(D_R)$ and $\calT^w_\frakg(w_k)$ using the ``in basis''
technique of Equation~\eqref{eq:Twb}. The outputs of these calculations lie
in $\calU(I\frakg)$; we display these results in a PBW basis in which the
elements of $\frakg^\ast$ precede the elements of $\frakg$:

\begin{eqnarray} 
  \calT^w_\frakg(D_L)
    &=& x_1\varphi^1+x_2\varphi^2 =
      \varphi^1x_1+\varphi^2x_2+[x_2,\varphi^2]
      = \varphi^1x_1+\varphi^2x_2+\varphi_1, \notag \\
  \calT^w_\frakg(D_R) &=& \varphi^1x_1+\varphi^2x_2, \label{eq:2DExample} \\
  \calT^w_\frakg(w_k) &=& (\varphi^1)^k. \notag
\end{eqnarray}

\parpic[r]{$\input figs/4wheel.pstex_t$}
For the last assertion above, note that all non-vanishing structure
constants $b_{ij}^k$ in our case have $k=2$, and therefore all indices
corresponding to edges that exit an internal vertex must be set equal to
$2$. This forces the ``hub'' of $w_k$ to be marked $2$ and therefore the
legs to be marked $1$, and therefore $w_k$ is mapped to $(\varphi^1)^k$.

Note that the calculations in~\eqref{eq:2DExample} are consistent with the
relation $D_L-D_R=w_1$ of Theorem~\ref{thm:Aw} and that they show that
other than that relation, the generators of ${\calA}^w$ are linearly
independent.
\wClipEnd{120314}

\draftcut \subsection{The Alexander Polynomial} \label{subsec:Alexander}
Let
\wClipStart{120404}{0-00-00}
$K$ be a long w-knot, let $Z(K)$ be the invariant of
Theorem~\ref{thm:ExpansionForKnots}. Theorem~\ref{thm:Alexander} below
asserts that apart from self-linking, $Z(K)$ contains precisely the
same information as the Alexander polynomial $A(K)$ of $K$ (defined
below). But we have to start with some definitions.

\begin{figure}
\begin{center}
  $\input figs/8-17.pstex_t$
  \qquad
  \raisebox{-18mm}{\includegraphics[height=40mm]{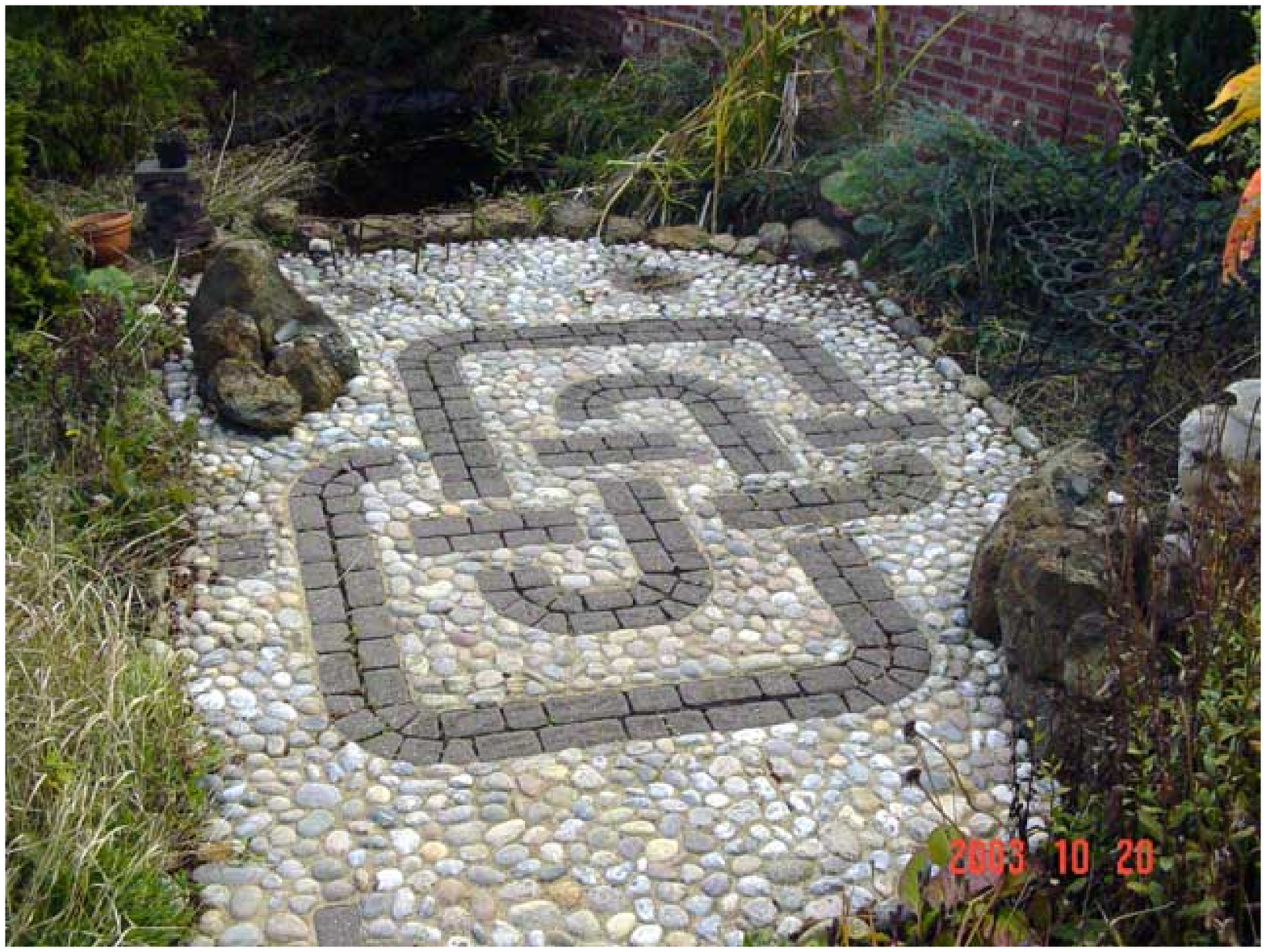}}
\end{center}
\caption{
  A long $8_{17}$, with the span of crossing $\#3$
  marked.  The projection is as in Brian Sanderson's garden.
  See~\cite{WKO}/\href{http://www.math.toronto.edu/~drorbn/papers/WKO/SandersonsGarden.html}{\tt
  SandersonsGarden.html}.
} \label{fig:817}
\end{figure}

\begin{definition} \label{def:STA} Enumerate the crossings of $K$
from $1$ to $n$ in some arbitrary order. For $1\leq j\leq n$, the
``span'' of crossing $\#i$ is the connected open interval along the line
parametrizing $K$ between the two times $K$ ``visits'' crossing $\#i$
(see Figure~\ref{fig:817}). Form a matrix $T=\glos{T(K)}$ with $T_{ij}$
the indicator function of ``the lower strand of crossing $\#j$ is within
the span of crossing $\#i$'' (so $T_{ij}$ is $1$ if for a given $i,j$
the quoted statement is true, and $0$ otherwise). Let $\glos{s_i}$ be the
sign of crossing $\#i$ ($(-,-,-,-,+,+,+,+)$ for Figure~\ref{fig:817}),
let $\glos{d_i}$ be $+1$ if $K$ visits the ``over'' strand of crossing
$\#i$ before visiting the ``under'' strand of that crossing, and let
$d_i=-1$ otherwise ($(-,+,-,+,-,+,-,+$) for Figure~\ref{fig:817}). Let
$S=\glos{S(K)}$ be the diagonal matrix with $S_{ii}=s_id_i$, and for
an indeterminate $\glos{X}$, let $X^{-S}$ denote the diagonal matrix
with diagonal entries $X^{-s_id_i}$.  Finally, let $\glos{A(K)}$ be the
Laurent polynomial in $\bbZ[X,X^{-1}]$ given by
\begin{equation} \label{eq:AKDef}
  A(K)(X) := \det\left(I+T(I-X^{-S})\right).
\end{equation}
\end{definition}

\begin{example} For the knot diagram in Figure~\ref{fig:817},
\[ \scriptstyle
  T = \left(\begin{smallmatrix}
    0 & 1 & 1 & 1 & 1 & 0 & 1 & 0 \\
    0 & 0 & 1 & 0 & 1 & 0 & 0 & 0 \\
    0 & 1 & 0 & 0 & 1 & 0 & 0 & 0 \\
    0 & 1 & 0 & 0 & 1 & 0 & 1 & 0 \\
    0 & 1 & 0 & 1 & 0 & 1 & 1 & 1 \\
    0 & 1 & 0 & 1 & 0 & 0 & 1 & 0 \\
    0 & 0 & 0 & 1 & 0 & 1 & 0 & 0 \\
    0 & 0 & 0 & 1 & 0 & 1 & 0 & 0
  \end{smallmatrix}\right),
  \quad
  S = \left(\begin{smallmatrix}
    1 & 0 & 0 & 0 & 0 & 0 & 0 & 0 \\
    0 & -1 & 0 & 0 & 0 & 0 & 0 & 0 \\
    0 & 0 & 1 & 0 & 0 & 0 & 0 & 0 \\
    0 & 0 & 0 & -1 & 0 & 0 & 0 & 0 \\
    0 & 0 & 0 & 0 & -1 & 0 & 0 & 0 \\
    0 & 0 & 0 & 0 & 0 & 1 & 0 & 0 \\
    0 & 0 & 0 & 0 & 0 & 0 & -1 & 0 \\
    0 & 0 & 0 & 0 & 0 & 0 & 0 & 1
  \end{smallmatrix}\right),
  \quad\text{and}\quad
  A = \left|\begin{smallmatrix}
    1 & 1-X & 1-X^{-1} & 1-X & 1-X & 0 & 1-X & 0 \\
    0 & 1 & 1-X^{-1} & 0 & 1-X & 0 & 0 & 0 \\
    0 & 1-X & 1 & 0 & 1-X & 0 & 0 & 0 \\
    0 & 1-X & 0 & 1 & 1-X & 0 & 1-X & 0 \\
    0 & 1-X & 0 & 1-X & 1 & 1-X^{-1} & 1-X & 1-X^{-1} \\
    0 & 1-X & 0 & 1-X & 0 & 1 & 1-X & 0 \\
    0 & 0 & 0 & 1-X & 0 & 1-X^{-1} & 1 & 0 \\
    0 & 0 & 0 & 1-X & 0 & 1-X^{-1} & 0 & 1
  \end{smallmatrix}\right|.
\]
The last determinant equals $-X^3+4X^2-8X+11-8X^{-1}+4X^{-2}-X^{-3}$,
the Alexander polynomial of the knot $8_{17}$
(e.g.~\cite{Rolfsen:KnotsAndLinks}).
\end{example}

\begin{theorem} \label{thm:AlexanderFormula}
(P.~Lee,~\cite{Lee:AlexanderInvariant}) For any (classical)
knot $K$, $A(K)$ is equal to the normalized Alexander
polynomial~\cite{Rolfsen:KnotsAndLinks} of $K$. \qed
\end{theorem}

The Mathematica notebook~\cite[``wA'']{WKO} verifies
Theorem~\ref{thm:AlexanderFormula} for all prime knots with up to 11
crossings.

The following theorem asserts that $Z(K)$ can be computed from $A(K)$
(Equation~\eqref{eq:AtoZ}) and that modulo a certain additional relation
and with the appropriate identifications in place, $Z(K)$ {\em is} $A(K)$
(Equation~\eqref{eq:ZisA}).

\begin{theorem} \label{thm:Alexander} (Proof in
Section~\ref{subsec:AlexanderProof}). Let $x$ be an indeterminate, let $\sl$
be self-linking as in Exercise~\ref{ex:sl}, let $D_A:=D_L=D_R$ and $w_k$
be as in Figure~\ref{fig:AwGenerators}, and let $\glos{w}\colon
\bbQ\llbracket x\rrbracket \to\calA^w$ be the linear map defined by
$x^k\mapsto w_k$. Then for a w-knot $K$,
\begin{equation} \label{eq:AtoZ}
  Z(K) = 
    \underbrace{
      \exp_{\calA^{sw}}\left(\sl(K)D_A\right)
    }_\text{$\sl$ coded in arrows} \cdot
    \underbrace{
      \exp_{\calA^{sw}}\left(-w\left(\log_{\bbQ\llbracket x\rrbracket}
        A(K)(e^x)
      \right)\right)
    }_\text{main part: Alexander coded in wheels},
\end{equation}
where the logarithm and inner exponentiation are computed by formal power
series in $\bbQ\llbracket x\rrbracket$ and the outer exponentiations
are likewise computed in $\calA^{sw}$.
\end{theorem}

\parpic[r]{$\input figs/wkl.pstex_t$}
Let $\calA^\text{reduced}$ be $\calA^{sw}$ modulo the additional
relations $D_A=0$ and $w_kw_l=w_{k+l}$ for $k,l\neq 1$. The quotient
$\calA^\text{reduced}$ can be identified with vector space of (infinite)
linear combinations of $w_k$'s (with $k\neq 1$).  Identifying the
$k$-wheel $w_k$ with $x^k$, we see that $\calA^\text{reduced}$ is
the space of power series in $x$ having no linear terms. Note by
inspecting~\eqref{eq:AKDef} that $A(K)(e^x)$ never has a term linear
in $x$, and that modulo $w_kw_l=w_{k+l}$, the exponential and the
logarithm in~\eqref{eq:AtoZ} cancel each other out. Hence within
$\calA^\text{reduced}$,

\begin{equation} \label{eq:ZisA} Z(K) = A^{-1}(K)(e^x). \end{equation}

\begin{remark} In~\cite{HabiroKanenobuShima:R2K} K.~Habiro, T.~Kanenobu,
and A.~Shima show that all coefficients of the Alexander polynomial are
finite type invariants of w-knots, and in~\cite{HabiroShima:R2KII}
K.~Habiro and A.~Shima show that all finite type invariants of w-knots are
polynomials in the coefficients of the Alexander polynomial. Thus
Theorem~\ref{thm:Alexander} is merely an ``explicit form'' of these earlier
results.
\end{remark}

\draftcut
\subsection{Proof of Theorem~\ref{thm:Alexander}}
\label{subsec:AlexanderProof}

We start with a sketch. The proof of Theorem~\ref{thm:Alexander} can be
divided in three parts: differentiation, bulk management, and computation.
\wClipEnd{120404}

\newpage
\noindent{\bf Differentiation.} \wClipComment{120418}{0-00-46}{has
further background on $E$, the differential of $\exp$, and the BCH formula.}
Both sides of our goal,
Equation~\eqref{eq:AtoZ}, are exponential in nature. When seeking to
show an equality of exponentials it is often beneficial to compare
their derivatives\footnote{Thanks, Dylan.}. In our case the useful
``derivatives'' to use are the ``Euler operator'' $\glos{E}$ (``multiply
every term by its degree'', an analogue of $f\mapsto xf'$, defined
in Section~\ref{subsubsec:Euler}), and the ``normalized Euler
operator'' $Z\mapsto\glos{\tilE} Z:=Z^{-1}EZ$, which is a variant of the
logarithmic derivative $f\mapsto x(\log f)'=xf'/f$. Since $\tilE$
is one to one (Section~\ref{subsubsec:Euler}) and since we know how
to apply $\tilE$ to the right hand side of Equation~\eqref{eq:AtoZ}
(Section~\ref{subsubsec:Euler}), it is enough to show that with
$\glos{B}:=T(\exp(-xS)-I)$ and suppressing the fixed w-knot $K$ from the
notation,
\begin{equation} \label{eq:EofAtoZ}
  EZ = Z\cdot\left(
    \sl\cdot D_A-w\!\left[x\tr\left( (I-B)^{-1}TS\exp(-xS) \right)\right]
  \right) \qquad \text{ in }\calA^{sw}.
\end{equation}

\noindent{\bf Bulk Management.}
\wClipStart{120425}{0-00-07}
Next we seek to understand the left hand
side of~\eqref{eq:EofAtoZ}. $Z$ is made up of ``quantities in bulk'':
arrows that come in exponential ``reservoirs''. As it turns out,
$EZ$ is made up of the same bulk quantities, but also allowing for a
single non-bulk ``red excitation'' (compare with $Ee^x=x\cdot e^x$; the
``bulk'' $e^x$ remains, and single ``excited red'' $x$ gets created). We
wish manipulate and simplify that red excitation. This is best done by
introducing a certain module, $\glos{\IAM_K}$, the ``Infinitesimal Alexander
Module'' of $K$ (see Section~\ref{subsubsec:IAM}). The elements of $\IAM_K$
can be thought of as names for ``bulk objects with a red excitation'',
and hence there is an ``interpretation map'' $\glos{\iota}\colon
\IAM_K\to\calA^{sw}$, which maps every ``name'' into the object it
represents. There are three special elements in $\IAM_K$: an element
$\glos{\lambda}$, which is the name of $EZ$ (that is, $\iota(\lambda)=EZ$),
the element $\glos{\delta_A}$ which is the name of $D_A\cdot Z$
(so $\iota(\delta_A)=D_A\cdot Z$), and an element $\glos{\omega_1}$
which is the name of a ``detached'' 1-wheel that is appended to
$Z$. The latter can take a coefficient which is a power of $x$,
with $\iota(x^k\omega_1)=w(x^{k+1})\cdot Z=(Z\text{ times a }
(k+1)\text{-wheel})$. Thus it is enough to show that in $\IAM_K$,
\begin{equation} \label{eq:GoalInIAM}
  \lambda = \sl\cdot\delta_A
    - \tr\left((I-B)^{-1}TSX^{-S}\right)\omega_1,
  \quad\text{with}\quad X=e^x.
\end{equation}
Indeed, applying $\iota$ to both sides of the above equation, we get
Equation~\eqref{eq:EofAtoZ} back again.

\noindent{\bf Computation.} Last, we show in
Section~\ref{sec:ComputeLambda} that~\eqref{eq:GoalInIAM} holds true. This
is a computation that happens entirely in $\IAM_K$ and does not mention
finite type invariants, expansions or arrow diagrams in any way.

\subsubsection{The Euler Operator} \label{subsubsec:Euler} Let $A$ be
a completed graded algebra with unit, in which all degrees are $\geq
0$. Define a continuous linear operator $E\colon A\to A$ by setting $Ea=(\deg
a)a$ for homogeneous $a\in A$. In the case $A=\bbQ\llbracket x\rrbracket$,
we have $Ef=xf'$, the standard ``Euler operator'', and hence we adopt
this name for $E$ in general.

We say that $Z\in A$ is a ``perturbation of the identity'' if its
degree 0 piece is 1. Such a $Z$ is always invertible. For such a $Z$,
set $\tilE Z:=Z^{-1}\cdot EZ$, and call the thus (partially) defined
operator $\tilE \colon A\to A$ the ``normalized Euler operator''. From this
point on when we write $\tilE Z$ for some $Z\in A$, we automatically
assume that $Z$ is a perturbation of the identity or that it is trivial
to show that $Z$ is a perturbation of the identity. Note that for
$f\in\bbQ\llbracket x\rrbracket$, we have $\tilE f=x(\log f)'$,
so $\tilE$ is a variant of the logarithmic derivative.

\begin{claim} $\tilE$ is one to one.
\end{claim}

\begin{proof} Assume $Z_1\neq Z_2$ and let $d$ be the smallest degree
in which they differ. Then $d>0$ and in degree $d$ the difference
$\tilE Z_1-\tilE Z_2$ is $d$ times the difference $Z_1-Z_2$, and
hence $\tilE Z_1\neq\tilE Z_2$. \qed
\end{proof}

Thus in order to prove our goal, Equation~\eqref{eq:AtoZ}, it is enough to
compute $\tilE$ of both sides and to show the equality then. We start
with the right hand side of~\eqref{eq:AtoZ}; but first, we need some
simple properties of $E$ and $\tilE$. The proofs of these properties are
routine and hence they are omitted.

\begin{proposition} The following hold true:
\begin{enumerate}
\item $E$ is a derivation: $E(fg)=(Ef)g+f(Eg)$.
\item If $Z_1$ commutes with $Z_2$, then $\tilE(Z_1Z_2)=\tilE Z_1+\tilE Z_2$.
\item If $z$ commutes with $Ez$, then $Ee^z=e^z(Ez)$ and $\tilE e^z=Ez$.
\item If $w\colon A\to\calA$ is a morphism of graded algebras,
then it commutes with $E$ and $\tilE$. \qed
\end{enumerate}
\end{proposition}

Let us denote the right hand side of~\eqref{eq:AtoZ} by $Z_1(K)$. Then by
the above proposition, remembering (Theorem~\ref{thm:Aw}) that $\calA^{sw}$ is
commutative and that $\deg D_A=1$, we have
\[ \tilE Z_1(K) = \sl\cdot D_A-w(E\log A(K)(e^x))
  = \sl\cdot D_A-w\left(x\frac{d}{dx}\log A(K)(e^x)\right).
\]
The rest is an exercise in matrices and
differentiation. $A(K)$ is a determinant~\eqref{eq:AKDef}, and in general,
$\frac{d}{dx}\log\det(M) = \tr\left(M^{-1}\frac{d}{dx}M\right)$. So with
$B=T(e^{-xS}-I)$ (so $M=I-B$), we have
\[ \tilE Z_1(K) =
  \sl\cdot D_A + w\left(x\tr\left((I-B)^{-1}\frac{d}{dx}B\right)\right)
  = \sl\cdot D_A - w\left(x\tr\left((I-B)^{-1}TSe^{-xS}\right)\right),
\]
as promised in Equation~\eqref{eq:EofAtoZ}.

\subsubsection{The Infinitesimal Alexander Module} \label{subsubsec:IAM}
Let $K$ be a w-knot diagram. The Infinitesimal Alexander Module $\IAM_K$
of $K$ is a certain module made from a certain space $\glos{\IAM^0_K}$ of
pictures ``annotating'' $K$ with ``red excitations'' modulo some pictorial
relations that indicate how the red excitations can be moved around. The
space $\IAM^0_K$ in itself is made of three pieces, or ``sectors''. The
``A sector'' in which the excitations are red arrows, the ``Y sector''
in which the excitations are ``red hairy Y-diagrams'', and a rank 1
``W sector'' for ``red hairy wheels''. There is an ``interpretation
map'' $\glos{\iota}\colon \IAM^0_K\to\calA^w$ which descends to a well defined
(and homonymous) $\iota\colon \IAM_K\to\calA^w$. Finally, there are some
special elements $\lambda$ and $\delta_A$ that live in the A sector of
$\IAM^0_K$ and $\omega_1$ that lives in the W sector.

In principle, the description of $\IAM^0_K$ and of $\IAM_K$ can be given
independently of the interpretation map $\iota$, and there are some good
questions to ask about $\IAM_K$ (and the special elements in it) that are
completely independent of the interpretation of the elements of $\IAM_K$ as
``perturbed bulk quantities'' within $\calA^{sw}$. Yet $\IAM_K$ is a
complicated object and we fear its definition will appear completely
artificial without its interpretation. Hence below the two definitions will
be woven together.

\begin{figure}
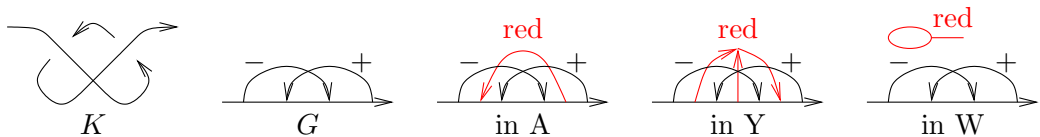

\[ \input figs/IAM0Def.pstex_t \]
\caption{
  A sample w-knot $K$, it's Gauss diagram $G$, and one generator from
  each of the A, Y, and W sectors of $\IAM^0_K$. Red parts are marked
  with the word ``red''.
} \label{fig:IAM0Def}
\end{figure}

$\IAM_K$ and $\iota$ may equally well be described in terms of $K$ or in
terms of the Gauss diagram of $K$ (Remark~\ref{rem:GD}). For pictorial
simplicity, we choose to use the latter; so let $G=G(K)$ be the Gauss
diagram of $K$. It is best to read the following definition while at the
same time studying Figure~\ref{fig:IAM0Def}.

\begin{definition} Let $\glos{R}$ be the ring $\bbZ[X,X^{-1}]$ of
Laurent polynomials in $X$, and let $\glos{R_1}$ be the subring of
polynomials that vanish at $X=1$ (i.e., whose sum of coefficients
is $0$)\footnote{$R_1$ is only very lightly needed, and only within
Definition~\ref{def:InterpretationMap}. In particular, all that we say
about $\IAM_K$ that does not concern the interpretation map $\iota$ is
equally valid with $R$ replacing $R_1$.}.  Let $\IAM^0_K$ be
the direct sum of the following three modules (which for the purpose of
taking the direct sum, are all regarded
as $\bbZ$-modules):
\begin{enumerate}
\item The ``A sector'' is the free $\bbZ$-module generated by all diagrams
made from $G$ by the addition of a single unmarked ``red excitation''
arrow, whose endpoints are on the skeleton of $G$ and are distinct from
each other and from all other endpoints of arrows in $G$. Such diagrams
are considered combinatorially --- so two are equivalent iff they differ
only by an orientation preserving diffeomorphism of the skeleton. Let
us count: if $K$ has $n$ crossings, then $G$ has $n$ arrows and the
skeleton of $G$ get subdivided into $m:=2n+1$ arcs. An A sector diagram
is specified by the choice of an arc for the tail of the red arrow and
an arc for the head ($m^2$ choices), except if the head and the tail
fall within the same arc, their relative ordering has to be specified
as well ($m$ further choices). So the rank of the A sector over $\bbZ$
is $m(m+1)$.
\item The ``Y sector'' is the free $R_1$-module generated by all
diagrams made from $G$ by the addition of a single ``red excitation''
$Y$-shape single-vertex graph, with two incoming edges (``tails'') and
one outgoing (``head''), modulo anti-symmetry for the two incoming edges
(again, considered combinatorially). Counting is more elaborate: when
the three edges of the $Y$ end in distinct arcs in the skeleton of $G$,
we have $\frac12m(m-1)(m-2)$ possibilities ($\frac12$ for the
antisymmetry). When the two tails of the $Y$ lie on the same arc, we get $0$
by anti-symmetry. The remaining possibility is to have the head and
one tail on one arc (order matters!) and the other tail on another,
at $2m(m-1)$ possibilities. So the rank of the Y sector over $R_1$
is $m(m-1)(\frac12m+1)$.
\item The ``W sector'' is the rank 1 free $R$-module with a single
generator $w_1$. It is not necessary for $w_1$ to have a pictorial
representation, yet one, involving a single ``red'' 1-wheel, is shown in
Figure~\ref{fig:IAM0Def}.
\end{enumerate}
\end{definition}

\begin{definition} \label{def:InterpretationMap}
The ``interpretation map'' $\iota\colon \IAM^0_K\to\calA^w$
is defined by sending the arrows (marked $+$ or $-$) of a diagram in
$\IAM^0_K$ to $e^{\pm a}$-exponential reservoirs of arrows, as in the
definition of $Z$ (see Remark~\ref{rem:ZwForGD}). In addition, the red
excitations of diagrams in $\IAM^0_K$ are interpreted as follows:
\begin{enumerate}
\item In the A sector, the red arrow is simply mapped to itself, with the
colour red suppressed.
\item In the $Y$ sector diagrams have red $Y$'s and coefficients $f\in
R_1$. Substitute $X=e^x$ in $f$, expand in powers of $x$,
and interpret $x^kY$ as a ``hairy $Y$ with $k-1$ hairs'' as in
Exercise~\ref{ex:Hair}. Note that $f(1)=0$, so only positive powers of $x$
occur, so we never need to worry about ``$Y$'s with $-1$ hairs''. This is
the only point where the condition $f\in R_1$ (as opposed to $f\in R$) is
needed.
\item In the $W$ sector treat the coefficients as above, but interpret
$x^kw_1$ as a detached $w_{k+1}$. I.e., as a detached wheel with $k+1$
spokes, as in Exercise~\ref{ex:Hair}.
\end{enumerate}
\end{definition}

As stated above, $\IAM_K$ is the quotient of $\IAM^0_K$ by some set of
relations. The best way to think of this set of relations is as
``everything that's obviously annihilated by $\iota$''. Here's the same
thing, in a more formal language:

\begin{figure}
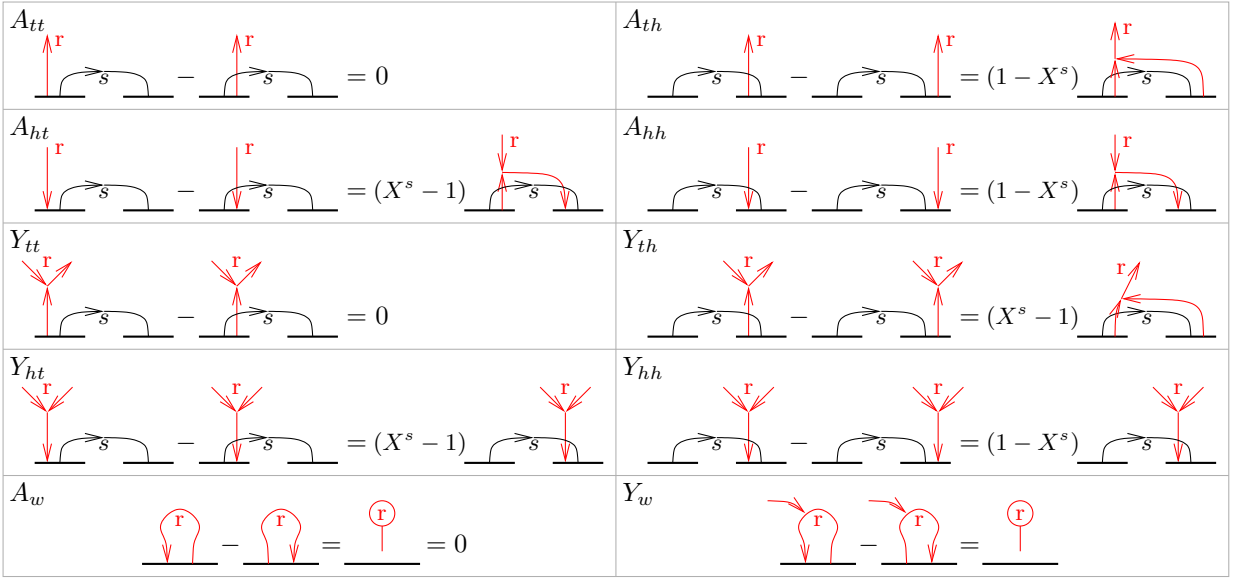

\[ \input figs/IAMRelations.pstex_t \]
\caption{The relations $\calR$ making $\IAM_K$.} \label{fig:IAMRelations}
\end{figure}

\begin{definition} Let $\glos{\IAM_K}:=\IAM^0_K/\calR$, where
$\glos{\calR}$ is the linear span of the relations depicted in
Figure~\ref{fig:IAMRelations}. The top 8 relations are about moving
a leg of the red excitation across an arrow head or an arrow tail in
$G$. Since the red excitation may be either an arrow ($A$) or a $Y$,
its leg in motion may be either a tail or a head, and it may be moving
either past a tail or past a head, there are 8 relations of that type. The
next relation corresponds to $D_L-D_R=w_1=0$. The last relation indicates
the ``price'' (always a red $w_1$) of commuting a red head across a red
tail. As per custom, in each case only the changing part of the diagrams
involved is shown. Further, the red excitations are marked with the
letter ``r'' and the sign of an arrow in $G$ is marked $s$; so always
$s\in\{\pm 1\}$. The relations in the left column may be multiplied
by a scalar in $\bbZ$, while the relations in the right column may be
multiplied by a scalar in $R$. Hence, for example, $x^0w_1=0$ by $A_w$,
yet $x^kw_1\neq 0$ for $k>0$.
\end{definition}

\begin{proposition} The interpretation map $\iota$ indeed annihilates all
the relations in $\calR$.
\end{proposition}

\begin{proof} $\iota A_{tt}$ and $\iota Y_{tt}$ follow immediately from
``Tails Commute''. The formal identity $e^{\ad b}(a)=e^bae^{-b}$ implies
$e^{\ad b}(a)e^b=e^ba$ and hence $ae^b-e^ba=(1-e^{\ad b})(a)e^b$. With
$a$ interpreted as ``red head'', $b$ as ``black head'', and $\ad b$
as ``hair'' (justified by the $\iota$-meaning of hair and by the
$\aSTU_1$ relation, Figure~\ref{fig:aSTU}), the last equality becomes
a proof of $\iota Y_{hh}$.  Further pushing that same equality, we get
$ae^b-e^ba=\frac{1-e^{\ad b}}{\ad b}([b,a])$, where $\frac{1-e^{\ad
b}}{\ad b}$ is first interpreted as a power series $\frac{1-e^y}{y}$
involving only non-negative powers of $y$, and then the substitution
$y=\ad b$ is made. But that's $\iota A_{hh}$, when one remembers
that $\iota$ on the Y sector automatically contains a single
``$\frac{1}{\text{hair}}$'' factor. Similar arguments, though using
$\aSTU_2$ instead of $\aSTU_1$, prove that $Y_{ht}$, $Y_{th}$, $A_{ht}$,
and $A_{th}$ are all in $\ker\iota$. Finally, $\iota A_w$ is RI,
and $\iota Y_w$ is a direct consequence of $\aSTU_2$. \qed
\end{proof}

Finally, we come to the special elements $\lambda$, $\delta_A$, and $\omega_1$.

\begin{figure}
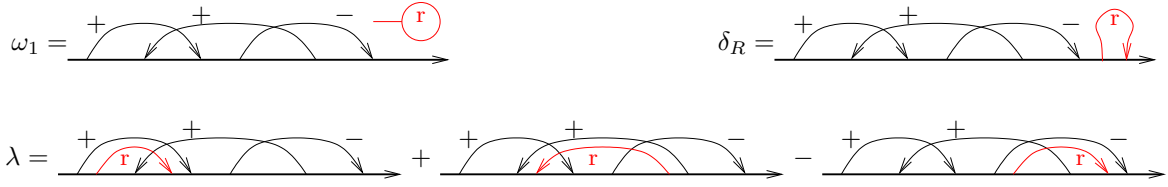

\[ \input figs/SpecialElements.pstex_t \]
\caption{The special elements $\omega_1$, $\delta_A$, and $\lambda$
  in $\IAM_G$, for a sample 3-arrow Gauss diagram $G$.
} \label{fig:SpecialElements}
\end{figure}

\begin{definition} Within $\IAM_G$, let $\omega_1$ be, as before, the
generator of the W sector. Let $\delta_A$ be a ``short''
red arrow, as in the $A_w$ relation (exercise:
modulo $\calR$, this is independent of the placement of the short
arrows within $G$). Finally, let $\lambda$ be the signed sum of exciting
each of the (black) arrows in $G$ in turn. The picture says all, and it is
Figure~\ref{fig:SpecialElements}.
\end{definition}

\begin{proposition} In $\calA^{sw}(\uparrow)$, the special elements of
$\IAM_G$ are interpreted as follows: $\iota(\omega_1)=Zw_1$,
$\iota(\delta_A)=ZD_A$, and most interesting, $\iota(\lambda)=EZ$.
Therefore, Equation~\eqref{eq:GoalInIAM} (if true) implies
Equation~\eqref{eq:EofAtoZ} and hence it implies our goal,
Theorem~\ref{thm:Alexander}.
\end{proposition}

\begin{proof} For the proof of this proposition, the only thing that isn't
done yet and isn't trivial is the assertion $\iota(\lambda)=EZ$. But this
assertion is a consequence of $Ee^{\pm a}=\pm ae^{\pm a}$ and of a
Leibniz law for the derivation $E$, appropriately generalized to a
context where $Z$ can be thought of as a ``product'' of ``arrow
reservoirs''. The details are left to the reader. \qed
\end{proof}

\subsubsection{The Computation of $\lambda$} \label{sec:ComputeLambda}

Naturally, our next task is to prove Equation~\eqref{eq:GoalInIAM}. This is
done entirely algebraically within the finite rank module $\IAM_G$. To
read this section one need not know about $\calA^{sw}(\uparrow)$, or $\iota$,
or $Z$, but we do need to lay out some notation. Start by marking the arrows
of $G$ with $a_1$ through $a_n$ in some order.

\begin{figure}
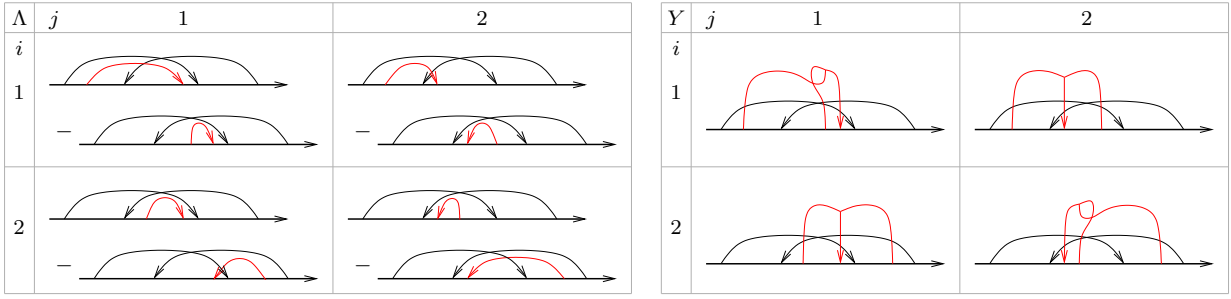

\[ \input figs/LambdaAndY.pstex_t \]
\caption{The matrices $\Lambda$ and $Y$ for a sample 2-arrow Gauss
  diagram (the signs on $a_1$ and $a_2$ are suppressed, and so are the $r$
  marks). The twists in $y_{11}$ and $y_{22}$ may be replaced by minus
  signs.
} \label{fig:LambdaAndY}
\end{figure}

Let $\epsilon$ stand for the informal yet useful quantity
``a little''. Let $\lambda_{ij}$ denote the difference
$\lambda'_{ij}-\lambda''_{ij}$ of red excitations in the A sector of
$\IAM_G$, where $\lambda'_{ij}$ is the diagram with a red arrow whose
tail is $\epsilon$ to the right of the left end of $a_i$ and whose head
is $\frac12\epsilon$ away from head of $a_j$ in the direction of the
tail of $a_j$, and where $\lambda''_{ij}$ has a red arrow whose tail
is $\epsilon$ to the left of the right end of $a_i$ and whose head is
as before, $\frac12\epsilon$ away from head of $a_j$ in the direction
of the tail of $a_j$.  Let $\Lambda=(\lambda_{ij})$ be the matrix whose
entries are the $\lambda_{ij}$'s, as shown in Figure~\ref{fig:LambdaAndY}.

Similarly, let $y_{ij}$ denote the element in the Y sector of $\IAM_G$
whose red Y has its head $\frac12\epsilon$ away from head of $a_j$
in the direction of the tail of $a_j$, its right tail (as seen
from the head) $\epsilon$ to the left of the right end of $a_i$ and
its left tail $\epsilon$ to the right of the left end of $a_i$. Let
$Y=(y_{ij})$ be the matrix whose entries are the $y_{ij}$'s, as shown
in Figure~\ref{fig:LambdaAndY}.

\begin{proposition} \label{prop:IAMStructure}
With $S$ and $T$ as in Definition~\ref{def:STA}, and with $B=T(X^{-S}-I)$
and $\lambda$ as above, the following identities between elements
of $\IAM_G$ and matrices with entries in $\IAM_G$ hold true:
\begin{eqnarray}
  \lambda-\sl\cdot D_A &=& \tr S\Lambda \label{eq:lambda} \\
  \Lambda &=& -BY-TX^{-S}w_1 \label{eq:Lambda} \\
  Y &=& BY + TX^{-S}w_1 \label{eq:Y}
\end{eqnarray}
\end{proposition}

\noindent{\em Proof of Equation~\eqref{eq:GoalInIAM} given
Proposition~\ref{prop:IAMStructure}.} The last of the equalities above
implies that $Y=(I-B)^{-1}TX^{-S}w_1$. Thus
\begin{align*}
  \lambda-\sl\cdot D_A = \tr S\Lambda = -\tr S(BY+TX^{-S}w_1) &=
    -\tr S(B(I-B)^{-1}TX^{-S}+TX^{-S})w_1 \\
  &= -\tr\left((I-B)^{-1}TSX^{-S}\right)w_1,
\end{align*}
and this is exactly Equation~\eqref{eq:GoalInIAM}. \qed

\noindent{\em Proof of Proposition~\ref{prop:IAMStructure}.}
Equation~\eqref{eq:lambda} is trivial. The proofs of
Equations~\eqref{eq:Lambda} and~\eqref{eq:Y} both have the same simple
cores, that have to be supplemented by highly unpleasant tracking of signs
and conventions and powers of $X$. Let us start from the cores.

To prove Equation~\eqref{eq:Lambda} we wish to ``compute''
$\lambda_{ik}=\lambda'_{ik}-\lambda''_{ik}$. As $\lambda'_{ik}$ and
$\lambda''_{ik}$ have their heads in the same place, we can compute their
difference by gradually sliding the tail of $\lambda'_{ik}$ from its
original position near the left end of $a_i$ towards the right end of
$a_i$, where it would be cancelled by $\lambda''_{ik}$. As the tail slides
we pick up a $y_{jk}$ term each time it crosses a head of an $a_j$ (relation
$A_{th}$), we pick up a vanishing term each time it crosses a tail
(relation $A_{tt}$), and we pick up a $w_1$ term if the tail needs to
cross over its own head (relation $A_w$). Ignoring signs and $(X^{\pm
1}-1)$ factors, the sum of the $y_{jk}$-terms should be proportional
to $TY$, for indeed, the matrix $T$ has non-zero entries precisely when
the head of an $a_j$ falls within the span of an $a_i$. Unignoring these
signs and factors, we get $-BY$ (recall that $B=T(X^{-S}-I)$ is just $T$
with added $(X^{\pm 1}-1)$ factors). Similarly, a $w_1$ term arises
in this process when a tail has to cross over its own head, that is,
when the head of $a_k$ is within the span of $a_i$. Thus the $w_1$
term should be proportional to $Tw_1$, and we claim it is $-TX^{-S}w_1$.

The core of the proof of Equation~\eqref{eq:Y} is more or less the
same. We wish to ``compute'' $y_{ik}$ by sliding its left leg, starting
near the left end of $a_i$, towards its right leg, which is stationary
near the right end of $a_i$. When the two legs come together, we get 0
because of the anti-symmetry of Y excitations. Along the way we pick up
further Y terms from the $Y_{th}$ relations, and sometimes a $w_1$ term
from the $Y_w$ relation. When all signs and $(X^{\pm 1}-1)$ factors are
accounted for, we get Equation~\eqref{eq:Y}.

We leave it to the reader to complete the details in the above proofs. It
is a major headache, and we would not have trusted ourselves had we not
written a computer program to manipulate quantities in $\IAM_G$ by a
brute force application of the relations in $\calR$. Everything checks;
see~\cite[``The Infinitesimal Alexander Module'']{WKO}. \qed

This concludes the proof of Theorem~\ref{thm:Alexander}. \qed
\wClipEnd{120425}

\begin{remark} We chose the name ``Infinitesimal Alexander Module'' as in
our mind there is some similarity between $\IAM_K$ and the ``Alexander
Module'' of $K$. Yet beyond the above, we did not embark on any serious
study of $\IAM_K$. In particular, we do not know if $\IAM_K$ in itself
is an invariant of $K$ (though we suspect it wouldn't be hard to show
that it is), we do not know if $\IAM_K$ contains any further information
beyond $\sl$ and the Alexander polynomial, and we do not know if there is any
formal relationship between $\IAM_K$ and the Alexander module of $K$.
\end{remark}

\begin{remark} The logarithmic derivative of the Alexander polynomial
also appears in Lescop's~\cite{Lescop:EquivariantLinking, Lescop:Cube}. We
don't know if its appearances there are related to its appearance here.
\end{remark}

\draftcut
\subsection{The Relationship with u-Knots} \label{subsec:RelWithKont}
Unlike in the case of braids, there is a canonical universal finite type
invariant of $u$-knots: the Kontsevich integral $\glos{Z^u}$. So it
makes sense to ask how it is related to the expansion $Z^w$.

\parpic[l]{$\xymatrix{
  \calK^u(\uparrow) \ar[r]^{Z^u} \ar[d]^a
    & \glos{\calA^u}(\uparrow) \ar[d]^\alpha \\
  \calK^w(\uparrow) \ar[r]^{Z^w}
    & \calA^w(\uparrow)
}$}
We claim that the square on the left commutes, where
$\glos{\calK^u}(\uparrow)$ stands for long $u$-knots (knottings of
an oriented line), and similarly $\calK^w(\uparrow)$ denotes long
$w$-knots. As before, $a$ is the composition of the maps $u$-knots
$\to$ $v$-knots $\to$ $w$-knots, and $\alpha$ is the induced map on
the projectivizations, mapping each chord to the sum of the two ways to
direct it.

Recall that $\alpha$ kills everything but wheels and arrows.
We are going to use the formula for the ``wheel part'' of the Kontsevich integral
as stated in \cite{Kricker:Kontsevich}. 
Let $K$ be a 0-framed long knot, and let $A(K)$ denote the Alexander polynomial. Then by \cite{Kricker:Kontsevich},
$$Z^u(K)= \exp_{\calA^u}\left(-\frac{1}{2} \log A(K)(e^h)|_{h^{2n}\to w^u_{2n}}\right)+\text{ ``loopy terms''},$$
where $w^u_{2n}$ stands for the unoriented wheel with $2n$ spokes; and ``loopy terms'' means terms that contain 
diagrams with more than one loop, which
are killed by $\alpha$. Note that by the symmetry $A(z)=A(z^{-1})$ of the Alexander polynomial,
$A(K)(e^h)$ contains only even powers of $h$, as suggested by the formula.

We need to understand how $\alpha$ acts on wheels. Due to the two-in-one-out
rule, a wheel is zero unless all the ``spokes'' are oriented inward, and the cycle oriented in
one direction. In other words, there are two ways to orient an unoriented wheel:
clockwise or counterclockwise. Due to 
the anti-symmetry of chord vertices, we get that for odd wheels $\alpha(w^u_{2h+1})=0$ and
for even wheels $\alpha(w^u_{2h})=2w^w_{2h}$. As a result,
$$\alpha Z^u(K)=\exp_{\calA^w}\left(-\frac{1}{2} \log A(K)(e^h)|_{h^{2n}\to 2w_{2n}}\right)
=\exp_{\calA^w}\left(-\log A(K)(e^h)|_{h^{2n}\to w_{2n}}\right)$$
which agrees with the formula (\ref{eq:AtoZ}) of Theorem
\ref{thm:Alexander}. Note that since $K$ is 0-framed, the first part
(``$\sl$ coded in arrows'') of (\ref{eq:AtoZ}) is trivial.

\clearpage\draftcut
\section{Algebraic Structures, Projectivizations, Expansions, Circuit Algebras}
\label{sec:generalities}

\begin{quote} \small {\bf Section Summary. }
  \summaryalg
\end{quote}

\subsection{Algebraic Structures} \label{subsec:AlgebraicStructures}

An
\wClipStart{120321}{0-03-41}
``algebraic structure'' $\glos{\calO}$ is some collection $(\calO_\alpha)$
of sets of objects of different kinds, where the subscript
$\alpha$ denotes the ``kind'' of the objects in $\calO_\alpha$,
along with some collection of ``operations'' $\glos{\psi_\beta}$, where
each $\psi_\beta$ is an arbitrary map with domain some product
$\calO_{\alpha_1}\times\dots\times\calO_{\alpha_k}$ of sets of objects,
and range a single set $\calO_{\alpha_0}$ (so operations may be unary or
binary or multinary, but they always return a value of some fixed kind).
We also allow some named ``constants'' within some $\calO_\alpha$'s
(or equivalently, allow some 0-nary operations).\footnote{
alternatively define ``algebraic structures'' using the theory of
``multicategories''~\cite{Leinster:Higher}. Using this language,
an algebraic structure is simply a functor from some ``structure''
multicategory $\calC$ into the multicategory {\bf Set} (or into {\bf
Vect}, if all $\calO_i$ are vector spaces and all operations are
multilinear). A ``morphism'' between two algebraic structures over the
same multicategory $\calC$ is a natural transformation between the two
functors representing those structures.} The operations may or may not
be subject to axioms --- an ``axiom'' is an identity asserting that
some composition of operations is equal to some other composition of
operations.

\begin{figure}[h]
\parbox[m]{7cm}{\caption{An algebraic structure $\calO$ with 4 kinds of
objects and one binary, 3 unary and two 0-nary operations (the constants
$1$ and $\sigma$).}
\label{fig:AlgebraicStructure}}
\parbox[m]{9cm}{\centering
\def\ObjTypeThree{{\raisebox{3mm}{\parbox{0.7in}{\footnotesize
  $\left\{\parbox{0.5in}{\centering objects of kind 3}\right\}=$
}}}}
\input figs/AlgebraicStructure.pstex_t%
}
\end{figure}

Figure~\ref{fig:AlgebraicStructure} illustrates the general notion of an
algebraic structure. Here are a few specific examples:
\begin{itemize}
\item Groups: one kind of objects, one binary ``multiplication'',
  one unary ``inverse'', one constant ``the identity'', and some axioms.
\item Group homomorphisms: Two kinds of objects, one for each
  group. 7 operations --- 3 for each of the two groups and the homomorphism
  itself, going between the two groups. Many axioms.
\item A group acting on a set, a group extension, a split group extension
  and many other examples from group theory.
\item A quandle. It is worthwhile to quote the abstract of the paper that
  introduced the definition (Joyce,~\cite{Joyce:TheKnotQuandle}):
  \begin{quote} \em
    The two operations of conjugation in a group, $x\rhd y=y^{-1}xy$ and
    $x\rhd^{-1}y=yxy^{-1}$ satisfy certain identities. A set with two
    operations satisfying these identities is called a quandle. The
    Wirtinger presentation of the knot group involves only relations of
    the form $y^{-1}xy=z$ and so may be construed as presenting a quandle
    rather than a group. This quandle, called the knot quandle, is not only
    an invariant of the knot, but in fact a classifying invariant of the
    knot.
  \end{quote}
  Also see Definition~\ref{def:quandle}.
\item Planar algebras as in~\cite{Jones:PlanarAlgebrasI} and circuit
  algebras as in Section~\ref{subsec:CircuitAlgebras}.
\item The algebra of knotted trivalent graphs as
  in~\cite{Bar-Natan:AKT-CFA, Dancso:KIforKTG}.
\item Let $\varsigma\colon B\to S$ be an arbitrary homomorphism of groups (though
  our notation suggests what we have in mind --- $B$ may well be braids,
  and $S$ may well be permutations). We can consider an algebraic structure
  $\calO$ whose kinds are the elements of $S$, for which the objects of
  kind $s\in S$ are the elements of $\calO_s:=\varsigma^{-1}(s)$, and with
  the product in $B$ defining operations
  $\calO_{s_1}\times\calO_{s_2}\to\calO_{s_1s_2}$.
\item Clearly, many more examples appear throughout mathematics.
\end{itemize}

\draftcut
\subsection{Projectivization} \label{subsec:Projectivization}

Any algebraic structure $\calO$ has a projectivization. First extend
$\calO$ to allow formal linear combinations of objects of the same kind
(extending the operations in a linear or multi-linear manner), then let
$\glos{\calI}$, the ``augmentation ideal'', be the sub-structure made out of
all such combinations in which the sum of coefficients is $0$, then let
$\calI^m$ be the set of all outputs of algebraic expressions (that is,
arbitrary compositions of the operations in $\calO$) that have at least
$m$ inputs in $\calI$ (and possibly, further inputs in $\calO$), and
finally, set
\begin{equation} \label{eq:projO}
  \glos{\proj}\calO:=\bigoplus_{m\geq 0} \calI^m/\calI^{m+1}.
\end{equation}
Clearly, with the operations inherited from $\calO$, the projectivization
$\proj\calO$ is again algebraic structure with the same multi-graph
of spaces and operations, but with new objects and with new operations
that may or may not satisfy the axioms satisfied by the operations of
$\calO$. The main new feature in $\proj\calO$ is that it is a ``graded''
structure; we denote the degree $m$ piece $\calI^m/\calI^{m+1}$ of
$\proj\calO$ by $\projs_m\calO$.

We believe that many of the most interesting graded structures that appear
in mathematics are the result of this construction, and
that many of the interesting graded equations that appear in mathematics
arise when one tries to find ``expansions'', or ``universal finite type
invariants'', which are also morphisms\footnote{Indeed, if $\calO$ is
finitely presented then finding such a morphism $Z\colon \calO\to\proj\calO$
amounts to finding its values on the generators of $\calO$, subject to
the relations of $\calO$. Thus it is equivalent to solving a system
of equations written in some graded spaces.} $Z\colon \calO\to\proj\calO$
(see Section~\ref{subsec:Expansions}) or when one studies
``automorphisms'' of such expansions\footnote{The Drinfel'd graded
Grothendieck-Teichmuller group $\mathit{GRT}$ is an example of such an
automorphism group. See~\cite{Drinfeld:GalQQ, Bar-Natan:Associators}.}.
Indeed, the paper you are reading now is really the study of the
projectivizations of various algebraic structures associated with
w-knotted objects. We would like to believe that much of the theory of
quantum groups (at ``generic'' $\hbar$) will eventually be shown to be a
study of the projectivizations of various algebraic structures associated
with v-knotted objects.

Thus we believe that the operation described in Equation~\eqref{eq:projO} is
truly fundamental and therefore worthy of a catchy name. So why
``projectivization''? Well, it reminds us of graded spaces, but really,
that's all. We simply found no better name. We're open to suggestions.

Let us end this section with two examples.

\begin{proposition} If $G$ is a group, $\proj G$ is a graded associative
algebra with unit. \qed
\end{proposition}

\begin{definition} \label{def:quandle}
A quandle is a set $Q$ with a binary operation $\glos{\up}\colon Q\times
Q\to Q$ satisfying the following axioms:
\begin{enumerate}
\item $\forall x\in Q,\,x\up x=x$.
\item For any fixed $y\in Q$, the map $x\mapsto x\up y$ is
invertible\footnote{\label{foot:upinv}This can alternatively be stated
as ``there exists a second binary operation $\up^{-1}$ so that $\forall
x,\,x=(x\up y)\up^{-1}y=(x\up^{-1}y)\up y$'', so this axiom can still
be phrased within the language of ``algebraic structures''. Yet note
that below we do not use this axiom at all.}.
\item Self-distributivity:
  $\forall x,y,x\in Q,\,(x\up y)\up z=(x\up z)\up(y\up z)$.
\end{enumerate}
We say that a quandle $Q$ has a unit, or is unital, if there is a
distinguished element $1\in Q$ satisfying the further axiom:
\begin{enumerate}
\setcounter{enumi}{3}
\item $\forall x\in Q,\,x\up 1=x \text{ and } 1\up x=1$.
\end{enumerate}
\end{definition}

If $G$ is a group, it is also a (unital) quandle by setting $x\up
y:=y^{-1}xy$, yet there are many quandles that do not arise from groups
in this way.

\begin{proposition} \label{prop:ProjQ} If $Q$ is a unital quandle,
$\projs_0 Q$ is one-dimensional and $\projs_{>0} Q$ is a graded 
right Leibniz 
algebra\footnote{A Leibniz algebra is a Lie algebra without anticommutativity, as
defined by Loday in \cite{Loday:LeibnizAlg}.}
generated by $\projs_1 Q$.
\end{proposition}

\begin{proof} For any algebraic structure $A$ with just one kind of
objects, $\projs_0 A$ is one-dimensional, generated by the equivalence
class $[x]$ of any single object $x$. In particular, $\projs_0 Q$ is
one-dimensional and generated by $[1]$. Let $\calI\subset\bbQ Q$ be the
augmentation ideal of $Q$. For any $x\in Q$ set $\bar{x}:=x-1\in\calI$.
Then $\calI$ is generated by the $\bar{x}$'s, and therefore $\calI^m$
is generated by expressions involving the operation $\up$ applied to some
$m$ elements of $\bar{Q}:=\{\bar{x}\colon x\in Q\}$ and possibly some further
elements $y_i\in Q$. When regarded in $\calI^m/\calI^{m+1}$, any $y_i$ in
such a generating expression can be replaced by $1$, for the difference
would be the same expression with $y_i$ replaced by $\bar{y}_i$, and
this is now a member of $\calI^{m+1}$. But for any element $z\in\calI$
we have $z\up 1=z$ and $1\up z=0$, so all the $1$'s can be eliminated
from the expressions generating $\calI^m$. Thus $\projs_{>0} Q$ is
generated by $\bar{Q}$ and hence by $\projs_1 Q$.

Let $\glos{\Delta}\colon \bbQ Q\to\bbQ Q\otimes\bbQ Q$ be the linear
extension of the operation $x\mapsto x\otimes x$ defined on $x\in Q$,
and extend $\up$ to a binary operator $\glos{\up_2}\colon (\bbQ
Q\otimes\bbQ Q)\otimes(\bbQ Q\otimes\bbQ Q)\to\bbQ Q\otimes\bbQ Q$ by
using $\up$ twice, to pair the first and third tensor factors and then
to pair the second and the fourth tensor factors. With this language in
place, the self-distributivity axiom becomes the following {\em linear}
statement, which holds for every $x,y,z\in\bbQ Q$:
\begin{equation} \label{eq:LinSelfDist}
  (x\up y)\up z = \up\circ\up_2(x\otimes y\otimes\Delta z).
\end{equation}

Clearly, we need to understand $\Delta$ better.
By direct computation, if $x\in Q$ then $\Delta\bar{x}=\bar{x}\otimes
1+1\otimes\bar{x}+\bar{x}\otimes\bar{x}$. We claim that in general, if $z$
is a generating expression of $\calI^m$ (that is, a formula made of $m$
elements of $\bar{Q}$ and $m-1$ applications of $\up$), then
\begin{equation} \label{eq:QuandleDelta}
  \Delta z=z\otimes 1 + 1\otimes z + \sum z_i'\otimes z_i'',
  \qquad\text{with}\qquad
  \sum z_i'\otimes z_i'' \in
    \sum_{m'+m''=m+1,\atop m',m''>0}\calI^{m'}\otimes\calI^{m''}.
\end{equation}
Indeed, for the generators of $\calI^1$ this had just been shown, and if
$z=z_1\up z_2$ is a generator of $\calI^m$, with $z_1$ and $z_2$ generators
of $\calI^{m_1}$ and $\calI^{m_2}$ with $1\leq m_1,m_2<m$ and $m_1+m_2=m$,
then (using $w\up 1=w$ and $1\up w=0$ for $w\in\calI$),
\begin{eqnarray*}
  \Delta z &=& \Delta(z_1\up z_2)=(\Delta z_1)\up_2(\Delta z_2) \\
  &=& (z_1\otimes 1 + 1\otimes z_1 + \sum z_{1j}'\otimes z_{1j}'')
    \up_2 (z_2\otimes 1 + 1\otimes z_2 + \sum z_{2k}'\otimes z_{2k}'') \\
  &=& (z_1\up z_2)\otimes 1 + 1\otimes(z_1\up z_2) \\
  && + \sum_j \left(
      (z_{1j}'\up z_2)\otimes z_{1j}''
      + z_{1j}'\otimes (z_{1j}''\up z_2)
      + \sum_k (z_{1j}'\up z_{2k}')\otimes (z_{1j}''\up z_{2k}'')
    \right),
\end{eqnarray*}
and it is easy to see that the last line agrees
with~\eqref{eq:QuandleDelta}.

We can now combine Equations~\eqref{eq:LinSelfDist}
and~\eqref{eq:QuandleDelta} to get that for any $x,y,z\in\bbQ Q$,
\[ (x\up y)\up z = (x\up z)\up y + x\up(y\up z)
  +\sum (x\up z_i')\up(y\up z_i'').
\]
If $x\in\calI^{m_1}$, $y\in\calI^{m_2}$, and $z\in\calI^{m_3}$, then
by~\eqref{eq:QuandleDelta} the last term above is in
$\calI^{m_1+m_2+m_3+1}$, and so the above identity becomes the Jacobi
identity $(x\up y)\up z = (x\up z)\up y + x\up(y\up z)$ in
$\projs_{m_1+m_2+m_3}Q$.
\end{proof}

Note that in the above proof neither axiom (1) nor axiom (2) of
Definition~\ref{def:quandle} was used. 

\begin{exercise}
Show that axiom (1) implies the antisymmetry of $\up$ on $\calI^1$. 
\end{exercise}

\draftcut
\subsection{Expansions and Homomorphic Expansions}
\label{subsec:Expansions}
We
\wClipStart{120502}{0-05-48}
start with the definition. Given an algebraic structure $\calO$ let
$\glos{\fil}\calO$ denote the filtered structure of linear combinations of
objects in $\calO$ (respecting kinds), filtered by the powers $(\calI^m)$
of the augmentation ideal $\calI$. Recall also that any graded space
$G=\bigoplus_mG_m$ is automatically filtered, by $\left(\bigoplus_{n\geq
m}G_n\right)_{m=0}^\infty$.

\begin{definition} An ``expansion'' $Z$ for $\calO$
is a map $Z\colon \calO\to\proj\calO$ that preserves the kinds of objects
and whose linear extension (also called $Z$) to $\fil\calO$ respects the
filtration of both sides, and for which $\left(\gr Z\right):
\left(\gr\fil\calO=\proj\calO\right) \to
\left(\gr\proj\calO=\proj\calO\right)$ is the identity map of
$\proj\calO$.
\end{definition}

In practical terms, this is equivalent to saying that $Z$ is a map
$\calO\to\proj\calO$ whose restriction to $\calI^m$ vanishes in degrees
less than $m$ (in $\proj\calO$) and whose degree $m$ piece is the
projection $\calI^m\to\calI^m/\calI^{m+1}$.

We come now to what is perhaps the most crucial definition in this paper.

\begin{definition} A ``homomorphic expansion'' is an expansion which
also commutes with all the algebraic operations defined on the algebraic
structure $\calO$.
\wClipEnd{120321}
\end{definition}

\noindent{\bf Why Bother with Homomorphic Expansions?} Primarily, for two
reasons:
\begin{itemize}
\item Often times $\proj\calO$ is simpler to work with than $\calO$;
for one, it is graded and so it allows for finite ``degree by degree''
computations, whereas often times, such as in many topological examples,
anything in $\calO$ is inherently infinite. Thus it can be beneficial to
translate questions about $\calO$ to questions about $\proj\calO$. A
simplistic example would be, ``is some element $a\in\calO$ the square
(relative to some fixed operation) of an element $b\in\calO$?''. Well, if
$Z$ is a homomorphic expansion and by a finite computation it can be shown
that $Z(a)$ is not a square already in degree $7$ in $\proj\calO$, then
we've given a conclusive negative answer to the example question. Some less
simplistic and more relevant examples appear in~\cite{Bar-Natan:AKT-CFA}.
\item Often times $\proj\calO$ is ``finitely presented'', meaning that it
is generated by some finitely many elements $g_1,\dots,g_k\in\calO$,
subject to some relations $R_1\dots R_{n}$ that can be written in terms
of $g_1,\dots,g_k$ and the operations of $\calO$. In this case, finding a
homomorphic expansion $Z$ is essentially equivalent to guessing the values
of $Z$ on $g_1,\dots,g_k$, in such a manner that these values
$Z(g_1),\dots,Z(g_k)$ would satisfy the $\proj\calO$ versions of the
relations $R_1\dots R_{n}$. So finding $Z$ amounts to solving equations in
graded spaces. It is often the case (as will be demonstrated in this paper;
see also~\cite{Bar-Natan:NAT, Bar-Natan:Associators})
that these equations are very interesting for their own algebraic sake, and
that viewing such equations as arising from an attempt to solve a problem
about $\calO$ sheds further light on their meaning. 
\end{itemize}

In practise, often times the first difficulty in searching for an
expansion (or a homomorphic expansion) $Z\colon \calO\to\proj\calO$ is that its
would-be target space $\proj\calO$ is hard to identify. It is typically
easy to make a suggestion $\calA$ for what $\proj\calO$ could be. It
is typically easy to come up with a reasonable generating set $\calD_m$
for $\calI^m$ (keep some knot theoretic examples in mind, or the case of
quandles as in Proposition~\ref{prop:ProjQ}). It is a bit harder but not
exceedingly difficult to discover some relations $\calR$ satisfied by the
elements of the image of $\calD$ in $\calI^m/\calI^{m+1}$ (4T, $\aft$, and
more in knot theory, the Jacobi relation in Proposition~\ref{prop:ProjQ}).
Thus we set $\calA:=\calD/\calR$; but it is often very hard to be
sure that we found everything that ought to go in $\calR$; so perhaps
our suggestion $\calA$ is still too big? Finding 4T, or Jacobi in
Proposition~\ref{prop:ProjQ} was actually not {\em that} easy. Perhaps we
missed some further relations that are hiding in $\proj Q$, for example?

The notion of an $\calA$-expansion, defined below, solves two problems are
once. Once we find an $\calA$-expansion we know that we've identified
$\proj\calO$ correctly, and we automatically get what we really wanted, a
($\proj\calO$)-valued expansion.

\parpic[r]{\raisebox{-9mm}{$\xymatrix{
  & \calA \ar@<-2pt>[d]_\pi \\
  \calO \ar[ur]^{Z_{\calA}} \ar[r]_<>(0.4)Z
  & \proj\calO \ar@<-2pt>[u]_{\gr Z_\calA}
}$}}
\begin{definition} \label{def:CanProj}
A ``candidate projectivization'' for an algebraic structure
$\calO$ is a graded structure $\glos{\calA}$ with the same operations as
$\calO$ along with a homomorphic surjective graded map $\pi\colon
\calA\to\proj\calO$. An ``$\calA$-expansion'' is a kind and filtration
respecting map $\glos{Z_\calA}\colon \calO\to\calA$ for which $(\gr
Z_\calA)\circ\pi\colon \calA\to\calA$ is the identity. There's no need
to define ``homomorphic $\calA$-expansions''.
\end{definition}

\begin{proposition} \label{prop:CanProj}
If $\calA$ is a candidate projectivization of $\calO$
and $Z_\calA\colon \calO\to\calA$ is a homomorphic $\calA$-expansion, then
$\pi:\calA\to\proj\calO$ is an isomorphism and $Z:=\pi\circ Z_\calA$ is a
homomorphic expansion. (Often in this case, $\calA$ is identified with
$\proj\calO$ and $Z_\calA$ is identified with $Z$).
\end{proposition}

\begin{proof} $\pi$ is surjective by birth. Since $(\gr Z_\calA)\circ\pi$
is the identity, $\pi$ it is also injective and hence it is an
isomorphism. The rest is immediate. \qed
\end{proof}

\draftcut
\subsection{Circuit Algebras} \label{subsec:CircuitAlgebras}

``Circuit algebras'' are so common and everyday, and they make
such a useful language (definitely for the purposes of this paper,
but also elsewhere), we find it hard to believe they haven't made
it into the standard mathematical vocabulary\footnote{Or have they,
and we've been looking the wrong way?}. People familiar with planar
algebras~\cite{Jones:PlanarAlgebrasI} may note that circuit algebras are
just the same as planar algebras, except with the planarity requirement
dropped from the ``connection diagrams'' (and all colourings are dropped
as well). For the rest, we'll start with an image and then move on to
the dry definition.

\begin{figure}
\parpic[r]{\hspace{-1cm}\raisebox{-18mm}{$\input figs/FlipFlop.pstex_t$}}
\caption{
  The J-K flip flop, a very basic memory cell, is an electronic circuit
  that can be realized using 9 components --- two triple-input ``and''
  gates, two standard ``nor'' gates, and 5 ``junctions'' in which 3 wires
  connect (many engineers would not consider the junctions to be real
  components, but we do). Note that the ``crossing'' in the middle of
  the figure is merely a projection artifact and does not indicate an
  electrical connection, and that electronically speaking, we need not
  specify how this crossing may be implemented in $\bbR^3$. The J-K flip
  flop has 5 external connections (labelled J, K, CP, Q, and Q') and hence
  in the circuit algebra of computer parts, it lives in $C_5$. In the
  directed circuit algebra of computer parts it would be in $C_{3,2}$
  as it has 3 incoming wires (J, CP, and K) and two outgoing wires
  (Q and Q').
} \label{fig:FlipFlop}
\end{figure}

\begin{figure}
\parpic[r]{\raisebox{-28mm}{\includegraphics[width=50mm]{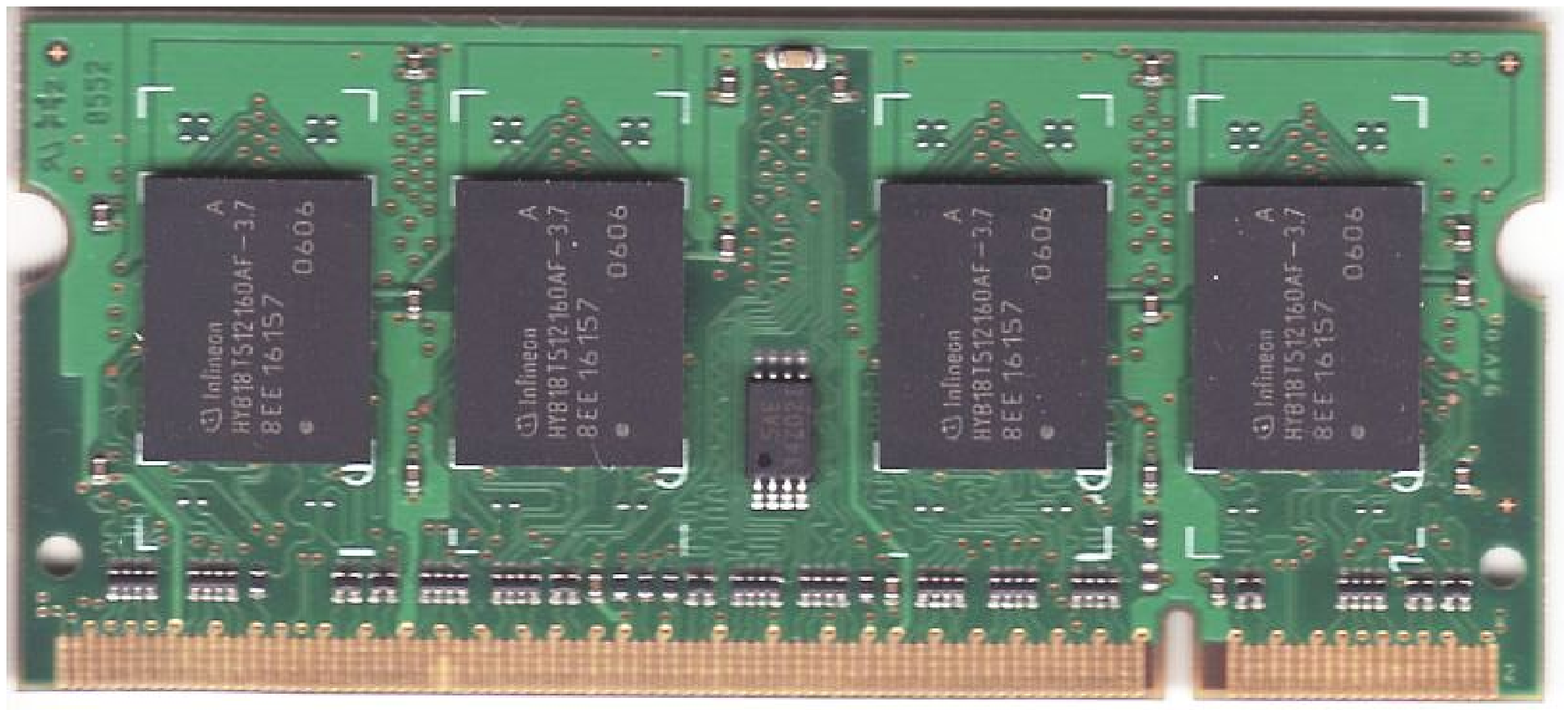}}}
\caption{
  The circuit algebra product of 4 big black components and 1 small black
  component carried out using a green wiring diagram, is an even bigger
  component that has many golden connections (at bottom). When plugged into
  a yet bigger circuit, the CPU board of a laptop, our circuit functions
  as 4,294,967,296 binary memory cells.
} \label{fig:Circuit}
\end{figure}

\begin{image} Electronic circuits are made of ``components'' that can
be wired together in many ways. On a logical level, we only care to
know which pin of which component is connected with which other pin of
the same or other component. On a logical level, we don't really need
to know how the wires between those pins are embedded in space (see
Figures~\ref{fig:FlipFlop} and~\ref{fig:Circuit}). ``Printed Circuit
Boards'' (PCBs) are operators that make smaller components (``chips'')
into bigger ones (``circuits'') --- logically speaking, a PCB is simply a
set of ``wiring instructions'', telling us which pins on which components
are made to connect (and again, we never care precisely how the wires
are routed provided they reach their intended destinations, and ever
since the invention of multi-layered PCBs, all conceivable topologies for
wiring are actually realizable). PCBs can be composed (think ``plugging
a graphics card onto a motherboard''); the result of a composition of
PCBs, logically speaking, is simply a larger PCB which takes a larger
number of components as inputs and outputs a larger circuit. Finally,
it doesn't matter if several PCB are connected together and then the
chips are placed on them, or if the chips are placed first and the PCBs
are connected later; the resulting overall circuit remains the same.
\end{image}

We start process of drying (formalizing) this image by defining
``wiring diagrams'', the abstract analogs of printed circuit boards.
Let $\bbN$ denote the set of natural numbers including $0$, and for
$n\in\bbN$ let $\underline{n}$ denote some fixed set with $n$ elements,
say $\{1,2,\dots,n\}$.

\begin{definition} Let $k, n, n_1,\dots,n_k\in\bbN$ be natural
numbers. A ``wiring diagram'' $D$ with inputs $\underline{n_1},\dots
\underline{n_k}$ and outputs $\underline{n}$ is an unoriented
compact 1-manifold whose boundary is $\underline{n}\amalg
\underline{n_1}\amalg\cdots\amalg\underline{n_k}$, regarded
up to homeomorphism. In strictly combinatorial terms, it is
a pairing of the elements of the set $\underline{n}\amalg
\underline{n_1}\amalg\cdots\amalg\underline{n_k}$ along
with a single further natural number that counts closed
circles. If $D_1;\dots;D_m$ are wiring diagrams with inputs
$\underline{n_{11}},\dots,\underline{n_{1k_1}}; \dots;
\underline{n_{m1}},\dots,\underline{n_{mk_m}}$ and outputs
$\underline{n_1};\dots;\underline{n_m}$ and $D$ is a wiring diagram
with inputs $\underline{n_1};\dots;\underline{n_m}$ and outputs
$\underline{n}$, there is an obvious ``composition'' $D(D_1,\dots,D_m)$
(obtained by gluing the corresponding 1-manifolds, and also describable
in completely combinatorial terms) which is a wiring diagram with
inputs $(\underline{n_{ij}})_{1\leq i\leq k_j,1\leq j\leq m}$ and
outputs $\underline{n}$ (note that closed circles may be created in
$D(D_1,\dots,D_m)$ even if none existed in $D$ and in $D_1;\dots;D_m$).
\end{definition}

A circuit algebra is an algebraic structure (in the sense of
Section~\ref{subsec:Projectivization}) whose operations are parametrized by
wiring diagrams. Here's a formal definition:

\begin{definition} A circuit algebra consists of the following data:
\begin{itemize}
\item For every natural number $n\geq 0$ a set (or a $\bbZ$-module) $C_n$ ``of
circuits with $n$ legs''.
\item For any wiring diagram $D$ with inputs $\underline{n_1},\dots
\underline{n_k}$ and outputs $\underline{n}$, an operation (denoted by the
same letter) $D\colon C_{n_1}\times\dots\times C_{n_k}\to C_n$ (or linear
$D\colon C_{n_1}\otimes\dots\otimes C_{n_k}\to C_n$ if we work with
$\bbZ$-modules).
\end{itemize}
We insist that the obvious ``identity'' wiring diagrams with
$\underline{n}$ inputs and $\underline{n}$ outputs act as the identity of
$C_n$, and that the actions of wiring diagrams be compatible in the obvious
sense with the composition operation on wiring diagrams.
\end{definition}

A silly but useful example of a circuit algebra is the circuit algebra
$\glos{\calS}$ of empty circuits, or in our context, of ``skeletons''. The
circuits with $n$ legs for $\calS$ are wiring diagrams with $n$ outputs
and no inputs; namely, they are 1-manifolds with boundary $\underline{n}$
(so $n$ must be even).

More generally one may pick some collection of ``basic components''
(perhaps some logic gates and junctions for electronic circuits as in
Figure~\ref{fig:FlipFlop}) and speak of the ``free circuit algebra''
generated by these components. Even more generally we can speak of
circuit algebras given in terms of ``generators and relations''; in the
case of electronics, our relations may include the likes of De Morgan's
law $\neg(p\vee q)=(\neg p)\wedge(\neg q)$ and the laws governing the
placement of resistors in parallel or in series. We feel there is no need
to present the details here, yet many examples of circuit algebras given
in terms of generators and relations appear in this paper, starting with
the next section. We will use the notation $C=\CA\langle \, G \mid R \, \rangle$
to denote the circuit algebra generated by a collection of elements
$G$ subject to some collection $R$ of relations. 

People familiar with electric circuits know very well that connectors
sometimes come in ``male'' and ``female'' versions, and that you
can't plug a USB cable into a headphone jack and expect your system to
cooperate. Thus one may define ``directed circuit algebras'' in which
the wiring diagrams are oriented, the circuit sets $C_n$ get replaced
by $C_{n_1n_2}$ for ``circuits with $n_1$ incoming wires and $n_2$
outgoing wires'' and only orientation preserving connections are ever
allowed. Likewise there is a ``coloured'' version of everything, in which
the wires may be coloured by the elements of some given set $X$ which may
include among its members the elements ``USB'' and ``audio'' and in which
connections are allowed only if the colour coding is respected. We will
not give formal definitions of directed and/or coloured circuit algebras
here, yet we will allow ourselves to freely use these notions. Likewise
for the obvious analogues of the skeletons algebra $\calS$ and for
algebras given in terms of generators and relations.

Note that there is an obvious notion of ``a morphism between
two circuit algebras'' and that circuit algebras (directed or not,
coloured or not) form a category. We feel that a precise definition
is not needed. Yet a lovely example is the ``implementation morphism''
of logic circuits in the style of Figure~\ref{fig:FlipFlop} into more
basic circuits made of transistors and resistors.

Perhaps the prime mathematical example of a circuit algebra is tensor
algebra. If $t_1$ is an element (a ``circuit'') in some tensor product of
vector spaces and their duals, and $t_2$ is the same except in a possibly
different tensor product of vector spaces and their duals, then once an
appropriate pairing $D$ (a ``wiring diagram'') of the relevant vector
spaces is chosen, $t_1$ and $t_2$ can be contracted (``wired together'')
to make a new tensor $D(t_1,t_2)$. The pairing $D$ must pair a vector
space with its own dual, and so this circuit algebra is coloured by the
set of vector spaces involved, and directed, by declaring (say) that
some vector spaces are of one gender and their duals are of the other. We
have in fact encountered this circuit algebra already, in
Section~\ref{subsec:LieAlgebras}.

Let $G$ be a group. A $G$-graded algebra $A$ is a collection $\{A_g\colon g\in
G\}$ of vector spaces, along with products $A_g\otimes A_h\to A_{gh}$ that
induce an overall structure of an algebra on $A:=\bigoplus_{g\in G}A_g$. In
a similar vein, we define the notion of an $\calS$-graded circuit algebra:

\begin{definition} An $\calS$-graded circuit algebra, 
or a ``circuit algebra with skeletons'', is an algebraic structure $C$ with
spaces $C_\beta$, one for each element $\beta$ of the circuit algebra of
skeletons $\calS$, along with composition operations
$D_{\beta_1,\dots,\beta_k}\colon C_{\beta_1}\times\dots\times C_{\beta_k}\to
C_\beta$, defined whenever $D$ is a wiring diagram and
$\beta=D(\beta_1,\dots,\beta_k)$, so that with the obvious induced
structure, $\coprod_\beta C_\beta$ is a circuit algebra. A similar
definition can be made if/when the skeletons are
taken to be directed or coloured.
\end{definition}

Loosely speaking, a circuit algebra with skeletons is a circuit
algebra in which every element $T$ has a well-defined skeleton
$\varsigma(T)\in\calS$. Yet note that as an algebraic structure a circuit
algebra with skeletons has more ``spaces'' than an ordinary circuit
algebra, for its spaces are enumerated by skeleta and not merely by
integers. The prime examples for circuit algebras with skeletons appear in
the next section.

\clearpage
\draftcut
\section{w-Tangles} \label{sec:w-tangles}

\begin{quote} \small {\bf Section Summary. }
  \summarytangles
\end{quote}

\subsection{v-Tangles and w-Tangles} \label{subsec:vw-tangles} With
\wClipStart{120510}{0-01-11}
The (surprisingly pleasant) task of defining circuit algebras completed
in Section~\ref{subsec:CircuitAlgebras}, the definition of v-tangles
and w-tangles is simple.

\begin{definition} The ($\calS$-graded) circuit algebra $\glos{\vT}$ of
v-tangles is the $\calS$-graded directed circuit algebra generated by
two generators in $C_{2,2}$ called the ``positive crossing'' and the
``negative crossing'', modulo the usual \Rs, R2 and R3 moves as depicted in
Figure~\ref{fig:VKnotRels} (these relations clearly make sense as circuit
algebra relations between our two generators), with the obvious meaning
for their skeleta. The circuit algebra $\wT$ of w-tangles is the same,
except we also mod out by the OC relation of Figure~\ref{fig:VKnotRels}
(note that each side in that relation involves only two generators,
with the apparent third crossing being merely a projection artifact).
In fewer words, $\vT:=$\raisebox{-1.5mm}{\input{figs/vTDef.pstex_t}}, and
$\glos{\wT}:=$\raisebox{-1.8mm}{\begin{picture}(0,0)%
\includegraphics{figs/wTDef.pstex}%
\end{picture}%
%
%
\setlength{\unitlength}{3158sp}%
\begingroup\makeatletter\ifx\SetFigFont\undefined%
\gdef\SetFigFont#1#2#3#4#5{%
  \reset@font\fontsize{#1}{#2pt}%
  \fontfamily{#3}\fontseries{#4}\fontshape{#5}%
  \selectfont}%
\fi\endgroup%
\begin{picture}(1121,271)(2011,-1670)
\put(2708,-1577){\makebox(0,0)[lb]{\smash{{\SetFigFont{8}{9.6}{\rmdefault}{\mddefault}{\updefault}{\color[rgb]{0,0,0}$=$}%
}}}}
\put(2026,-1561){\makebox(0,0)[lb]{\smash{{\SetFigFont{10}{12.0}{\rmdefault}{\mddefault}{\updefault}{\color[rgb]{0,0,0}$\vT$}%
}}}}
\end{picture}%
}.
\end{definition}
\wClipEnd{120502}

\parpic[r]{\raisebox{-17mm}{\begin{picture}(0,0)%
\includegraphics{figs/TangleExample.pstex}%
\end{picture}%
%
%
\setlength{\unitlength}{3158sp}%
\begingroup\makeatletter\ifx\SetFigFont\undefined%
\gdef\SetFigFont#1#2#3#4#5{%
  \reset@font\fontsize{#1}{#2pt}%
  \fontfamily{#3}\fontseries{#4}\fontshape{#5}%
  \selectfont}%
\fi\endgroup%
\begin{picture}(775,775)(1189,-2174)
\end{picture}%
}}
\begin{remark} One may also define v-tangles and w-tangles using the
language of planar algebras, except then another generator is required
(the ``virtual crossing'') and also a few further relations (VR1--VR3,
M), and some of the operations (non-planar wirings) become less
elegant to define.
\end{remark}

Our next task is to study the projectivizations $\proj\vT$ and $\proj\wT$
of $\vT$ and $\wT$. Again, the language of circuit algebras makes it
exceedingly simple.

\parpic[r]{\raisebox{-8mm}{$\input figs/arrows.pstex_t$}}
\begin{definition} The ($\calS$-graded) circuit algebra
$\glos{\calD^v}=\glos{\calD^w}$ of arrow diagrams is the graded and
$\calS$-graded directed circuit algebra generated by a single degree
1 generator $a$ in $C_{2,2}$ called ``the arrow'' as shown on the
right, with the obvious meaning for its skeleton. There are morphisms
$\pi\colon \calD^v\to\vT$ and $\pi\colon \calD^w\to\wT$ defined by
mapping the arrow to an overcrossing minus a no-crossing. (On the
right some virtual crossings were added to make the skeleta match). Let
$\glos{\calA^v}$ be $\calD^v/6T$, let
$\glos{\calA^w}:=\calA^v/TC=\calD^w/(\aft,TC)$, and let
$\glos{\calA^{sv}}:=\calA^v/RI$ and $\glos{\calA^{sw}}:=\calA^w/RI$
as usual, with RI, $6T$, $\aft$, and $TC$ being the same relations as
in Figures~\ref{fig:ADand6T} and~\ref{fig:TCand4TForKnots} (allowing
skeleta parts that are not explicitly connected to really lie on separate
skeleton components).
\end{definition}

\begin{proposition} The maps $\pi$ above induce surjections
$\pi\colon \calA^{sv}\to\proj\vT$ and $\pi\colon
\calA^{sw}\to\proj\wT$. Hence in the language of
Definition~\ref{def:CanProj}, $\calA^{sv}$ and $\calA^{sw}$ are candidate
projectivizations of $\vT$ and $\wT$.
\end{proposition}

\begin{proof} Proving that $\pi$ is well-defined amounts to checking
directly that the RI and 6T or RI, $\aft$ and TC relations are in the
kernel of $\pi$. (Just like in the finite type theory of virtual knots and
braids.) Thanks to the circuit algebra structure, it is enough to verify
the surjectivity of $\pi$ in degree 1. We leave this as an exercise for
the reader. \qed
\end{proof}

We do not know if $\calA^{sv}$ is indeed the projectivizations of $\vT$ (also
see~\cite{Bar-NatanHalachevaLeungRoukema:v-Dims}). Yet in the w case, the
picture is simple:

\begin{theorem} The assignment $\overcrossing\mapsto e^a$ (with $e^a$
denoting the exponential of a single arrow from the over strand to the
under strand) extends to a well defined $Z\colon \wT\to\calA^{sw}$. The
resulting map $Z$ is a homomorphic $\calA^{sw}$-expansion, and in particular,
$\calA^{sw}\cong\proj\wT$ and $Z$ is a homomorphic expansion.
\end{theorem}

\begin{proof} There is nothing new here. $Z$ is satisfies the Reidemeister
moves for the same reasons as in Theorem~\ref{thm:RInvariance} and
Theorem~\ref{thm:ExpansionForKnots} and as there it also
satisfies the universality property. The rest follows from
Proposition~\ref{prop:CanProj}. \qed
\end{proof}

In a similar spirit to Definition~\ref{def:wJac}, one may define a
``w-Jacobi diagram'' (often shorts to ``arrow diagram'') on an arbitrary
skeleton. Denote the circuit algebra of formal linear combinations of arrow
diagrams modulo $\aSTU_1$, $\aSTU_2$, and TC relations by $\calA^{wt}$. We
have the following bracket-rise theorem:

\begin{theorem} The obvious inclusion of diagrams induces a circuit
algebra isomorphism $\calA^w\cong\calA^{wt}$. Furthermore, the $\aAS$
and $\aIHX$ relations of Figure~\ref{fig:aIHX} hold in $\calA^{wt}$.
Similarly, $\calA^{sw}\cong\calA^{swt}$, with the expected definition for
$\calA^{swt}$.
\end{theorem}

\begin{proof} The proof of Theorem~\ref{thm:BracketRise} can be repeated
verbatim. Note that that proof does not make use of the connectivity of the
skeleton. \qed
\end{proof}

Given the above theorem, we no longer keep the distinction between
$\calA^w$ and $\calA^{wt}$ and between $\calA^{sw}$ and $\calA^{swt}$.

\begin{remark} \label{rem:HeadInvariance}
Note that if $T$ is an arbitrary $w$ tangle, then the equality on the
left side of the figure below always holds, while the one on the right
generally doesn't:
\begin{equation} \label{eq:TangleLassoMove}
  \begin{array}{c}\input{figs/TangleLassoMove.pstex_t}\end{array}
\end{equation}
The
\wClipComment{120510}{0-19-41}{
  shows a direct proof of~\eqref{eq:HeadInvariance}
}
arrow diagram version of this statement is that if $D$ is an arbitrary
arrow diagram in $\calA^w$, then the left side equality in the
figure below always holds (we will sometimes refer to this as the
``head-invariance'' of arrow diagrams), while the right side equality
(``tail-invariance'') generally fails.
\begin{equation} \label{eq:HeadInvariance}
  \begin{array}{c}\input{figs/HeadInvariance.pstex_t}\end{array}
\end{equation}
We leave it to the reader to ascertain that
Equation~\eqref{eq:TangleLassoMove} implies
Equation~\eqref{eq:HeadInvariance}. There is also a direct
proof of Equation~\eqref{eq:HeadInvariance} which we also leave
to the reader, though see an analogous statement and proof in
\cite[Lemma~3.4]{Bar-Natan:NAT}. Finally note that a restricted version of
tail-invariance does hold --- see Section~\ref{subsec:sder}.
\end{remark}

\draftcut
\subsection{$\calA^w(\uparrow_n)$ and the Alekseev-Torossian Spaces}
\label{subsec:ATSpaces}

\begin{definition} Let $\glos{\vT(\uparrow_n)}$ (likewise
$\glos{\wT(\uparrow_n)}$) be the set of v-tangles (w-tangles) whose
skeleton is the disjoint union of $n$ directed lines. Likewise
let $\glos{\calA^v(\uparrow_n)}$ be the part of $\calA^v$ in
which the skeleton is the disjoint union of $n$ directed lines,
with similar definitions for $\glos{\calA^w(\uparrow_n)}$,
$\glos{\calA^{sv}(\uparrow_n)}$, and $\glos{\calA^{sw}(\uparrow_n)}$.
\end{definition}

\begin{figure}
\input{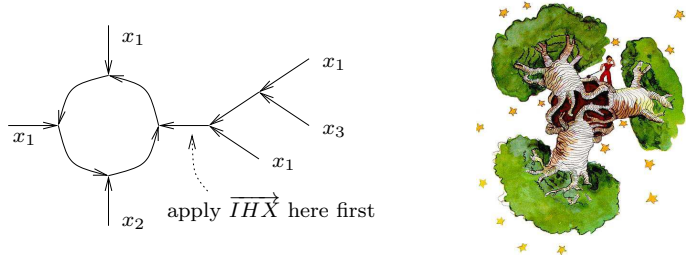}
\caption{A wheel of trees can be reduced to a combination of wheels, and a wheel of trees with 
a Little Prince.}\label{fig:WheelOfTreesAndPrince}
\end{figure}
In the same manner as in the case of knots (Theorem~\ref{thm:Aw}),
$\calA^w(\uparrow_n)$ is a bi-algebra isomorphic (via a diagrammatic
PBW theorem, applied independently on each component of the
skeleton) to a space $\glos{\calB^w_n}$ of unitrivalent diagrams
with symmetrized ends coloured with colours in some $n$-element set
(say $\{x_1,\ldots,x_n\}$), modulo $\aAS$ and $\aIHX$.
Note that the RI relation becomes $w_1=0$, where $w_1$
denotes the 1-wheel of any colour.

The primitives $\glos{\calP^w_n}$ of $\calB^w_n$ are the
connected diagrams (and hence the primitives of $\calA^w(\uparrow_n)$
are the diagrams that remain connected even when the skeleton is
removed). Given the ``two in one out'' rule for internal vertices,
the diagrams in $\calP^w_n$ can only be trees or wheels (``wheels of
trees'' can be reduced to simple wheels by repeatedly using $\aIHX$,
as in Figure~\ref{fig:WheelOfTreesAndPrince}).

Thus as a vector space $\calP^w_n$ is easy to identify. It is a direct sum
$\calP^w_n=\langle\text{trees}\rangle\oplus\langle\text{wheels}\rangle$.
The wheels part is simply the graded vector space generated by
all cyclic words in the letters $x_1,\ldots,x_n$.  Alekseev and
Torossian~\cite{AlekseevTorossian:KashiwaraVergne} denote the
space of cyclic words by $\glos{\attr_n}$, and so shall we. The trees in
$\calP^w_n$ have leafs coloured $x_1,\ldots,x_n$. Modulo $\aAS$ and
$\aIHX$, they correspond to elements of the free Lie algebra $\glos{\lie_n}$
on the generators $x_1,\ldots,x_n$. But the root of each such tree
also carries a label in $\{x_1,\ldots,x_n\}$, hence there are $n$
types of such trees as separated by their roots, and so $\calP^w_n$
is isomorphic to the direct sum $\attr_n\oplus\bigoplus_{i=1}^n\lie_n$.
With $\calB_n^{sw}$ and $\calP_n^{sw}$ defined in the analogous manner,
we can also conclude that
$\calP^{sw}_n\cong\attr_n/(\text{deg }1)\oplus\bigoplus_{i=1}^n\lie_n$.

By the Milnor-Moore theorem~\cite{MilnorMoore:Hopf}, $\calA^w(\uparrow_n)$
is isomorphic to the universal enveloping algebra $\calU(\calP^w_n)$,
with $\calP^w_n$ identified as the subspace $\glos{\calP^w(\uparrow_n)}$
of primitives of $\calA^w(\uparrow_n)$ using the PBW symmetrization
map $\chi\colon \calB^w_n\to\calA^w(\uparrow_n)$. Thus in order to
understand $\calA^w(\uparrow_n)$ as an associative algebra, it is enough
to understand the Lie algebra structure induced on $\calP^w_n$ via the
commutator bracket of $\calA^w(\uparrow_n)$.

We now wish to identify $\calP^w(\uparrow_n)$ as the Lie algebra
$\attr_n\rtimes(\fraka_n\oplus\tder_n)$,
which in itself is a combination of the Lie algebras
$\fraka_n$, $\tder_n$ and $\attr_n$ studied by Alekseev and
Torossian~\cite{AlekseevTorossian:KashiwaraVergne}. Here are the relevant
definitions:

\begin{definition} Let $\glos{\fraka_n}$ denote the vector space with basis
$x_1,\ldots,x_n$, also regarded as an Abelian Lie algebra of dimension $n$.
As before, let $\lie_n=\lie(\fraka_n)$ denote the free Lie algebra on $n$
generators, now identified as the basis elements of $\fraka_n$. Let
$\glos{\der_n}=\der(\lie_n)$ be the (graded) Lie algebra of derivations
acting on $\lie_n$, and let
\[ \glos{\tder_n}=\left\{D\in\der_n\colon \forall i\ \exists a_i\text{ s.t.{}
  }D(x_i)=[x_i,a_i]\right\}
\]
denote the subalgebra of ``tangential derivations''. A tangential
derivation $D$ is determined by the $a_i$'s for which $D(x_i)=[x_i,a_i]$,
and determines them up to the ambiguity $a_i\mapsto a_i+\alpha_ix_i$, where
the $\alpha_i$'s are scalars. Thus as vector spaces,
$\fraka_n\oplus\tder_n\cong\bigoplus_{i=1}^n\lie_n$.
\end{definition}

\begin{definition} Let $\glos{\Ass_n}=\calU(\lie_n)$ be the free associative
algebra ``of words'', and let $\glos{\Ass_n^+}$ be the degree $>0$ part of
$\Ass_n$. As before, we let $\attr_n=\Ass^+_n/(x_{i_1}x_{i_2}\cdots
x_{i_m}=x_{i_2}\cdots x_{i_m}x_{i_1})$ denote ``cyclic words'' or
``(coloured) wheels''.  $\Ass_n$, $\Ass_n^+$, and $\attr_n$ are
$\tder_n$-modules and there is an obvious equivariant ``trace''
$\tr\colon \Ass^+_n\to\attr_n$.
\end{definition}
\wClipEnd{120510}

\begin{proposition}\label{prop:Pnses}
There\wClipStart{120523}{0-00-07}
is a split short exact sequence of Lie algebras 
\[ 0 \longrightarrow \attr_n
  \stackrel{\glos{\iota}}{\longrightarrow} \calP^w(\uparrow_n) 
  \stackrel{\glos{\pi}}{\longrightarrow} \fraka_n \oplus \tder_n 
  \longrightarrow 0.
\]
\end{proposition}

\begin{proof}
The inclusion $\iota$ is defined the natural way: $\attr_n$ is
spanned by coloured ``floating'' wheels, and such a wheel is mapped
into $\calP^w(\uparrow_n)$ by attaching its ends to their assigned strands in
arbitrary order. Note that this is well-defined: wheels have only tails,
and tails commute.

 As vector spaces, the statement is already proven: $\calP^w(\uparrow_n)$ 
is generated by trees 
and wheels (with the all arrow endings fixed on $n$ strands). When factoring out by the wheels,
only trees remain. Trees have one head and many tails. All the tails commute with 
each other, and commuting a tail with a head on a strand costs a wheel (by $\aSTU$), 
thus in the quotient the head also commutes with the tails. Therefore, the quotient
is the space of floating (coloured) trees, which we have previously identified with
$\bigoplus_{i=1}^{n} \lie_n \cong \fraka_n\oplus\tder_n$.

It remains to show that the maps $\iota$ and $\pi$ are Lie algebra maps as well. For $\iota$ this
is easy: the Lie algebra $\attr_n$ is commutative, and is mapped to the commutative
(due to $TC$)
subalgebra of $\calP^w(\uparrow_n)$ generated by wheels.

To show that $\pi$ is a map of Lie algebras we give two proofs,
first a ``hands-on'' one, then a ``conceptual'' one.

{\bf Hands-on argument.} $\fraka_n$ is the image of single arrows on one strand.
These commute with everything in $\calP^w(\uparrow_n)$, and so does $\fraka_n$
in the direct sum $\fraka_n \oplus \tder_n$.

It remains to show that the bracket of $\tder_n$ works the same way as
commuting trees in $\calP^w(\uparrow_n)$. Let $D$ and $D'$ be elements of
$\tder_n$ represented by $(a_1,\ldots ,a_n)$ and $(a_1',\ldots ,a_n')$, meaning
that $D(x_i)=[x_i,a_i]$ and $D'(x_i)=[x_i,a_i']$ for $i=1,\ldots ,n$. Let
us compute the commutator of these elements:
\begin{multline*}
  [D,D'](x_i)=(DD'-D'D)(x_i)=D[x_i,a_i']-D'[x_i,a_i]= \\
  =[[x_i,a_i],a_i']+[x_i,Da_i']-[[x_i,a_i'],a_i]-[x_i,D'a_i]
    = [x_i,Da_i'-D'a_i+[a_i,a_i']].
\end{multline*}

Now let $T$ and $T'$ be two trees in $\calP^w(\uparrow_n)/\attr_n$,
their heads on strands $i$ and $j$, respectively ($i$ may or may not
equal $j$).  Let us denote by $a_i$ (resp. $a_j'$) the element in $\lie_n$ given by forming
the appropriate commutator of the colours of the tails of $T$'s (resp. $T'$). 
In $\tder_n$, let $D=\pi(T)$ and
$D'=\pi(T')$.  $D$ and $D'$ are determined by $(0,\ldots,a_i,\ldots,0)$,
and $(0,\ldots,a_j',\ldots0)$, respectively. (In each case, the $i$-th or
the $j$-th is the only non-zero component.) The commutator of these
elements is given by $[D,D'](x_i)=[Da_i'-D'a_i+[a_i,a_i'],x_i]$, and
$[D,D'](x_j)=[Da_j'-D'a_j+[a_j,a_j'],x_j].$ Note that unless $i=j$,
$a_j=a_i'=0$.

In $\calP^w(\uparrow_n)/\attr_n$, all tails commute, as well as a head of a tree with its
own tails. Therefore, commuting two trees only incurs a cost when commuting a head of
one tree over the tails of the other on the same strand, and the two heads over each other,
if they are on the same strand.

If $i \neq j$, then commuting the head of $T$ over the tails of $T'$ by $\aSTU$ 
costs a sum of trees given by $Da_j'$, with heads on strand $j$, while moving
the head of $T'$ over the tails of $T$ costs exactly $-D'a_i$, with heads on strand $i$,
as needed.

If $i=j$, then everything happens on strand $i$, and the cost is 
$(Da_i'-D'a_i+[a_i,a_i'])$, where the last term happens when
the two heads cross each other.

{\bf Conceptual argument.}
There is an action of $\calP^w(\uparrow_n)$ on $\lie_n$, as follows: introduce
and extra strand on the right. An element $L$ of $\lie_n$ corresponds to a tree with 
its head on the extra strand. Its commutator with an element of $\calP^w(\uparrow_n)$ 
(considered as an element of $\calP^w(\uparrow_{n+1})$ by the obvious inclusion)
is again a tree with head on strand $(n+1)$, defined to be the result of the action.

Since $L$ has only tails on the first $n$ strands, 
elements of $\attr_n$, which
also only have tails, act trivially. So do single (local) arrows on one strand
($\fraka_n$). It remains to show that trees act as $\tder_n$, and it is enough
to check this on the generators of $\lie_n$ (as the Leibniz rule is obviously
satisfied). The generators of $\lie_n$ are arrows pointing from one of the first 
$n$ strands, say strand $i$, to strand $(n+1)$. A tree with head on strand $i$
acts on this element, according $\aSTU$, by forming the commutator, which
is exactly the action of $\tder_n$.
\end{proof}

To identify $\calP^w(\uparrow_n)$ as the semidirect product
$\attr_n\rtimes(\fraka_n\oplus\tder_n)$, it remains to show that
the short exact sequence of the Proposition splits. This is indeed the case,
although not canonically.  Two ---of the many--- splitting maps
$\glos{u},\glos{l}\colon \tder_n\oplus\fraka_n \to \calP^w(\uparrow_n)$
are described as follows: $\tder_n\oplus\fraka_n$ is identified with
$\bigoplus_{i=1}^n\lie_n$, which in turn is identified with floating
(coloured) trees. A map to $\calP^w(\uparrow_n)$ can
be given by specifying how to place the legs on their specified strands.
A tree may have many tails but has only one head, and due to $TC$, only
the positioning of the head matters. Let $u$ (for {\it upper}) be the map
placing the head of each tree above all its tails on the same strand,
while $l$ (for {\it lower}) places the head below all the tails. It is
obvious that these are both Lie algebra maps and that $\pi \circ u$ and
$\pi \circ l$ are both the identity of $\tder_n \oplus \fraka_n$. This
makes $\calP^w(\uparrow_n)$ a semidirect product. \qed

\begin{remark} Let $\glos{\attr_n^s}$ denote $\attr_n$ mod out by its
degree one part (one-wheels). Since the RI relation is in the kernel of
$\pi$, there is a similar split exact sequence
\[ 0\to \attr_n^s \stackrel{\overline{\iota}}{\rightarrow} \calP^{sw}
  \stackrel{\overline{\pi}}{\rightarrow} \fraka_n \oplus \tder_n.
\]
\end{remark}

\begin{definition}\label{div} 
For any $D \in \tder_n$, $(l-u)D$ is in the kernel of $\pi$, therefore
is in the image of $\iota$, so $\iota^{-1}(l-u)D$ makes sense. We call
this element $\glos{\divop}D$.
\end{definition}

\begin{definition}
In \cite{AlekseevTorossian:KashiwaraVergne} 
div is defined as follows: div$(a_1,\ldots,a_n):=\sum_{k=1}^n \tr((\partial_k a_k)x_k)$,
where $\partial_k$ picks out the words of a sum which end in $x_k$ and deletes their last letter
$x_k$, and deletes all other words (the ones which do not end in $x_k$).
\end{definition}

\begin{proposition}
The div of Definition \ref{div} and the div of \cite{AlekseevTorossian:KashiwaraVergne} are 
the same.  
\end{proposition}

\parpic[r]{\input{figs/combtree.pstex_t}}
{\it Proof.}
It is enough to verify the claim for the linear generators of $\tder_n$, namely, elements
of the form $(0,\ldots,a_j,\ldots,0)$, where $a_j \in \lie_n$ or equivalently, single (floating, 
coloured) trees, where the colour of
the head is $j$. By the Jacobi identity, each $a_j$ can be written 
in a form $a_j=[x_{i_1},[x_{i_2},[\ldots,x_{i_k}]\ldots]$. 
Equivalently, by $\aIHX$, each tree has a 
standard ``comb'' form, as shown on the picture on the right.

For an associative word $Y=y_1y_2\ldots y_l \in \Ass_n^+$, 
we introduce the notation $[Y]:=[y_1,[y_2,[\ldots,y_l]\ldots]$.
The div of \cite{AlekseevTorossian:KashiwaraVergne} picks out the
words that end in $x_j$, forgets the rest, and considers these as
cyclic words. Therefore, by interpreting the Lie brackets as commutators,
one can easily check that for $a_j$ written as above,
\begin{equation}\label{divformula}
{\rm div}((0,\ldots,a_j,\ldots,0))=\sum_{\alpha\colon  i_{\alpha}=x_j} 
-x_{i_1}\ldots x_{i_{\alpha-1}}[x_{i_{\alpha+1}}\ldots x_{i_k}]x_j.
\end{equation}

\parpic[r]{\input{figs/divproof.pstex_t}}
In Definition \ref{div}, div of a tree is the difference between attaching its
head on the appropriate strand (here, strand $j$) below all of its tails and above.
As shown in the figure on the right, moving the head across each of the tails on 
strand $j$ requires an $\aSTU$ relation,
which ``costs'' a wheel (of trees, which is equivalent to a sum of honest wheels). 
Namely, the head gets connected to the tail in question.
So div of the tree represented by $a_j$ is given by
\begin{center}
$\sum_{\alpha\colon  x_{i_{\alpha}}=j}$``\rm connect the head to the $\alpha$ leaf''.
\end{center}

\noindent
This in turn gets mapped to the formula above via the correspondence between 
wheels and cyclic words. \qed
\wClipEnd{120523}
\wClipComment{120530}{0-00-07}{has extra material on the relationship
of all this with differential operators}

\parpic[r]{\input{figs/treeactonwheel.pstex_t}}
\begin{remark}\label{rem:tderontr}
There is an action of $\tder_n$ on $\attr_n$ as
follows. Represent a cyclic word $w \in \attr_n$ as a
wheel in $\calP^w(\uparrow_n)$ via the map $\iota$. Given
an element $D \in \tder_n$, $u(D)$, as defined above, is a tree
in $\calP^w(\uparrow_n)$ whose head is above all of its tails. We
define $D \cdot w:=\iota^{-1}(u(D)\iota(w)-\iota(w)u(D))$. Note that
$u(D)\iota(w)-\iota(w)u(D)$ is in the image of $\iota$, i.e., a linear
combination of wheels, for the following reason. The wheel $\iota(w)$ has only tails. As we commute
the tree $u(D)$ across the wheel, the head of the tree is commuted
across tails of the wheel on the same strand. Each time this happens
the cost, by the $\aSTU$ relation, is a wheel with the tree attached
to it, as shown on the right, which in turn (by $\aIHX$ relations,
as Figure~\ref{fig:WheelOfTreesAndPrince} shows) is a sum of wheels.
Once the head of the tree has been moved to the top, the tails of the
tree commute up for free by $TC$. Note that the alternative definition,
$D \cdot w:=\iota^{-1}(l(D)\iota(w)-\iota(w)l(D))$ is in fact equal to
the definition above.
\end{remark}

\begin{definition}
In \cite{AlekseevTorossian:KashiwaraVergne}, the group $\glos{\TAut_n}$
is defined as $\exp(\tder_n)$. Note that $\tder_n$ is positively
graded, hence it integrates to a group. Note also that $\TAut_n$ is
the group of ``basis-conjugating'' automorphisms of $\lie_n$, i.e.,
for $g \in \TAut_n$, and any $x_i$, $i=1,\ldots ,n$ generator of
$\lie_n$, there exists an element $g_i \in \exp(\lie_n)$ such that
$g(x_i)=g_i^{-1}x_ig_i$.
\end{definition}

The action of $\tder_n$ on $\attr_n$ lifts to an action of $\TAut_n$ on $\attr_n$,
by interpreting exponentials formally, in other words $e^D$ acts as 
$\sum_{n=0}^\infty\frac{D^n}{n!}$. The lifted action is by conjugation:
for $w \in \attr_n$ and $e^D \in \TAut_n$, 
$e^D \cdot w=\iota^{-1}(e^{uD} \iota(w) e^{-uD})$.

Recall that in Section 5.1 of \cite{AlekseevTorossian:KashiwaraVergne}
Alekseev and Torossian construct a map $\glos{j}\colon \TAut_n \to
\attr_n$ which is characterized by two properties: the cocycle property
\begin{equation}\label{eq:jcocycle}
 j(gh)=j(g)+g\cdot j(h),
\end{equation}
where in the second term multiplication by $g$ denotes the action described above;
and the condition
\begin{equation}\label{eq:jderiv}
\frac{d}{ds}j(\exp(sD))|_{s=0}=\divop(D). 
\end{equation}

Now let us interpret $j$ in our context.
\begin{definition}\label{def:Adjoint}
The adjoint map $\glos{*}\colon \calA^w(\uparrow_n) \to
\calA^w(\uparrow_n)$ acts by ``flipping over diagrams and negating arrow
heads on the skeleton''. In other words, for an arrow diagram $D$,
\[ D^*:=(-1)^{\#\{\text{tails on skeleton}\}}S(D), \]
where $S$ denotes the map which switches the orientation of the skeleton
strands (i.e. flips the diagram over), and multiplies by $(-1)^{\#
\text{skeleton vertices}}$.
\end{definition}

\begin{proposition}\label{prop:Jandj}For $D \in \tder_n$,
define a map $\glos{J}\colon \TAut_n \to \exp(\attr_n)$ by
$J(e^D):=e^{uD}(e^{uD})^*$. Then
$$\exp(j(e^D))=J(e^D).$$
\end{proposition}

\begin{proof}
Note that $(e^{uD})^*=e^{-lD}$, due to ``Tails Commute'' and the fact that a 
tree has only one head.

Let us check that $\log J$ satisfies properties \eqref{eq:jcocycle} and 
\eqref{eq:jderiv}. Namely, with $g=e^{D_1}$ and $h=e^{D_2}$, and 
using that $\attr_n$ is commutative, we need to show that
\begin{equation}
 J(e^{D_1}e^{D_2})=J(e^{D_1})\big(e^{uD_1}\cdot J(e^{D_2})\big),
\end{equation}
where $\cdot$ denotes the action of $\tder_n$ on $\attr_n$; and that
\begin{equation}
 \frac{d}{ds}J(e^{sD})|_{s=0}=\divop D.
\end{equation}

Indeed, with $\operatorname{BCH}(D_1,D_2)=\log e^{D_1}e^{D_2}$ being the 
standard Baker--Campbell--Hausdorff formula,
\begin{multline*}
  J(e^{D_1}e^{D_2})=J(e^{\operatorname{BCH}(D_1,D_2)})
  =e^{u(\operatorname{BCH}(D_1,D_2)}
  e^{-l(\operatorname{BCH}(D_1,D_2)}=
  e^{\operatorname{BCH}(uD_1,uD_2)}
  e^{-\operatorname{BCH}(lD_1,lD_2)} \\
  =e^{uD_1}e^{uD_2}e^{-lD_2}e^{-lD_1}=
  e^{uD_1}(e^{uD_2}e^{-lD_2})e^{-uD_1}e^{uD_1}e^{lD_1}
  =(e^{uD_1}\cdot J(D_2))J(D_1),
\end{multline*}
as needed.

As for condition~\eqref{eq:jderiv}, a direct computation of the derivative
yields
$$\frac{d}{ds}J(e^{sD})|_{s=0}=uD-lD=\divop D,$$
as desired. \qed
\end{proof}

\draftcut
\subsection{The Relationship with u-Tangles} \label{subsec:sder} Let
$\glos{\uT}$ be the planar algebra of classical, or ``{\it u}sual''
tangles.  There is a map $a\colon \uT \to \wT$ of $u$-tangles into
$w$-tangles: algebraically, it is defined in the obvious way on the planar
algebra generators of $\uT$. (It can also be interpreted topologically
as Satoh's tubing map, as in Section~\ref{subsubsec:TopTube},
where a u-tangle is a tangle drawn on a sphere. However, it is only
conjectured that the circuit algebra presented here is a Reidemeister
theory for ``tangled ribbon tubes in $\bbR^4$''.)  The map $a$ induces a
corresponding map $\alpha\colon  \calA^u \to \calA^{sw}$, which maps an
ordinary Jacobi diagram (i.e., unoriented chords with internal trivalent
vertices modulo the usual $AS$, $IHX$ and $STU$ relations) to the sum
of all possible orientations of its chords (many of which are zero in
$\calA^{sw}$ due to the ``two in one out'' rule).

\parpic[l]{$\xymatrix{
  \uT \ar@{.>}[r]^{Z^u} \ar[d]^a & \calA^u \ar[d]^\alpha \\
  \wT \ar[r]^{Z^w} & \calA^{sw}
}$}
It is tempting to ask whether the square on the left
commutes. Unfortunately, this question hardly makes sense, as there
is no canonical choice for the dotted line in it. Similarly to the
braid case in Section~\ref{subsubsec:RelWithu}, the definition of the
Kontsevich integral for $u$-tangles typically depends on various choices
of ``parenthesizations''. Choosing parenthesizations, this square becomes
commutative up to some fixed corrections. The details are in
Proposition~\ref{prop:uwBT}.

Yet already at this point we can recover something from the existence of
the map $a\colon\uT\to\wT$, namely an interpretation of the
Alekseev-Torossian~\cite{AlekseevTorossian:KashiwaraVergne} space of
special derivations, $$\glos{\sder_n}:=\{ D\in\tder_n\colon D(\sum_{i=1}^n
x_i)=0\}.$$ Recall from Remark \ref{rem:HeadInvariance} that
in general it is not possible to slide a strand under an arbitrary $w$-tangle.
However, it is possible to slide strands freely under
tangles {\em in the image of $a$}, and thus by reasoning similar to the
reasoning in Remark~\ref{rem:HeadInvariance}, diagrams $D$ in the image
of $\alpha$ respect ``tail-invariance'':
\begin{equation} \label{eq:TailInvariance}
  \begin{array}{c}\input{figs/TailInvariance.pstex_t}\end{array}
\end{equation}

Let $\calP^u(\uparrow_n)$ denote the primitives of $\calA^u(\uparrow_n)$,
that is, Jacobi diagrams that remain connected when the skeleton is
removed. Remember that $\calP^{w}(\uparrow_n)$ stands for the primitives
of $\calA^{w}(\uparrow_n)$. Equation~\eqref{eq:TailInvariance} readily
implies that the image of the composition
\[ \xymatrix{
  \calP^u(\uparrow_n) \ar[r]^(0.48){\alpha}
  & \calP^w(\uparrow_n) \ar[r]^(0.45)\pi
  & \fraka_n \oplus \tder_n
} \]
is contained in $\fraka_n \oplus \sder_n$. Even better is true.

\begin{theorem}\label{thm:sder}
The image of $\pi\alpha$ is precisely $\fraka_n \oplus \sder_n$. 
\end{theorem}

This theorem was first proven by Drinfel'd (Lemma after Proposition 6.1
in \cite{Drinfeld:GalQQ}), but the proof we give here is due to Levine
\cite{Levine:Addendum}.

\begin{proof}
Let $\lie_n^d$ denote the degree $d$ piece of $\lie_n$. Let $V_n$ be
the vector space with basis $x_1, x_2, \ldots , x_n$.  Note that
$$V_n \otimes \lie_n^d \cong \bigoplus_{i=1}^n \lie_n^d \cong
(\tder_n \oplus \fraka_n)^d,$$
where $\tder_n$ is graded by the number of tails of a tree, and $\fraka_n$ 
is contained in degree 1.  

The bracket defines a map $\beta\colon  V_n \otimes \lie_n^d \to \lie_n^{d+1}$:
for $a_i \in \lie_n^d$ where $i=1,\ldots ,n$, the ``tree'' 
$D=(a_1,a_2,\ldots ,a_n) \in (\tder_n \oplus \fraka_n)^d$ is mapped to 
$$\beta(D)=\sum_{i=1}^n[x_i,a_i]=D\left(\sum_{i=1}^n x_i\right),$$
where the first equality is by the definition of tensor product and the bracket,
and the second is by the definition of the action of $\tder_n$ on $\lie_n$.

Since $\fraka_n$ is contained in degree 1, by definition 
$\sder_n^d=(\operatorname{ker}\beta)^d$ for $d\geq2$. In degree 
1, $\fraka_n$ is obviously in the kernel, hence 
$(\operatorname{ker}\beta)^1= \fraka_n \oplus \sder_n^1$. So overall,
$\operatorname{ker}\beta=\fraka_n\oplus\sder_n$.

We want to study the image of the map $\calP^u(\uparrow^n)
\stackrel{\pi\alpha}{\longrightarrow} \fraka_n \oplus \tder_n$.
Under $\alpha$, all connected Jacobi diagrams that are not trees or
wheels go to zero, and under $\pi$ so do all wheels. Furthermore, $\pi$
maps trees that live on $n$ strands to ``floating'' trees with univalent
vertices coloured by the strand they used to end on. So for determining
the image, we may replace $\calP^u(\uparrow^n)$ by the space $\calT_n$
of connected {\em un}oriented ``floating trees'' (uni-trivalent graphs), the ends (univalent vertices)
of which are coloured by the $\{x_i\}_{i=1,..,n}$. We denote the degree
$d$ piece of $\calT_n$, i.e., the space of trees with $d+1$ ends,
by $\calT_n^{d}$. Abusing notation, we shall denote the map induced by
$\pi\alpha$ on $\calT_n$ by $\alpha\colon  \calT_n \to \fraka_n \oplus
\tder_n$. Since choosing a ``head'' determines the entire orientation of
a tree by the two-in-one-out rule, $\alpha$ maps a tree in $\calT_n^d$
to the sum of $d+1$ ways of choosing one of the ends to be the ``head''.

We want to show that $\operatorname{ker}\beta=\operatorname{im}\alpha$.
This is equivalent to saying that $\bar{\beta}$ is injective, where
$\bar{\beta}\colon V_n\otimes\lie_n/\operatorname{im}\alpha
\to \lie_n$ is map induced by $\beta$ on the quotient by
$\operatorname{im}\alpha$.

\parpic[r]{\input{figs/beta.pstex_t}}
The degree $d$ piece of $V_n \otimes \lie_n$, in the pictorial
description, is generated by floating trees with $d$ tails and one head,
all coloured by $x_i$, $i=1,\ldots ,n$. This is mapped to $\lie_n^{d+1}$,
which is isomorphic to the space of floating trees with $d+1$ tails and
one head, where only the tails are coloured by the $x_i$. The map $\beta$
acts as shown on the picture on the right.

\parpic[r]{\begin{picture}(0,0)%
\includegraphics{figs/taudef.pstex}%
\end{picture}%
%
%
\setlength{\unitlength}{3158sp}%
\begingroup\makeatletter\ifx\SetFigFont\undefined%
\gdef\SetFigFont#1#2#3#4#5{%
  \reset@font\fontsize{#1}{#2pt}%
  \fontfamily{#3}\fontseries{#4}\fontshape{#5}%
  \selectfont}%
\fi\endgroup%
\begin{picture}(2799,849)(1489,-448)
\put(3376,-361){\makebox(0,0)[b]{\smash{{\SetFigFont{10}{12.0}{\rmdefault}{\mddefault}{\updefault}{\color[rgb]{0,0,0}$+$}%
}}}}
\put(2626,239){\makebox(0,0)[b]{\smash{{\SetFigFont{10}{12.0}{\rmdefault}{\mddefault}{\updefault}{\color[rgb]{0,0,0}$\tau$}%
}}}}
\end{picture}%
}
We show that $\bar{\beta}$ is injective by exhibiting a map $\tau\colon
\lie_n^{d+1} \to V_n\otimes\lie_n^d/\operatorname{im}\alpha$ so that
$\tau\bar{\beta}=I$. $\tau$ is defined as follows: given a tree with
one head and $d+1$ tails $\tau$ acts by deleting the head and the
arc connecting it to the rest of the tree and summing over all ways of
choosing a new head from one of the tails on the left half of the tree relative to the
original placement of the head (see the
picture on the right). As long as we show that $\tau$ is well-defined,
it follows from the definition and the pictorial description of $\beta$
that $\tau\bar{\beta}=I$.

For well-definedness we need to check that the images of $\aAS$ and
$\aIHX$ relations under $\tau$ are in the image of $\alpha$. This we do 
in the picture below. In both cases it is enough to check the
case when the ``head'' of the relation is the head of the tree
itself, as otherwise an $\aAS$ or $\aIHX$ relation in the domain is mapped
to an $\aAS$ or $\aIHX$ relation, thus zero, in the image.
\[ \input figs/tauproof.pstex_t \]
\[ \input figs/tauproof2.pstex_t \]
In the $\aIHX$ picture, in higher degrees $A$, $B$ and $C$ may denote
an entire tree. In this case, the arrow at $A$ (for example) means the
sum of all head choices from the tree $A$.
\qed
\end{proof}

\begin{comment} In view of the relation between the right half of
Equation~\eqref{eq:TailInvariance} and the special derivations $\sder$,
it makes sense to call w-tangles that satisfy the condition in the left
half of Equation~\eqref{eq:TailInvariance} ``special''. The $a$ images
of u-tangles are thus special. We do not know if the global version of
Theorem~\ref{thm:sder} holds true. Namely, we do not know whether every
special w-tangle is the $a$-image of a u-tangle.
\end{comment}

\draftcut
\subsection{The local topology of w-tangles}\label{subsec:TangleTopology}
So far throughout this section we have presented $w$-tangles as a Reidemeister theory: 
a circuit algebra given by generators and relations. Note that Satoh's tubing map (see Sections \ref{subsubsec:ribbon} and \ref{subsubsec:TopTube}) 
does extend to w-tangles in the obvious way, although it is not known whether it is an isomorphism between
the circuit algebra described here and tangled tubes in $\bbR^4$.
Nonetheless, this intuition explains the local relations (Reidemeister moves). The purpose of this subsection is
to explain the local topology of crossings and understand orientations, signs and orientation
reversals.

\parpic[r]{\begin{picture}(0,0)%
\includegraphics{figs/TubeOrientation.pstex}%
\end{picture}%
%
%
\setlength{\unitlength}{4934sp}%
\begingroup\makeatletter\ifx\SetFigFont\undefined%
\gdef\SetFigFont#1#2#3#4#5{%
  \reset@font\fontsize{#1}{#2pt}%
  \fontfamily{#3}\fontseries{#4}\fontshape{#5}%
  \selectfont}%
\fi\endgroup%
\begin{picture}(557,1046)(3731,-217)
\put(3810,433){\makebox(0,0)[lb]{\smash{{\SetFigFont{11}{13.2}{\rmdefault}{\mddefault}{\updefault}{\color[rgb]{0,0,0}1D}%
}}}}
\put(4079,185){\makebox(0,0)[lb]{\smash{{\SetFigFont{11}{13.2}{\rmdefault}{\mddefault}{\updefault}{\color[rgb]{0,0,0}2D}%
}}}}
\end{picture}%
}
The tubes we consider are endowed with two orientations, we will call these the 1- and 2-dimensional orientations. The one
dimensional orientation is the direction of the tube as a ``strand'' of the tangle. In other 
words, each tube has a ``core''\footnote{The core of Lord Voldemort's wand was made of a phoenix feather.}: 
a distinguished line along the tube,  
which is oriented as a 1-dimensional manifold. Furthermore, the tube as a 
2-dimensional surface is oriented as given by the tubing map. An example is shown on the right.

Note that a tube in $\bbR^4$ has a ``filling'': 
a solid (3-dimensional) cylinder
embedded in $\bbR^4$, with boundary the tube, and the 2D orientation of the tube induces an orientation
of its filling as a 3-dimensional manifold. A (non-virtual) crossing is when the core of one tube intersects the
filling of another transversely. Due to the complementary dimensions, the intersection is a single point,
and the 1D orientation of the core along with the 3D orientation of the filling it passes through determines
an orientation of the ambient space. We say that the crossing is positive if this agrees with the standard orientation
of $\bbR^4$, and negative otherwise. Hence, there are four types of crossings, given by whether the core
of tube A intersects the filling of B or vice versa, and two possible signs in each case. 

As discussed in Section \ref{subsec:wBraids}, braided tubes in $\bbR^4$ can be thought
of as movies of flying rings in $\bbR^3$, and in particular a crossing
represents a ring flying through another ring. In this interpretation, the 1D orientation
of the tube is given by time moving forward. The 2D and 1D orientations of the tube together induce an orientation
of the flying ring which is a cross-section of the tube at each moment. Hence, saying ``below'' and ``above'' the ring 
makes sense, and as mentioned in 
Exercise \ref{ex:swBn} there are four types of crossings:
ring A flies through ring B from below or from above; and ring B flies through ring A from below
or from above. A crossing is positive if the inner ring comes from below, and negative otherwise.

\parpic[r]{\input{figs/PushMembranes.pstex_t}}
In Sections \ref{subsubsec:ribbon} and \ref{subsubsec:TopTube} we have discussed the tubing map from v- or w-diagrams
of braids or knots to ribbon tubes in $\bbR^4$: the under-strand of a crossing is interpreted as a thinner tube (or a ring flying through another). 
This generalizes to tangles easily. We take the opportunity here to introduce another notation, 
to be called the ``band notation'', which is more suggestive of the 4D topology than the strand notation. We represent a tube in $\bbR^4$
by a picture of an oriented band in $\bbR^3$.
By ``oriented band'' we mean that it has two orientations: a 1D direction (for example an orientation of one of the edges),
and a 2D orientation as a surface. To interpret the 3D picture
of a band as an tube in $\bbR^4$, we add an extra coordinate. Let us refer to the $\bbR^3$ coordinates as $x, y$ and $t$,
and to the extra coordinate as $z$. Think of $\bbR^3$ as being embedded in $\bbR^4$ as the hyperplane $z=0$, and think of the 
band as being made of a thin double membrane. Push the membrane up and down
in the $z$ direction at each point as far as the distance of that point from the boundary of the band, as shown on the right. 
Furthermore, keep the 2D orientation of the top membrane (the one being pushed up), but reverse it on the bottom. This produces 
an oriented tube embedded in $\bbR^4$.

In band notation, the four possible crossings appear as follows, where underneath each crossing we indicate the corresponding
strand picture, as mentioned in Exercise \ref{ex:swBn}:
\begin{center}
\input{figs/BandCrossings.pstex_t}
\end{center}
The signs for each type of crossing are shown in the figure above. Note that the sign of a crossing depends of the 2D orientation of the
over-strand, as well as the 1D direction of the under-strand. Hence, switching only
the direction (1D orientation) of a strand changes the sign of the crossing if and only if the strand of changing direction is the under
strand. However, fully changing the orientation (both 1D and 2D) always switches the 
sign of the crossing. Note that switching the strand orientation in the strand notation corresponds to the total (both 1D and 2D)
orientation switch.

\draftcut
\subsection{Good properties and uniqueness of the homomorphic expansion}
\label{subsec:UniquenessForTangles}

In much the same way as in Section \ref{subsubsec:BraidCompatibility}, $Z$
has a number of good properties with respect to various tangle operations:
it is group-like; it commutes with adding an inert strand (note that
this is a circuit algebra operation, hence it doesn't add anything beyond
homomorphicity); and it commutes with deleting a strand and with strand
orientation reversals. All but the last of these were explained in the
context of braids and the explanations still hold. Orientation reversal
$\glos{S_k}\colon\wT\to\wT$ is the operation which reverses the
orientation of the $k$-th component. Note that in the world of topology (via Satoh's 
tubing map) this means reversing both the 1D and the 2D orientations.
The induced diagrammatic operation
$S_k\colon  \calA^w(T) \to \calA^w(S_k(T))$, where $T$ denotes the
skeleton of a given w-tangle, acts by multiplying each arrow diagram by
$(-1)$ raised to the power the number of arrow endings (both heads and
tails) on the $k$-th strand, as well as reversing the strand orientation.  
Saying that ``$Z$ commutes with $S_k$'' means that
the appropriate square commutes.

The following theorem asserts that a well-behaved homomorphic expansion of 
$w$-tangles is unique:
\begin{theorem}\label{thm:Tangleuniqueness}
The only homomorphic expansion satisfying the good properties described
above is the $Z$ defined in Section \ref{subsec:vw-tangles}.
\end{theorem}

\parpic[r]{\begin{picture}(0,0)%
\includegraphics{figs/rho.pstex}%
\end{picture}%
%
%
\setlength{\unitlength}{2368sp}%
\begingroup\makeatletter\ifx\SetFigFont\undefined%
\gdef\SetFigFont#1#2#3#4#5{%
  \reset@font\fontsize{#1}{#2pt}%
  \fontfamily{#3}\fontseries{#4}\fontshape{#5}%
  \selectfont}%
\fi\endgroup%
\begin{picture}(2389,944)(214,-83)
\put(676,314){\makebox(0,0)[rb]{\smash{{\SetFigFont{10}{12.0}{\rmdefault}{\mddefault}{\updefault}{\color[rgb]{0,0,0}$\rho=$}%
}}}}
\put(1726,314){\makebox(0,0)[b]{\smash{{\SetFigFont{10}{12.0}{\rmdefault}{\mddefault}{\updefault}{\color[rgb]{0,0,0}$+$}%
}}}}
\end{picture}%
}
\begin{proof}
We first prove the following claim: Assume, by contradiction, that $Z'$ is a different 
homomorphic expansion
of $w$-tangles with the good properties described above. Let $R'=Z'(\overcrossing)$ and
$R=Z(\overcrossing)$, and denote by $\rho$ the lowest degree homogeneous
non-vanishing term of $R'-R$. (Note that $R'$ determines $Z'$, so if $Z'\neq Z$, then
$R' \neq R$.) Suppose $\rho$ is of degree $k$. 
Then we claim that $\rho=\alpha_1 w_k^1+\alpha_2 w_k^2$ is a linear combination of $w_k^1$ and $w_k^2$, 
where $w_k^i$ denotes a $k$-wheel 
living on strand $i$, as shown on the right.

Before proving the claim, note that it leads to a contradiction.
Let $d_i$ denote the operation ``delete strand $i$''.
Then up to degree $k$, we have $d_1(R')=\alpha_2 w_k^1$ and $d_2(R')=\alpha_1 w_k^2$, but
$Z'$ is compatible with strand deletions, so $\alpha_1=\alpha_2=0$. Hence
$Z$ is unique, as stated.

On to the proof of the claim, note that $Z'$ being an expansion determines the degree 1 term of $R'$ 
(namely, the single arrow 
$a^{12}$ from strand 1 to strand 2, with coefficient 1). So we can assume that $k \geq 2$. Note also that since both $R'$ and $R$ are 
group-like, $\rho$ is primitive. Hence $\rho$ is a linear combination of connected diagrams,
namely trees and wheels. 

Both $R$ and $R'$ satisfy the Reidemeister 3 relation:
$$R^{12}R^{13}R^{23}=R^{23}R^{13}R^{12}, \qquad R'^{12}R'^{13}R'^{23}=R'^{23}R'^{13}R'^{12}$$
where the superscripts denote the strands on which $R$ is placed
(compare with Remark \ref{rem:YangBaxter}).
We focus our attention on the degree $k+1$ part of the equation for $R'$,
and use that up to degree $k+1$. We can write $R'=R+\rho+\mu$, where $\mu$ denotes the degree
$k+1$ homogeneous part of $R'-R$. Thus, up to degree $k+1$, we have
$$(R^{12}\!+\!\rho^{12}\!+\!\mu^{12})(R^{13}\!+\!\rho^{13}\!+\!\mu^{13})(R^{23}\!+\!\rho^{23}\!+\!\mu^{23})=
(R^{23}\!+\!\rho^{23}\!+\!\mu^{23})(R^{13}\!+\!\rho^{13}\!+\!\mu^{13})(R^{12}\!+\!\rho^{12}\!+\!\mu^{12}).$$
The homogeneous degree $k+1$ part of this equation is a sum of some terms which contain $\rho$
and some which don't. The diligent reader can check that those which don't involve $\rho$ 
cancel on both sides, either due to the
fact that $R$ satisfies the Reidemeister 3 relation, or by simple degree counting. 
Rearranging all the terms which do involve $\rho$ to the left side, we get the following equation,
where $a^{ij}$ denotes an arrow pointing from strand $i$ to strand $j$:
\begin{equation}\label{eq:Reid3forrho}
[a^{12}, \rho^{13}]+[\rho^{12},a^{13}]+[a^{12},\rho^{23}]+[\rho^{12},a^{23}]+
[a^{13},\rho^{23}]+[\rho^{13},a^{23}]=0. 
\end{equation}

The third and fifth terms sum to $[a^{12}+a^{13},\rho^{23}]$,
which is zero due to the ``head-invariance'' of diagrams, as in Remark
\ref{rem:HeadInvariance}.

We treat the tree and wheel components of $\rho$ separately.
Let us first assume that $\rho$ is a linear combination of trees. Recall that the
space of trees on two strands is isomorphic to $\lie_2 \oplus \lie_2$, the
first component given by trees whose head is on the first strand, and the second 
component by trees with their head on the second strand.
Let $\rho=\rho_1 +\rho_2$, where $\rho_i$ is the projection to the $i$-th component
for $i=1,2$.

Note that due to $TC$, we have $[a^{12}, \rho^{13}_2]=[\rho^{12}_2,a^{13}]=
[\rho^{12}_1,a^{23}]=0$. So Equation (\ref{eq:Reid3forrho}) reduces to
$$[a^{12},\rho^{13}_1]+[\rho^{12}_1,a^{13}]+[\rho^{12}_2,a^{23}]+[\rho^{13}_1,a^{23}]+[\rho^{13}_2,a^{23}]=0$$
The left side of this equation lives in $\bigoplus_{i=1}^3 \lie_3$. Notice that only the
first term lies in the second direct sum component, while the second, third and last terms live in the third one,
and the fourth term lives in the first.
This in particular means that the first term is itself zero. By $\aSTU$, this implies 
$$0=[a^{12},\rho^{13}_1]=-[\rho_1, x_1]^{13}_2,$$
where $[\rho_1, x_1]^{13}_2$ means the tree defined by the element $[\rho_1,x_1] \in \lie_2$,
with its tails on strands 1 and 3, and head on strand 2. Hence, $[\rho_1, x_1]=0$, so $\rho_1$
is a multiple of $x_1$. The tree given by $\rho_1=x_1$ is a degree 1 element, a possibility we have eliminated, so
$\rho_1=0$.

Equation (\ref{eq:Reid3forrho}) is now reduced to  
$$[\rho^{12}_2,a^{23}]+[\rho^{13}_2,a^{23}]=0.$$
Both terms are words in $\lie_3$, but notice that the first term does not involve
the letter $x_3$. This means that if the second term involves $x_3$ at all, i.e., if
$\rho_2$ has tails on the second strand, then both terms have to be zero individually.
Assuming this and looking at the first term, $\rho^{12}_2$ is a Lie word in $x_1$ and $x_2$,
which does involve $x_2$ by assumption. We have
$[\rho^{12}_2,a^{23}]=[x_2, \rho^{12}_2]=0$, which implies $\rho^{12}_2$ is a multiple of $x_2$, in
other words, $\rho$ is a single arrow on the second strand. This is ruled out by the 
assumption that $k \geq 2$.

On the other hand if the second term does not involve $x_3$ at all, then $\rho_2$ has no tails on the second
strand, hence it is of degree 1, but again $k \geq 2$. We have proven that the ``tree part''
of $\rho$ is zero.

So $\rho$ is a linear combination of wheels. 
Wheels have only tails, so the
first, second and fourth terms of (\ref{eq:Reid3forrho}) are zero due to the tails commute relation.
What remains is $[\rho^{13}, a^{23}]=0$. We assert that this is true if and only if each
linear component of $\rho$ has all of its tails on one strand. 

To prove this, recall each wheel of $\rho^{13}$ represents a cyclic word in letters $x_1$ and $x_3$.
The map $r\colon  \rho^{13} \mapsto [\rho^{13}, a^{23}]$ is a map $\attr_2 \to \attr_3$, which sends each
cyclic word in letters $x_1$ and $x_3$ to the sum of all ways of substituting $[x_2,x_3]$ for one 
of the $x_3$'s in the word.
Note that if we expand the commutators, then all terms that have $x_2$
between two $x_3$'s cancel. Hence all remaining terms will be cyclic words in $x_1$ and $x_3$ with
a single occurrence of $x_2$ in between an $x_1$ and an $x_3$. 

We construct an almost-inverse $r'$ to $r$: for a cyclic word $w$ in $\attr_3$ with one occurrence of $x_2$,
let $r'$ be the map that deletes $x_2$ from $w$ and maps it to the resulting word in 
$\attr_2$ if $x_2$ is followed by $x_3$ in $w$, and maps it to 0 otherwise. On the rest of $\attr_3$
the map $r'$ may be defined to be 0.

The composition $r'r$ takes a cyclic word in $x_1$ and $x_3$ to itself multiplied by the number of times
a letter $x_3$ follows a letter $x_1$ in it. The kernel of this map can consist only of cyclic words 
that do not contain the sub-word $x_3x_1$, namely, these are the words of the form $x_3^k$ or $x_1^k$.
Such words are indeed in the kernel of $r$, so these make up exactly the kernel of $r$. This is exactly what 
needed to be proven: all wheels in $\rho$ have all their tails on one strand.

This concludes the proof of the claim, and the proof of the theorem. \qed
\end{proof}

\clearpage\draftcut
\section{w-Tangled Foams} \label{sec:w-foams}

\begin{quote} \small {\bf Section Summary. }
  \summaryfoams
\end{quote}

\subsection{The Circuit Algebra of w-Tangled Foams} \label{subsec:wTFo}
For reasons we will reluctantly acknowledge later in this section (see
Comment~\ref{com:MissingTopology}), we will present the circuit algebra
of w-tangled foams via its Reidemeister-style diagrammatic description
(accompanied by a local topological interpretation) rather than as an
entirely topological construct.

\begin{definition}\label{def:wTFo}
Let $\glos{\wTFo}$ (where $o$ stands for ``orientable'', to be explained
in Section~\ref{subsec:TheWen}) be the algebraic structure
\[
  \wTFo=\CA\!\left.\left.\left\langle
  \raisebox{-2mm}{\input{figs/wTFgens.pstex_t}}
  \right|
  \parbox{1.2in}{\centering w-relations as in Section~\ref{subsubsec:wrels}}
  \right|
  \parbox{1.2in}{\centering w-operations as in Section~\ref{subsubsec:wops}}
  \right\rangle.
\]
Hence $\wTFo$ is the circuit algebra generated by the generators listed
above and described below, modulo the relations described in
Section~\ref{subsubsec:wrels}, and augmented with several ``auxiliary
operations'', which are a part of the algebraic structure of $\wTFo$ but are
not a part of its structure as a circuit algebra, as described in
Section~\ref{subsubsec:wops}.

\parpic[r]{\begin{picture}(0,0)%
\includegraphics{figs/VertexExamples.pstex}%
\end{picture}%
%
%
\setlength{\unitlength}{3947sp}%
\begingroup\makeatletter\ifx\SetFigFont\undefined%
\gdef\SetFigFont#1#2#3#4#5{%
  \reset@font\fontsize{#1}{#2pt}%
  \fontfamily{#3}\fontseries{#4}\fontshape{#5}%
  \selectfont}%
\fi\endgroup%
\begin{picture}(1228,325)(3028,-1274)
\end{picture}%
}
To be completely precise, we have to admit that $\wTFo$ as a circuit
algebra has more generators than shown above. The last two generators
are ``foam vertices'', as will be explained shortly, and exist in all
possible orientations of the three strands. Some examples are shown on
the right. However, in Section~\ref{subsubsec:wops} we will describe the
operation ``orientation switch'' which allows switching the orientation of
any given strand.  In the algebraic structure which includes this extra
operation in addition to the circuit algebra structure, the generators
of the definition above are enough.
\end{definition}

\subsubsection{The generators of $\wTFo$}\label{subsubsec:wTFgens}
There is topological meaning to each of the generators of $\wTFo$:
they each stand for a certain local feature of framed knotted ribbon
tubes in $\bbR^4$. As in Section \ref{subsec:TangleTopology}, 
the tubes are oriented as 2-dimensional surfaces, and also have a
distinguished core with a 1-dimensional orientation (direction).

The crossings are as explained in Section \ref{subsubsec:ribbon} and Section
\ref{subsec:TangleTopology}: the under-strand denotes the ring flying through,
or the ``thin'' tube. Remember that there really are four kinds of crossings,
but in the circuit algebra the two not shown are obtained from the two that are shown
by adding virtual crossings.

\parpic[r]{\begin{picture}(0,0)%
\includegraphics{figs/Cap.pstex}%
\end{picture}%
%
%
\setlength{\unitlength}{3158sp}%
\begingroup\makeatletter\ifx\SetFigFont\undefined%
\gdef\SetFigFont#1#2#3#4#5{%
  \reset@font\fontsize{#1}{#2pt}%
  \fontfamily{#3}\fontseries{#4}\fontshape{#5}%
  \selectfont}%
\fi\endgroup%
\begin{picture}(1187,771)(176,-1123)
\end{picture}%
}
The bulleted end denotes a cap on the tube, or a flying ring that shrinks to a point, 
as in the figure on the right. In terms of 
Satoh's tubing map, the cap means that ``the string is attached to the bottom of the thickened
surface'', as shown in the figure below. Recall from Section \ref{subsubsec:TopTube} that the tubing map is the
composition 
$$\gamma\times S^1 \hookrightarrow \Sigma \times [-\epsilon,\epsilon] \hookrightarrow \bbR^4.$$
Here $\gamma$ is a trivalent tangle with ``drawn on the virtual surface $\Sigma$'', with caps ending on 
$\Sigma \times [-\epsilon, \epsilon]$. The first embedding above is the product of this ``drawing'' with an
$S^1$, while the second arises from the unit normal bundle of $\Sigma$ in $\bbR^4$. For each cap $(c, -\epsilon)$ the
tube resulting from Satoh's map has a boundary component $\partial_c=(c,-\epsilon)\times S^1$.
Follow the tubing map by gluing a disc to this boundary component to obtain the capped tube mentioned above.  

\begin{center}
 \input figs/SatohCap.pstex_t
\end{center}

\parpic[l]{\includegraphics[height=5cm]{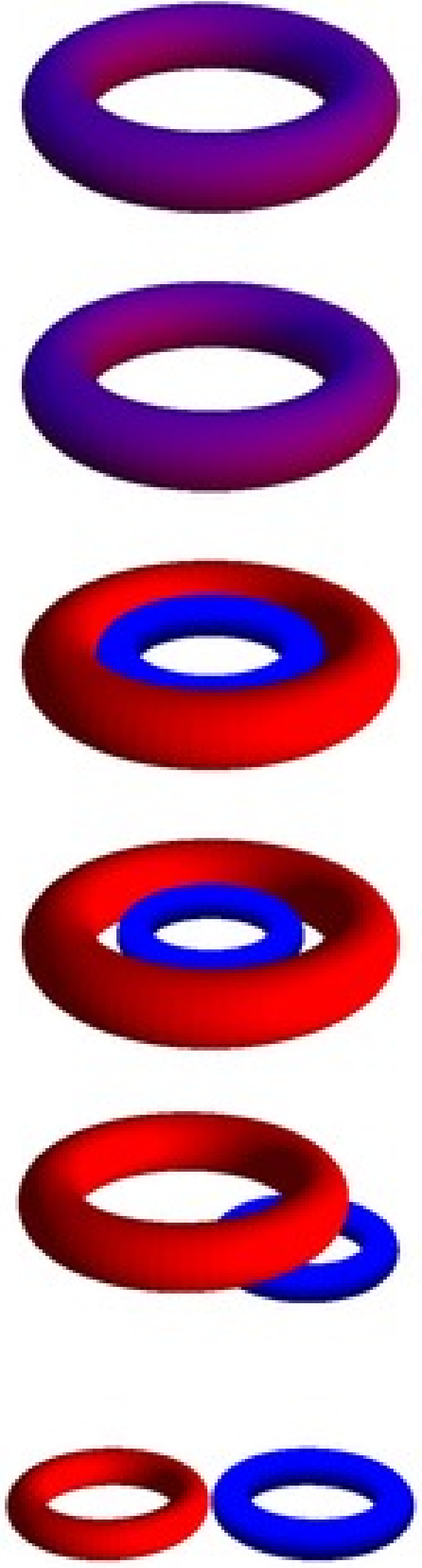}}
The last two generators denote singular ``foam vertices''. As the notation
suggests, a vertex can be thought of as ``half of a crossing''. To make
this precise using the flying rings interpretation, the first singular
vertex represents the movie shown on the left: the ring corresponding
to the right strand approaches the ring represented by the left strand
from below, flies inside it, and then the two rings fuse (as opposed to
a crossing where the ring coming from the right would continue to fly
out to above and to the left of the other one).  The second vertex is
the movie where a ring splits radially into a smaller and a larger ring,
and the small one flies out to the right and below the big one.

\vspace{5mm}

\parpic[r]{\begin{picture}(0,0)%
\includegraphics{figs/VertexInSurface.pstex}%
\end{picture}%
%
%
\setlength{\unitlength}{3158sp}%
\begingroup\makeatletter\ifx\SetFigFont\undefined%
\gdef\SetFigFont#1#2#3#4#5{%
  \reset@font\fontsize{#1}{#2pt}%
  \fontfamily{#3}\fontseries{#4}\fontshape{#5}%
  \selectfont}%
\fi\endgroup%
\begin{picture}(2949,1149)(2539,-1948)
\put(4426,-1111){\makebox(0,0)[lb]{\smash{{\SetFigFont{10}{12.0}{\rmdefault}{\mddefault}{\updefault}{\color[rgb]{0,0,0}$\Sigma \times [-\epsilon,\epsilon]$}%
}}}}
\end{picture}%
}
The vertices can also be interpreted topologically via a natural extension of Satoh's tubing map. For
the first generating vertex, imagine the broken right strand approaching the continuous left strand directly from below
in a thickened surface, as shown.

\begin{figure}[h!]
\input{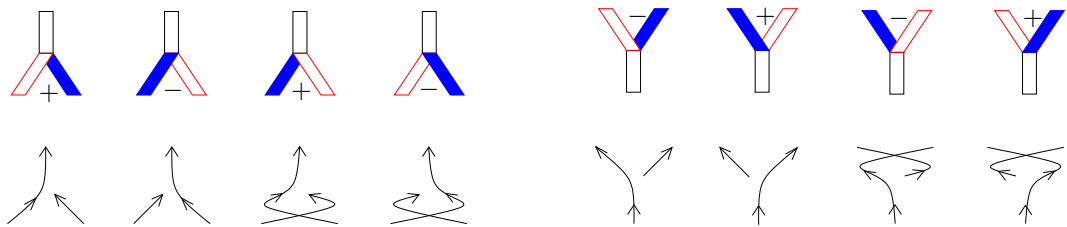}
\caption{Vertex types in $\wTFo$.}\label{fig:VertexTypes}
\end{figure}

The reader might object that there really are four types of vertices (as
there are four types of crossings), and each of these can be viewed as
a ``fuse'' or a ``split'' depending on the strand directions, as shown
in Figure~\ref{fig:VertexTypes}. However, looking at the fuse vertices
for example, observe that the last two of these can be obtained from the
first two by composing with virtual crossings, which always exist in a
circuit algebra.

The sign of a vertex can be defined the same way as the sign of a crossing
(see Section~\ref{subsec:TangleTopology}).  We will sometimes refer to
the first generator vertex as ``the positive vertex'' and to the second
one as ``the negative vertex''.  We use the band notation for vertices
the same way we do for crossings: the fully coloured band stands for
the thin (inner) ring.

\subsubsection{The relations of $\wTFo$} \label{subsubsec:wrels} 
In addition to the usual \Rs, R2, R3, and OC moves of
Figure~\ref{fig:VKnotRels}, we need more relations to describe the
behaviour of the additional features.

\begin{comment}\label{com:MissingTopology}
As before, the relations have local topological explanations, and
we conjecture that together they provide a Reidemeister theory for
``w-tangled foams'', that is, knotted ribbon tubes with foam vertices
in $\bbR^4$.  In this section we list the relations along with the
topological reasoning behind them. However, for any rigorous purposes below,
$\wTFo$ is studied as a circuit algebra given by the declared generators and relations,
regardless of their topological meaning.
\end{comment}

Recall that topologically, a cap represents a capped tube or equivalently, 
flying ring shrinking to a point. Hence, a cap
on the thin (or under) strand can be ``pulled out'' from a crossing,
but the same is not true for a cap on the thick (or over) strand, as
shown below. This is the case for any orientation of the strands. We
denote this relation by \glost{CP}, for Cap Pull-out.
\[ \input{figs/CapRel.pstex_t} \]

The Reidemeister 4 relations assert that a strand can be moved under
or over a crossing, as shown in the picture below. The ambiguously
drawn vertices in the picture denote a vertex of any kind (as described
in Section~\ref{subsubsec:wTFgens}), and the strands can be oriented
arbitrarily. The local topological (tube or flying ring) interpretations
can be read from the pictures below. These relations will be denoted
\glost{R4}.
\begin{center}
 \begin{picture}(0,0)%
\includegraphics{figs/R4.pstex}%
\end{picture}%
%
%
\setlength{\unitlength}{3158sp}%
\begingroup\makeatletter\ifx\SetFigFont\undefined%
\gdef\SetFigFont#1#2#3#4#5{%
  \reset@font\fontsize{#1}{#2pt}%
  \fontfamily{#3}\fontseries{#4}\fontshape{#5}%
  \selectfont}%
\fi\endgroup%
\begin{picture}(7077,2124)(7561,-1273)
\put(7576,-286){\makebox(0,0)[lb]{\smash{{\SetFigFont{12}{14.4}{\rmdefault}{\mddefault}{\updefault}{\color[rgb]{0,0,0}$R4:$}%
}}}}
\end{picture}%

\end{center}

\subsubsection{The auxiliary operations of $\wTFo$} \label{subsubsec:wops}
The circuit algebra $\wTFo$ is equipped with several extra operations.

The first of these is the familiar orientation switch. We will, as mentioned in
Section~\ref{subsec:TangleTopology}, distinguish between switching both the 2D and 1D orientations, or just the 
strand (1D) direction. 

Topologically {\it orientation switch}, denoted $\glos{S_e}$, is the
switch of both orientations of the strand $e$.  Diagrammatically (and this
is the definition) $S_e$ is the operation which reverses the orientation
of a strand in a $\wTFo$ diagram. The reader can check that when applying
Satoh's tubing map, this amounts to reversing both the direction and
the 2D orientation of the tube arising from the strand.

\parpic[r]{\begin{picture}(0,0)%
\includegraphics{figs/Adjoint.pstex}%
\end{picture}%
%
%
\setlength{\unitlength}{3947sp}%
\begingroup\makeatletter\ifx\SetFigFont\undefined%
\gdef\SetFigFont#1#2#3#4#5{%
  \reset@font\fontsize{#1}{#2pt}%
  \fontfamily{#3}\fontseries{#4}\fontshape{#5}%
  \selectfont}%
\fi\endgroup%
\begin{picture}(1333,1299)(4260,-1048)
\put(4275,-979){\makebox(0,0)[lb]{\smash{{\SetFigFont{12}{14.4}{\rmdefault}{\mddefault}{\updefault}{\color[rgb]{0,0,0}$e$}%
}}}}
\put(5177,-979){\makebox(0,0)[lb]{\smash{{\SetFigFont{12}{14.4}{\rmdefault}{\mddefault}{\updefault}{\color[rgb]{0,0,0}$e$}%
}}}}
\end{picture}%
}
The operation which, in topology world, reverses a tube's direction
but not its 2D orientation is called {\it ``adjoint''}, and denoted by
$\glos{A_e}$. This is slightly more intricate to define rigorously in terms
of diagrams. In addition to reversing the direction of the strand $e$
of the $\wTFo$ diagram, $A_e$ also locally changes each crossing of $e$
{\em over} another strand by adding two virtual crossings, as shown
on the right. We recommend for the reader to convince themselves that
this indeed represents a direction switch in topology after reading
Section~\ref{subsec:TheWen}.

\begin{figure}
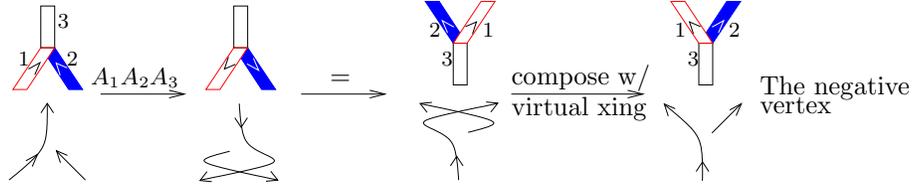

\input figs/VertexSwitch.pstex_t
\caption{Switching strand orientations at a vertex. The adjoint operation
only switches the tube direction, hence in the \emph{band picture}
only the arrows change. To express this vertex in terms of the negative
generating vertex in strand notation, we use a virtual crossing (see
Figure~\ref{fig:VertexTypes}).}\label{fig:VertexSwitch}
\end{figure}

\begin{remark}\label{rem:SwitchingVertices}
As an example, let us observe how the negative generator vertex
can be obtained from the positive generator vertex by adjoint
operations and composition with virtual crossings, as shown in
Figure~\ref{fig:VertexSwitch}. Note that also all other vertices can
be obtained from the positive vertex via orientation switch and adjoint
operations and composition by virtual crossings.

As a small exercise, it is worthwhile to convince ourselves of the effect
of orientation switch operations on the {\em band picture}. For example,
replace $A_1A_2A_3$ by $S_1S_2S_3$ in figure \ref{fig:VertexSwitch}.
In the strand diagram, this will only reverse the direction of the
strands. The reader can check that in the band picture not only the
arrows will reverse but also the blue band will switch to be on top of
the red band.
\end{remark}

\begin{comment}\label{com:wTFFraming}
Framings were discussed in Remark~\ref{rem:Framing}, but have not played
a significant role so far, except to explain the lack of a Reidemeister 1
relation. Now we will need to discuss framings in order to provide a topological
explanation for the unzip (tube doubling) operation.

In the local topological interpretation of $\wTFo$, strands represent ribbon-knotted tubes with foam vertices, 
which are also equipped with a framing, arising from the blackboard framing of the strand diagrams via Satoh's tubing map. 
Strand doubling is the operation of doubling a tube by ``pushing it off itself slightly'' in the framing direction,
as shown in Figure \ref{fig:BandUnzip}.

Recall that ribbon knotted tubes 
have a ``filling'', with only ``ribbon'' self-intersections. When we double a tube, we want this ribbon property to be preserved.
This is equivalent to saying that the ring obtained by pushing off any given girth of the tube in the
framing direction is not linked with the original tube, which is indeed the case, as mentioned in Remark \ref{rem:Framing}.

\parpic[r]{$\begin{picture}(0,0)%
\includegraphics{figs/DetailedVertex.pstex}%
\end{picture}%
\setlength{\unitlength}{3947sp}%
\begingroup\makeatletter\ifx\SetFigFont\undefined%
\gdef\SetFigFont#1#2#3#4#5{%
  \reset@font\fontsize{#1}{#2pt}%
  \fontfamily{#3}\fontseries{#4}\fontshape{#5}%
  \selectfont}%
\fi\endgroup%
\begin{picture}(1265,1555)(5561,-3576)
\end{picture}%
$} 
Framings arising from the blackboard framing of strand diagrams via Satoh's tubing map 
always match at the vertices, with the normal vectors pointing either directly towards or
away from the centre of the singular ring. Note that the orientations of the three tubes may or may not match.
An example of a vertex with the orientations and framings shown is on the right. Note that the
framings on the two sides of each band are mirror images of each other, as they should be.
\end{comment}

\begin{figure}
\input{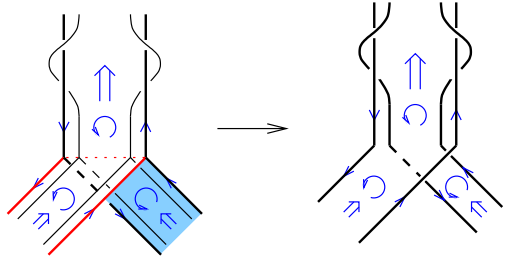}
\caption{Unzipping a tube, in band notation with orientations and framing marked.}\label{fig:BandUnzip}
\end{figure}

Unzip, or tube doubling is perhaps the most interesting of the auxiliary $\wTFo$ operations. As mentioned above, topologically 
this means pushing the tube off
itself slightly in the framing direction. At each of the vertices at the two ends of the doubled tube there are two tubes
to be attached to the doubled tube. At each end, the normal vectors pointed either directly towards or away from the centre,
so there is an ``inside'' and an ``outside'' ending ring. The two tubes to be attached also come as an ``inside'' and an
``outside'' one, which defines which one to attach to which. An example is shown in Figure \ref{fig:BandUnzip}. Unzip
can only be done if the 1D and 2D orientations match at both ends.

\parpic[r]{\input{figs/StringUnzip.pstex_t}}
To define unzip rigorously, we must talk only of strand diagrams. The
natural definition is to let $\glos{u_e}$ double the strand $e$ using the
blackboard framing, and then attach the ends of the doubled strand to the
connecting ones, as shown on the right.  We restrict unzip to strands
whose two ending vertices are of different signs. This is a somewhat
artificial condition which we impose to get equations equivalent to
the~\cite{AlekseevTorossian:KashiwaraVergne} equations.

A related operation, {\it disk unzip}, is unzip done on a capped strand, pushing the tube off in the direction of the framing
(in diagrammatic world, in the direction of the blackboard framing), as 
before. An example in the line and band notations (with the framing suppressed) is shown below.

\begin{center}
\input{figs/CapUnzip.pstex_t}
\end{center}

Finally, we allow the deletion of ``long linear'' strands, meaning strands that do not end in a vertex on either side.

The goal, as before, is to construct a homomorphic expansion for $\wTFo$. However, first we need to understand its target space, 
the projectivization $\proj \wTFo$.

\draftcut
\subsection{The projectivization} \label{subsec:fproj}
Mirroring the previous section, we describe the projectivization
$\glos{\calA^{sw}}$ of $\wTFo$ and its ``full version'' $\glos{\calA^w}$
as circuit algebras on certain generators modulo a number of
relations. From now on we will write $\glos{A^{(s)w}}$ to mean ``$\calA^{w}$
and/or $\calA^{sw}$''.
\[ \calA^{(s)w}=\CA\!\left.\left.\left\langle
  \raisebox{-2mm}{\input{figs/wTFprojgens.pstex_t}}
  \right|
  \parbox{1.2in}{\centering relations as in Section~\ref{subsubsec:wTFProjRels}}
  \right|
  \parbox{1.2in}{\centering operations as in Section~\ref{subsubsec:wTFProjOps}}
  \right\rangle.
\]
In other words, $\calA^{(s)w}$ are the circuit algebras of arrow diagrams on trivalent (or foam) skeletons with 
caps. Note that all but the first of the generators are skeleton features (of degree 0), and that the single arrow is
the only generator of degree $1$. As for the generating vertices, the same remark applies as in Definition \ref{def:wTFo},
that is, there are more vertices with all possible strand orientations needed to generate $\calA^{(s)w}$ as circuit algebras.

\subsubsection{The relations of $\calA^{(s)w}$}\label{subsubsec:wTFProjRels}
In addition to the usual $\aft$ and TC 
relations (see Section \ref{subsec:FT4Braids}), as well as RI in the case of $\calA^{sw}=\calA^w/RI$, 
diagrams in $\calA^{(s)w}$ satisfy the following additional relations:

{\it Vertex invariance}, denoted by \glost{VI}, are relations arising
the same way as $\aft$ does, but with the participation of a vertex as
opposed to a crossing:
\begin{center}
 \input figs/VI.pstex_t
\end{center}
The other end of the arrow is in the same place throughout the relation, somewhere outside the picture
shown. The signs are positive whenever the strand on which the arrow ends is directed towards the vertex,
and negative when directed away. The ambiguously drawn vertex means any kind of vertex, but the same one throughout.

\parpic[r]{\begin{picture}(0,0)%
\includegraphics{figs/CapHeads.pstex}%
\end{picture}%
%
%
\setlength{\unitlength}{3158sp}%
\begingroup\makeatletter\ifx\SetFigFont\undefined%
\gdef\SetFigFont#1#2#3#4#5{%
  \reset@font\fontsize{#1}{#2pt}%
  \fontfamily{#3}\fontseries{#4}\fontshape{#5}%
  \selectfont}%
\fi\endgroup%
\begin{picture}(925,632)(288,-983)
\put(826,-736){\makebox(0,0)[lb]{\smash{{\SetFigFont{10}{12.0}{\rmdefault}{\mddefault}{\updefault}{\color[rgb]{0,0,0}$=0$}%
}}}}
\end{picture}%
}
The CP relation (a cap can be pulled out from under a strand but not from
over, Section \ref{subsubsec:wrels}) implies that arrow heads near a cap
are zero, as shown on the right. Denote this relation also by
\glost{CP}. (Also note that a tail near a cap is not set to zero.)

As in the previous sections, and in particular in Definition~\ref{def:wJac}, we define a
``w-Jacobi diagram'' (or just ``arrow diagram'') on a foam
skeleton by allowing trivalent chord vertices. Denote the circuit algebra of formal 
linear combinations of arrow diagrams by $\calA^{(s)wt}$. We have the following bracket-rise theorem:

\begin{theorem} The obvious inclusion of diagrams induces a circuit
algebra isomorphism $\calA^{(s)w}\cong\calA^{(s)wt}$. Furthermore, the $\aAS$
and $\aIHX$ relations of Figure~\ref{fig:aIHX} hold in $\calA^{(s)wt}$.
\end{theorem}

\begin{proof} Same as the proof of Theorem~\ref{thm:BracketRise}. \qed
\end{proof}

As in Section~\ref{subsec:vw-tangles}, the primitive elements of
$\calA^{(s)w}$ are connected diagrams, namely trees and wheels. Before
moving on to the auxiliary operations of $\calA^{(s)w}$, let us make
two useful observations:

\begin{lemma}\label{lem:CapIsWheels}
$\calA^w(\raisebox{-1mm}{\begin{picture}(0,0)%
\includegraphics{figs/SmallCap.pstex}%
\end{picture}%
%
%
\setlength{\unitlength}{3947sp}%
\begingroup\makeatletter\ifx\SetFigFont\undefined%
\gdef\SetFigFont#1#2#3#4#5{%
  \reset@font\fontsize{#1}{#2pt}%
  \fontfamily{#3}\fontseries{#4}\fontshape{#5}%
  \selectfont}%
\fi\endgroup%
\begin{picture}(64,201)(151,-1123)
\end{picture}%
})$, the part of
$\calA^w$ with skeleton $\raisebox{-1mm}{\begin{picture}(0,0)%
\includegraphics{figs/SmallCap.pstex}%
\end{picture}%
%
%
\setlength{\unitlength}{3947sp}%
\begingroup\makeatletter\ifx\SetFigFont\undefined%
\gdef\SetFigFont#1#2#3#4#5{%
  \reset@font\fontsize{#1}{#2pt}%
  \fontfamily{#3}\fontseries{#4}\fontshape{#5}%
  \selectfont}%
\fi\endgroup%
\begin{picture}(64,201)(151,-1123)
\end{picture}%
}$,
is isomorphic as a vector space to the completed polynomial 
algebra freely generated by wheels $w_k$ with $k \geq 1$. Likewise
$\calA^{sw}(\raisebox{-1mm}{\begin{picture}(0,0)%
\includegraphics{figs/SmallCap.pstex}%
\end{picture}%
%
%
\setlength{\unitlength}{3947sp}%
\begingroup\makeatletter\ifx\SetFigFont\undefined%
\gdef\SetFigFont#1#2#3#4#5{%
  \reset@font\fontsize{#1}{#2pt}%
  \fontfamily{#3}\fontseries{#4}\fontshape{#5}%
  \selectfont}%
\fi\endgroup%
\begin{picture}(64,201)(151,-1123)
\end{picture}%
})$, except here
$k \geq 2$.
\end{lemma}

\begin{proof}
 Any arrow diagram with an arrow head at its top is zero by the Cap Pull-out (CP) relation. If $D$ is an arrow
diagram that has a head somewhere on the skeleton but not at the top, then one can use repeated $\aSTU$ relations 
to commute the head to the top at the cost of diagrams with one fewer skeleton head. 

Iterating
this procedure, we can get rid of all arrow heads, and hence write $D$ as a linear combination of 
diagrams having no heads on the skeleton. All connected components of such diagrams are wheels. 

To prove that there are no relations between wheels in $\calA^{(s)w}(\raisebox{-1mm}{\begin{picture}(0,0)%
\includegraphics{figs/SmallCap.pstex}%
\end{picture}%
%
%
\setlength{\unitlength}{3947sp}%
\begingroup\makeatletter\ifx\SetFigFont\undefined%
\gdef\SetFigFont#1#2#3#4#5{%
  \reset@font\fontsize{#1}{#2pt}%
  \fontfamily{#3}\fontseries{#4}\fontshape{#5}%
  \selectfont}%
\fi\endgroup%
\begin{picture}(64,201)(151,-1123)
\end{picture}%
})$, 
let $S_L\colon \calA^{(s)w}(\uparrow_1) \to \calA^{(s)w}(\uparrow_1)$ 
(resp. $S_R$) be the map that sends an arrow diagram to the sum of all ways of dropping one left (resp. right) arrow 
(on a vertical strand, left means down and right means up). Define
$$\glos{F}:=\sum_{k=0}^{\infty}\frac{(-1)^k}{k!}D_R^k(S_L+S_R)^k,$$
where $D_R$ is the short right arrow as shown in Figure~\ref{fig:AwGenerators}.
We leave it as an exercise for the reader to check that $F$ is a bi-algebra homomorphism that kills diagrams with an arrow head at the top
(i.e., CP is in the kernel of $F$), and $F$ is injective on wheels. This concludes the proof.
\qed
\end{proof}

\begin{lemma}
$\calA^{(s)w}(Y)=\calA^{(s)w}(\uparrow_2)$, where $\calA^{(s)w}(Y)$
stands for the space of arrow diagrams whose skeleton is a $Y$-graph with
any orientation of the strands, and as before $\calA^{(s)w}(\uparrow_2)$
is the space of arrow diagrams on two strands.
\end{lemma}

\begin{proof}
 We can use the vertex invariance (VI) relation to push all arrow heads and tails from the ``trunk'' of the vertex to the other two strands.
\qed
\end{proof}

\subsubsection{The auxiliary operations of $\calA^{(s)w}$}\label{subsubsec:wTFProjOps}
Recall from Section \ref{subsec:TangleTopology} that the orientation switch $S_e$ (i.e. changing both the $1D$ and $2D$ orientations of a strand) always changes
the sign of a crossing involving the strand $e$. Hence, letting $S$ denote any foam (trivalent) skeleton, the induced arrow diagrammatic
operation is a map $S_e\colon \calA^{(s)w}(S) \to \calA^{(s)w}(S_e(S))$ which acts by multiplying  
each arrow diagram by $(-1)$ raised to the number of arrow endings on $e$ (counting both heads and tails).

The adjoint operation $A_e$ (i.e. switching only the strand direction), on the other hand, only changes the sign of a crossing when the
strand being switched is the under- (or through) strand. (See section \ref{subsec:TangleTopology} for pictures and explanation.) 
Therefore, the arrow diagrammatic $A_e$ acts by switching the direction
of $e$ and multiplying each arrow diagram by $(-1)$ raised to the number
of {\it arrow heads} on $e$. Note that in $\calA^{(s)w}(\uparrow_n)$
taking the adjoint on every strand gives the adjoint map of Definition
\ref{def:Adjoint}.

\parpic[r]{\input{figs/Unzip.pstex_t}}
The arrow diagram operations induced by unzip and disc unzip (both to be denoted $u_e$, and interpreted appropriately according to whether the 
strand $e$ is capped) are maps $u_e\colon  \calA^{(s)w}(S) \to \calA^{(s)w}(u_e(S))$, where each arrow ending (head or tail) on $e$ is mapped to
a sum of two arrows, one ending on each of the new strands, as shown on the right. In other words, if in an arrow diagram $D$ there are $k$ arrow
ends on $e$, then $u_e(D)$ is a sum of $2^k$ arrow diagrams.

The operation induced by deleting the long linear strand $e$ is the map $d_e\colon  \calA^{(s)w}(S) \to \calA^{(s)w}(d_e(S))$ which kills arrow diagrams with
any arrow ending (head or tail) on $e$, and leaves all else unchanged, except with $e$ removed.

\draftcut
\subsection{The homomorphic expansion}\label{subsec:wTFExpansion}
\begin{theorem}\label{thm:wTFExpansionExists}
There exists a group-like\footnote{The formal definition of the group-like property
is along the lines of \ref{par:Delta}. In practise, it means that the $Z$-values of the vertices, crossings, and cap 
(denoted $V$, $R$ and $C$ below) are exponentials of linear combinations of
connected diagrams.}
homomorphic expansion for $\wTFo$, i.e. a group-like expansion $Z\colon
\wTFo \to \calA^{sw}$ which is a map of circuit algebras and also
intertwines the auxiliary operations of $\wTFo$ with their arrow
diagrammatic counterparts.
\end{theorem}

Since both $\wTFo$ and $\calA^{sw}$ are circuit algebras defined by generators and relations,
when looking for a suitable $Z$ all we have to do is to find values for each of the generators of $\wTFo$ so that these satisfy (in $\calA^{sw}$) the equations
which arise from the relations in $\wTFo$ and the homomorphicity requirement. In this section we will derive these equations and show that they are equivalent to
the Alekseev-Torossian version of the Kashiwara-Vergne equations \cite{AlekseevTorossian:KashiwaraVergne}. In \cite{AlekseevEnriquezTorossian:ExplicitSolutions}
Alekseev Enriquez and Torossian construct explicit solutions to these equations using associators. In a later paper we will interpret these results in
our context of homomorphic expansions for w-tangled foams.

Let $\glos{R}:=Z(\overcrossing) \in \calA^{sw}(\uparrow_2)$. It
follows from the Reidemeister 2 relation that
$Z(\undercrossing)=(R^{-1})^{21}$. As discussed in Sections
\ref{subsec:vw-tangles} and \ref{subsec:UniquenessForTangles},
Reidemeister 3 with group-likeness and homomorphicity
implies that $R=e^a$, where $a$ is a single arrow
pointing from the over to the under strand.  Let
$\glos{C}:=Z(\raisebox{-1mm}{\begin{picture}(0,0)%
\includegraphics{figs/SmallCap.pstex}%
\end{picture}%
%
%
\setlength{\unitlength}{3947sp}%
\begingroup\makeatletter\ifx\SetFigFont\undefined%
\gdef\SetFigFont#1#2#3#4#5{%
  \reset@font\fontsize{#1}{#2pt}%
  \fontfamily{#3}\fontseries{#4}\fontshape{#5}%
  \selectfont}%
\fi\endgroup%
\begin{picture}(64,201)(151,-1123)
\end{picture}%
})\in
\calA^{sw}(\raisebox{-1mm}{\begin{picture}(0,0)%
\includegraphics{figs/SmallCap.pstex}%
\end{picture}%
%
%
\setlength{\unitlength}{3947sp}%
\begingroup\makeatletter\ifx\SetFigFont\undefined%
\gdef\SetFigFont#1#2#3#4#5{%
  \reset@font\fontsize{#1}{#2pt}%
  \fontfamily{#3}\fontseries{#4}\fontshape{#5}%
  \selectfont}%
\fi\endgroup%
\begin{picture}(64,201)(151,-1123)
\end{picture}%
})$.
By Lemma \ref{lem:CapIsWheels}, we know that
$C$ is made up of wheels only.  Finally, let
$\glos{V}=\glos{V^+}:=Z(\raisebox{-1mm}{\begin{picture}(0,0)%
\includegraphics{figs/PlusVertex.pstex}%
\end{picture}%
%
%
\setlength{\unitlength}{3947sp}%
\begingroup\makeatletter\ifx\SetFigFont\undefined%
\gdef\SetFigFont#1#2#3#4#5{%
  \reset@font\fontsize{#1}{#2pt}%
  \fontfamily{#3}\fontseries{#4}\fontshape{#5}%
  \selectfont}%
\fi\endgroup%
\begin{picture}(211,211)(3439,-1160)
\end{picture}%
})\in
\calA^{sw}(\raisebox{-1mm}{\begin{picture}(0,0)%
\includegraphics{figs/PlusVertex.pstex}%
\end{picture}%
%
%
\setlength{\unitlength}{3947sp}%
\begingroup\makeatletter\ifx\SetFigFont\undefined%
\gdef\SetFigFont#1#2#3#4#5{%
  \reset@font\fontsize{#1}{#2pt}%
  \fontfamily{#3}\fontseries{#4}\fontshape{#5}%
  \selectfont}%
\fi\endgroup%
\begin{picture}(211,211)(3439,-1160)
\end{picture}%
})\cong
\calA^{sw}(\uparrow_2)$, and
$\glos{V^-}:=Z(\raisebox{-1mm}{\begin{picture}(0,0)%
\includegraphics{figs/MinusVertex.pstex}%
\end{picture}%
%
%
\setlength{\unitlength}{3158sp}%
\begingroup\makeatletter\ifx\SetFigFont\undefined%
\gdef\SetFigFont#1#2#3#4#5{%
  \reset@font\fontsize{#1}{#2pt}%
  \fontfamily{#3}\fontseries{#4}\fontshape{#5}%
  \selectfont}%
\fi\endgroup%
\begin{picture}(211,211)(3065,-1347)
\end{picture}%
})\in
\calA^{sw}(\raisebox{-1mm}{\begin{picture}(0,0)%
\includegraphics{figs/MinusVertex.pstex}%
\end{picture}%
%
%
\setlength{\unitlength}{3158sp}%
\begingroup\makeatletter\ifx\SetFigFont\undefined%
\gdef\SetFigFont#1#2#3#4#5{%
  \reset@font\fontsize{#1}{#2pt}%
  \fontfamily{#3}\fontseries{#4}\fontshape{#5}%
  \selectfont}%
\fi\endgroup%
\begin{picture}(211,211)(3065,-1347)
\end{picture}%
})\cong
\calA^{sw}(\uparrow_2)$.

Before we translate each of the relations of Section \ref{subsubsec:wrels} to equations let us slightly extend the notation used in 
Section \ref{subsec:UniquenessForTangles}. Recall that $R^{23}$, for instance, meant ``$R$ placed on strands 2 and 3''. In this section 
we also need notation such as $R^{(23)1}$, which means ``$R$ with its first strand doubled, placed on strands 2, 3 and 1''.

Now on to the relations, note that Reidemeister 2 and 3 and Overcrossings Commute have already been dealt with. Of the two Reidemeister 4
relations, the first one induces an equation that is automatically satisfied. Pictorially, the equation 
looks as follows:
\begin{center}
 \input{figs/R4ToEquation.pstex_t}
\end{center}
In other words, we obtained the equation
$$V^{12}R^{3(12)}=R^{32}R^{31}V^{12}.$$
However, observe that by the ``head-invariance'' property of arrow diagrams (Remark \ref{rem:HeadInvariance})
$V^{12}$ and $R^{3(12)}$ commute on the left hand side. Hence the left hand side equals $R^{3(12)}V^{12}=R^{32}R^{31}V^{12}$.
Also, $R^{3(12)}=e^{a^{31}+a^{32}}=e^{a^{32}}e^{a^{31}}=R^{32}R^{31}$, where the second step is due to the fact that 
$a^{31}$ and $a^{32}$ commute. Therefore, the equation is true independently of the choice of $V$.

We have no such luck with the second Reidemeister 4 relation, which, in the same manner as in the paragraph above, 
translates to the equation
\begin{equation}\label{eq:HardR4}
 V^{12}R^{(12)3}=R^{23}R^{13}V^{12}.
\end{equation}
There is no ``tail invariance'' of arrow diagrams, so $V$ and $R$ do not commute on the left hand side; also, $R^{(12)3}\neq R^{23}R^{13}$. 
As a result, this equation puts a genuine restriction on the choice of $V$. 

The Cap Pull-out (CP) relation translates to the equation $R^{12}C^2=C^2$. This is true independently of the choice of $C$: by head-invariance,
$R^{12}C^2=C^2R^{12}$. Now $R^{12}$ is just below the cap on strand $2$, and the cap ``kills heads'', in other words, every term of $R^{12}$
with an arrow head at the top of strand $2$ is zero. Hence, the only surviving term of $R^{12}$ is $1$ (the empty diagram), which makes the
equation true.

The homomorphicity of the orientation switch operation was used to prove the uniqueness of $R$ in 
Theorem \ref{thm:Tangleuniqueness}. The homomorphicity of the adjoint leads to the equation
$V_-=A_1A_2(V)$ (see Figure \ref{fig:VertexSwitch}), eliminating $V_-$ as an unknown.
Note that we also silently assumed these homomorphicity properties when we did not introduce
32 different values of the vertex depending on the strand orientations.

Homomorphicity of the (annular) unzip operation leads to an equation for $V$, which we are going to refer to as ``unitarity''. This 
is illustrated in the figure below. Recall that $A_1$ and $A_2$ denote the adjoint (direction switch) operation on strand 1
and 2, respectively. 
\begin{center}
 \input figs/Unitarity.pstex_t
\end{center}
Reading off the equation, we have 
\begin{equation}\label{eq:unitarity}
 V\cdot A_1A_2(V)=1.
\end{equation}

\parpic[r]{\input{figs/CapEqn.pstex_t}}\picskip{5}
Homomorphicity of the disk unzip leads to an equation for $C$ which
we will refer to as the ``cap equation''.  The translation from
homomorphicity to equation is shown in the figure on the right.  $C$,
as we introduced before, denotes the $Z$-value of the cap. Hence, the
cap equation reads

\begin{equation}\label{eq:CapEqn}
  V^{12}C^{(12)}=C^1C^2 \qquad\text{in}\quad
  \calA^{sw}(\raisebox{-1mm}{\begin{picture}(0,0)%
\includegraphics{figs/SmallCap.pstex}%
\end{picture}%
%
%
\setlength{\unitlength}{3947sp}%
\begingroup\makeatletter\ifx\SetFigFont\undefined%
\gdef\SetFigFont#1#2#3#4#5{%
  \reset@font\fontsize{#1}{#2pt}%
  \fontfamily{#3}\fontseries{#4}\fontshape{#5}%
  \selectfont}%
\fi\endgroup%
\begin{picture}(64,201)(151,-1123)
\end{picture}%
}_2)
\end{equation}

The homomorphicity of deleting long strands does not lead to an equation on its own, however it was
used to prove the uniqueness of $R$ (Theorem \ref{thm:Tangleuniqueness}).

To summarize, we have reduced the problem of finding a homomorphic
expansion $Z$ to finding the $Z$-values of the (positive) vertex and
the cap, denoted $V$ and $C$, subject to three equations: the ``hard
Reidemeister 4'' equation (\ref{eq:HardR4}); ``unitarity of V'' equation
(\ref{eq:unitarity}); and the ``cap equation'' (\ref{eq:CapEqn}).

\draftcut
\subsection{The equivalence with the Alekseev-Torossian equations}\label{subsec:EqWithAT}
First let us recall Alekseev and Torossian's
formulation of the generalized Kashiwara-Vergne problem
(see~\cite[Section~5.3]{AlekseevTorossian:KashiwaraVergne}):

{\bf Generalized KV problem:} Find an element $\glos{F}\in \TAut_2$ with the properties
\begin{equation}\label{eq:ATKVEqns}
 F(x+y)=\log(e^xe^y), \text{ and } j(F)\in \im(\tilde{\delta}).
\end{equation}
Here $\tilde{\delta}\colon  \attr_1 \to \attr_2$ is defined by $(\tilde{\delta}a)(x,y)=a(x)+a(y)-a(\log(e^{x}e^{y}))$,
where elements of $\attr_2$ are cyclic words in the letters $x$ and $y$. (See
\cite{AlekseevTorossian:KashiwaraVergne}, Equation (8)). Note that an element of $\attr_1$ is a polynomial 
with no constant term in one variable. 
In other words, the second condition says that there exists 
$a \in \attr_1$ such that $jF=a(x)+a(y)-a(\log(e^{x}e^{y}))$.

\begin{theorem}\label{thm:ATEquivalence}
Theorem~\ref{thm:wTFExpansionExists}, namely the existence of a
group-like homomorphic expansion for $\wTFo$, is equivalent to the
generalized Kashiwara-Vergne problem.
\end{theorem}

\begin{proof}
We have reduced the problem of finding a homomorphic expansion to finding group-like solutions $V$ and $C$ to the hard Reidemeister 4 equation (\ref{eq:HardR4}),
the unitarity equation (\ref{eq:unitarity}), and the cap equation (\ref{eq:CapEqn}).

Suppose we have found such solutions and write $V=e^be^{uD}$, where $b \in \tr_2^s$, $D \in  \tder_2\oplus \fraka_2$, 
and where $u$ is the map $u\colon  \tder_2 \to \calA^{sw}(\uparrow_2)$ which plants the head of a tree
above all of its tails, as introduced in Section \ref{subsec:ATSpaces}. $V$ can be written in this form without loss of generality because wheels can always be brought to the
bottom of a diagram (at the possible cost of more wheels). Furthermore, $V$ is group-like and hence it can be written in exponential form. Similarly, write
$C=e^c$ with $c \in \attr_1^s$.

Note that $u(\fraka_2)$ is central in $\calA^{sw}(\uparrow_2)$ and
that replacing a solution $(V,C)$ by $(e^{u(a)}V, C)$ for any $a \in
\fraka_2$ does not interfere with any of the equations (\ref{eq:HardR4}),
(\ref{eq:unitarity}) or (\ref{eq:CapEqn}). Hence we may assume that $D$
does not contain any single arrows, that is, $D \in \tder_2$. Also, a
solution $(V,C)$ in $\calA^{sw}$ can be lifted to a solution in $\calA^w$
by simply setting the degree one terms of $b$ and $c$ to be zero. It is
easy to check that this $b \in \attr_2$ and $c \in \attr_1$ along with $D$
still satisfy the equations. (In fact, in $\calA^w$ (\ref{eq:unitarity})
and (\ref{eq:CapEqn}) respectively imply that $b$ is zero in degree 1,
and that the degree 1 term of $c$ is arbitrary, so we may as well assume
it to be zero.)  In light of this we declare that $b\in \attr_2$ and $c
\in \attr_1$.

The hard Reidemeister 4 equation (\ref{eq:HardR4}) reads $V^{12}R^{(12)3}=R^{23}R^{13}V^{12}$. Denote the arrow from strand 1 to strand 3 by $x$, and the
arrow from strand 2 to strand 3 by $y$. Substituting the known value for $R$ and rearranging, we get 
$$e^be^{uD}e^{x+y}e^{-uD}e^{-b}=e^ye^x.$$ Equivalently, $e^{uD}e^{x+y}e^{-uD}=e^{-b}e^ye^xe^b.$ Now on the right side there are only tails on the
first two strands, hence $e^b$ commutes with $e^ye^x$, so $e^{-b}e^b$ cancels. Taking logarithm of both sides we obtain 
$e^{uD}(x+y)e^{-uD}=\log e^ye^x$. Now for notational alignment with \cite{AlekseevTorossian:KashiwaraVergne} we switch strands 1 and 2, which exchanges 
$x$ and $y$ so we obtain:
\begin{equation}\label{eq:HardR4Translated}
e^{uD^{21}}(x+y)e^{-uD^{21}}=\log e^xe^y.
\end{equation}

The unitarity of $V$ (Equation (\ref{eq:unitarity})) translates to $1=e^be^{uD}(e^be^{uD})^*,$ where $*$ denotes the adjoint map (Definition \ref{def:Adjoint}). Note that the adjoint switches 
the order of a product and acts trivially on wheels. Also, $e^{uD}(e^{uD})^*=J(e^D)=e^{j(e^D)}$, by Proposition \ref{prop:Jandj}. 
So we have $1=e^be^{j(e^D)}e^b$. Multiplying by $e^{-b}$ on the right and by $e^b$ on the left, we get $1=e^{2b}e^{j(e^D)}$, and again by switching strand 1 and 2 we arrive at
\begin{equation}\label{eq:UnitarityTranslated}
1=e^{2b^{21}}e^{j(e^{D^{21}})}.
\end{equation}

As for the cap equation, if $C^1=e^{c(x)}$ and $C^2=e^{c(y)}$, then $C^{12}=e^{c(x+y)}$. Note that wheels
on different strands commute, hence $e^{c(x)}e^{c(y)}=e^{c(x)+c(y)}$, so the cap equation reads $$e^be^{uD}e^{c(x+y)}=e^{c(x)+c(y)}.$$ As this equation lives
in the space of arrow diagrams on two \emph{capped} strands, we can multiply the left side on the right by $e^{-uD}$: $uD$ has its head at the top, so it
is 0 by the Cap relation, hence $e^{uD}=1$ near the cap. Hence, $$e^be^{uD}e^{c(x+y)}e^{-uD}=e^{c(x)+c(y)}.$$ 

\parpic[r]{\begin{picture}(0,0)%
\includegraphics{figs/Sigma.pstex}%
\end{picture}%
%
%
\setlength{\unitlength}{3947sp}%
\begingroup\makeatletter\ifx\SetFigFont\undefined%
\gdef\SetFigFont#1#2#3#4#5{%
  \reset@font\fontsize{#1}{#2pt}%
  \fontfamily{#3}\fontseries{#4}\fontshape{#5}%
  \selectfont}%
\fi\endgroup%
\begin{picture}(3207,1265)(10299,-683)
\put(11688,-30){\makebox(0,0)[lb]{\smash{{\SetFigFont{12}{14.4}{\rmdefault}{\mddefault}{\updefault}{\color[rgb]{0,0,0}$\sigma$}%
}}}}
\end{picture}%
}
On the right side of the equation above \linebreak $e^{uD}e^{c(x+y)}e^{-uD}$ reminds us of Equation (\ref{eq:HardR4Translated}), however we cannot use (\ref{eq:HardR4Translated})
directly as we live in a different space now. In particular, $x$ there meant an arrow from strand 1 to strand 3, while here it means a one-wheel on (capped) 
strand 1, and similarly for $y$. Fortunately, there is a map $\sigma\colon  \calA^{sw}(\uparrow_3) \to \calA^{sw}(\raisebox{-1mm}{\begin{picture}(0,0)%
\includegraphics{figs/SmallCap.pstex}%
\end{picture}%
%
%
\setlength{\unitlength}{3947sp}%
\begingroup\makeatletter\ifx\SetFigFont\undefined%
\gdef\SetFigFont#1#2#3#4#5{%
  \reset@font\fontsize{#1}{#2pt}%
  \fontfamily{#3}\fontseries{#4}\fontshape{#5}%
  \selectfont}%
\fi\endgroup%
\begin{picture}(64,201)(151,-1123)
\end{picture}%
}_2)$,
where $\sigma$ ``closes the third strand and turns it into a chord (or internal) strand, and caps the first two strands'', as shown on the right. This map is
well defined (in fact, it kills almost all relations, and turns one $\aSTU$ into an $\aIHX$). Under this map, using our abusive notation, $\sigma(x)=x$ and 
$\sigma(y)=y$.

Now we can apply Equation (\ref{eq:HardR4Translated}) and get $e^be^{c(\log e^y e^x)}=e^{c(x)+c(y)}$, which, using that tails commute, implies
$b=c(x)+c(y)-c(\log e^y e^x)$. Switching strands 1 and 2, we obtain
\begin{equation}\label{eq:CapEqnTranslated}
b^{21}=c(x)+c(y)-c(\log e^x e^y)
\end{equation}

In summary, we can use $(V,C)$ to produce $F:=e^{D^{21}}$ (sorry\footnote{%
  We apologize for the annoying $2\leftrightarrow 1$ transposition in this equation,
  which makes some later equations, especially~\eqref{eq:ATPhiandV},
  uglier than they could have been. There is no depth here, just
  mis-matching conventions between us and Alekseev-Torossian.
})
which satisfies the Alekseev-Torossian equations
(\ref{eq:ATKVEqns}): $e^{D^{21}}$ acts on $\lie_2$ by conjugation
by $e^{uD^{21}}$, so the first part of (\ref{eq:ATKVEqns})
is implied by (\ref{eq:HardR4Translated}). The second half of
(\ref{eq:ATKVEqns}) is true due to (\ref{eq:UnitarityTranslated}) and
(\ref{eq:CapEqnTranslated}).

On the other hand, suppose that we have found $F\in \TAut_2$ and $a \in \tr_1$ satisfying (\ref{eq:ATKVEqns}). 
Then set $D^{21}:=\log F$, $b^{21}:=\frac{-j(e^{D^{21}})}{2}$,
and $c \in \tilde{\delta}^{-1}(b^{21})$, in particular $c=-\frac{a}{2}$ works. Then $V=e^be^{uD}$ and $C=e^c$ satisfy the equations for 
homomorphic expansions (\ref{eq:HardR4}), (\ref{eq:unitarity})
and (\ref{eq:CapEqn}).\qed

\end{proof}

\draftcut
\subsection{The wen}\label{subsec:TheWen}
A topological feature of w-tangled foams which we excluded from the theory
so far is the wen $\glos{w}$. The wen was introduced
in~\ref{subsubsec:NonHorRings} as a Klein bottle cut apart; it amounts
to changing the 2D orientation of a tube, as shown in the picture below:
\begin{center}
 \begin{picture}(0,0)%
\includegraphics{figs/Wen2.pstex}%
\end{picture}%
%
%
\setlength{\unitlength}{4934sp}%
\begingroup\makeatletter\ifx\SetFigFont\undefined%
\gdef\SetFigFont#1#2#3#4#5{%
  \reset@font\fontsize{#1}{#2pt}%
  \fontfamily{#3}\fontseries{#4}\fontshape{#5}%
  \selectfont}%
\fi\endgroup%
\begin{picture}(4737,1339)(-2861,-1744)
\put(-1488,-1180){\makebox(0,0)[lb]{\smash{{\SetFigFont{14}{16.8}{\rmdefault}{\mddefault}{\updefault}{\color[rgb]{0,0,0}$w$}%
}}}}
\put(301,-1186){\makebox(0,0)[lb]{\smash{{\SetFigFont{14}{16.8}{\rmdefault}{\mddefault}{\updefault}{\color[rgb]{0,0,0}$=$}%
}}}}
\end{picture}%

\end{center}
In this section we study the circuit algebra of w-Tangled Foams with the wen rightfully included as a generator,
and denote this space by $\glos{\wTF}$. 

\subsubsection{The relations and auxiliary operations of $\wTF$.}
\label{subsubsec:WenRels}
Adding the wen as a generator means we have to impose additional relations
involving the wen to keep our topological heuristics intact, as follows:

The interaction of a wen and a crossing has already been mentioned in Section \ref{subsubsec:NonHorRings},
and is described by Equation (\ref{eq:FlipRels}), which we repeat here for convenience:
\begin{center}
 \input figs/FlipRels2.pstex_t
\end{center}
Recall that in flying ring language, a wen is a ring flipping over. It does not matter whether ring B flips first and
then flies through ring A or vice versa. However, the movies in which ring A first flips and then ring B flies through it, 
or B flies through A first and then A flips differ in the fly-through direction, which is cancelled by virtual crossings,
as in the figure above. We will refer to these relations as the Flip Relations, and abbreviate them by \glost{FR}. 

A double flip is homotopic to no flip, in other words two consecutive wens
equal no wen. Let us denote this relation by $\glos{W^2}$, for Wen
squared. Note that this relation explains why there are no ``left and
right wens''.

\parpic[r]{$\begin{picture}(0,0)%
\includegraphics{figs/CapWen.pstex}%
\end{picture}%
%
%
\setlength{\unitlength}{3158sp}%
\begingroup\makeatletter\ifx\SetFigFont\undefined%
\gdef\SetFigFont#1#2#3#4#5{%
  \reset@font\fontsize{#1}{#2pt}%
  \fontfamily{#3}\fontseries{#4}\fontshape{#5}%
  \selectfont}%
\fi\endgroup%
\begin{picture}(2666,771)(347,-1123)
\put(526,-811){\makebox(0,0)[lb]{\smash{{\SetFigFont{10}{12.0}{\rmdefault}{\mddefault}{\updefault}{\color[rgb]{0,0,0}$w$}%
}}}}
\end{picture}%
$}
A cap can slide through a wen, hence a capped wen disappears, as shown
on the right, to be denoted \glost{CW}.

\vspace{3mm}

\parpic[l]{\parbox{1.6in}{
  \begin{picture}(0,0)%
\includegraphics{figs/LongUparrow.pstex}%
\end{picture}%
%
%
\setlength{\unitlength}{3158sp}%
\begingroup\makeatletter\ifx\SetFigFont\undefined%
\gdef\SetFigFont#1#2#3#4#5{%
  \reset@font\fontsize{#1}{#2pt}%
  \fontfamily{#3}\fontseries{#4}\fontshape{#5}%
  \selectfont}%
\fi\endgroup%
\begin{picture}(84,2874)(2059,-2023)
\end{picture}%

  \hspace{3mm}
  \includegraphics[angle=90, height=5 cm]
    {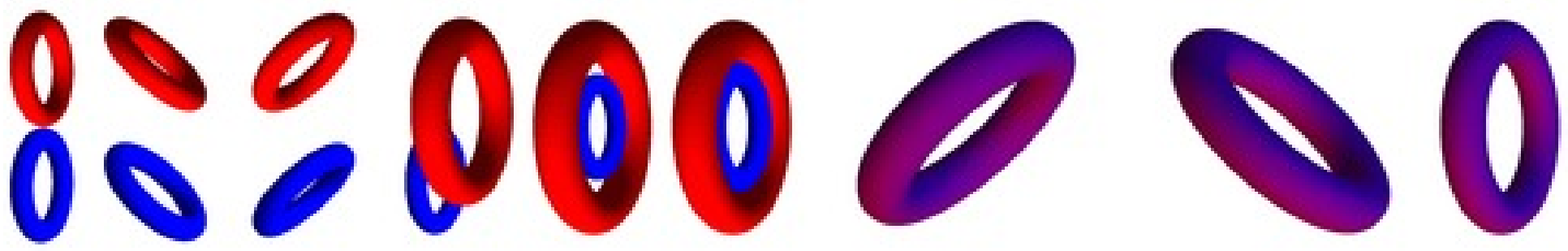}
  \hspace{3mm}\raisebox{2.3cm}{$\leftrightarrow$}
  \raisebox{-2mm}{\includegraphics[angle=90, height=5 cm]
    {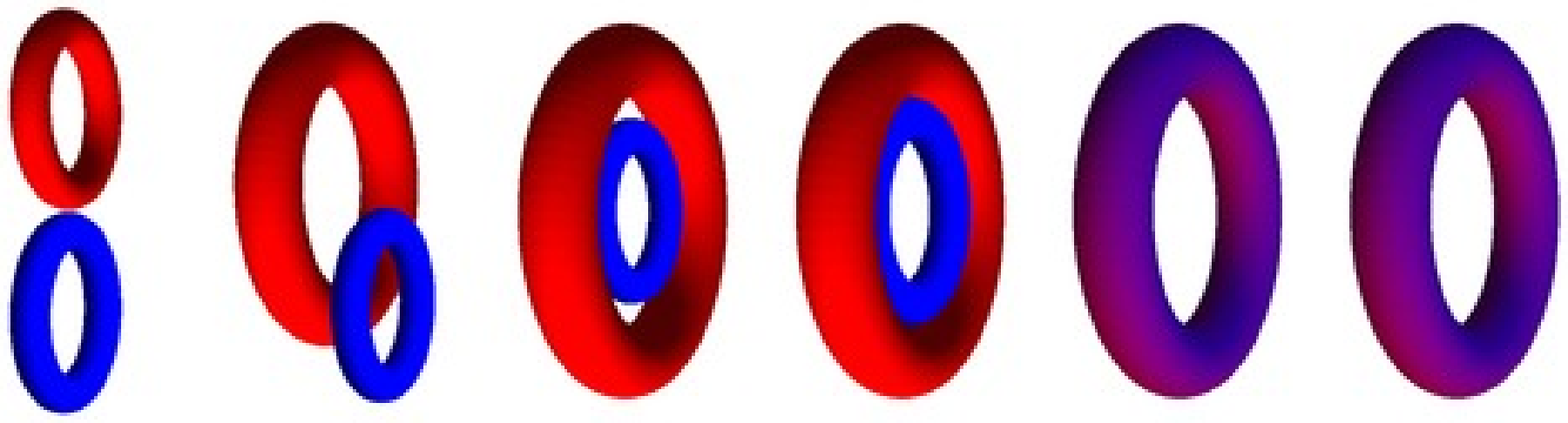}}
  \newline
  \begin{picture}(0,0)%
\includegraphics{figs/VertexWen.pstex}%
\end{picture}%
%
%
\setlength{\unitlength}{3158sp}%
\begingroup\makeatletter\ifx\SetFigFont\undefined%
\gdef\SetFigFont#1#2#3#4#5{%
  \reset@font\fontsize{#1}{#2pt}%
  \fontfamily{#3}\fontseries{#4}\fontshape{#5}%
  \selectfont}%
\fi\endgroup%
\begin{picture}(2096,1013)(2464,-673)
\put(3451,-211){\makebox(0,0)[lb]{\smash{{\SetFigFont{10}{12.0}{\rmdefault}{\mddefault}{\updefault}{\color[rgb]{0,0,0}$\leftrightarrow$}%
}}}}
\end{picture}%

}}
The last wen relation describes the interaction of wens and
vertices. Recall that there are four types of vertices with the same
strand orientation: among the bottom two bands (in the pictures on
the left) there is a non-filled and a filled band (corresponding to
over/under in the strand diagrams), meaning the ``large'' ring and the
``small'' one which flies into it before they merge.  Furthermore, there
is a top and a bottom band (among these bottom two, with apologies for the
ambiguity in overusing the word bottom): this denotes the fly-in direction
(flying in from below or from above). Conjugating a vertex by three wens
switches the top and bottom bands, as shown in the figure on the left:
if both rings flip, then merge, and then the merged ring flips again,
this is homotopic to no flips, except the fly-in direction (from below
or from above) has changed. We are going to denote this relation by
\glost{TV}, for ``twisted vertex''.

The auxiliary operations are the same as for $\wTFo$: orientation
switches, adjoints, deletion of long linear strands, cap unzips
and unzips\footnote{%
We need not specify how to unzip an edge $e$ that carries a wen. To unzip
such $e$, first use the TV relation to slide the wen off $e$.}.
Thus, informally we can say that \linebreak $\wTF=(\wTFo +
\text{wens})/\text{FR},W^2,\text{CW},\text{TV}$.

\subsubsection{The projectivization}\label{subsubsec:AwWen}
The projectivization of $\wTF$ (still denoted $\calA^{sw}$) is the same as the projectivization
for $\wTFo$ but with the wen added as a generator (a degree 0 skeleton feature), and with extra relations
describing the behaviour of the wen. Of course, the relations describing the interaction of wens
with the other skeleton features ($W^2$, TV, and CW) still apply, as well as the old RI, $\aft$, and TC relations.

In addition, the Flip Relations FR imply that wens ``commute'' with
arrow heads, but ``anti-commute'' with tails. We also call these \glost{FR}
relations:
\begin{center}
  \input figs/WenRel.pstex_t
\end{center}

\subsubsection{The homomorphic expansion}\label{subsubsec:ZwithWen}
The goal of this section is to prove that there exists a homomorphic expansion $Z$ for $\wTF$. This involves solving a similar system of
equations to Section \ref{subsec:wTFExpansion}, but with an added unknown for the value of the wen, as well as added equations arising from
the wen relations.
Let $\glos{W} \in \calA(\uparrow_1)$ denote the $Z$-value of the wen, and let us agree that the arrow diagram $W$ always appears just above the
wen on the skeleton. In fact, we are going to show that there exists a homomorphic expansion with $W=1$.

As two consecutive wens on the skeleton cancel, we obtain the equation shown in the picture and explained below:
\begin{center}
 \input{figs/WenSquare.pstex_t}
\end{center}
The $Z$-value of two consecutive wens on a strand is a skeleton wen followed by $W$ followed by a skeleton wen and another $W$. Sliding the
bottom $W$ through the skeleton wen ``multiplies each tail by $(-1)$''. Let us denote this operation by ``bar'', i.e. for an arrow diagram $D$,
$\overline{D}=D\cdot(-1)^{\#\text{ of tails in }D}$. Cancelling the two skeleton wens, we obtain $\overline{W}W=1$. If $W=1$ then this 
equation is certainly satisfied.

Now recall the Twisted Vertex relation of Section \ref{subsubsec:WenRels}. 
Note that the negative the $Z$-value of the vertex on the right hand side of the relation can be 
written as $S_1S_2A_1A_2(V)=\overline(V)$.
(Compare with Remark~\ref{rem:SwitchingVertices}.)
On the other hand, applying $Z$ to the left hand side of the relation, assuming $W=1$, we get:
\begin{center}
 \input figs/VbarEquation.pstex_t
\end{center}
Thus, the equation arising from the twisted vertex relation with $W=1$ is automatically satisfied.

The CW (Capped Wen) relation says that a cap can slide through a wen. The value of the wen is $1$, but the wen as a skeleton feature anti-commutes
with tails (this is the Flip Relation of Section \ref{subsubsec:wTFProjRels}). The value of the cap $C$ is made up of only wheels (Lemma \ref{lem:CapIsWheels}),
hence the CW relation translates to the equation $\overline{C}=C$, which is equivalent to saying that $C$ consists only of even wheels. 

The fact that this
is possible follows from Proposition 6.2 of \cite{AlekseevTorossian:KashiwaraVergne}: the value of the cap is $C=e^c$, where can be set to $c=-\frac{a}{2}$, as 
explained in the proof of Theorem \ref{thm:ATEquivalence}. Here $a$ is such that $\tilde{\delta}(a)=jF$ as in Equation \ref{eq:ATKVEqns}. A power series $f$
so that $a=\tr f$ (where $\tr$ is the trace which turns words into cyclic words) is called the Duflo function of $F$. In Proposition 6.2 Alekseev and Torossian 
show that the even part of $f$ is $\frac{1}{2}\frac{ln(e^{x/2}-e^{-x/2})}{x}$, and that for any $f$ with this even part there is a corresponding solution $F$
of the generalized $KV$ problem. In particular, $f$ can be assumed to be even, namely the power series above, and hence it can be guaranteed that $C$ consists
of even wheels only. Thus we have proven the following:

\begin{theorem}
There exists a group-like homomorphic expansion $Z: \wTF \to \calA^{sw}$. \qed
\end{theorem}

\draftcut
\subsection{The relationship with $u$-Knotted Trivalent Graphs}
\label{subsec:KTG}
The ``$u$sual'', or classical topological objects corresponding to
$\wTF$ are loosely speaking Knotted Trivalent Graphs, or \glost{KTGs}.
A trivalent graph is a graph with three edges meeting at each vertex,
equipped with a cyclic orientation of the three half-edges at each
vertex. KTGs are framed embeddings of trivalent graphs into $\bbR^3$,
regarded up to isotopies. The skeleton of a
KTG is the trivalent graph (as a combinatorial object) behind it.  For a
detailed introduction to KTGs see for example \cite{Bar-NatanDancso:KTG}.
Here we only recall the most important facts. The reader might recall
that in Section \ref{sec:w-knots} we only dealt with long $w$-knots,
as the $w$-theory of round knots is essentially trivial (see Theorem
\ref{thm:AwCirc}). A similar issue arises with ``$w$-knotted trivalent
graphs''. Hence, the space we are really interested in is ``long KTGs'',
in other words trivalent $(1,1)$-tangles whose ``top end'' is connected
to the ``bottom end'' by some path along the tangle.

\parpic[r]{\input{figs/UnzipAndInsertion.pstex_t}}
Long KTGs form an algebraic structure with the operations orientation
\linebreak switch; edge unzip (as shown on the right); and tangle
insertion (I.e., inserting a small copy of a $(1,1)$-tangle $S$ into
the middle of some strand of a $(1,1)$-tangle $T$, also shown on
the right. It is a slightly weaker operation than the connected sum
of~\cite{Bar-NatanDancso:KTG}).  The projectivization of the space
of long KTGs is the space $\glos{\calA^u}$ of chord diagrams on long
trivalent graph skeleta, modulo the $\glos{4T}$ and vertex invariance
(VI) relations. The induced operations on $\calA^u$ are as expected:
orientation switch multiplies a chord diagram by $(-1)$ to the number
of chord endings on the edge.  Edge unzips $u_e$ maps a chord diagram
with $k$ chord endings on the edge $e$ to a sum of $2^k$ diagrams where
each chord ending has a choice between the two daughter edges. Finally,
tangle insertion induces the insertion of chord diagrams.

\parpic[r]{\input{figs/glitch.pstex_t}}
In \cite{Bar-NatanDancso:KTG} the authors prove that there is no
\emph{homomorphic} expansion for KTGs. This theorem, as well as the proof,
applies to long KTGs with slight modifications. There is a well-known
expansion constructed by extending the Kontsevich integral to KTGs and
renormalizing at the vertices. There are several constructions that do
this (\cite{MurakamiOhtsuki:KTGs}, \cite{ChepteaLe:EvenAssociator},
\cite{Dancso:KIforKTG}), and not all of these are ``compatible''
with a corresponding $Z^w$. For now, let us choose one (any) such
expansion and following the notation of \cite{Bar-NatanDancso:KTG}
denote it by $Z^{old}$.  It turns out that any of the above $Z^{old}$
is almost homomorphic but not quite: they all intertwine the orientation
switch, strand delete and tangle composition operations with their
chord-diagrammatic counterparts, but commutativity with unzip fails by
a controlled amount, as shown on the right. Here $\glos{\nu}$ denotes
the ``invariant of the unknot'', the value of which was conjectured in
\cite{Bar-NatanGaroufalidisRozanskyThurston:WheelsWheeling} and proven
in \cite{Bar-NatanLeThurston:TwoApplications}.

In \cite{Bar-NatanDancso:KTG} the authors fix this anomaly by slightly
changing the space of KTGs and adding some extra combinatorics (``dots'' on the edges), and construct a
homomorphic expansion for this new space by a slight adjustment of $Z^{old}$. Here we are going
to use a similar but different adjustment of the space of long KTGs, namely breaking the symmetry of the vertices
and restricting the domain of unzip. 

\begin{figure}[h!]
\input{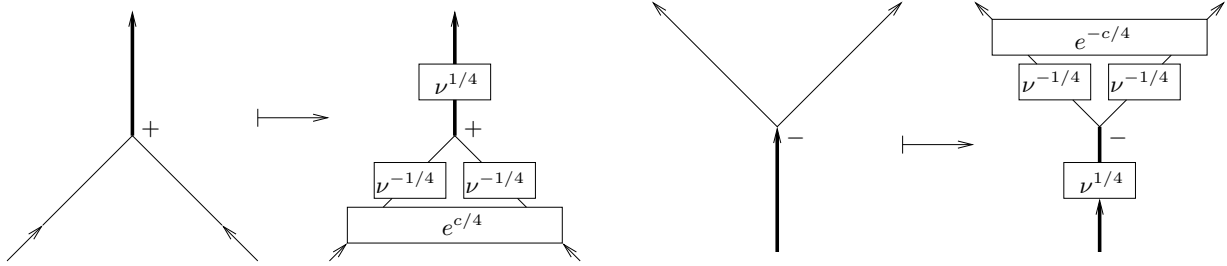}
\caption{Normalizations for $Z^u$ at the vertices.}\label{fig:ZatVertices}
\end{figure}

In this model, denoted by $\glos{\sKTG}$ for ``signed long KTGs'', each
vertex has a distinguished edge coming out of it (denoted by a thick
line in Figure~\ref{fig:ZatVertices}), as well as a sign. Our
pictorial convention will be that a vertex drawn in a ``$\lambda$''
shape with all strands oriented up and the top strand distinguished
is always positive and a vertex drawn in a ``$Y$'' shape with strands
oriented up and the bottom strand distinguished is always negative,
as in Figure~\ref{fig:ZatVertices}.

Orientation switch of either of the non-distinguished strands changes the sign of the vertex, switching the orientation of 
the distinguished strand does not. Unzip of an edge
is only allowed if the edge is distinguished at both of its ends and the vertices at either end are of opposite signs. 

\parpic[r]{\input{figs/ZoldOfTangle.pstex_t}}
The homomorphic expansion $Z^u \colon  \sKTG \to \calA^u$ is computed from $Z^{old}$ as follows. First of all we need to interpret
$Z^{old}$ as an invariant of $(1,1)$-tangles. This is done by connecting the top and bottom ends by a non-interacting long
strand followed by a normalization, as shown on the right. By ``multiplying by $\nu^{-1}$'' we mean that after computing $Z^{old}$
we insert $\nu^{-1}$ on the long strand.

To compute $Z^u$ from $Z^{old}$ the following normalizations are added
near the vertices, as in Figure~\ref{fig:ZatVertices}. Note that in that
figure the symbol $\glos{c}$ denotes a horizontal chord going from left to
right, and $e^{\pm c/4}$ denotes the exponential of $\pm c/4$ in a sense
similar to the exponentiation of arrows in Equation~\eqref{eq:reservoir}.

Checking that $Z^u$ is a homomorphic expansion is a simple calculation
using the almost homomorphicity of $Z^{old}$, which we leave to the
reader. Now let us move on the the question of compatibility between
$Z^u$ and $Z^w$ (from now on we are going to refer to the homomorphic
expansion of $\wTF$ --called $Z$ in the previous section-- as $Z^w$
to avoid confusion).

\def\uwsquare{{\xymatrix{
  \sKTG \ar[d]^{Z^u} \ar[r]^a & \wTF \ar[d]^{Z^w} \\
  \calA^u \ar[r]^\alpha & \calA^{sw}
}}}
\parpic[l]{$\uwsquare$}
There is a map $a\colon  \sKTG \to \wTF$, given by interpreting $\sKTG$
diagrams as $\wTF$ diagrams. In particular, positive vertices (of edge
orientations shown above) are interpreted as the positive $\wTF$ vertex
\begin{picture}(0,0)%
\includegraphics{figs/PlusVertex.pstex}%
\end{picture}%
%
%
\setlength{\unitlength}{3947sp}%
\begingroup\makeatletter\ifx\SetFigFont\undefined%
\gdef\SetFigFont#1#2#3#4#5{%
  \reset@font\fontsize{#1}{#2pt}%
  \fontfamily{#3}\fontseries{#4}\fontshape{#5}%
  \selectfont}%
\fi\endgroup%
\begin{picture}(211,211)(3439,-1160)
\end{picture}%
 and negative vertices as the negative
\begin{picture}(0,0)%
\includegraphics{figs/NegVertex.pstex}%
\end{picture}%
%
%
\setlength{\unitlength}{3158sp}%
\begingroup\makeatletter\ifx\SetFigFont\undefined%
\gdef\SetFigFont#1#2#3#4#5{%
  \reset@font\fontsize{#1}{#2pt}%
  \fontfamily{#3}\fontseries{#4}\fontshape{#5}%
  \selectfont}%
\fi\endgroup%
\begin{picture}(249,249)(3514,-1498)
\end{picture}%
.  The induced map $\alpha\colon  \calA^u
\to \calA^{sw}$ is defined as in Section \ref{subsec:sder}, that is,
$\alpha$ maps each chord to the sum of its two possible orientations. Now
we can ask the question whether the square on the left commutes, or more
precisely, whether we can choose $Z^u$ and $Z^w$ so that it does.

As a first step to answering this question, we prove that $\sKTG$ is
finitely generated (and therefore $Z^u$ is determined by its values on
finitely many generators, and these values will later be compared with the
values $V$ and $C$ that determine $Z^w$):

\begin{proposition}\label{prop:sKTGgens}
The algebraic structure $\sKTG$ is finitely generated by the following
list of elements:
\begin{center} \input{figs/sKTGgens.pstex_t} \end{center}
\end{proposition}

\begin{proof}
First of all note that throughout this proof (in fact even in the
statement of the proposition) we are ignoring the issue of strand
orientations. We can do this as orientation switches are allowed in
$\sKTG$ without restriction. We are also going to omit vertex signs from
the pictures given the pictorial convention stated before.

We need to prove that any $\sKTG$ (call it $G$) can be built from the
generators above using $\sKTG$ operations. To show this, consider a Morse
drawing of $G$, that is, a planar projection of $G$ with a height function
so that all singularities along the strands are Morse and so that every
``feature'' of the projection (local minima and maxima, crossings and
vertices) occurs at a different height.

The idea in short is to decompose $G$ into levels of this Morse drawing where at each level only one ``feature'' occurs. The levels themselves
are not $\sKTG$'s, but we show that the composition of the levels can be achieved by composing their ``closed-up'' $\sKTG$ versions followed
by some unzips. Each feature gives rise to a generator by ``closing up'' extra ends at its top and bottom. We then show that we can
construct each level using the generators and the tangle insert operation.

So let us decompose $G$ into a composition of trivalent tangles, each of which has one ``feature'' and (possibly) some 
straight vertical strands. An example is shown in the figure below. Note that these tangles are not necessarily $(1,1)$-tangles,
and hence need not be elements of $\sKTG$. However, we can turn each of them into a $(1,1)$-tangle by ``closing up'' their tops and bottoms
by arbitrary trees. In the example below we show this for one level of the Morse-drawn $\sKTG$ containing a crossing and two vertical strands.
\begin{center}
\input{figs/MorseTangle.pstex_t}
\end{center}

Now we can compose the $\sKTG$'s obtained from closing up each level, as tangle composition is a special case of tangle insertion. Each tree
that we used to close up the tops and bottoms of levels determines a ``parenthesization'' of the strand endings. If these parenthesizations
match on the top of each level with the bottom of the next, then we can recreate tangle composition of the levels by composing their closed
versions followed by a number of unzips performed on the connecting trees. This is illustrated in the example below, for two consecutive levels
of the $\sKTG$ of the previous example.
\begin{center}
 \input{figs/CombineLevels.pstex_t}
\end{center}

If the trees used to close up consecutive levels correspond to different parenthesizations, then we can use insertion of the left and right associators 
(the last two pictures of the list of generators in the statement of the theorem) to change one parenthesization to match the other. This is 
illustrated in the figure below.
\begin{center}
 \begin{picture}(0,0)%
\includegraphics{figs/Reassociate.pstex}%
\end{picture}%
\setlength{\unitlength}{4144sp}%
\begingroup\makeatletter\ifx\SetFigFont\undefined%
\gdef\SetFigFont#1#2#3#4#5{%
  \reset@font\fontsize{#1}{#2pt}%
  \fontfamily{#3}\fontseries{#4}\fontshape{#5}%
  \selectfont}%
\fi\endgroup%
\begin{picture}(6729,2364)(7144,-2863)
\put(12016,-1636){\makebox(0,0)[lb]{\smash{{\SetFigFont{12}{14.4}{\rmdefault}{\mddefault}{\updefault}{\color[rgb]{0,0,0}unzips}%
}}}}
\end{picture}%

\end{center}

So far we have shown that $G$ can be assembled from closed versions of the levels in its Morse drawing. The closed versions of the levels of $G$ are
simpler $\sKTG$'s, and it remains to show that these can be obtained from the generators using $\sKTG$ operations. 

\parpic[r]{\input{figs/SrtandsToBubble.pstex_t}}
Let us examine what each level
might look like. First of all, in the absence of any ``features'' a level might be a single strand, in which case it is the first generator itself. 
Two parallel strands when closed up become the ``bubble'', as shown on the right.

Now suppose that a level consists of $n$ parallel strands, and that the trees used to close it up on the top and bottom are horizontal mirror images 
of each other, as shown below (if not, then this can be achieved by associator insertions and unzips). We want to show that this $\sKTG$ can
be obtained from the generators using $\sKTG$ operations. Indeed, this can be achieved by repeatedly inserting bubbles into a bubble, as shown:
\begin{center}
\input{figs/BubbleInsertions.pstex_t}
\end{center}

A level consisting of a single crossing becomes a left or right twist
when closed up (depending on the sign of the crossing). Similarly, a
single vertex becomes a bubble. A level can not contain a single minimum
or maximum by itself, since we required that the top end of an $\sKTG$
be connected to its bottom end via a path. Hence, any minimum or maximum
must be accompanied by at least one through strand. A maximum with one
through strand becomes the balloon after closing up, and a minimum with
one through strand becomes the noose.

It remains to see that the $\sKTG$'s obtained when closing up simple features accompanied by more through strands can be built from the generators.
This is achieved by inserting the corresponding generators into nested bubbles (bubbles inserted into bubbles), as in the example shown below.
Recall that the trees (parenthesizations) used for the closing up process can be changed arbitrarily by inserting associators and unzipping, and
hence we are free to use the most convenient tree in the example below. This completes the proof.
\begin{center}
 \input{figs/MoreStrands.pstex_t}
\end{center}
\qed 
\end{proof}

We are now equipped to answer the main question of this section:

\noindent\parbox[t]{4.25in}{
  \begin{theorem}\label{thm:ZuwCompatible}
  There exists a homomorphic expansion for the combined algebraic structure
  $\left(\sKTG\overset{a}{\longrightarrow}\wTF\right)$. In other words,
  there exist homomorphic expansions $Z^u$ and $Z^w$ for which the square
  on the right commutes.
  \end{theorem}
}\hfill\begin{minipage}[t]{2in}
  \begin{equation} \label{eq:uwcompatibility}
    \begin{array}{c} \uwsquare \end{array}
  \end{equation}
\end{minipage}


Before moving on to the proof let us state and prove the following Lemma,
to be used repeatedly in the proof of the theorem.

\begin{lemma}\label{lem:TreesAndUnitarity}
If $a$ and $b$ are group-like elements in $\calA^{sw}(\uparrow_n)$, then $a=b$ if and only if $\pi(a)=\pi(b)$ and $aa^*=bb^*$. Here $\pi$
is the projection induced by $\pi\colon  \calP^w(\uparrow_n) \to \tder_n \oplus \fraka_n$ (see Section \ref{subsec:ATSpaces}),
and $*$ refers to the adjoint map of Definition \ref{def:Adjoint}. In the notation of this section $*$ is applying the adjoint $A$ on all strands.
\end{lemma}

\begin{proof}
Write $a=e^we^{uD}$ and $b=e^{w'}e^{uD'}$, where $w\in \attr_n$, $D\in \tder_n\oplus \fraka_n$ 
and $u\colon  \tder_n\oplus \fraka_n \to \calP_n$ is the ``upper'' map of 
Section \ref{subsec:ATSpaces}. Assume that $\pi(a)=\pi(b)$ and $aa^*=bb^*$. Since $\pi(a)=e^D$ and $\pi(b)=e^{D'}$, we conclude
that $D=D'$. Now we compute $aa^*=e^we^{uD}e^{-lD}e^w=e^we^{j(D)}e^w,$ where $j\colon  \tder_n \to \attr_n$ is the map defined in Section 5.1 of 
\cite{AlekseevTorossian:KashiwaraVergne} and discussed in \ref{prop:Jandj} of this paper. Now note that both $w$ and $j(D)$ are elements of 
$\attr_n$, hence they commute, so $aa^*=e^{2w+j(D)}$. Thus, $aa^*=bb^*$ means that $e^{2w+j(D)}=e^{2w'+j(D)}$, which implies that $w=w'$ and
$a=b$. \qed
\end{proof}

\noindent{\em Proof of Theorem \ref{thm:ZuwCompatible}.}
Since $\sKTG$ is finitely generated, we only need to check that the
square~\eqref{eq:uwcompatibility} commutes for each of the generators.

\noindent{\em Proof of commutativity of~\eqref{eq:uwcompatibility}
for the strand and the bubble.} For the single strand commutativity is
obvious: both the $Z^u$ and $Z^w$ values are trivial.

\parpic[r]{\input{figs/ValueOfTheBubble.pstex_t}}
We claim that the $Z^u$ value of the bubble is also trivial. By connecting
the top and bottom of the bubble we obtain a ``theta-graph'', and
$Z^{old}$ of a theta graph has $\nu^{1/2}$ on each strand, as shown
on the right (for a computation see \cite{Bar-NatanDancso:KTG}
for example). After applying the vertex normalizations of
Figure~\ref{fig:ZatVertices}, everything cancels, so the $Z^u$-value
of the bubble is trivial.  As for $Z^w$, the value of the bubble is
$V_-V$, as shown, which equals to 1 by the Unitarity property of $V$,
Equation~\eqref{eq:unitarity}. This proves the commutativity of the
square for bubbles. \qed

\noindent{\em Proof of commutativity of~\eqref{eq:uwcompatibility}
for the twists.} First note that the $Z^u$-value of the right twist
is $R^u=e^{c/2}$, where $c$ denotes a single chord between the two
twisted strands (see \cite{Bar-NatanDancso:KTG} for details). Hence the
commutativity of $Z^u$ and $Z^w$ for the right twist is equivalent to the
``Twist Equation'' $\alpha(R^u)=V^{-1}RV^{21}$,
 where $R=e^{a_{12}}$ is the $Z^w$-value of the crossing, that is, the exponential of a single arrow pointing from strand 1 to strand 2. By definition
of $\alpha$, $\alpha(R^u)=e^{\frac{1}{2}(a_{12}+a_{21})}$, where $a_{12}$ and $a_{21}$ are single arrows pointing from strand 1 to 2 and 2 to 1, 
respectively. So the Twist Equation becomes
\begin{equation}\label{eq:twist}
 e^{\frac{1}{2}(a_{12}+a_{21})}=V^{-1}RV^{21}.
\end{equation}
If $V$ is to give rise to a homomorphic expansion $Z^w$ that is compatible with $Z^u$, then $V$ has to satisfy the Twist Equation in addition to 
the previous equations (\ref{eq:HardR4}),(\ref{eq:unitarity}) and (\ref{eq:CapEqn}). 
To prove that such a $V$ exists, we use Lemma \ref{lem:TreesAndUnitarity}.
Lemma \ref{lem:TreesAndUnitarity} implies that it is enough
to find a $V$ which satisfies the Twist Equation ``on tree level'' (i.e., after applying $\pi$), and for which the adjoint condition of the Lemma holds.

Let us start with the adjoint condition. Multiplying the left hand side of the Twist Equation by its adjoint, we get 
$$e^{\frac{1}{2}(a_{12}+a_{21})}(e^{\frac{1}{2}(a_{12}+a_{21})})^*=e^{\frac{1}{2}(a_{12}+a_{21})}e^{-\frac{1}{2}(a_{12}+a_{21})}=1.$$
As for the right hand side, we have to compute $V^{-1}RV^{21}(V^{21})^*R^*(V^{-1})^*$. Since $V$ is unitary (Equation (\ref{eq:unitarity})), $VV^*=V\cdot A_1A_2(V)=1$.
Now $R=e^{a_{12}}$, so $R^*=e^{-a_{12}}=R^{-1}$, hence the expression on the right hand side also simplifies to 1, as needed.  

As for the ``tree level'' of the Twist Equation, recall that in Section \ref{subsec:wTFExpansion} we deduced the existence of a solution to all the previous
equations from Alekseev and Torossian's solution $F\in \TAut_2$ to the Kashiwara--Vergne
equations \cite{AlekseevTorossian:KashiwaraVergne}. 
We produced $V$ from $F$ by setting $F=e^{D^{21}}$ with $D\in \tder_2^s$, $b:=\frac{-j(F)}{2}\in \attr_2$ and $V:=e^be^{uD}$, so $F$ is ``the tree part'' of $V$, 
up to re-numbering strands. 
Substituting this into the Twist Equation we obtain:
\begin{equation}\label{eq:TwistWithF}
 e^{\frac{1}{2}(a_{12}+a_{21})}=e^{-uD}e^{-b}e^{a_{12}}e^{b^{21}}e^{uD^{21}}.
\end{equation}
Applying $\pi$, we get
$$
 e^{\frac{1}{2}(a_{12}+a_{21})}=e^{-uD}e^{a_{12}}e^{uD^{21}}=(F^{21})^{-1}e^{a_{12}}F.
$$  
The existence of a solution $F$ of the KV equations which also satisfies
the above is equivalent to the existence of ``symmetric solutions of
the Kashiwara-Vergne problem'' discussed and proven in Sections 8.2
and 8.3 of \cite{AlekseevTorossian:KashiwaraVergne} (note that in
\cite{AlekseevTorossian:KashiwaraVergne} $R$ denotes $e^{a_{21}}$). \qed

\noindent{\em Proof of commutativity of~\eqref{eq:uwcompatibility} for the
associators.} Let us recall that a Drinfel'd associator is a group-like
element of $\calA^u(\uparrow_3)$ satisfying the so-called pentagon and
positive and negative hexagon equations, as well as a non-degeneracy
and mirror skew-symmetry property.  For a detailed explanation see
Section 4 of \cite{Bar-NatanDancso:KTG}; associators were first defined
in \cite{Drinfeld:QuasiHopf}.  The $Z^u$-value of the generator shown
in the statement of Proposition~\ref{prop:sKTGgens} called ``right
associator'' is a Drinfel'd associator. The proof of this statement
is the same as the proof of Theorem 4.2 of \cite{Bar-NatanDancso:KTG},
with minor modifications. (I.e., the graphs have positive and negative
vertices as opposed to ``dots and crosses'' on the edges. Note that the
vertex re-normalizations for the four vertices of an associator cancel
each other out). Let us call this associator $\glos{\Phi}$.

What we need to show is that there exists a $V$ satisfying all previous equations including the Twist Equation (\ref{eq:twist}),
so that 
\begin{equation}\label{eq:AssociatorAndV}
  \alpha(\Phi)=V_-^{(12)3}V_-^{12}V^{23}V^{1(23)}
  \qquad\text{in}\qquad
  \calA^{sw}(\uparrow_3),
\end{equation}
where $\alpha\colon \calA^u \to \calA^{sw}$ is the map defined in Section
\ref{subsec:sder}, and keeping in mind that $V_-=V^{-1}$. The reasoning
behind this equation is shown in the figure below.
\begin{center}
 \input{figs/Associator.pstex_t}
\end{center}

We proceed in a similar manner as we did for the Twist Equation, treating the ``tree and wheel parts'' separately
using Lemma \ref{lem:TreesAndUnitarity}. As $\Phi$ is by
definition group-like, let us denote $\Phi=:\glos{e^\phi}$. 

First we verify that the ``wheel level'' adjoint condition holds. Starting
with the right hand side of Equation~\eqref{eq:AssociatorAndV}, the
unitarity $VV^*=1$
of $V$ implies that
\[ V_-^{(12)3}V_-^{12}V^{23}V^{1(23)}
  (V^{1(23)})^*(V^{23})^*(V_-^{12})^*(V_-^{(12)3})^*=1.
\]
For the left hand side of~\eqref{eq:AssociatorAndV} we need to study
$e^{\alpha(\phi)}(e^{\alpha(\phi)})^*$ and show that it equals
1 as well. This is assured if we pick a $Z^u$ for which $\Phi$
is a group-like horizontal chord associator (possible for example
using~\cite{ChepteaLe:EvenAssociator}, as mentioned at the beginning of
this section). Indeed restricted to the $\alpha$-images of horizontal
chords $*$ is multiplication by $-1$, and as it is an anti-Lie morphism,
this fact extends to the Lie algebra generated by $\alpha$-images
of horizontal chords. Hence $e^{\alpha(\phi)}(e^{\alpha(\phi)})^*
=e^{\alpha(\phi)}e^{\alpha(\phi)^*}=e^{\alpha(\phi)}e^{-\alpha(\phi)}=1$.

On to the tree part. Applying $\pi$ to Equation (\ref{eq:AssociatorAndV})
we obtain
\begin{multline}\label{eq:ATPhiandV}
  e^{\pi\alpha(\phi)}
  = (F^{3(12)})^{-1} (F^{21})^{-1} F^{32} F^{(23)1}
  =e^{-D^{(12)3}}e^{-D^{12}}e^{D^{23}}e^{D^{1(23)}} \\
  \text{in }
  \glos{\SAut_3}:=\exp(\sder_3)\subset\TAut_3.
\end{multline}
This is Equation (26) of \cite{AlekseevTorossian:KashiwaraVergne},
up to re-numbering strands 1 and 2 as 2 and 1\footnote{Note that
in \cite{AlekseevTorossian:KashiwaraVergne} ``$\Phi'$ is an
associator'' means that $\Phi'$ satisfies the pentagon equation,
mirror skew-symmetry, and positive and negative hexagon equations
in the space $\SAut_3$. These equations are stated in
\cite{AlekseevTorossian:KashiwaraVergne} as equations (25), (29),
(30), and (31), and the hexagon equations are stated with strands 1
and 2 re-named to 2 and 1 as compared to \cite{Drinfeld:QuasiHopf}
and \cite{Bar-NatanDancso:KTG}. This is consistent with $F=e^{D^{21}}$.}.
To prove it in our context, we need the following fact from
\cite{AlekseevTorossian:KashiwaraVergne} (their Theorem 7.5, Propositions
9.2 and 9.3 combined):

\begin{fact}
If $\Phi'=e^{\phi'}$ is an associator in $\SAut_3$ so that $j(\Phi')=0$\footnote{The condition
$j(\phi')=0$ is equivalent to the condition $\Phi\in KRV^0_3$ in \cite{AlekseevTorossian:KashiwaraVergne}.
The relevant definitions in \cite{AlekseevTorossian:KashiwaraVergne} can be found in Remark 4.2 and at the bottom of 
page 434 (before Section 5.2).}
then Equation~(\ref{eq:ATPhiandV}) has a solution $F=e^{D^{21}}$ which is
also a solution to the KV equations, and all such solutions are symmetric
(i.e. verify the Twist Equation (\ref{eq:TwistWithF})). \qed
\end{fact}

To use this fact, we need to show that $\Phi':=\pi\alpha(\Phi)$ is an
associator in $\SAut_3$ and that $j(\Phi')=j(\pi\alpha(\Phi))=0$. The
latter is the unitarity of $\Phi$ which is already proven. The
former follows from the fact that $\Phi$ is an associator and the fact
(Theorem~\ref{thm:sder}) that the image of $\pi\alpha$ is contained in
$\sder$ (ignoring degree 1 terms, which are not present in an associator
anyway). \qed

\parpic[r]{\input{figs/uValueOfTheNoose.pstex_t}}
\noindent{\em Proof of commutativity of~\eqref{eq:uwcompatibility}
for the balloon and the noose.}
Connecting the top and bottom end of the noose picture creates a
``dumbbell graph'', and $Z^{old}$ of the dumbbell is a $\nu$ placed on
each of the circles with nothing on the line connecting them. Applying
the vertex normalizations and the $\nu^{-1}$ normalization on the
long strand, we obtain that $Z^u$ of the noose has chords only on the
circle, namely $e^{-c/4}\nu^{1/2}$, as shown on the right. We leave it
to the reader to check this, keeping in mind the fact that in $\calA^u$,
any chord diagram with chord endings on a bridge in the graph (i.e.,
an edge whose deletion increases the number of connected components)
is zero. Also keep in mind that the bottom vertex is not a positive
vertex: the orientation of the left strand is switched, so we have
to apply an orientation switch operation of that strand to the value
of the normalization. As $S(\nu)=\nu$, this only affects the sign of
the exponent.  A similar computation can be done for the balloon, 
where the result is $e^{c/4}\nu^{1/2}$ on the circle.

\parpic[r]{\input{figs/wValueOfTheNoose.pstex_t}}
$Z^w$ on the other hand assigns a $V$ value to each vertex, one of which
has its first strand orientation switched as shown in the figure on the
right. The top copy of $V$ appearing there cancels: pushing arrow heads
and tails onto the noose using VI results in two terms that have opposite
signs but are otherwise equal (we can slide arrow heads/tails across the
$S_1(V)$ term as anything concentrated on one strand is a combination of
wheels and $D_A$ arrows, and we can slide across these using $\aft$/TC).

Hence, what we need to show is that the two equations 
below hold, arising from the noose and the balloon, respectively.
\begin{center}
 \input{figs/NooseEquations.pstex_t}
\end{center}

\begin{figure}
  \input{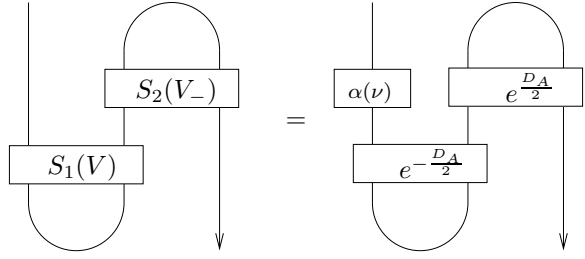}
\caption{The product equation.}\label{fig:NooseBalloonProduct}
\end{figure}

We will start by proving that the product of these two equations,
shown in Figure~\ref{fig:NooseBalloonProduct}, is satisfied. Note
that any local (small) arrow diagram on a single strand is central in
$\calA^{sw}(\uparrow_n)$: a diagram on one strand can be written in
terms of only wheels and isolated arrows, both of which commute with
both arrow heads and tails by $\aft$ and $TC$. Hence we can slide and
merge the $\alpha(\nu)$ terms as we wish.

\begin{figure}[h]
  \input{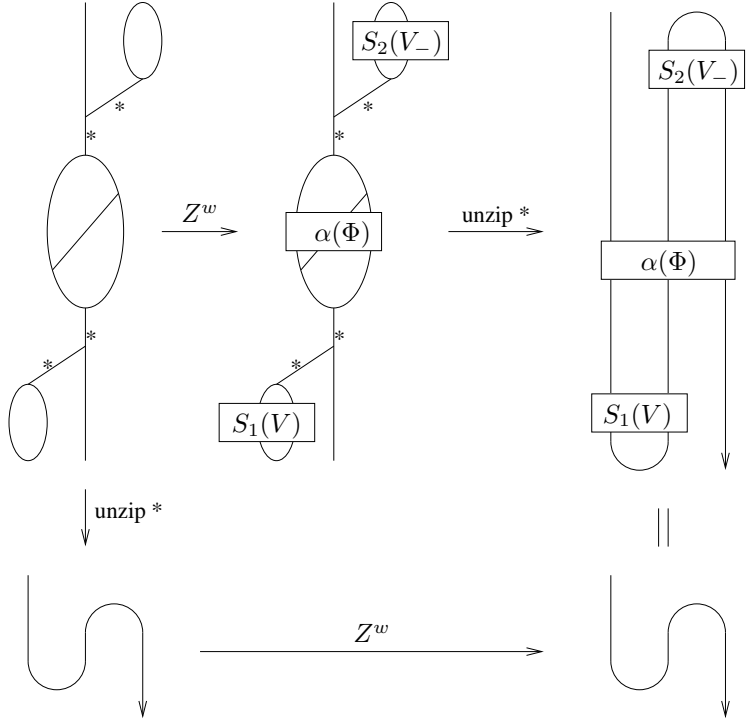}
\caption{Unzipping a noose and a balloon to a squiggle.}
\label{fig:NooseBalloonWProof}
\end{figure}

To show that the product equation is satisfied, consider 
Figure~\ref{fig:NooseBalloonWProof}. We start with the $\wTF$ on the top left and either apply $Z^w$ followed by unzipping the 
edges marked 
by stars, or first unzip the same edges and then apply $Z^w$. We 
use that by the compatibility with associators, $Z^w$ of an associator is $\alpha(\Phi)$. 
Since $Z^w$ is homomorphic, the two results in the bottom right corner must agree. 
(Note that two of the four unzips we perform are ``illegal'',
as the strand directions don't match. However, it is easy to get around this issue by inserting small bubbles at the top of the balloon and the bottom 
of the noose, and switching the appropriate edge orientations before and after the unzips. The $Z^w$-value of a bubble is 1, hence this will not effect 
the computation and so we ignore the issue for simplicity.)

We conclude that to prove that the product equation of Figure
\ref{fig:NooseBalloonProduct} is satisfied, it is enough to show that
the left equality of Figure~\ref{fig:NooseBalloonReduced} holds. Note
that in Figure \ref{fig:NooseBalloonReduced} the inverse is taken in
$\calA^{sw}(\uparrow_1)$. As both sides of this equation are in the image
of $\alpha$, it is enough to prove the pre-image of the equation in $\calA^u$,
as shown on the right of Figure~\ref{fig:NooseBalloonReduced}.
That equation in turn follows from an argument identical to that
of Figure~\ref{fig:NooseBalloonWProof} but carried out in $\sKTG$
and $\calA^u$, using that $Z^u$ is homomorphic with respect to tangle
insertion. This finishes the proof that the product of the noose and
balloon equations holds.

\begin{figure}
  \input figs/ReducedNooseEq.pstex_t
\caption{The reduced noose and balloon equation.}
\label{fig:NooseBalloonReduced}
\end{figure}

What remains is to show that the noose and balloon equations hold individually. In light of the results so far, it is sufficient to show that
\begin{equation}\label{eq:NooseSymmetry}
 \raisebox{-6mm}{\input figs/NooseSymmetry.pstex_t}
\end{equation}
As stated in Theorem \ref{thm:Aw}, $\calA^{sw}(\uparrow_1)$ is the
polynomial algebra freely generated by the arrow $D_A$ 
and wheels of degrees 2 and higher. Since $V$
is group-like, the ``one-strand version'' of $S_1(V)$ (resp. $S_2(V_-)$)
shown in Equation (\ref{eq:NooseSymmetry}) is an exponential $e^{A_1}$
(resp. $e^{A_2}$) with \linebreak $A_1, A_2 \in \calA^{sw}(\uparrow_1)$. We want to
show that $e^{A_1}=e^{A_2}\cdot e^{-D_A}$, equivalently that $A_1=A_2-D_A$.

\begin{figure}
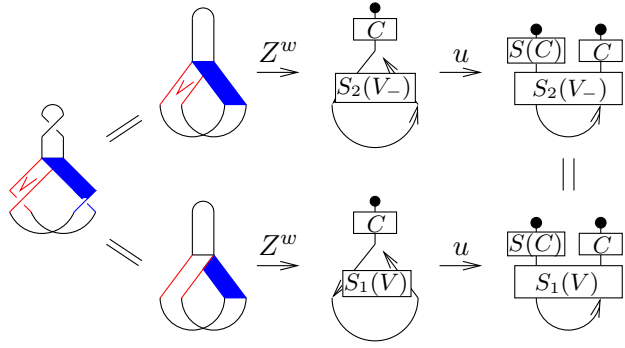

  \input figs/NooseCappedProof.pstex_t
\caption{The proof of Equation (\ref{eq:NooseCapped}). Note that the unzips are ``illegal'', as the strand directions don't match. This can be fixed
by inserting a small bubble at the bottom of the noose and doing a number of orientation switches. As this doesn't change the result or the main argument,
we suppress the issue for simplicity. Equation (\ref{eq:NooseCapped}) is obtained from this result by multiplying by $S(C)^{-1}$ on the bottom and by $C^{-1}$
on the top.}\label{fig:NooseCappedProof}
\end{figure}

In degree 1, this can be done by explicit verification. Let $A_1^{\geq
2}$ and $A_2^{\geq 2}$ denote the degree 2 and higher parts of $A_1$ and
$A_2$, respectively. We claim that capping the strand at both its top
and its bottom takes $e^{A_1}$ to $e^{A_1^{\geq 2}}$, and similarly $e^{A_2}$
to $e^{A_2^{\geq 2}}$. (In other words, capping kills arrows but leaves
wheels un-changed.) This can be proven similarly to the proof of
Lemma~\ref{lem:CapIsWheels}, but using
\[
  F' := \sum_{k_1,k_2=0}^{\infty}
    \frac{(-1)^{k_1+k_2}}{k_1!k_2!}D_A^{k_1+k_2}S_L^{k_1}S_R^{k_2}
\]
in place of $F$ in the proof. What we want to show, then, is that
\begin{equation}\label{eq:NooseCapped}
  \raisebox{-6mm}{\input figs/NooseCapped.pstex_t}
\end{equation}
The proof of this is shown in Figure~\ref{fig:NooseCappedProof}. \qed

Having verified the commutativity of~\eqref{eq:uwcompatibility} for all the
generators of $\sKTG$ appearing in Proposition~\ref{prop:sKTGgens}, we have
concluded the proof of Theorem~\ref{thm:ZuwCompatible}. \qed

Recall from Section~\ref{subsec:sder} that there is no commutative
square linking $Z^u\colon\uT\to\calA^u$ and $Z^w\colon\wT\to\calA^{sw}$,
for the simple reason that the Kontsevich integral for tangles $Z^u$
is not canonical, but depends on a choice of parenthesizations for
the ``bottom'' and the ``top'' strands of a tangle $T$. Yet given
such choices, a tangle $T$ can be ``closed'' as within the proof of
Proposition~\ref{prop:sKTGgens} into an $\sKTG$ which we will denote
$G$. For $G$ a commutativity statement does hold as we have just
proven. The $Z^u$ and $Z^w$ invariants of $T$ and of $G$ differ only
by a number of vertex-normalizations and vertex-values on skeleton-trees
at the bottom or at the top of $G$, and using VI, these values can slide
so they are placed on the original skeleton of $T$. This is summarized
as the following proposition:

\begin{proposition} \label{prop:uwBT} Let $n$ and $n'$ be natural numbers.
Given choices $c$ and and $c'$ of parenthesizations of $n$ and $n'$
strands respectively, there exists invertible elements
$C\in\calA^{sw}(\uparrow_n)$ and $C'\in\calA^{sw}(\uparrow_{n'})$ so
that for any u-tangle $T$ with $n$ ``bottom'' ends and  $n'$ ``top'' ends
we have
\[ \alpha Z^u_{c,c'}(T)=C^{-1}Z^w(aT)C', \]
where $Z^u_{c,c'}$ denotes the usual Kontsevich integral of $T$ with
bottom and top parenthesizations $c$ and $c'$.
\end{proposition}

For u-braids the above proposition may be stated with $c=c'$ and then $C$
and $C'$ are the same.

\clearpage\draftcut
\section{Odds and Ends} \label{sec:OddsAndEnds}

\subsection{What means ``closed form''?} \label{subsec:ClosedForm}

As stated earlier, one of my hopes for this paper is that it will
lead to closed-form formulae for tree-level associators. The
notion ``closed-form'' in itself requires an explanation (see
footnote~\ref{foot:ClosedForm}). Is $e^x$ a closed form expression for
$\sum_{n=0}^\infty\frac{x^n}{n!}$, or is it just an artificial name given
for a transcendental function we cannot otherwise reduce? Likewise,
why not call some tree-level associator $\Phi^\text{tree}$ and now it is
``in closed form''?

For us, ``closed-form'' should mean ``useful for computations''. More
precisely, it means that the quantity in question is an element of some
space $\calAcf$ of ``useful closed-form thingies'' whose elements have
finite descriptions (hopefully, finite and short) and on which some
operations are defined by algorithms which terminate in finite time
(hopefully, finite and short). Furthermore, there should be a finite-time
algorithm to decide whether two descriptions of elements of $\calAcf$
describe the same element\footnote{In our context, if it is hard to
decide within the target space of an invariant whether two elements
are equal or not, the invariant is not too useful in deciding whether
two knotted objects are equal or not.}. It is even better if the said
decision algorithm takes the form ``bring each of the two elements in question
to a canonical form by means of some finite (and hopefully short)
procedure, and then compare the canonical forms verbatim''; if this is the
case, many algorithms that involve managing a large number of elements
become simpler and faster.

Thus for example, polynomials in a variable $x$ are always of closed form,
for they are simply described by finite sequences of integers (which in
themselves are finite sequences of digits), the standard operations on
polynomials ($+$, $\times$, and, say, $\frac{d}{dx}$) are algorithmically
computable, and it is easy to write the ``polynomial equality'' computer
program. Likewise for rational functions and even for rational functions
of $x$ and $e^x$.

On the other hand, general elements $\Phi$ of the space
$\calA^\text{tree}(\uparrow_3)$ of potential tree-level associators
are not closed-form, for they are determined by infinitely many
coefficients. Thus iterative constructions of associators, such
as the one in~\cite{Bar-Natan:NAT} are computationally useful only
within bounded-degree quotients of $\calA^\text{tree}(\uparrow_3)$
and not as all-degree closed-form formulae. Likewise, ``explicit''
formulae for an associator $\Phi$ in terms of multiple $\zeta$-values
(e.g.~\cite{LeMurakami:HOMFLY}) are not useful for computations as it
is not clear how to apply tangle-theoretic operations to $\Phi$ (such as
$\Phi\mapsto\Phi^{1342}$ or $\Phi\mapsto(1\otimes\Delta\otimes 1)\Phi$)
while staying within some space of ``objects with finite description in
terms of multiple $\zeta$-values''. And even if a reasonable space of such
objects could be defined, it remains an open problem to decide whether
a given rational linear combination of multiple $\zeta$-values is equal
to $0$.

\draftcut
\subsection{Arrow Diagrams to Degree 2} \label{subsec:ToTwo} Just as
an example, in this section we study the spaces $\calA^-(\uparrow)$,
$\calA^{s-}(\uparrow)$, $\calA^{r-}(\uparrow)$, $\calP^-(\uparrow)$,
$\calA^-(\bigcirc)$, $\calA^{s-}(\bigcirc)$, and $\calA^{r-}(\bigcirc)$
in degrees $m\leq 2$ in detail, both in the ``v'' case and in the ``w''
case (the ``u'' case has been known since long).

\subsubsection{Arrow Diagrams in Degree 0} There is only
one degree 0 arrow diagram, the empty diagram $D_0$ (see
Figure~\ref{fig:Deg0-2Diagrams}).  There are no relations, and thus
$\{D_0\}$ is the basis of all $\calG_0\calA^-(\uparrow)$ spaces
and its obvious closure, the empty circle, is the basis of all
$\calG_0\calA^-(\bigcirc)$ spaces. $D_0$ is the unit $1$, yet $\Delta
D_0=D_0\otimes D_0=1\otimes 1\neq D_0\otimes 1+1\otimes D_0$, so $D_0$
is not primitive and $\dim\calG_0\calP^-(\uparrow)=0$.

\subsubsection{Arrow Diagrams in Degree 1} \label{subsubsec:DegreeOne}
\begin{figure}
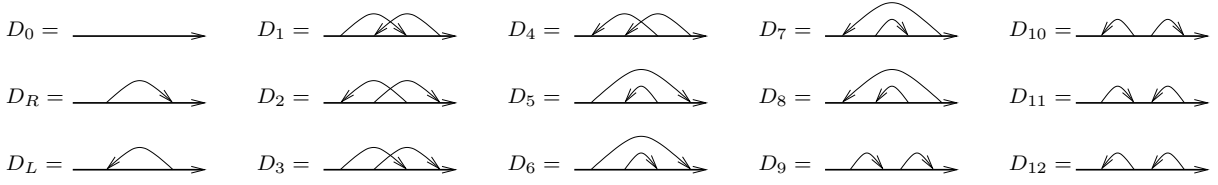

\[ \input figs/Deg0-2Diagrams.pstex_t \]
\caption{The 15 arrow diagrams of degree at most 2.}
\label{fig:Deg0-2Diagrams}
\end{figure}

There is only two degree 1 arrow diagrams, the ``right
arrow'' diagram $D_R$ and the ``left arrow'' diagram $D_L$ (see
Figure~\ref{fig:Deg0-2Diagrams}).  There are no $6T$ relations, and
thus $\{D_R, D_L\}$ is the basis of $\calG_1\calA^-(\uparrow)$. Modulo
RI, $D_L=D_R$ and hence $D_A:=D_L=D_R$ is the single basis element of
$\calG_1\calA^{s-}(\uparrow)$. Both $D_R$ and $D_L$ vanish modulo FI, so
$\dim\calG_1\calA^{r-}(\uparrow)=\dim\calG_1\calA^{r-}(\bigcirc)=0$. Both
$D_R$ and $D_L$ are primitive, so $\dim\calG_1\calP^-(\uparrow)=2$.
Finally, the closures ${\bar D}_R$ and ${\bar D}_L$ of $D_R$ and $D_L$ are
equal, so $\calG_1\calA^{s-}(\bigcirc)=\calG_1\calA^-(\bigcirc)=\langle
{\bar D}_R\rangle=\langle {\bar D}_L\rangle=\langle {\bar D}_A\rangle$.

\subsubsection{Arrow Diagrams in Degree 2} There are 12 degree
2 arrow diagrams, which we denote $D_1,\dots,D_{12}$ (see
Figure~\ref{fig:Deg0-2Diagrams}). There are six $6T$ relations,
corresponding to the 6 ways of ordering the 3 vertical strands that
appear in a $6T$ relation (see Figure~\ref{fig:6T}) along a long line. The
ordering $(ijk)$ becomes the relation $D_3+D_9+D_3=D_6+D_3+D_6$. Likewise,
$(ikj)\mapsto D_6+D_1+D_{11}=D_3+D_5+D_1$, $(jik)\mapsto
D_{10}+D_2+D_6=D_2+D_5+D_3$, $(jki)\mapsto D_4+D_7+D_1=D_8+D_1+D_{11}$,
$(kij)\mapsto D_2+D_7+D_4=D_{10}+D_2+D_8$, and $(kji)\mapsto
D_8+D_4+D_8=D_4+D_{12}+D_4$. After some linear algebra, we find
that $\{D_1, D_2, D_6, D_8, D_9, D_{11}, D_{12}\}$ form a basis
of $\calG_2\calA^v(\uparrow)$, and that the remaining diagrams
reduce to the basis as follows: $D_3=2D_6-D_9$, $D_4=2D_8-D_{12}$,
$D_5=D_9+D_{11}-D_6$, $D_7=D_{11}+D_{12}-D_8$, and $D_{10}=D_{11}$.  In
$\calG_2\calA^{sv}(\uparrow)$ we further have that $D_5=D_6$, $D_7=D_8$,
and $D_9=D_{10}=D_{11}=D_{12}$, and so $\calG_2\calA^{sv}(\uparrow)$
is 3-dimensional with basis $D_1$, $D_2$, and $D_3=\ldots=D_{12}$.
In $\calG_2\calA^{rv}(\uparrow)$ we further have that $D_{5-12}=0$.
Thus $\{D_1, D_2\}$ is a basis of $\calG_2\calA^{rv}(\uparrow)$.

There are 3 OC relations to write for $\calG_2\calA^w(\uparrow)$:
$D_2=D_{10}$, $D_3=D_6$, and $D_4=D_8$. Along with the $6T$ relations,
we find that $\{D_1, D_3=D_6=D_9, D_2=D_5=D_7=D_{10}=D_{11},
D_4=D_8=D_{12}\}$ is a basis of $\calG_2\calA^w(\uparrow)$. Similarly
$\{D_1,D_2=\ldots=D_{12}\}$ is a basis of the two-dimensional
$\calG_2\calA^{sw}(\uparrow)$.  When we mod out by FI, only one
diagram remains non-zero in $\calG_2\calA^{rw}(\uparrow)$ and it is $D_1$.

We leave the determination of the primitives and the spaces with a circle
skeleton as an exercise to the reader.

\clearpage\draftcut \section{Glossary of notation} \label{sec:glossary}
Greek letters, then Latin, then symbols:

\noindent
{\small \begin{multicols}{2}
\begin{list}{}{
  \renewcommand{\makelabel}[1]{#1\hfil}
}

\item[{$\alpha$}] maps $\calA^u\to \calA^v$ or $\calA^u\to \calA^w$
  \hfill\ \ref{subsubsec:RelWithu}
\item[{$\Delta$}] cloning, co-product\hfill\ \ref{par:Delta},
  \ref{subsec:Projectivization}
\item[{$\delta$}] Satoh's tube map\hfill\ \ref{subsubsec:TopTube}
\item[{$\delta_A$}] a formal $D_A$\hfill\ \ref{subsec:AlexanderProof}
\item[{$\theta$}] inversion, antipode\hfill\ \ref{par:theta}
\item[{$\iota$}] an inclusion $\wB_n\to\wB_{n+1}$\hfill\
  \ref{subsubsec:McCool}
\item[{$\iota$}] interpretation map\hfill\ \ref{subsec:AlexanderProof},
  \ref{subsubsec:IAM}
\item[{$\iota$}] inclusion $\attr_n\to\calP^w(\uparrow_n)$\hfill\ 
  \ref{subsec:ATSpaces}
\item[{$\lambda$}] a formal $EZ$\hfill\ \ref{subsec:AlexanderProof}
\item[{$\nu$}] the invariant of the unknot\hfill\ \ref{subsec:KTG}
\item[{$\xi_i$}] the generators of $F_n$\hfill\ \ref{subsubsec:McCool}
\item[{$\pi$}] the projection $\calP^w(\uparrow_n)\to\fraka_n\oplus\tder_n$
  \hfill\ \ref{subsec:ATSpaces}
\item[{$\Sigma$}] a virtual surface\hfill\ \ref{subsubsec:TopTube}
\item[{$\sigma_i$}] a crossing between adjacent strands\hfill\ 
  \ref{subsubsec:Planar}
\item[{$\sigma_{ij}$}] strand $i$ crosses over strand $j$\hfill\ 
  \ref{subsubsec:Abstract}
\item[{$\varsigma$}] the skeleton morphism\hfill\ \ref{subsubsec:Planar}
\item[{$\phi$}] log of an associator\hfill\ \ref{subsec:KTG}
\item[{$\Phi$}] an associator\hfill\ \ref{subsec:KTG}
\item[{$(\varphi^i)$}] a basis of $\frakg^\ast$\hfill\ \ref{subsec:LieAlgebras}
\item[{$\psi_\beta$}] ``operations''\hfill\ \ref{subsec:AlgebraicStructures}
\item[{$\omega_1$}] a formal 1-wheel\hfill\ \ref{subsec:AlexanderProof}

\item
\item[{$\fraka_n$}] $n$-dimensional Abelian Lie algebra\hfill\ 
  \ref{subsec:ATSpaces}
\item[{$\calA$}] a candidate projectivization\hfill\ \ref{subsec:Expansions}
\item[{$\calA(G)$}] associated graded of $G$\hfill\
  \ref{subsubsec:FTAlgebraic}
\item[{$\calA^{sv}$}] $\calD^v$ mod 6T, RI\hfill\ \ref{subsec:vw-tangles}
\item[{$\calA^{sw}$}] $\calD^w$ mod $\aft$, TC, RI\hfill\ 
  \ref{subsec:vw-tangles}
\item[{$\calA^{sw}$}] $\proj\wTFo$\hfill\ \ref{subsec:fproj}
\item[{$\calA^{sw}$}] $\proj\wTF$\hfill\ \ref{subsubsec:AwWen}
\item[{$\calA^{(s)w}$}] $\calA^{w}$ and/or $\calA^{sw}$\hfill\ 
  \ref{subsec:fproj}
\item[{$\calA^u$}] chord diagrams mod rels for KTGs\hfill\ \ref{subsec:KTG}
\item[{$\calA^v$}] $\calD^v$ mod 6T\hfill\ \ref{subsec:vw-tangles}
\item[{$\calA^w$}] $\calD^w$ mod $\aft$, TC\hfill\ \ref{subsec:vw-tangles}
\item[{$\calA^w$}] $\proj\wTFo$ without RI\hfill\ \ref{subsec:fproj}
\item[{$\calA^-(\uparrow_n)$}] $\calA^-$ for pure $n$-tangles\hfill\ 
  \ref{subsec:ATSpaces}
\item[{$\calA^-_n$}] $\calD^v_n$ mod relations\hfill\
  \ref{subsubsec:FTPictorial}
\item[{$\calA^{-t}$}] $\calA^-$ allowing trivalent vertices\hfill\
  \ref{subsec:Jacobi}
\item[{$\calA^-(\uparrow)$}] $\calD^v(\uparrow)$ mod relations\hfill\
  \ref{subsec:FTforvwKnots}
\item[{$\calA^-(\bigcirc)$}] $\calA^-(\uparrow)$ for round skeletons\hfill\
  \ref{subsec:SomeDimensions}
\item[{$\calA^u$}] usual chord diagrams\hfill\ \ref{subsec:RelWithKont}
\item[{$A(K)$}] the Alexander polynomial\hfill\ \ref{subsec:Alexander}
\item[{$A_e$}] 1D orientation reversal\hfill\ \ref{subsubsec:wops}
\item[{$\aAS$}] arrow-AS relations\hfill\ \ref{subsec:Jacobi}
\item[{$\Ass$}] associative words\hfill\ \ref{subsec:ATSpaces}
\item[{$\Ass^+$}] non-empty associative words\hfill\ \ref{subsec:ATSpaces}
\item[{$a$}] maps $u\to v$ or $u\to w$\hfill\ \ref{subsubsec:RelWithu}
\item[{$a_{ij}$}] an arrow from $i$ to $j$\hfill\
  \ref{subsubsec:FTPictorial}
\item[{$\calB^w$}] unitrivalent arrow diagrams\hfill\ \ref{subsec:Jacobi}
\item[{$\calB^w_n$}] $n$-coloured unitrivalent arrow \newline
  diagrams\hfill\ \ref{subsec:ATSpaces}
\item[{$B$}] the matrix $T(\exp(-xS)-I)$\hfill\ \ref{subsec:AlexanderProof}
\item[{$b_{ij}^k$}] structure constants of $\frakg^\ast$\hfill\ 
  \ref{subsec:LieAlgebras}
\item[{$C$}] the invariant of a cap\hfill\ \ref{subsec:wTFExpansion}
\item[{CC}] the Commutators Commute relation\hfill\ \ref{subsec:Jacobi}
\item[{CP}] the Cap-Pull relation\hfill\ \ref{subsubsec:wrels},
  \ref{subsec:fproj}
\item[{CW}] Cap-Wen relations\hfill\ \ref{subsubsec:WenRels}
\item[{$c$}] a chord in $\calA^u$\hfill\ \ref{subsec:KTG}
\item[{$\der$}] Lie-algebra derivations\hfill\ \ref{subsec:ATSpaces}
\item[{$\calD^v$, $\calD^w$}] arrow diagrams for v/w-tangles\hfill\
  \ref{subsec:vw-tangles}
\item[{$\calD^v_n$}] arrow diagrams for braids\hfill\
  \ref{subsubsec:FTPictorial}
\item[{$\calD^{-t}$}] $\calD^-$ allowing trivalent vertices\hfill\
  \ref{subsec:Jacobi}
\item[{$\calD^v(\uparrow)$}] arrow diagrams long knots\hfill\ 
  \ref{subsec:FTforvwKnots}
\item[{$D_A$}] either $D_L$ or $D_R$\hfill\ \ref{subsec:Jacobi}
\item[{$D_L$}] left-going isolated arrow\hfill\ \ref{subsec:Jacobi}
\item[{$D_R$}] right-going isolated arrow\hfill\ \ref{subsec:Jacobi}
\item[{$\divop$}] the ``divergence''\hfill\ \ref{subsec:ATSpaces}
\item[{$d_k$}] strand deletion\hfill\ \ref{par:deletions}
\item[{$d_i$}] the direction of a crossing\hfill\ \ref{subsec:Alexander}
\item[{$E$}] the Euler operator\hfill\ \ref{subsec:AlexanderProof}
\item[{$\tilE$}] the normalized Euler operator\hfill\ 
  \ref{subsec:AlexanderProof}
\item[{$F$}] a map $\calA^w\to\calA^w$\hfill\ \ref{subsec:fproj}
\item[{$F$}] the main~\cite{AlekseevTorossian:KashiwaraVergne} unknown
  \hfill\ \ref{subsec:EqWithAT}
\item[{FI}] Framing Independence\hfill\ \ref{subsec:SomeDimensions}
\item[{FR}] Flip Relations\hfill\ \ref{subsubsec:WenRels},
  \ref{subsubsec:AwWen}
\item[{$F_n$}] the free group\hfill\ \ref{subsubsec:McCool}
\item[{$\FA_n$}] the free associative algebra\hfill\ \ref{par:action}
\item[{$\fil$}] a filtered structure\hfill\ \ref{subsec:Expansions}
\item[{$\frakg$}] a finite-dimensional Lie algebra\hfill\
  \ref{subsec:LieAlgebras}
\item[{$\calG_m$}] degree $m$ piece\hfill\ \ref{subsubsec:FTPictorial}
\item[{$\calI$}] augmentation ideal\hfill\ \ref{subsubsec:FTAlgebraic}, 
  \ref{subsec:Projectivization}
\item[{$I\frakg$}] $\frakg^\ast\rtimes\frakg$\hfill\ \ref{subsec:LieAlgebras}
\item[{$\IAM$}] Infinitesimal Alexander\newline Module\hfill\ 
  \ref{subsec:AlexanderProof}, \ref{subsubsec:IAM}
\item[{$\IAM^0$}] $\IAM$, before relations\hfill\ \ref{subsubsec:IAM}
\item[{$\aIHX$}] arrow-IHX relations\hfill\ \ref{subsec:Jacobi}
\item[{$i_u$}] an inclusion $F_n\to\wB_{n+1}$\hfill\ \ref{subsubsec:McCool}
\item[{$J$}] a map $\TAut_n\to\exp(\attr_n)$\hfill\ \ref{subsec:ATSpaces}
\item[{$j$}] a map $\TAut_n\to\attr_n$\hfill\ \ref{subsec:ATSpaces}
\item[{$\calK^u$}] usual knots\hfill\ \ref{subsec:RelWithKont}
\item[{KTG}] Knotted Trivalent Graphs\hfill\ \ref{subsec:KTG}
\item[{$\lie_n$}] free Lie algebra\hfill\ \ref{subsec:ATSpaces}
\item[{$l$}] a map $\tder_n\to\calP^w(\uparrow_n)$\hfill\ \ref{subsec:ATSpaces}
\item[{M}] the ``mixed'' move\hfill\ \ref{subsec:VirtualKnots}
\item[{$\calO$}] an ``algebraic structure''\hfill\ 
  \ref{subsec:AlgebraicStructures}
\item[{OC}] the Overcrossings Commute relation\hfill\ \ref{subsec:wBraids}
\item[{$\calP^w_n$}] primitives of $\calB^w_n$\hfill\ \ref{subsec:ATSpaces}
\item[{$\calP^-(\uparrow)$}] primitives of $\calA^-(\uparrow)$\hfill\ 
  \ref{subsec:FTforvwKnots}
\item[{$\calP^-(\uparrow_n)$}] primitives of $\calA^-(\uparrow_n)$\hfill\ 
  \ref{subsec:ATSpaces}
\item[{$\PvB_n$}] the group of pure v-braids\hfill\ \ref{subsubsec:Planar}
\item[{$\PwB_n$}] the group of pure w-braids\hfill\ \ref{subsec:wBraids}
\item[{$\proj$}] projectivization \hfill\ \ref{subsec:Projectivization}
\item[{$\calR$}] the relations in $\IAM$\hfill\ \ref{subsubsec:IAM}
\item[{$R$}] $Z(\overcrossing)$\hfill\ \ref{subsec:wBraidExpansion}
\item[{$R$}] the ring $\bbZ[X,X^{-1}]$\hfill\ \ref{subsubsec:IAM}
\item[{$R_1$}] the augmentation ideal of $R$\hfill\ \ref{subsubsec:IAM}
\item[{$R$}] the invariant of a crossing\hfill\ \ref{subsec:wTFExpansion}
\item[{RI}] Rotation number Independence\hfill\ \ref{subsec:FTforvwKnots}
\item[{R123}] Reidemeister moves\hfill\ \ref{subsec:VirtualKnots}
\item[{R4}] a Reidemeister move for foams/graphs\hfill\ \ref{subsubsec:wrels}
\item[{\Rs}] the ``spun'' R1 move\hfill\ \ref{subsec:VirtualKnots}
\item[{$\sder$}] special derivations\hfill\ \ref{subsec:sder}
\item[{$\calS$}] the circuit algebra of skeletons\hfill\ 
  \ref{subsec:CircuitAlgebras}
\item[{$\SAut_n$}] the group $\exp(\sder_n)$\hfill\ \ref{subsec:KTG}
\item[{$S(K)$}] a matrix of signs\hfill\ \ref{subsec:Alexander}
\item[{$S_k$}] complete orientation reversal\hfill\ 
  \ref{subsec:UniquenessForTangles}
\item[{$S_e$}] complete orientation reversal\hfill\ \ref{subsubsec:wops}
\item[{$S_n$}] the symmetric group\hfill\ \ref{subsubsec:Planar}
\item[{$\aSTU$}] arrow-STU relations\hfill\ \ref{subsec:Jacobi}
\item[{$s_i$}] a virtual crossing between adjacent strands\hfill\ 
  \ref{subsubsec:Planar}
\item[{$s_i$}] the sign of a crossing\hfill\ \ref{subsec:Alexander}
\item[{$\sKTG$}] signed long KTGs\hfill\ \ref{subsec:KTG}
\item[{$\sl$}] self-linking\hfill\ \ref{subsec:VirtualKnots}
\item[{TV}] Twisted Vertex relations\hfill\ \ref{subsubsec:WenRels}
\item[{$\tder$}] tangential derivations\hfill\ \ref{subsec:ATSpaces}
\item[{$\attr_n$}] cyclic words\hfill\ \ref{subsec:ATSpaces}
\item[{$\attr^s_n$}] cyclic words mod degree 1\hfill\ \ref{subsec:ATSpaces}
\item[{$\calT^w_\frakg$}] a map ${\calA}^w\to\calU(I\frakg)$\hfill\ 
  \ref{subsec:LieAlgebras}
\item[{$\TAut_n$}] the group $\exp(\tder_n)$\hfill\ \ref{subsec:ATSpaces}
\item[{TC}] Tails Commute\hfill\ \ref{subsubsec:FTPictorial}
\item[{$T(K)$}] the ``trapping'' matrix\hfill\ \ref{subsec:Alexander}
\item[{$\calU$}] universal enveloping algebra\hfill\ \ref{subsec:LieAlgebras}
\item[{UC}] Undercrossings Commute\hfill\ \ref{subsec:wBraids}
\item[{$u$}] a map $\tder_n\to\calP^w(\uparrow_n)$\hfill\ \ref{subsec:ATSpaces}
\item[{$u_e$}] strand unzips\hfill\ \ref{subsubsec:wops}
\item[{$u_k$}] strand unzips\hfill\ \ref{par:unzip}
\item[{$\uB_n$}] the (usual) braid group\hfill\ \ref{subsubsec:Planar}
\item[{$\uT$}] u-tangles\hfill\ \ref{subsec:sder}
\item[{$V$}] a finite-type invariant\hfill\ \ref{subsubsec:FTPictorial}
\item[{$V$, $V^+$}] the invariant of a (positive) vertex\hfill\ 
  \ref{subsec:wTFExpansion}
\item[{$V^-$}] the invariant of a negative vertex\hfill\ 
  \ref{subsec:wTFExpansion}
\item[{II}] Vertex Invariance\hfill\ \ref{subsec:fproj}
\item[{VR123}] virtual Reidemeister moves\hfill\ \ref{subsec:VirtualKnots}
\item[{$\vB_n$}] the virtual braid group\hfill\ \ref{subsubsec:Planar}
\item[{$\vT$}] v-tangles\hfill\ \ref{subsec:vw-tangles}
\item[{$\vT(\uparrow_n)$}] pure $n$-component v-tangles\hfill\ 
  \ref{subsec:ATSpaces}
\item[{$W$}] $Z(w)$\hfill\ \ref{subsubsec:ZwithWen}
\item[{$W_m$}] weight system\hfill\ \ref{subsubsec:FTPictorial}
\item[{$W^2$}] Wen squared\hfill\ \ref{subsubsec:WenRels}
\item[{$w$}] the map $x^k\mapsto w_k$\hfill\ \ref{subsec:Alexander}
\item[{$w$}] the wen\hfill\ \ref{subsec:TheWen}
\item[{$w_i$}] flip ring $\#i$\hfill\ \ref{subsubsec:FlyingRings}
\item[{$w_k$}] the $k$-wheel\hfill\ \ref{subsec:Jacobi}
\item[{$\wB_n$}] the group of w-braids\hfill\ \ref{subsec:wBraids}
\item[{$\wT$}] w-tangles\hfill\ \ref{subsec:vw-tangles}
\item[{$\wT(\uparrow_n)$}] pure $n$-component w-tangles\hfill\ 
  \ref{subsec:ATSpaces}
\item[{$\wTF$}] w-tangled foams with wens\hfill\ \ref{subsec:TheWen}
\item[{$\wTFo$}] orientable w-tangled foams\hfill\ \ref{subsec:wTFo}
\item[{$X$}] an indeterminate\hfill\ \ref{subsec:Alexander}
\item[{$X_n,\, \tilde{X}_n$}] moduli of horizontal rings\hfill\ 
  \ref{subsubsec:FlyingRings}
\item[{$x_i$}] the generators of $FA_n$\hfill\ \ref{par:action}
\item[{$(x_j)$}] a basis of $\frakg$\hfill\ \ref{subsec:LieAlgebras}
\item[{$Y_n,\, \tilde{Y}_n$}] moduli of rings\hfill\ \ref{subsubsec:NonHorRings}
\item[{$Z$}] expansions \hfill\ throughout
\item[{$Z_\calA$}] an $\calA$-expansion\hfill\ \ref{subsec:Expansions}
\item[{$Z^u$}] the Kontsevich integral\hfill\ \ref{subsec:RelWithKont}

\item
\item[{4T}] $4T$ relations\hfill\ \ref{subsec:KTG}
\item[{$\aft$}] $\aft$ relations\hfill\ \ref{subsubsec:FTPictorial}
\item[{$6T$}] $6T$ relations\hfill\ \ref{subsubsec:FTPictorial}
\item[{$\semivirtualover,\,\semivirtualunder$}] semi-virtual
  crossings\hfill\ \ref{subsubsec:FTPictorial}
\item[{$\sslash$}] right action\hfill\ \ref{subsubsec:McCool}
\item[{$\uparrow$}] a ``long'' strand\hfill\ throughout
\item[{$\up$}] the quandle operation\hfill\ \ref{subsec:Projectivization}
\item[{$\up_2$}] doubled $\up$\hfill\ \ref{subsec:Projectivization}
\item[{$*$}] the adjoint on $\calA^w(\uparrow_n)$\hfill\ \ref{subsec:ATSpaces}

\end{list}
\end{multicols}}

\clearpage\draftcut

\if\draft y
  \clearpage
  Everything below is to be blanked out before the completion of this paper.
  \section*{Recycling}

\begin{exercise} Do the same for the obviously-defined ``w-links'',
excluding the material about the Alexander polynomial. Note that the wheels
that are obtained in the case of w-links have legs coloured by the
components of the w-link in question. Hence if there is more than one
component, the number of such wheels grows exponentially in the degree and
thus $Z$ contains more information than can be coded in a polynomial of
even a multi-variable polynomial.
\end{exercise}

\vskip -5mm
\parpic[r]{\raisebox{-14mm}{$\input figs/CC.pstex_t$}}
\begin{conjecture} In the case of ordinary links seen as w-links,
if we mod out the target space of $Z$ by the ``Commutators Commute''
relation shown on the right, what remains of the wheels part of $Z$
is precisely the multi-variable Alexander polynomial.
\end{conjecture}

Note that $D \in \tder_n$ is never an arrow on a single strand (these are
elements of $\fraka_n$),
and hence $\operatorname{div}D$ is never a 1-wheel, more precisely it never
has a degree 1 component. Thus, even
though the target space of div is $\attr_n/(\text{deg }1)$, we can just as
well think of it as a
map to $\attr_n$ itself.

\subsection{The Injectivity of $i_u\colon F_n\to\wB_{n+1}$}
\label{subsec:FreeInW}

Just for completeness, we sketch here an algebraic proof of the
injectivity of the map $i_u\colon F_n\to\wB_{n+1}$ discussed in
Section~\ref{subsubsec:McCool}. There's some circularity in our argument
--- we need this injectivity in order to motivate the definition of the
map $\Psi\colon \wB_n\to\Aut(F_n)$, and in the proof below we use $\Psi$
to prove the injectivity of $i_u$. But $\Psi$ exists regardless of how
its definition is motivated, and it can be shown to be well defined by
explicitly verifying that it respects the relations defining $\wB_n$. So
our proof is logically valid.

\begin{claim} The map $i_u\colon F_n\to\wB_{n+1}$ is injective.
\end{claim}

\begin{proof} (sketch). Let $H$ be the subgroup of $\wB_{n+1}$ MORE
\end{proof}

\subsection{Finite Type Invariants of v-Braids and w-Braids, in some
Detail} \label{subsec:FTDetails}

As mentioned in Section~\ref{subsec:wBraids}, w-braids are v-braids modulo
an additional relation. So we start with a discussion of finite type
invariants of v-braids. For simplicity we take our base ring to be $\bbQ$;
everywhere we could replace it by an arbitrary field of characteristic
$0$\footnote{Or using the variation of constants method, we can simply
declare that $\bbQ$ is an arbitrary field of characteristic $0$.},
and many definitions make sense also over $\bbZ$ or even with $\bbQ$
replaced by an arbitrary Abelian group.

\subsubsection{Basic Definitions.}
Let $\bbQ\vB_n$ denote group ring of $\vB_n$, the algebra of formal linear
combinations of elements of $\vB_n$, and let $\bbQ S_n$ be the group
ring of $S_n$. The skeleton homomorphism of Remark~\ref{rem:Skeleton}
extends to a homomorphism $\varsigma\colon \bbQ\vB_n\to\bbQ S_n$. Let $\calI$
(or $\calI_n$ when we need to be more precise) denote the kernel of
the skeleton homomorphism; it is the ideal in $\bbQ\vB_n$ generated
by formal differences of v-braids having the same skeleton. One may
easily check that $\calI$ is generated by differences of the form
$\overcrossing-\virtualcrossing$ and $\virtualcrossing-\undercrossing$.
Following~\cite{GoussarovPolyakViro:VirtualKnots} we call such differences
``semi-virtual crossings'' and denote them by $\semivirtualover$
and $\semivirtualunder$, respectively\footnote{The signs in
$\semivirtualover\leftrightarrow\overcrossing-\virtualcrossing$ and
$\semivirtualunder\leftrightarrow\virtualcrossing-\undercrossing$ are
``crossings come with their sign and their virtual counterparts come with
the opposite sign''.}. In a similar manner, for any natural number $m$
the $m$th power $\calI^m$ of $\calI$ is generated by ``$m$-fold iterated
differences'' of v-braids, or equally well, by ``$m$-singular v-braids'',
which are v-braids that also have exactly $m$ semi-virtual crossings
(subject to relations which we don't need to specify).

Let $V\colon \vB_n\to A$ be an invariant of v-braids with values in some vector
space $A$. We say that $V$ is ``of type $m$'' (for some $m\in\bbZ_{\geq
0}$) if its linear extension to $\bbQ\vB_n$ vanishes on $\calI^{m+1}$
(alternatively, on all $(m+1)$-singular v-braids, in clear analogy with the
standard definition of finite type invariants). If $V$ is of type $m$
for some unspecified $m$, we say that $V$ is ``of finite type''. Given
a type $m$ invariant $V$, we can restrict it to $\calI^m$ and as it
vanishes on $\calI^{m+1}$, this restriction can be regarded as an element
of $\left(\calI^m/\calI^{m+1}\right)^\star$. If two type $m$ invariants
define the same element of $\left(\calI^m/\calI^{m+1}\right)^\star$ then
their difference vanishes on $\calI^m$, and so it is an invariant of type
$m-1$. Thus it is clear that an understanding of $\calI^m/\calI^{m+1}$
will be instrumental to an inductive understanding of finite type
invariants. Hence the following definition.

\begin{definition} The projectivization\footnote{Why ``projectivization''?
See Section~\ref{subsec:Projectivization}.} $\proj\vB_n$ is the direct sum
\[ \proj\vB_n:=\bigoplus_{m\geq 0} \calI^m/\calI^{m+1}. \]
Note that throughout this paper, whenever we write an infinite direct sum,
we automatically complete it. Therefore an element in $\proj\vB_n$ is an
infinite sequence $D=(D_0, D_1,\dots)$, where $D_m\in\calI^m/\calI^{m+1}$.
The projectivization $\proj\vB_n$ is a graded space, with the degree $m$
piece being $\calI^m/\calI^{m+1}$.
\end{definition}

We proceed with the study of $\proj\vB_n$ (and thus of finite type
invariants of v-braids) in three steps. In
Section~\ref{subsubsec:ArrowDiagrams} we introduce a space $\calD^v_n$ and
a surjection $\rho_0\colon \calD^v_n\to\proj\vB_n$. In Section~\ref{subsubsec:6T}
we find some relations in $\ker\rho_0$, most notably the $6T$ relation, and
introduce the quotient $\calA^v_n:=\calD^v_n/6T$. And then in
Section~\ref{subsubsec:UFTI} we introduce the notion of a ``universal
finite type invariant'' and explain how the existence of such a gadget
proves that $\proj\vB_n$ is isomorphic to $\calA^v_n$ (in a more
traditional language this is the statement that every weight system
integrates to an invariant).

Unfortunately, we do not know if there is a universal finite type invariant
of v-braids. Thus in Section~\ref{subsec:wbraids} we return to the subject
of w-braids and prove the weaker statement that there exists a universal
finite type invariant of w-braids.

\subsubsection{Arrow Diagrams.} \label{subsubsec:ArrowDiagrams}

We are looking for a space that will surject on $\calI^m/\calI^{m+1}$. In
other words, we are looking for a set of generators for $\calI^m$, and
we are willing to call two such generators the same if their difference
is in $\calI^{m+1}$. But that's easy. Left and right multiples of the
formal differences $\semivirtualover=\overcrossing-\virtualcrossing$
and $\semivirtualunder=\virtualcrossing-\undercrossing$
generate $I$, so products of the schematic form
\begin{equation} \label{eq:GeneratingProduct}
  B_0 (\semivirtualover|\semivirtualunder) B_1
  (\semivirtualover|\semivirtualunder) B_2 \cdots
  B_{m-1} (\semivirtualover|\semivirtualunder) B_m
\end{equation}
\parpic[r]{\begin{picture}(0,0)%
\includegraphics{figs/SemiVirtRels.pstex}%
\end{picture}%
%
%
\setlength{\unitlength}{3158sp}%
\begingroup\makeatletter\ifx\SetFigFont\undefined%
\gdef\SetFigFont#1#2#3#4#5{%
  \reset@font\fontsize{#1}{#2pt}%
  \fontfamily{#3}\fontseries{#4}\fontshape{#5}%
  \selectfont}%
\fi\endgroup%
\begin{picture}(2184,1824)(1309,-1573)
\put(2401,-286){\makebox(0,0)[b]{\smash{{\SetFigFont{10}{12.0}{\rmdefault}{\mddefault}{\updefault}{\color[rgb]{0,0,0}$=$}%
}}}}
\put(2401,-1336){\makebox(0,0)[b]{\smash{{\SetFigFont{10}{12.0}{\rmdefault}{\mddefault}{\updefault}{\color[rgb]{0,0,0}$=$}%
}}}}
\end{picture}%
 }
\noindent generate $\calI^m$ (here $(\semivirtualover|\semivirtualunder)$
means ``either a $\semivirtualover$ or a $\semivirtualunder$'', and there
are exactly $m$ of those in any product). Furthermore, inside such a
product any $B_k$ can be replaced by any other v-braid $B'_k$ having the
same skeleton (e.g., with $\varsigma(B_k)$), for then $B_k-B'_k\in\calI$
and the whole product changes by something in $\calI^{m+1}$.  Also,
the relations in~\eqref{eq:R3} and in~\eqref{eq:MixedRelations} imply
the relations shown on the right for $\semivirtualover$, and similar
relations for $\semivirtualunder$. With this freedom, a product as
in~\eqref{eq:GeneratingProduct} is determined by the strand-placements
of the $\semivirtualover$'s and the $\semivirtualunder$'s. That is,
for each semi-virtual crossing in such a product, we only need to
know which strand number is the ``over'' strand, which strand number
is the ``under'' strand, and a sign that determines whether it is the
positive semi-virtual $\semivirtualover$ or the negative semi-virtual
$\semivirtualunder$. This motivates the following definition.

\begin{definition} A ``horizontal $m$-arrow diagrams'' (analogues to the
``chord diagrams'' of, say, \cite{Bar-Natan:OnVassiliev}) is an ordered
pair $(D,\beta)$ in which $D$ is a word of length $m$ in the alphabet
$\{a^+_{ij},a^-_{ij}\colon i,j\in\{1,\ldots,n\},\,i\neq j\}$ and $\beta$ is a
permutation in $S_n$. Let $\calD_m^{vh}$ be the space of formal linear
combinations of horizontal $m$-arrow diagrams. We usually use a pictorial
notation for horizontal arrow diagram, as demonstrated in
Figure~\ref{fig:Dvh}.
\end{definition}

\begin{figure}
\parpic[r]{\hspace{-5mm}\raisebox{-29mm}{$\input figs/Dvh.pstex_t$}}
\caption{
  The horizontal $3$-arrow diagram
  $(D,\beta)=$ $(a^+_{12}a^-_{41}a^+_{23},\,3421)$ and its image via 
  $\rho_0$. The first arrow, $a^+_{12}$ starts at strand $1$, ends 
  at strand $2$ and carries a $+$ sign, so it is mapped to a positive
  semi-virtual crossing of strand $1$ over strand $2$. Likewise the second
  arrow $a^-_{41}$ maps to a negative semi-virtual crossing of strand 
  $4$ over strand $1$, and $a^+_{23}$ to a positive semi-virtual
  crossing of strand $2$ over strand $3$. We also show one possible
  choice for a representative of the image of $\rho_0(D,\beta)$ in
  $\calI^m/\calI^{m+1}$: it is a v-braid with semi-virtual crossings as
  specified by $D$ and whose overall skeleton is $3421$.
} \label{fig:Dvh}
\end{figure}

There is a surjection $\rho_0\colon \calD_m^{vh}\to\calI^m/\calI^{m+1}$. The
definition of $\rho_0$ is suggested by the first paragraph of this section
and an example is shown in Figure~\ref{fig:Dvh}; we will skip the formal
definition here. We also skip the formal proof of the surjectivity of
$\rho_0$.

Finally, consider the product $\semivirtualover\cdot\semivirtualunder$ and
use the second Reidemeister move for both virtual and non-virtual
crossings:
\[
  \semivirtualover\semivirtualunder
  = (\overcrossing-\virtualcrossing)(\virtualcrossing-\undercrossing)
  = \overcrossing\virtualcrossing+\virtualcrossing\undercrossing
    - \overcrossing\undercrossing - \virtualcrossing\virtualcrossing
  = (\overcrossing\virtualcrossing-1) + (\virtualcrossing\undercrossing)
  = \semivirtualover\virtualcrossing - \virtualcrossing\semivirtualunder.
\]
If a total of $m-1$ further semi-virtual crossings are multiplied into this
equality on the left and on the right, along with arbitrary further
crossings and virtual crossings, the left hand side of the equality becomes
a member of $\calI^{m+1}$, and therefore, as a member of
$\calI^m/\calI^{m+1}$, it is $0$. Thus with ``$\ldots$'' standing for
extras added on the left and on the right, we have that in
$\calI^m/\calI^{m+1}$,
\[ 0 = 
  \ldots(
    \semivirtualover\virtualcrossing-\virtualcrossing\semivirtualunder
  )\ldots
  = \rho_0(\ldots??\ldots)
\]

MORE.

\subsubsection{The $6T$ Relations.} \label{subsubsec:6T}

MORE.

\subsubsection{The Notion of a Universal Finite Type Invariant.}
\label{subsubsec:UFTI}

MORE.

\subsection{Finite type invariants of w-braids} \label{subsec:wbraids}

MORE.
 
  \section*{To Do}

\par\noindent{\bf Sorted.}
\begin{itemize}
\item Finish the paper.
\item Freeze Mathematica notebooks.
\end{itemize}

\fi

\end{document}